\documentclass[a4paper,12pt]{article}

\usepackage {amssymb,latexsym}
\usepackage {amsmath}
\usepackage[french]{babel}

\newtheorem{theorem}{Th\'{e}or\`{e}me}[section]
\newtheorem{prop}{Proposition}[section]
\newtheorem{lemma}{Lemme}[section]
\newtheorem{cor}{Corollaire}[section]

\newtheorem{defi}{D\'{e}finition}[section]
\newtheorem{exer}{Exercice}[section]
\newtheorem{exem}{Exemple}[section]
\newtheorem{prob}{Probl\`{e}me}[section]
\newtheorem{rem}{Remarque}[section]
\newenvironment{prooof}{
                        \noindent{\bf\small D\'{e}monstration: }\small}
                                       {\hfill {$\mathbf \Box$}\\}
\newcommand{\bdefi}{\begin{defi}}
\newcommand{\edefi}{\end{defi}}
\newcommand{\bexer}{\begin{exer}\small\rm}
\newcommand{\eexer}{\end{exer}}
\newcommand{\bexem}{\begin{exem}\rm}
\newcommand{\eexem}{\end{exem}}
\newcommand{\bsat}{\begin{theorem}}
\newcommand{\esat}{\end{theorem}}
\newcommand{\bprop}{\begin{prop}}
\newcommand{\eprop}{\end{prop}}
\newcommand{\bcor}{\begin{cor}}
\newcommand{\ecor}{\end{cor}}
\newcommand{\blem}{\begin{lemma}}
\newcommand{\elem}{\end{lemma}}
\newcommand{\brem}{\begin{rem}\rm}
\newcommand{\erem}{\end{rem}}
\newcommand{\bbew}{\begin{prooof}}
\newcommand{\ebew}{\end{prooof}}
\newcommand{\bprob}{\begin{prob}}
\newcommand{\eprob}{\end{prob}}

\newcommand{\beq}{\begin{equation}}
\newcommand{\bea}{\begin{eqnarray}}
\newcommand{\eea}{\end{eqnarray}}
\newcommand{\beas}{\begin{eqnarray*}}
\newcommand{\eeas}{\end{eqnarray*}}

\newcommand{\real}{\mathbb R}
\newcommand{\korps}{\mathbb K}

\newcommand{\complex}{\mathbb C}
\newcommand{\nat}{\mathbb N}
\newcommand{\entier}{\mathbb Z}

\newcommand{\ben}{\begin{enumerate}}
\newcommand{\een}{\end{enumerate}}
\newcommand{\ra}{\rightarrow}
\newcommand{\Cinf}{\mathcal{C}^\infty}
\newcommand{\CinfR}[1]{\mathcal{C}^\infty(#1,\real)}
\newcommand{\CinfC}[1]{\mathcal{C}^\infty(#1,\complex)}
\newcommand{\CinfK}[1]{\mathcal{C}^\infty(#1,\korps)}
\newcommand{\Ginf}{{\Gamma}^\infty}
\newcommand {\prol} {\mathsf{prol}}
\newcommand {\id} {{\mathrm{id}}}
\newcommand{\Dop}{\mathbf{D}}
\makeatletter
\@addtoreset{equation}{section}

\makeatother

\begin{document}

\pagestyle{empty}
~~~~~~~~~~~~
%\vspace{2cm}
\begin{center}
  {\LARGE \bf (Bi)Modules, morphismes et r\'{e}duction des star-produits:
     le cas symplectique, feuilletages et obstructions}\\[2mm]

\vspace{0.5cm}

D\'{e}di\'{e} au soixanti\`{e}me anniversaire d'Alan Weinstein.

\vspace{0.5cm}

Mars 2004

\vspace{1cm}

       {\bf Martin Bordemann} \\[2mm]
        Laboratoire des Math\'{e}matiques et Applications\\
%        Facult\'{e} des Sciences et Techniques\\
        Universit\'{e} de Haute Alsace, Mulhouse \\
%        4, rue des Fr\`{e}res Lumi\`{e}re\\
%        68093 Mulhouse, France \\
        e--mail: M.Bordemann@uha.fr \\

\vspace{1cm}

\footnotesize
{\bf English Abstract}\\[5mm]

\begin{minipage}{13cm}
We give a differential geometric framework for the description of
(bi)modules, morphisms and reduction of star-products in deformation
quantization in terms of multidifferential operators along maps.
We show that algebra morphisms deform Poisson maps and left modules
on function spaces deform coisotropic maps. Let $C$ be a closed coisotropic
submanifold of a Poisson manifold $M$: a star-product on $M$ is representable
(as differential operators on the space of smooth functions) on $C$ iff it
is equivalent to a star-product
for which the vanishing ideal $\mathcal{I}[[\nu]]$ of $C$ is a left ideal
in the deformed algebra (``coisotropic creed'' by Lu and Weinstein).
Moreover, for symplectic reduction the fact is important that
deformation of the vanishing ideal $\mathcal{I}$ as a subalgebra
automatically is either a left or right ideal in case $M$ is symplectic.
We show that a symplectic star-product is always representable on a
coisotropic submanifold in case the reduced space exists, and we classify
all bimodule structures on the function space on $C$ (with respect to the
deformed and the reduced algebra) in terms of its first de Rham cohomology.
We show that in the symplectic case there are obstructions to the
representability of a star-product which are related to a differential
topological invariant of the foliation of the coisotropic submanifold, the
so-called Atiyah-Molino class. If the latter vanishes, we prove that
star-products with suitable Deligne classes can be represented.
Finally, by a Fedosovian analysis we show that a Poisson map $\phi:M\ra M'$
between two symplectic manifolds is not always deformable into an algebra
morphism: the obstructions are once again related to the Atiyah-Molino
class of the symplectic orthogonal bundle of the kernel bundle of the
tangent map of $\phi$: if this class vanishes, we show that $\phi$ can
be quantized subject to suitable conditions on the Deligne classes.
The example of the cotangent bundle of the two-dimensional torus and
complex projective space are discussed.
\end{minipage}
\normalsize
\end{center}

\newpage
\footnotesize
\tableofcontents
\normalsize
\newpage

\pagestyle{plain}

\section*{Introduction}
  \addcontentsline{toc}{section}{Introduction}

Depuis l'article fondateur de Bayen, Flato, Fr{\o}nsdal, Lichnerowicz et
Sternheimer
en 1978 \cite{BFFLS78}, la quantification par d\'{e}formation est devenu un
domaine
de recherche assez vaste qui couvre plusieur domaines alg\'{e}briques comme
la th\'{e}orie des d\'{e}formations formelles des alg\`{e}bres associatives
voir \cite{Ger63}, \cite{Ger64}, \cite{Ger66}, \cite{Ger68}, \cite{Ger74})
et,
plus r\'{e}cemment, la th\'{e}orie des op\'{e}rades et le domaine des vari\'{e}t\'{e}s
de Poisson qui inclut celui des vari\'{e}t\'{e}s symplectiques, voir \'{e}galement
\cite{Wei95}, \cite{Ste98} et \cite{DS00}.
Plus pr\'{e}cis\'{e}ment, on consid\`{e}re
dans cette th\'{e}orie
une vari\'{e}t\'{e} de Poisson $(M,P)$ et l'alg\`{e}bre commutative associative
$\CinfK{M}$ de toutes les fonctions de classe
$\Cinf$ sur $M$ \`{a} valeurs dans $\korps=\real$ ou $\complex$.
Un star-produit
$*=\sum_{r=0}^\infty \nu^r\mathsf{C}_r$ sur $M$ est une d\'{e}formation formelle
associative de la multiplication point-par-point $\mathsf{C}_0$ telle
que $\mathsf{C}_1(f,g)-\mathsf{C}_1(g,f)$ est \'{e}gal \`{a} deux fois le crochet
de Poisson $P(df,dg)$ quels que soient $f,g\in\CinfK{M}$. Les applications
$\mathsf{C}_r$ sont des op\'{e}rateurs bidiff\'{e}rentiels, et le param\`{e}tre formel
$\nu$ correpond \`{a} $\frac{i\hbar}{2}$ dans des situations convergentes.
L'alg\`{e}bre d\'{e}form\'{e}e $\mathcal{A}=\big(\CinfK{M}[[\nu]],*\big)$ sur l'anneau
$\korps[[\nu]]$ est interpr\'{e}t\'{e}e comme alg\`{e}bre d'observables quantiques.
Deux star-produits $*$ et $*'$ sont dits \'{e}quivalents lorsque les alg\`{e}bres
d\'{e}form\'{e}es sont isomorphes par rapport \`{a} une s\'{e}rie formelle d'op\'{e}rateurs
diff\'{e}rentiels dont le coefficient d'ordre $0$ est l'application identique.
L'existence des star-produits
a \'{e}t\'{e} d\'{e}montr\'{e}e pour le cas symplectique par De Wilde et Lecomte (1983
dans \cite{DL83}), par Fedosov (1985 dans \cite{Fed85}) et par
Omori, Maeda, Yoshioka (1991 dans \cite{OMY91}). La classification \`{a}
\'{e}quivalence pr\`{e}s --dans le cas symplectique-- a \'{e}t\'{e} obtenue par Deligne
\cite{Del95}, Nest/Tsygan \cite{NT95a}, \cite{NT95b} et
Bertelsson/Cahen/Gutt \cite{BCG97} en termes
des s\'{e}ries formelles \`{a} coefficients dans le deuxi\`{e}me
groupe de cohomologie de de Rham $H^2_{dR}(M,\korps)$. En 1997
Kontsevitch a montr\'{e} l'existence et la classification dans le cas d'une
vari\'{e}t\'{e} de Poisson g\'{e}n\'{e}rale, \cite{Kon97b} (voir \'{e}galement
\cite{AMM99}).\\
Le grand avantage de cette approche est premi\`{e}rement
l'absence de toute obstruction,
deuxi\`{e}mement la `limite classique $\lim_{\nu\to 0}$' et donc la transition
\`{a} la m\'{e}canique classique automatiquement install\'{e}es, et troisi\`{e}mement
une description en termes g\'{e}om\'{e}triques intuitives. La quantification par
d\'{e}formation d\'{e}crit d'abord l'alg\`{e}bre des observables sans faire r\'{e}f\'{e}rence
\`{a} un espace (pr\'{e})hilbertien, ce qui est tr\`{e}s raisonnable pour la comparaison
de
la situation classique avec la situation quantique: les ph\'{e}nom\`{e}nes
`bizarres' de la m\'{e}canique quantique, comme les {\em verschr\"{a}nkten Zust\"{a}nde}
(les \'{e}tats enchev\^{e}t\-r\'{e}s) de Schr\"{o}din\-ger dans la description d'un syst\`{e}me
compos\'{e} de deux sys\-t\`{e}mes, ne se montrent pas de mani\`{e}re ensembliste
au niveau des observables: dans le cas classique on a $\CinfK{M_1\times M_2}=
\CinfK{M_1}\hat{\otimes}\CinfK{M_2}$, ce qui est un produit tensoriel
topologique ainsi que
dans le cas quantique usuel
o\`{u} l'alg\`{e}bre (stellaire) des observables $\mathcal{B}$ est donn\'{e}e par
$\mathcal{B}_1\hat{\otimes}\mathcal{B}_2$.

Face \`{a} la classe d'alg\`{e}bres associatives fournies par la quantification
par d\'{e}formation, les questions suivantes sont imm\'{e}diates d'un point de
vue alg\`{e}brique:
\begin{itemize}
 \item Quels sont les {\bf modules} $\mathcal{M}$ de $\mathcal{A}$, i.e. les
  $\korps[[\nu]]$-modules $\mathcal{M}$ et les applications
  $\korps[[\nu]]$-bilin\'{e}aires
  \[
      \rho:\mathcal{A}\times \mathcal{M}\ra \mathcal{M}
  \]
  avec
  \[
     \rho(f*g)(\varphi)=\rho(f)\rho(g)(\varphi)~~~~
      \mathrm{et}~~~\rho(1)\varphi=\varphi
  \]
  quels que soient $f,g\in\mathcal{A}$ et $\varphi\in\mathcal{M}$ ?
 \item Quels sont les {\bf morphismes} de l'alg\`{e}bre $\mathcal{A}$ dans
    $\mathcal{A}'$,
 i.e. les applications $\korps[[\nu]]$-lin\'{e}aires
 \[
    \Phi:\mathcal{A}\ra \mathcal{A}'
 \]
 avec
 \[
     \Phi(f*g)=\big(\Phi(f)\big)*'\big(\Phi(g)\big)
       ~~~\mathrm{et}~~~\Phi(1)=1
 \]
 quels que soient $f,g\in\mathcal{A}$ ?
\end{itemize}
N\'{e}ansmoins, l'\'{e}tude de ces
questions a commenc\'{e} beaucoup plus tard que celle de l'existence et
l'\`{e}quiva\-lence des star-produits,
et on n'a pas encore beaucoup de r\'{e}sultats:\\
Quant aux modules --qui sont \`{a} premi\`{e}re vue tr\`{e}s int\'{e}ressants pour un
physicien qui veut savoir comment les espaces hilbertiens de la m\'{e}canique
quantique trouvent leur place dans cette th\'{e}orie-- je vois
une des raisons dans le fait
qu'une des principales motivations de
la m\'{e}thode de la quantification par d\'{e}formation, con\c{c}ue en
$1978$, a \'{e}t\'{e} l'objectif de se d\'{e}barasser de l'espace de Hilbert utilis\'{e}
en m\'{e}canique quantique usuelle
et de tirer tout renseignement sur le
spectre et, plus g\'{e}n\'{e}ralement, sur les probabilit\'{e}s de transition de
physique, directement de l'alg\`{e}bre d\'{e}form\'{e}e. Pour cette fin, l'outil
clef de l'\'{e}cole de Bayen, Flato, Fr{\o}nsdal, Lichnerowicz et
Sternheimer a \'{e}t\'{e} la {\em star-exponentielle} qui
--interpr\'{e}t\'{e}e dans un cadre distributionnel d\'{e}pendant de la situation--
a donn\'{e} le spectre correct des principaux exemples de
la m\'{e}canique quantique, voir \cite{BFFLS78}. \\
Quant aux morphismes: on ne voit pas beaucoup de morphismes d'alg\`{e}bres
d'observables en m\'{e}canique quantique entre deux syst\`{e}mes de physique r\'{e}els,
mis \`{a} part les inclusions d'alg\`{e}bres
$\mathcal{A}_1\ra\mathcal{A}$ qu'on rencontre entre autres dans la th\'{e}orie
des champs quantiques locale de Haag/Kastler (voir par exemple
\cite{Haa92}) ou du type
$\mathcal{A}_1\ra\mathcal{A}_1\hat{\otimes}\mathcal{A}_2$ provenant
de la m\'{e}canique quantique de plusieurs particules et correspondant
\`{a} des projections
\beq \label{EqProjApplPoisEinl}
 \begin{array}{c}
    M_1 \times  M_2 \\
    pr_1\downarrow  \\
    M_1
 \end{array}
\end{equation}

Pourtant, les modules ont apparu plusieurs fois dans un cadre de
quantification par d\'{e}formation: Fedosov a \'{e}tudi\'{e} des {\em modules
projectifs}
\`{a} un nombre fini de g\'{e}n\'{e}rateurs pour obtenir des r\'{e}sultats en $K$-th\'{e}orie
et th\'{e}orie d'indice, voir \cite{Fed85}, \cite{Fed96}. D'apr\`{e}s le th\'{e}or\`{e}me
de Serre/Swan, il s'y agit d'une d\'{e}formation des espaces de sections de
classe $\Cinf$ d'un fibr\'{e} vectoriel sur $M$. En outre, motiv\'{e}s par des
questions
du physicien math\'{e}maticien Klaus Fredenhagen, avec S.Waldmann on a \'{e}tudi\'{e}
les repr\'{e}sentations sur un espace pr\'{e}hilbertien
du type {\em Gel'fand-Naimark-Segal (GNS)} pour les alg\`{e}bres
d\'{e}form\'{e}es, \cite{BW98}, tout en utilisant l'ordre nonarchim\'{e}dien de
l'anneau $\real[[\nu]]$ pour d\'{e}finir les fonctions lin\'{e}aires positives.
On avait ainsi retrouv\'{e} les repr\'{e}\-sen\-ta\-tions de Schr\"{o}dinger, de Wick et
WKB. De plus, il n'est pas \'{e}tonnant que pour certains
{\em calculs symboliques}
des op\'{e}rateurs diff\'{e}rentiels sur une vari\'{e}t\'{e} $Q$ (voir par exemple
\cite{widom:1980a} pour des d\'{e}finitions)
l'espace $\CinfK{Q}[[\nu]]$
apparaisse comme module pour l'alg\`{e}bre $\big(\CinfK{T^*Q}[[\nu]],*\big)$
avec le star-produit $*$ bien-choisi, voir par exemple \cite{CFS92}
(o\`{u} une trace (voir \'{e}galement \cite{OMY92}) et la cohomologie cyclique ont
\'{e}t\'{e} \'{e}tudi\'{e}es), \cite{BNW98},
\cite{BNW99} et \cite{BNPW03}. Mais il y a aussi des calculs symboliques
qui ne d\'{e}finissent pas un star-produit sur $T^*Q$, voir les calculs
invariants par l'alg\`{e}bre de Lie $\mathfrak{sl}(n+1,\real)$ ou
$\mathfrak{so}(p+1,n-p+1)$ pour $Q=\real^n$, voir \cite{LO99} et
\cite{DLO99} et le calcul invariant par un changement projectif de la
connexion, \cite{Bor02p}. Finalement, S.~Waldmann a continu\'{e} l'\'{e}tude
des repr\'{e}sentations GNS et de l'induction de Rieffel, \cite{Wal00},
\cite{waldmann:2002b}, arrivant (avec H.~Bursztyn)
\`{a} l'\'{e}tude de l'\'{e}quivalence de Morita
de deux alg\`{e}bres d\'{e}form\'{e}es non isomorphes,
\cite{bursztyn.waldmann:2000b}, \cite{bursztyn.waldmann:2000a}
\cite{bursztyn.waldmann:2001b}, \cite{bursztyn.waldmann:2001a} et
\cite{bursztyn.waldmann:2002a}.\\
Les morphismes on \'{e}t\'{e} \'{e}tudi\'{e}s pour le cas particulier d'une sym\'{e}trie, i.e.
o\`{u} $\phi:\mathcal{U}\ra \mathcal{A}$ avec $\mathcal{U}$ l'alg\`{e}bre
enveloppante
d'une alg\`{e}bre de Lie $\mathfrak{g}$ ou, plus g\'{e}n\'{e}ralement, un groupe
quantique,
c.-\`{a}-d. une alg\`{e}bre de Hopf d\'{e}formant l'alg\`{e}bre enveloppante de
$\mathfrak{g}$.
L'analogue classique pour le premier cas est une application moment
\beq\label{EqApplMomentErsteMal}
  J:M\ra \mathfrak{g}^*,
\end{equation}
c.-\`{a}-d. une application de Poisson dans $\mathfrak{g}^*$ muni de la
structure de Poisson lin\'{e}aire. La d\'{e}formation $\mathbb{J}:M\ra
\mathfrak{g}^*[[\nu]]$ de $J$ est appel\'{e}e
{\em application moment quantique} par Xu, \cite{Xu98}, quand
\[
 \langle \mathbb{J},\xi\rangle * \langle \mathbb{J},\eta\rangle
    - \langle \mathbb{J},\eta\rangle * \langle \mathbb{J},\xi\rangle
    = 2\nu \langle \mathbb{J},[\xi,\eta]\rangle
\]
quels que soient $\xi,\eta\in\mathfrak{g}$.
 Un star-produit
$*$ sur $\CinfK{M}[[\nu]]$ qui admet une telle application est dit
{\em covariant}, \cite{ACMP83}. Fedosov a montr\'{e} l'existence des star-produits
covariants au cas o\`{u} les champs de vecteurs hamiltoniens
$X_{\langle J,\xi\rangle}$ (avec $\xi\in\mathfrak{g}$) pr\'{e}servent une
connexion $\nabla$ sur $M$ (par exemple dans le cas d'un groupe de Lie
compact ou plus g\'{e}n\'{e}ralement pour une action propre d'un groupe de Lie),
\cite{Fed96}, \cite{GR03}. \\
Un cas particulier de cette situation est
donn\'{e} par
$\mathfrak{g}$ ab\'{e}lienne (et $M$ symplectique, $J$ une submersion presque
partout et
$2\dim \mathfrak{g}=\dim M$): dans ce cas on trouve des {\em syst\`{e}mes
int\'{e}grables
quantiques}, voir par exemple \cite{CRM75}, \cite{OP83}, \cite{BW99}.

Une situation de la g\'{e}om\'{e}trie symplectique, pour laquelle on ne verrait pas
--\`{a} premi\`{e}re vue-- une liaison avec les deux questions mentionn\'{e}es ci-dessus,
est donn\'{e}e par la {\em r\'{e}duction symplectique}, voir \cite{MW74},
\cite{Wei77b},
\cite{AM85}: on consid\`{e}re une {\em sous-vari\'{e}t\'{e} co\"{\i}sotrope} $i:C\ra M$
d'une vari\'{e}t\'{e} symplectique $(M,\omega)$ (pour laquelle le fibr\'{e} $E$ des
sous-espaces orthogonaux \`{a} tous les espaces tangents fait partie du fibr\'{e}
tangent
$TC$). Puisque $E$ est int\'{e}grable, $C$ est muni d'un
{\em feuilletage r\'{e}gulier},
et on
peut consid\'{e}rer l'espace des feuilles, $\pi:C\ra M_{\mathrm{red}}$. Au cas o\`{u}
$M_{\mathrm{red}}$ est muni de la structure d'une vari\'{e}t\'{e} diff\'{e}rentiable
telle que la projection canonique $\pi$ soit une submersion,
{\em l'espace r\'{e}duit}
$M_{\mathrm{red}}$ est muni d'une structure symplectique
$\omega_{\mathrm{red}}$
telle que
\beq\label{EqFormesSympReduites}
   i^*\omega = \pi^*\omega_{\mathrm{red}}
\end{equation}
D'un point de vue de physique, la r\'{e}duction symplectique est tr\`{e}s importante:
la vari\'{e}t\'{e} $M$ correspond \`{a} une description `\`{a} trop de degr\'{e}s de
libert\'{e}', mais qui est dans la plupart des cas `plus simple pour le calcul'.
$C$ est le r\'{e}sultat des `contraintes de physique'
qui --s'ils sont `first class' (c.-\`{a}-d. $C$ co\"{\i}sotropes)--
`agissent' sur $C$ et y produisent des `orbites de jauge'.
L'espace r\'{e}duit y est
d\'{e}crit par des choix de repr\'{e}sentants locaux dans les feuilles, on
`fixe la jauge'. \footnote{Le point de vue de r\'{e}duction adopt\'{e} dans ce
rapport est {\em purement g\'{e}om\'{e}trique} dans le sens que c'est seulement
la sous-vari\'{e}t\'{e} feuillet\'{e}e $C$ de $M$ qui est \'{e}tudi\'{e}e et en g\'{e}n\'{e}ral pas les
\'{e}ventuelles contraintes qui la d\'{e}finissent. J.~\'{S}niatycki
m'a signal\'{e} qu'en physique on \'{e}tudie \'{e}galement la r\'{e}duction \`{a} l'aide
de l'id\'{e}al engendr\'{e} par un ensemble de contraintes sp\'{e}cifiques qui n'\'{e}puise
en g\'{e}n\'{e}ral pas l'id\'{e}al annulateur $\mathcal{I}$ de $C$.}\\
La {\em r\'{e}duction de Marsden-Weinstein} \cite{MW74} est un cas
particulier
de cette construction o\`{u} $C=J^{-1}(0)$ pour la valeur r\'{e}guli\`{e}re $0$
d'une application moment $J$. Ici les restrictions des fonctions invariantes
par l'action de
l'alg\`{e}bre de Lie se projettent sur les fonctions sur l'espace r\'{e}duit.
En physique, on rencontre souvent des situations singuli\`{e}res
o\`{u} $C$ n'est plus une
vari\'{e}t\'{e}, mais plut\^{o}t une vari\'{e}t\'{e} stratifi\'{e}e (voir \cite{Pfl01} pour
des d\'{e}finitions) ainsi que l'espace r\'{e}duit, \cite{SL91} et \cite{LPS01}.
Malgr\'{e} leur importance, je ne vais pas aborder ces sujets dans ce
rapport.
\\
Dans les premiers travaux sur la quantification par d\'{e}formation de la
r\'{e}duction
symplectique, on n'a consid\'{e}r\'{e} que les deux vari\'{e}t\'{e}s symplectiques
$M$ et $M_{\mathrm{red}}$: plus pr\'{e}cis\'{e}ment, on a voulu calculer un
star-produit
sur l'espace r\'{e}duit en termes d'un star-produit sur $M$: Fedosov a appliqu\'{e}
sa construction
\cite{Fed85} au
cas d'une application moment d'une action du cercle, voir \cite{Fed94b} et
ensuite au cas d'une
action d'un groupe de Lie compact, voir \cite{Fed98}; une formule explicite
pour un star-produit sur l'espace projectif complexe a \'{e}t\'{e} trouv\'{e}e par
r\'{e}duction
dans \cite{BBEW96a}, et pour les vari\'{e}t\'{e}s de Grassmann dans \cite{Sch97}.
Finalement, la r\'{e}duction symplectique a \'{e}t\'{e} formul\'{e}e --dans le cadre de la
r\'{e}duction de Marsden-Weinstein et la cohomologie BRST-- entre autres pour
le cas de l'action propre d'un groupe de Lie
dans \cite{BHW00}.\\
On peut se poser la question de savoir {\em quel est le r\^{o}le alg\'{e}brique de la
sous-vari\'{e}t\'{e} co\"{\i}sotrope $C$}: \`{a} premi\`{e}re vue,
au diagramme de vari\'{e}t\'{e}s diff\'{e}rentiables correspond un diagramme
d'anneaux commutatifs:
\beq\label{EqDiagReduction}
     \begin{array}{ccc}
         &                         &   M \\
         &\stackrel{i}{\nearrow}   &     \\
       C &                         &     \\
         &\stackrel{\searrow}{\pi} &     \\
         &                         &  M_{\mathrm{red}}
     \end{array} ~~~~~~~~~~~
     \begin{array}{ccc}
         &                         &   \CinfK{M} \\
         &\stackrel{i^*}{\swarrow}   &     \\
       \CinfK{C} &                         &     \\
         &\stackrel{\nwarrow}{\pi^*} &     \\
         &                         &  \CinfK{M_{\mathrm{red}}}
     \end{array}
\end{equation}
mais la vari\'{e}t\'{e} $C$ n'est en g\'{e}n\'{e}ral pas une vari\'{e}t\'{e} de Poisson
telle que la projection $\pi$ soit une application de Poisson: par exemple,
pour la sph\`{e}re
$S^3$ --vue comme sous-vari\'{e}t\'{e} co\"{\i}sotrope de $\real^4$-- la projection
$\pi$
est la fibration de Hopf sur $M_{\mathrm{red}}=S^2=\complex P(1)$; s'il
existait une structure
de Poisson sur $S^3$ qui se projettait sur celle de $S^2$ (provenant de
la forme de Fubini-Study) on aurait un
feuilletage de codimension
$1$ sur $S^3$ qui serait transverse aux fibres de la fibration de Hopf,
qui serait donc donn\'{e}
par un rel\`{e}vement de l'action de $SU(2)$ sur $S^2$ et dont les feuilles
seraient toutes
diff\'{e}omorphes \`{a} $S^2$ contrairement \`{a} un th\'{e}or\`{e}me de Novikov (1965)
que les feuilles
compactes ont la topologie du $2$-tore. Alors on ne peut pas
esp\'{e}rer que la quantification de $C$ consiste \'{e}galement en une d\'{e}formation
associative
formelle de l'alg\`{e}bre de fonctions $\CinfK{C}$. On obtient une indication en
regardant
{\em l'id\'{e}al annulateur de $C$},
\beq \label{EqIdAnnEinleitung}
  \mathcal{I}:=\{f\in \CinfK{M}~|~f(c)=0~~\forall~c\in C\},
\end{equation}
qui est un id\'{e}al par rapport \`{a} la multiplication point-par-point, mais
--gr\^{a}ce \`{a} la co\"{\i}sotropie-- seulement
une sous-alg\`{e}bre par rapport au crochet de Poisson. Alors il n'est pas dur
\`{a} deviner
que la quantification de $\mathcal{I}$ devrait \^{e}tre
`une moyenne entre id\'{e}al et sous-alg\`{e}bre', alors un {\em id\'{e}al  unilat\'{e}ral}
dans l'alg\`{e}bre d\'{e}form\'{e}e $(\CinfK{M}[[\nu]],*)$. On a \'{e}t\'{e} tr\`{e}s fier
\`{a} Freiburg de cette id\'{e}e en 1997 en \'{e}tudiant la r\'{e}duction symplectique
quantique, jusqu'\`{a} notre d\'{e}couverte de l'article de Jiang Hua Lu \cite{Lu93}
de 1993
(sur l'action des alg\`{e}bres de Hopf) dans lequel elle a d\'{e}j\`{a} parl\'{e} de la
quantification de
$\mathcal{I}$ en tant qu'id\'{e}al unilat\'{e}ral comme `{\em coisotropic creed}'
sans toucher \`{a} la quantification par d\'{e}formation. Si l'id\'{e}al annulateur
d\'{e}form\'{e} $\hat{\mathcal{I}}$ est un id\'{e}al \`{a} gauche et en plus isomorphe
--en tant que $\korps[[\nu]]$-module--
\`{a} $\mathcal{I}[[\nu]]$, alors le quotient $\CinfK{M}[[\nu]]/\hat{\mathcal{I}}$
serait
isomorphe \`{a} $\CinfK{C}[[\nu]]$ et muni de la structure canonique d'un
$\CinfK{M}[[\nu]]$-module \`{a} gauche. Nous avons \'{e}tudi\'{e} la d\'{e}formation
de l'id\'{e}al annulateur --dans le cas particulier de la r\'{e}duction
de Marsden-Weinstein-- dans \cite{BHW00} en 1999.\\
On est donc men\'{e} \`{a} interpr\'{e}ter le diagramme \ref{EqDiagReduction} comme
l'action de l'alg\`{e}bre $\CinfK{M}[[\nu]]$ et \'{e}galement de
$\CinfK{M_{\mathrm{red}}}[[\nu]]$
sur $\CinfK{C}[[\nu]]$, alors comme un {\em bimodule} par rapport \`{a} ces
alg\`{e}bres,
chose que j'ai apprise par Alan Weinstein.
En outre, il est bien connu que l'alg\`{e}bre de Poisson des fonctions
$\Cinf$ sur l'espace
r\'{e}duit
se d\'{e}crit comme quotient $\mathcal{N}(\mathcal{I})/\mathcal{I}$ o\`{u}
\beq\label{EqNIEinleitung}
 \mathcal{N}(\mathcal{I}):=\{f\in\CinfK{M}~|~\{f,g\}\in\mathcal{I}~\forall
                               ~g\in\mathcal{I}\}
\end{equation}
d\'{e}signe {\em l'id\'{e}alisateur de Lie de $\mathcal{I}$} qui est bien connu
dans la th\'{e}orie de BRST, voir par exemple \cite{BHW00}.
Une notion plus faible de la r\'{e}duction serait de d\'{e}former
$\mathcal{N}(\mathcal{I})$
et l'id\'{e}al annulateur
seulement comme {\em sous-alg\`{e}bres} de l'alg\`{e}bre d\'{e}forme\'{e}e telles
que la d\'{e}formation de $\mathcal{I}$ ne soit un id\'{e}al bilat\`{e}re que de
la d\'{e}formation de $\mathcal{N}(\mathcal{I})$. Ainsi on aurait aussi une
alg\`{e}bre quotient
bien d\'{e}finie qui correspondrait \`{a} l'alg\`{e}bre r\'{e}duite.

Toutes ces consid\'{e}rations montrent que la th\'{e}orie des (bi)modules en
quantification par d\'{e}formation n'est pas negligeable
si on essaie de comprendre la r\'{e}duction symplectique. De plus, si la
codimension de la sous-vari\'{e}t\'{e} co\"{\i}sotrope $C$ est strictement positive,
{\em l'espace
$\CinfK{C}$ en tant que $\CinfK{M}$-module n'est plus projectif}, donc on
ne peut en g\'{e}n\'{e}ral plus utiliser les r\'{e}sultats positifs de la quantification
des modules projectifs qui sont cit\'{e}s ci-dessus. A priori, l'existence
des repr\'{e}sentations n'est pas \'{e}vidente, et il faut \'{e}tudier
{\em `die Bedingungen der M\"{o}glichkeit von Darstellung \"{u}berhaupt'}, i.e.
les {\em obstructions} possibles.\\

Dans ce rapport, je voudrais bien \'{e}tudier les morphismes, les modules et la
r\'{e}duction dans le cadre de la quantification par d\'{e}formation.
Il y a plusieurs objectifs:
\ben
 \item Etablissement d'un cadre de g\'{e}om\'{e}trie diff\'{e}rentielle qui permet
   de sp\'{e}cifier la classe des morphismes et des modules qu'on va regarder.
   Ici surtout les {\em op\'{e}rateurs multidiff\'{e}entiels le long des applications}
   se montrent tr\`{e}s utils.
  \item Identification des `limites classiques' des morphismes et des
  re\-pr\'{e}\-sen\-ta\-tions
     (sur l'espace de fonctions $\Cinf$ \`{a} valeurs dans $\korps$ sur une
     vari\'{e}t\'{e}
     $C$): {\em applications de Poisson} et
     {\em applications co\"{\i}sotropes}.
  \item R\'{e}duction des trois probl\`{e}mes des morphismes, des modules et de la
    r\'{e}duction symplectique au seul probl\`{e}me des modules.
  \item Classification des bimodules apparaissant dans la version quantique de la
      r\'{e}duction.
  \item Etude des {\em obstructions} (r\'{e}currentes et totales) pour la
       repr\'{e}sentabilit\'{e} d'un star-produit sur une sous-vari\'{e}t\'{e}
       co\"{\i}sotrope: ici un {\em invariant de feuilletage}
       (la {\em classe d'Atiyah-Molino})
       jouera un r\^{o}le
      d\'{e}cisif.
  \item Description d'une classe d'exemples pour lesquels la repr\'{e}sentation
       est possible.
  \item Etude des {\em obstructions totales pour le probl\`{e}me de la
         quantification des
      morphismes de Poisson} entre deux vari\'{e}t\'{e}s {\em symplectiques}
  \item Description d'une classe d'exemples pour des morphismes de
         Poisson quantifiables.
\een

\noindent On a obtenu les r\'{e}sultats suivants:

Le paragraphe \ref{SecDefPrelim} sert de rappel
de plusieurs notions dont j'aurai besoin par la suite:
les vari\'{e}t\'{e}s et les applications de Poisson
(paragraphe \ref{SubSecPois}), la th\'{e}orie (alg\'{e}brique) des sous-vari\'{e}t\'{e}s
  (paragraphe \ref{SubSecSousvar}), les sous-vari\'{e}t\'{e}s et applications
  co\"{\i}sotropes (paragraphe \ref{SubSecCoiso}), les vari\'{e}t\'{e}s
  feuillet\'{e}es (paragraphe \ref{SubSecVarFeuil}), la r\'{e}duction symplectique
  (pargraphe \ref{SubSecRed1}), les op\'{e}rateurs multidiff\'{e}rentiels le long
  des applications (paragraphe \ref{SubSecOpMultidiffAppl})
  et les star-produits (paragraphe \ref{SubSecStarprod}).

   En paragraphe \ref{SubSecMor1}),
  on montre d'abord le r\'{e}sultat sans doute bien connu que
   tout {\em morphisme} $\Phi:\big(\CinfK{M}[[\nu]],*\big)\ra
   \big(\CinfK{M'}[[\nu]],*'\big)$ entre deux al\-g\`{e}b\-res d\'{e}form\'{e}es
   est la d\'{e}formation du pull-back avec un {\em morphisme de Poisson}
   $\phi:M'\ra M$ (voir le lemme \ref{LHomoEgalPoissonOrdreZero}).
   Alors il est naturel
  d'introduire la notion d'un {\em morphisme diff\'{e}rentiel}
  (voir la d\'{e}finition \ref{DDefHomStarGenDiff}) qui est donn\'{e}
  par une s\'{e}rie d'op\'{e}rateurs diff\'{e}rentiels le long de $\phi$
  (voir le paragraphe \ref{SubSecOpMultidiffAppl}) et de poser le
  {\em probl\`{e}me de la quantification des applications de Poisson}.
  De plus, on \'{e}tablit
  le lien entre les morphismes et les applications moment quantiques
  de Xu \cite{Xu98}, voir la proposition \ref{PMorphApplMomentQ}.

  En paragraphe \ref{SubSecRep1} on donne une d\'{e}finition d'une
  {\em repr\'{e}sentation $\rho$ de l'alg\`{e}bre d\'{e}form\'{e}e} $\big(\CinfK{M},*\big)$
  comme un homomorphisme dans l'alg\`{e}bre associative de toutes les
  s\'{e}ries formelles \`{a} coefficients dans {\em l'alg\`{e}bre des op\'{e}rateurs
  diff\'{e}rentiels des sections d'une fibr\'{e} vectoriel $E$ sur une vari\'{e}t\'{e}
  $C$} (d\'{e}finition \ref{DefRep1}) et les d\'{e}finitions usuelles
  des bimodules (l'\'{e}qn (\ref{EVraiBimoduleQuant})), des \'{e}quivalences
  (d\'{e}finition \ref{DefModEquiv}) de repr\'{e}senta\-tions et des liaisons
  entre deux star-produits et leurs repr\'{e}sentations
  (o\`{u} on admet des morphismes (l'\'{e}qn (\ref{EqDeuxRepLiees})).
   En particulier, pour le fibr\'{e}
  trivial $C\times \korps$ on peut montrer (proposition
  \ref{PRepCois}) que `la limite classique',
  i.e. le terme d'ordre $0$ de la repr\'{e}sentation est donn\'{e} par
  $\rho_0(f)\varphi=(i^*f)\varphi$ o\`{u} $i:C\ra M$ est une {\em application
  co\"{\i}sotrope} (d\'{e}finition \ref{DefAplCoi}), i.e. pour
  laquelle $\mathrm{Ker} i^*$ est une sous-alg\`{e}bre de Poisson de
  $\CinfK{M}$. Ceci provoque la d\'{e}finition d'une {\em repr\'{e}sentation
  diff\'{e}rentielle} o\`{u} les coefficients $\rho_r$ de la repr\'{e}sentation
  sont des op\'{e}ra\-teurs bidiff\'{e}rentiels le long des applications $i:C\ra M$
  et $id:C\ra C$ (d\'{e}finition \ref{DRepDiff}). En particulier,
  pour une {\em sous-vari\'{e}t\'{e} co\"{\i}sotrope $C$ d'une vari\'{e}t\'{e} de Poisson}
  le {\em probl\`{e}me de la repr\'{e}sentabilit\'{e} d'un star-produit donn\'{e}} se pose,
  i.e. si on peut trouver une repr\'{e}sentation diff\'{e}ren\-ti\-elle sur les
  fonctions sur $C$ qui d\'{e}forme l'application de restriction $i^*$
  (d\'{e}finition \ref{DRepresentabilite}). Cette d\'{e}formation de
  l'application de restriction a d\'{e}j\`{a} \'{e}t\'{e} \'{e}tudi\'{e}e dans le contexte
  BRST de
  la r\'{e}duction de Marsden Weinstein dans \cite{BHW00}.  De plus, il n'est pas
  dur \`{a} voir qu'un star-produit est repr\'{e}sentable sur $C$ si et seulement
  s'il est repr\'{e}sentable dans un voisinage ouvert (par exemple un
  voisinage tubulaire) de $C$ (corollaire
  \ref{CRepresentabiliteLocale}). On rappelle quelques calculs symboliques
  des op\'{e}rateurs
  diff\'{e}rentiels sur une vari\'{e}t\'{e} diff\'{e}rentiable $Q$ qui
  forment des exemples importants des repr\'{e}sentations diff\'{e}rentielles
  des star-produits
  sur $T^*Q$, voir \cite{widom:1980a}, \cite{BNPW03}, et on en
  tire la conclusion que dans une carte contractile de sous-vari\'{e}te
  la restriction de tout star-produit symplectique est repr\'{e}sentable
  (corollaire \ref{CRepresentCarteSympl}).

 On donne deux d\'{e}finitions possibles de la r\'{e}duction
   symplectique quantique par rapport \`{a} une sous-vari\'{e}t\'{e} co\"{\i}sotrope
   ferm\'{e}e $C$ en paragraphe \ref{SubSecReducStarProduit}:\\
   On appelle un star-produit {\em projetable} si les
   $\korps[[\nu]]$-modules $\mathcal{I}[[\nu]]$ et
   $\mathcal{N}(\mathcal{I})[[\nu]]$ sont des sous-alg\`{e}bres de
   $\big(\CinfK{M}[[\nu]],*\big)$ avec $\mathcal{I}[[\nu]]$ un id\'{e}al
   bilat\`{e}re dans $\mathcal{N}(\mathcal{I})[[\nu]]$ (d\'{e}finition
   \ref{DStarProjet}). Il est imm\'{e}diat que le $\korps[[\nu]]$-module quotient
   $\mathcal{N}(\mathcal{I})[[\nu]]/\mathcal{I}[[\nu]]$, qui
   est isomorphe \`{a} $\CinfK{M_{\mathrm{red}}}[[\nu]]$, est muni d'un
   star-produit (proposition \ref{PStarProjet}) pour lequel il existe
   une formule explicite (\ref{EqDefStarRed}) en termes du star-produit sur
   $M$. On appelle un
   star-produit {\em r\'{e}ductible} lorsqu'il est \'{e}quivalent \`{a} un
   star-produit projetable (d\'{e}finition \ref{DStarProjet}).
   On peut en d\'{e}finir une
   version locale (pour des quotients locaux) (d\'{e}finition
   \ref{DStarProjetLocal}).\\
   La deuxi\`{e}me version est celle du `{\em coisotropic creed}' de
   Lu et Weinstein, \cite{Lu93}, o\`{u} on essaie de r\'{e}aliser
   $\CinfK{C}[[\nu]]$ comme
   $\CinfK{M}[[\nu]]-\CinfK{M_{\mathrm{red}}}[[\nu]]$-bimodule
   (d\'{e}finition \ref{DStarStarredBimod}). En
   particulier, dans ce cas le star-produit sur $M$ doit \^{e}tre
   repr\'{e}sentable sur $C$.

 D'apr\`{e}s un th\'{e}or\`{e}me de Weinstein, \cite{Wei88}, le graphe
   d'une application de Poisson $\phi:M'\ra M$ est une sous-vari\'{e}t\'{e}
   co\"{\i}sotrope de $M'\times M$ diff\'{e}omorphe \`{a} $M'$. En paragraphe
   \ref{SubSecMorRep} on montre que la d\'{e}formation de $\phi$ en tant
   que morphisme de star-produits (par rapport \`{a} $*$ et $*'$)
   est \'{e}quivalent \`{a} la d\'{e}formation
   de $\CinfK{M'}[[\nu]]$ en tant que $*$-$*^{\prime\mathrm{opp}}$
   bimodule (proposition \ref{PDefoMorphDefoBimod}): donc, si
   l'on \'{e}large les ensembles de morphismes $Hom(\mathcal{A},\mathcal{B})$
   de la cat\'{e}gorie   des alg\`{e}bres
   associatives sur $\korps[[\nu]]$ par les classes d'isomorphie
   des $\mathcal{A}$-$\mathcal{B}$-bimodules \'{e}tant projectifs et \`{a} un
   nombre fini de g\'{e}n\'{e}rateurs en tant que $\mathcal{B}$-modules \`{a} droite
   (voir par exemple \cite{Lod92}), on n'aura pas plus d'espace pour
   d\'{e}former les applications de Poisson.

 En paragraphe \ref{SubSecRepDeformIdealAnnul}, on montre d'abord
   le r\'{e}sultat assez utile que tout star-produit $*$ repr\'{e}sentable sur une
   sous-vari\'{e}t\'{e} co\"{\i}sotrope est \'{e}qui\-val\-ent \`{a} un star-produit $*'$
   {\em adapt\'{e} \`{a} la sous-vari\'{e}t\'{e}}, i.e. pour lequel l'id\'{e}al annulateur
   $\mathcal{I}[[\nu]]$ est directement un id\'{e}al \`{a} gauche de $*'$,
   voir la proposition \ref{PKerRho1congId}, le th\'{e}or\`{e}me
   \ref{TIdAnIdGauche} et la d\'{e}finition \ref{DDefStarAdapte}.
   Il est clair que les star-produits adapt\'{e}s ont des repr\'{e}sentations
   canoniques sur les fonctions sur la sous-vari\'{e}t\'{e}.
   La g\'{e}n\'{e}ralisation de cette d\'{e}finition aux op\'{e}rateurs
   multidiff\'{e}rentiels adapt\'{e}s, i.e. $\mathsf{C}(f_1,\ldots,f_{k-1},g)\in
   \mathcal{I}$ quels que soient $g\in\mathcal{I}$,
   $f_1,\ldots,f_{k-1}\in \CinfK{M}$, est une op\'{e}rade (donc ferm\'{e}e par le
   crochet de Gerstenhaber) et joue un r\^{o}le d\'{e}cisif pour la formulation
   d'une `{\em formalit\'{e} adapt\'{e}e}', voir \cite{BGHHW03p} pour plus de
   d\'{e}tails. Une autre \'{e}tude de formalit\'{e} li\'{e}e aux sous-vari\'{e}t\'{e}s
   co\"{\i}sotropes a \'{e}t\'{e} faite par
   Cattaneo
   et Felder dans \cite{CF03p} en utilisant des graphes de Kontsevitch
   avec deux types d'ar\^{e}tes. Si deux star-produits adapt\'{e}s sont
   \'{e}quivalents \`{a} l'aide d'une transformation d'\'{e}quivalence adapt\'{e}e, alors
   leurs repr\'{e}sentations canoniques sont \'{e}quivalentes (proposition
   \ref{PCouplesEquiSiTrafoAdap}). D'un autre c\^{o}t\'{e}, il est possible que
   deux star-produits adapt\'{e}s soient \'{e}quivalents, mais pas via une
   transformation d'\'{e}quivalence adapt\'{e}e, voir l'exemple de la r\'{e}duction
   de l'espace projectif complexe en paragraphe \ref{SubSecEspProj}.
   Mais dans le cas symplectique, si le premier groupe de cohomologie
   verticale de $C$
   (voir l'\'{e}qn (\ref{EqDefdverticale}) du paragraphe \ref{SubSecVarFeuil})
   s'annule, la r\'{e}ciproque de la proposition \ref{PCouplesEquiSiTrafoAdap}
   est vraie (th\'{e}or\`{e}me \ref{Tdaserstemalfalsch}).
   Le {\em commutant} d'une repr\'{e}sentation d'un star-produit $*$ sur une
   sous-vari\'{e}t\'{e}
   co\"{\i}sotrope, i.e. l'alg\`{e}bre de tous les homomorphismes de modules
   est \'{e}tudi\'{e} dans la d\'{e}finition \ref{DefCommutant} et
   montr\'{e} d'\^{e}tre \'{e}gal \`{a} {\em l'id\'{e}alisateur \`{a} droite}
   $\mathcal{N}_*(\mathcal{I})$ de $\mathcal{I}[[\nu]]$ modulo
   $\mathcal{I}[[\nu]]$ en proposition \ref{PCommutantNdeI}. Pour un
   star-produit projetable et adapt\'{e}, on montre dans la m\^{e}me
   proposition \ref{PCommutantNdeI} que le commutant est
   ainsi anti-isomorphe \`{a} l'alg\`{e}bre r\'{e}duite.

  Les sous-vari\'{e}t\'{e}s co\"{\i}sotropes qui ne posent aucun
   probl\`{e}me de repr\'{e}sen\-ta\-tion sont les {\em sous-vari\'{e}t\'{e}s lagrangiennes},
   i.e. des sous-vari\'{e}t\'{e}s co\"{\i}so\-tropes de dimension minimale des
   vari\'{e}t\'{e}s symplectiques: en utilisant un th\'{e}or\`{e}me de Weinstein
   \cite{Wei71} disant qu'il existe un voisinage tubulaire autour de toute
   sous-vari\'{e}t\'{e} lagrangienne ferm\'{e}e $i:L\ra M$ qui est symplectomorphe \`{a} un
   voisinage
   ouvert de la z\'{e}ro-section de $T^*L$, on montre qu'un star-produit
   $*$ est repr\'{e}sentable sur les fonctions sur $L$ ssi la restriction de la
   classe de Deligne de $*$ \`{a} $L$, $i^*[*]$, s'annule
   (corollaire \ref{CorRepLagArbi} en paragraphe
   \ref{SubSecRepSousvarLag}). De plus, les classes d'isomorphie se
   classifient par les s\'{e}ries formelles \`{a} coefficients dans le
   {\em premier groupe de cohomologie de de Rham de $L$}. Le th\'{e}or\`{e}me
   clef pour ces \'{e}nonc\'{e}s est une \'{e}tude du calcul symbolique sur $L$,
   voir le th\'{e}or\`{e}me \ref{TRepLagTStarQ},
   qui a \'{e}t\'{e} faite presque enti\`{e}rement d\'{e}j\`{a} dans \cite{BNPW03}.

   Le cas d'une {\em sous-vari\'{e}t\'{e} co\"{\i}sotrope $i:C\ra M$
    d'une vari\'{e}t\'{e} symplectique $(M,\omega)$ pour laquelle l'espace r\'{e}duit
    $M_{\mathrm{red}}$ existe} est trait\'{e} en paragraphe
    \ref{SubSecEspRedReduction}. On r\'{e}sout compl\`{e}tement
    les probl\`{e}mes de (bi)modules et de
    r\'{e}duction en donnant ainsi une r\'{e}ponse partielle \`{a} une question --pos\'{e}e
    par Weinstein-- de savoir quelles sont les condition d'existence et
    une classification des bimodules,
   voir le th\'{e}or\`{e}me
    \ref{TRepUndAllesEspRedExist}: en utilisant le fait montr\'{e} par
    Weinstein que la
    sous-vari\'{e}t\'{e} $C$ se plonge comme sous-vari\'{e}t\'{e} lagrangienne de
    $M\times M_{\mathrm{red}}$ (voir le diagramme (\ref{EqDiagReduction}))
    et le corollaire \ref{CorRepLagArbi} dans le paragraphe
   \ref{SubSecRepSousvarLag} qui pr\'{e}c\`{e}de, on montre qu'il existe une
   repr\'{e}sentation d'un star-produit $[*]$ sur les fonctions sur $C$
   ssi la restriction de la classe de Deligne de $*$ \`{a} $C$, $i^*[*]$,
    est une {\em classe basique}, i.e. $i^*[*]=\pi^*\beta$ avec $\beta$ une
    s\'{e}rie formelle \`{a} coefficients dans
    $H^2_{dR}(M_{\mathrm{red}},\korps)$. Il s'ensuit que le star-produit
    $*$ et un star-produit $*_r$ sur $M_{\mathrm{red}}$ admettent un
    bimodule (qui d\'{e}forme $i$ et $\pi$) si et seulement si l'analogue
    suivant de la relation (\ref{EqFormesSympReduites}) entre les formes
    symplectiques sur $M$ et sur $M_{\mathrm{red}}$ est vrai:
    \beq \label{EqRelationPrincipaleClassesDel}
       i^*[*] = \pi^*[*_r].
    \end{equation}
    Les classes d'isomorphie de $*$-$*_r$-bimodules sur $C$ sont class\'{e}s
    par les s\'{e}ries formelles \`{a} coefficients dans $H_{dR}^1(C,\korps)$.
    En outre, on a une {\em \'{e}quival\-ence entre l'existence d'un star-produit
    projetable sur $M$ et le coisotropic creed}: le th\'{e}or\`{e}me clef --pour le
    cas symplectique-- est
    le fait que {\em la d\'{e}formation de l'id\'{e}al annulateur $\mathcal{I}$ en
    tant que sous-alg\`{e}bre est automatiquement un id\'{e}al unilat\'{e}ral}
    (th\'{e}or\`{e}me \ref{TIAlgId} en paragraphe \ref{SubSecISousAlg}).
    Finalement, {\em l'alg\`{e}bre r\'{e}duite s'av\`{e}re toujours comme le commutant
    de la repr\'{e}sentation de $*$ sur $C$}.
    Le cas de la r\'{e}duction Marsden Weinstein d'un fibr\'{e} cotangent a recemment
   \'{e}t\'{e} trait\'{e} de fa\c{c}on d\'{e}taill\'{e}e par Kowalzig, Neumaier et Pflaum dans
   \cite{KNP04p}.

  On \'{e}tudie les {\em obstructions \`{a} la repr\'{e}sentabilit\'{e}} d'un
    star-produit
    sur une sous-vari\'{e}t\'{e} co\"{\i}sotrope en paragraphe
    \ref{SecObstrucFeuilletage}: les {\em obstructions r\'{e}currentes},
    i.e. celles qui se posent de nouveau \`{a} chaque ordre de la
    construction, sont faciles \`{a} trouver: puisque tout star-produit
    repr\'{e}sentable est \'{e}quivalent \`{a} un star-produit adapt\'{e}, il suffit
    d'\'{e}tudier les derniers, et on montre dans le
    th\'{e}or\`{e}me \ref{TObstructionsRecurrentes}
    que les obstructions r\'{e}currentes
    se trouvent dans $H^2_P(C,\korps)$, le deuxi\`{e}me groupe de cohomologie
    BRST --pour lequel un cocomplexe est donn\'{e} par les sections du fibr\'{e}
    normal $\Lambda TM|_C/TC$, voir le lemme \ref{LCohdP}.
    Dans le cas
    symplectique, ce groupe co\"{\i}ncide avec le deuxi\`{e}me groupe de
    cohomologie verticale. Un corollaire imm\'{e}diat est un r\'{e}sultat de
    P.Gl\"{o}{\ss}ner de 1998 \cite{Glo98}: si $C$ est de {\em codimension 1},
    {\em alors $*$ est toujours repr\'{e}sentable}, voir le th\'{e}or\`{e}me
    \ref{TGloessnerCodimUn}. Ceci n'est pas aussi \'{e}tonnant: dans le cas o\`{u}
    $C$ est l'image r\'{e}ciproque pour la valeur r\'{e}guli\`{e}re $0$
    d'une fonction $\Cinf$ \`{a} valeurs r\'{e}elles $J$, l'id\'{e}al
    annulateur est l'ensemble
    des multiples de $J$, donc un module projectif de $\CinfK{M}$
    qui est quantifiable. \\
    Mais comme dans le cas du premier th\'{e}or\`{e}me d'existence d'un
    star-produit symplectique de Neroslavski et Vlassov \cite{NV81}
    (en 1981) pour lequel on avait besoin de supposer que $H^3_{dR}(M,\korps)$
    s'annule, les obstructions r\'{e}currentes ne donnent pas en g\'{e}n\'{e}ral
    de vraies conditions n\'{e}cessaires (car il y a des star-produits
    sur $T^*S^3$, voir \cite{BFFLS78}).
    En paragraphe \ref{SubSecObstrRepOrdreTrois} j'ai calcul\'{e} {\em les
    obstructions `totales' jusqu'\`{a} l'ordre $3$} pour la repr\'{e}sentabilit\'{e}
    d'un star-produit symplectique sur une sous-vari\'{e}t\'{e} co\"{\i}sotrope:
    pour le faire, il faut v\'{e}rifier {\em tous les star-produits} jusqu'\`{a}
    l'ordre $3$. Cette t\^{a}che assez fastidieuse se simplifie un tout petit
    peu (mais vraiment pas trop) en utilisant le th\'{e}or\`{e}me de
    Weinstein/Gotay \cite{Got82} sur un voisinage tubulaire de $C$ et
    le fait qu'il existent des connexions symplectiques adapt\'{e}es
    \`{a} cette situation, voir le paragraphe \ref{SubSecConnSympAdaptCois}.
    Jusqu'\`{a} l'ordre 2 on peut calculer le star-produit symplectique
    g\'{e}n\'{e}ral $*$ adapt\'{e} (th\'{e}or\`{e}me \ref{TObsRepOrdreDeux}): la seule
   obstruction est
   \beq \label{EqObstructionDeuxIntro}
       p_v i^*[*]_0=0
   \end{equation}
   i.e. la restriction $p_v$ du coefficient d'ordre $0$ de la classe de
   Deligne
   $[*]$ de $*$ aux feuilles doit s'annuler, condition qui ne pose aucun
   probl\`{e}me. En th\'{e}or\`{e}me \ref{TObstructionsOrdreTrois}, on calcule
   les obstructions \`{a} l'ordre $3$ et on obtient
   \bea
      \lefteqn{\frac{1}{12}\overline{P}^{(3)}
              \big(\kappa_{AM}(C,E),\kappa_{AM}(C,E)\big)} \nonumber \\
        & & ~~~~~~~~~~~+\frac{1}{2}
          \overline{P}^{(1)}([i^*\alpha_0]_{(1,1)},[i^*\alpha_0]_{(1,1)})
        -p_vi^*[*]_1=0, \label{EqObstructionImportanteIntro}
    \eea
    o\`{u}
    \begin{itemize}
     \item $E$ est le sous-fibr\'{e} de $TC$ tangent aux feuilles,
     \item $\alpha_0$ est un repr\'{e}sentant de la classe $[*]_0$,
     \item $\overline{P}^{(k)}$ (avec $k=1,2,3$) d\'{e}signe un accouplement
    de deux $1$-formes le long des feuilles \`{a} valeur dans la $k$me puissance
    sym\'{e}trique du fibr\'{e} $TC/E$ qui est form\'{e} de $k$ copies de la structure
    de Poisson $P$,
     \item $[i^*\alpha_0]_{(1,1)}$ est la classe de $\alpha_0$ vue comme
      $1$-forme le long des feuilles \`{a} valeurs dans $TC/E$, et
     \item $\kappa_{AM}(C,E)$ est la {\em classe
    d'Atiyah-Molino de la vari\'{e}t\'{e} feuillet\'{e}e $C$},
    un invariant de feuilletages
    qui est repr\'{e}sent\'{e}e
    par une $1$-forme le long des feuilles \`{a} valeurs dans $S^3(TC/E)$
    donn\'{e}e par une partie du tenseur de courbure d'une connexion
    symplectique adapt\'{e}e (th\'{e}or\`{e}me \ref{TMolinoPresymp}).
    \end{itemize}
   Je ne vois pas comment en g\'{e}n\'{e}ral le terme contenant la classe
   d'Atiyah-Molino dans la condition n\'{e}cessaire
   (\ref{EqObstructionImportanteIntro}) puisse \^{e}tre compens\'{e} par
   $\alpha_0$ et $[*]_1$.

  En paragraphe \ref{SecMorphRepAtiyahMolinoNulle} on \'{e}tudie d'abord
  le probl\`{e}me de la {\em quantification d'une application de Poisson $\phi$
  entre deux vari\'{e}t\'{e}s symplectiques $(M,\omega)\ra (M',\omega')$}: une telle
  application
  a localement la structure d'une projection d'un produit cart\'{e}sien
  sur un des facteurs, (voir l'\'{e}qn (\ref{EqProjApplPoisEinl})),
  mais globalement la situation
  \[
   \begin{array}{c}
    M'\times M'' \\
      \downarrow ~pr_1 \\
      M'
   \end{array} \mathrm{~est~remplac\acute{e}e~par~}
   \begin{array}{c}
      M \\
      \downarrow \phi \\
      M'
   \end{array} \mathrm{~~avec~~}TM=E\oplus Ker~T\phi
  \]
  o\`{u}
  $\phi:M\ra M'$ est une vari\'{e}t\'{e} fibr\'{e}e sur l'ouvert $\phi(M)$ de $M'$
  dont les fibres sont
  symplectiques et le fibr\'{e} orthogonal $E$ (par rapport \`{a} $\omega$) est
  int\'{e}grable (voir la proposition \ref{PMorphPoissonSymp} et
  \cite[p.76/77]{CHW98}).
  Alors $M$ est une vari\'{e}t\'{e} munie de {\em deux feuilletages
  symplectiques transverses} dont le feuilletage vertical n'a pas de classes
  caract\'{e}ristiques
  int\'{e}ressantes, tandis que le feuilletage horizontal peut a priori \^{e}tre
  nontrivial. On d\'{e}montre d'abord dans le th\'{e}or\`{e}me
  \ref{TFedQuantApplPois} en paragraphe \ref{SubSecQuantifMorphVarSymp}
  que le probl\`{e}me admet une {\em \'{e}tude
  fedosovienne}, qui est rappel\'{e}e dans le th\'{e}or\`{e}me
  \ref{TFedosovConstructionPartielle}.
  Grosso modo, il faut et il suffit d'arranger les constructions de Fedosov
  sur $M$ et
  $M'$ de mani\`{e}re \`{a} ce que {\em les deux d\'{e}riv\'{e}es de Fedosov soient
  `$\phi$-li\'{e}es'}. La proposition \ref{PFedosovEinfuehrungGamma}
  donne une condition suffisante et n\'{e}cessaire pour la quantification
  de $\phi$: on obtient toujours une d\'{e}formation de la multiplication
  ext\'{e}rieure des formes diff\'{e}rentielles le long du fibr\'{e} horizontal et une
  d\'{e}formation de la diff\'{e}rentielle
  $d_E$ le long de $E$ qui reste une d\'{e}rivation gradu\'{e}e, mais n'est pas de
  carr\'{e} nul. Si on calcule les deux premiers termes de cette condition
  on obtient exactement les m\^{e}mes obstructions que pour la
  repr\'{e}sentabilit\'{e} (\ref{EqObstructionDeuxIntro}) et
  (\ref{EqObstructionImportanteIntro}): ici c'est la {\em classe
  d'Atiyah-Molino du fibr\'{e} horizontal} qui entre dans la condition
  n\'{e}cessaire (th\'{e}or\`{e}me \ref{TQuantApplPoissCondNec}).
  D'un autre c\^{o}t\'{e}: {\em si la classe d'Atiyah-Molino
  du fibr\'{e} horizontal s'annule et les deux classes de Deligne sont
  `en quelque sorte $\phi$-li\'{e}es', il existe une quantification
  de $\phi$}, (th\'{e}or\`{e}me \ref{TQuantApplPoissonAtiyahMolinoNulle}).\\
  Le paragraphe \ref{SubSecFibEspRedLoc} sert \`{a} g\'{e}n\'{e}raliser l'astuce de
  Weinstein de plonger une sous-vari\'{e}t\'{e} co\"{\i}sotrope $C$, pour laquelle
  l'espace r\'{e}duit $M_{\mathrm{red}}$ existe, dans $M\times
  M_{\mathrm{red}}$ comme sous-vari\'{e}t\'{e} lagrangienne (utilis\'{e}e pour le
  th\'{e}or\`{e}me \ref{TRepUndAllesEspRedExist}): la projection
  \[
   \begin{array}{ccc}
      &                                 & M\times M_{\mathrm{red}}  \\
      & \stackrel{i\times\pi}{\nearrow} &                           \\
    C &                                 & \downarrow ~pr_1          \\
      & \stackrel{i}{\searrow}          &                           \\
      &                                 &        M
   \end{array}~~\mathrm{est~remplac\acute{e}e~par~~~}
   \begin{array}{ccc}
      &                                 & E^C                       \\
      & \stackrel{j}{\nearrow}          &                           \\
    C &                                 & \downarrow \phi          \\
      & \stackrel{i}{\searrow}          &                           \\
      &                                 &        M
   \end{array}
  \]
  o\`{u} $E^C$ est une vari\'{e}t\'{e} symplectique et $\phi$ est une application de
  Poisson (donc une submersion) {\em dont les fibres
  sont des espaces r\'{e}duits locaux} (th\'{e}or\`{e}me
  \ref{TVarieteEspacesReduitsLocaux}). De plus, $C$ se plonge de fa\c{c}on
  naturelle dans $E^C$
  en tant que sous-vari\'{e}t\'{e} lagrangienne. Finalement, si la classe
  d'Atiyah-Molino de la vari\'{e}t\'{e} feuillet\'{e}e $C$ s'annule, il en est de m\^{e}me
  avec la classe d'Atiyah-Molino du fibr\'{e} horizontal de $E^C$.
  \\
  En paragraphe \ref{SubSecQuantSousVarCoisClasAMNulle} on utilise la
  construction pr\'{e}c\'{e}dente de $E^C$ pour {\em contruire une repr\'{e}sentation
  d'un star-produit sur $M$ sur les fonctions sur $C$ si la classe
  d'Atiyah-Molino de $C$ s'annule} (th\'{e}or\`{e}me
  \ref{TQuantSousVarCoisClasAMNulle}): puisque la classe d'Atiyah-Molino
  du fibr\'{e} horizontal de $E^C$ s'annule, alors
  on peut quantifier l'application de Poisson $\phi$ d'apr\`{e}s le th\'{e}or\`{e}me
  \ref{TQuantApplPoissonAtiyahMolinoNulle}. Par cons\'{e}quent, on peut
  consid\'{e}rer  $\CinfK{M}[[\nu]]$ comme sous-alg\`{e}bre de
  $\CinfK{E^C}[[\nu]]$, et cette derni\`{e}re alg\`{e}bre admet une repr\'{e}sentation
  sur les fonctions sur la sous-vari\'{e}t\'{e} lagrangienne $C$ pourvu que
  les classes de Deligne soient bien choisies. Un cas particulier de ce
  th\'{e}or\`{e}me est donn\'{e} par $C$ \'{e}gale \`{a} l'image r\'{e}ciproque $J^{-1}(0)$
  d'une application moment associ\'{e}e \`{a} {\em l'action propre hamiltonienne}
  d'un groupe de Lie, voir le corollaire
  \ref{CActionsPropresClasseAMNulle}, r\'{e}sultat qui a d\'{e}j\`{a} \'{e}t\'{e} trait\'{e} dans
  \cite{BHW00}.

  Comme illustration, j'ai inclus deux exemples, celui du fibr\'{e} cotangent
  du $2$-tore en paragraphe \ref{SubSecFibCotDeuxTore} (qui provient de
  \cite{BHW00}) et de la r\'{e}duction symplectique
  $\complex^{n+1}\setminus \{0\}\leadsto \complex P(n)$
  en paragraphe
  \ref{SubSecEspProj} qui a \'{e}galement \'{e}t\'{e} trait\'{e} avant dans
  \cite{BBEW96a}, \cite{BBEW96b} et \cite{Wal98}: ici on construit des
  star-produits adapt\'{e}s explicites, et on montre qu'un star-produit
  adpat\'{e} projetable peut \^{e}tre \'{e}quivalent \`{a} un star-produit adapt\'{e} non
  projetable (pour le cas de $T^*T^2$) et que deux star-produits adapt\'{e}s
  projetables \'{e}quivalents peuvent avoir des alg\`{e}bres r\'{e}duites non
  isomorphes (pour le cas de l'espace projectif).
  \\

Ce rapport couvre une partie de mes activit\'{e}s de recherche
entre 2000 et 2003 dont j'ai donn\'{e} plusieurs expos\'{e}s: j'ai pr\'{e}sent\'{e}
la proposition \ref{PDefoMorphDefoBimod} le 23 novembre 2000 \`{a} Strasbourg.
Vivant pendant une certaine p\'{e}riode `{\em im dogmatischen Schlummer}'
qu'il n'y ait pas d'obstructions du tout \`{a} la repr\'{e}sentabilit\'{e} des
star-produits, j'ai d\'{e}couvert l'obstruction \`{a} l'ordre 3
(l'\'{e}qn (\ref{EqObstructionImportanteIntro})) contenant la classe
d'Atiyah-Molino
--mais sans le terme quadratique en $\alpha_0$-- en \'{e}t\'{e} 2001 en faisant une
\'{e}tude fedosovienne, et je l'ai expos\'{e}e la premi\`{e}re fois le 27 septembre 2001
lors d'un atelier organis\'{e} par l'Institut Max Planck des math\'{e}matiques
appliqu\'{e}es aux sciences \`{a} Leipzig. Pendant le
`Joint seminar Bruxelles-Warwick' \`{a} l'ULB Bruxelles, le 15 mars 2002,
j'ai pr\'{e}sent\'{e} les obstructions r\'{e}currentes (proposition
\ref{TObstructionsRecurrentes}) et le th\'{e}or\`{e}me de r\'{e}duction
\ref{TRepUndAllesEspRedExist} avec la relation
(\ref{EqRelationPrincipaleClassesDel}) des deux classes
de Deligne. Finalement, j'ai rapport\'{e} sur l'ensemble de tous les
r\'{e}sultats de ce rapport --sauf la d\'{e}monstration tr\`{e}s technique de la formule
(\ref{EqObstructionImportanteIntro}), le th\'{e}or\`{e}me
\ref{TIAlgId} et le th\'{e}or\`{e}me \ref{Tdaserstemalfalsch} dont la premi\`{e}re
version \'{e}tait fausse comme j'ai remarqu\'{e} apr\`{e}s une conversation avec
Nikolai Neumaier au mois de janvier 2004-- lors de l'atelier
`Quantisation of Poisson spaces with singularities' \`{a} Oberwolfach,
le 22 janvier 2003.

\subsubsection*{Remerciements}

Je voudrais remercier les personnes suivantes pour de nombreuses
discussions et propositions utiles: Didier Arnal, Pierre Bieliavsky,
Philippe Bonneau, Michel Cahen,
Alberto Cattaneo, Aziz El Kacimi, Giovanni Felder,
Gr\'{e}gory Ginot, Simone Gutt, Gilles Halbout, Hans-Christian Herbig,
Jens Hoppe, Johannes Huebschmann,
Franz Kamber, Pierre Lecomte, Sophie Lef\`{e}\-vre, Nikolai Neumaier,
Markus Pflaum,
Martin Schlichenmaier, Pierre Sleewaegen (pour m'avoir appris le th\'{e}or\`{e}me
de Weinstein-Gotay), Jerzy \'{S}niatyc\-ki, Daniel Sternheimer,
Thomas Strobl (pour m'avoir appris la cohomologie
relative), Vladimir Turaev (pour m'avoir incit\'{e} de r\'{e}fl\'{e}chir
sur des applications de Poisson entre vari\'{e}t\'{e}s symplectiques),
Stefan Waldmann, Alan Weinstein, Julius Wess et Eberhard Zeidler.
De plus, dans cette premi\`{e}re version du rapport, je prie le lecteur de me
pardonner
des omissions dans la bibliographie et je lui remercie d'avance
de me les signaler.

\subsubsection*{Notations}

Le symbole $\korps$ d\'{e}signera le corps de tous les nombres r\'{e}els, $\real$,
ou le
corps des tous les nombres complexes, $\complex$. Pour deux vari\'{e}t\'{e}s
diff\'{e}rentiables $M$ et $M'$, le symbole $\Cinf(M,M')$ d\'{e}signe l'ensemble
de toutes les applications de classe $\Cinf$ de $M$ dans $M'$. Pour un
fibr\'{e} vectoriel $E$ sur une vari\'{e}t\'{e} diff\'{e}rentiable $M$, on \'{e}crira
$\Ginf(M,E)$ pour l'espace de toutes les sections de classe $\Cinf$ du
fibr\'{e} $E$. On ne notera pas la diff\'{e}rence entre $\Ginf(M,E)$ et sa
complexification quand l'anneau $\CinfC{M}$ agit sur $\Ginf(M,E)$.
Le symbole $pr_1$ (resp. $pr_2$) repr\'{e}sentera toujours la projection
sur le premier facteur $M_1\times M_2\ra M_1$
(resp. sur le deuxi\`{e}me facteur $M_1\times M_2\ra M_2$) d'un produit
cart\'{e}sien de deux ensembles.

\newpage

\section{D\'{e}finitions pr\'{e}liminaires}
 \label{SecDefPrelim}

Dans ce paragraphe, on rappelle quelques notions autour de la g\'{e}om\'{e}trie de
Poisson qui sont plus ou moins connues, voir par exemple
\cite{Vai94} et \cite{CHW98}.

\subsection{Vari\'{e}t\'{e}s et applications de Poisson}
          \label{SubSecPois}

Soit $M$ une vari\'{e}t\'{e} diff\'{e}rentiable et $P$ un champs de bivecteurs
appartenant \`{a} $\Ginf(M,\Lambda^2 TM)\cong \Ginf\big(M,(\Lambda^2 T^*M)^*)$.
Le champ s'appelle {\em structure de Poisson} lorsque
le crochet de Poisson associ\'{e} \`{a} $P$ sur deux fonctions
\beq \label{EqCroPoi}
  \{f,g\}_P:=\{f,g\}:= P(df,dg) ~~\mathrm{quels~que~soient~}f,g\in\CinfK{M}
\end{equation}
satisfait l'identit\'{e} de Jacobi
\[
   \{\{f,g\},h\}+\{\{h,f\},g\}+\{\{g,h\},f\}=0
\]
quelles que soient $f,g,h\in\CinfK{M}$. Le couple $(M,P)$ est dite
{\em vari\'{e}t\'{e} de Poisson}. L'espace $\CinfK{M}$ muni de la multiplication
point-par-point et du crochet de
Poisson (\ref{EqCroPoi}) est une {\em alg\`{e}bre de Poisson}: ceci est
une alg\`{e}bre associative commutative unitaire sur $\korps$ munie d'une
structure d'alg\`{e}bre de Lie $\{~,~\}$ telle que les \'{e}quations de
compatibilit\'{e} $\{fg,h\}=f\{g,h\}+\{f,h\}g$ soient satisfaites quelles
que soient $f,g,h\in\CinfK{M}$.

Chaque structure de Poisson $P$ induit le morphisme de fibr\'{e}s vectoriels
suivant:
\beq \label{EPoissonSharp}
  P^\sharp:T^*M\ra TM:\alpha\mapsto P(~,\alpha).
\end{equation}
Le {\em champ hamiltonien} $X_f$ de $f\in\CinfK{M}$ est d\'{e}finie par
$X_f:=P^\sharp(df)$, et il vient $[X_f,X_g]=-X_{\{f,g\}}$ quelles
que soient $f,g\in\CinfK{M}$.

Soit $(M',P')$ une autre vari\'{e}t\'{e} de Poisson et $s,t\in\real$. Pour
$(m,m')\in M\times M'$ on note
$i_{(m,m')}:T_mM\ra T_mM\times T_{m'}M':v\mapsto (v,0)$ et
 $i'_{(m,m')}:T_{m'}M'\ra T_mM\times T_{m'}M':w\mapsto (0,w)$. On \'{e}crit
 $P_{(1)~(m,m')}:=(i_{(m,m')}\otimes i_{(m,m')})P_m$ et
 $P'_{(2)~(m,m')}:=(i'_{(m,m')}\otimes i'_{(m,m')})(P'_{m'})$. Il s'ensuit
que $sP_{(1)}+tP'_{(2)}$ est une structure de Poisson sur la vari\'{e}t\'{e}
 produit $M\times M'$.

On rappelle qu'une vari\'{e}t\'{e} symplectique $(M,\omega)$ est une vari\'{e}t\'{e}
diff\'{e}rentia\-ble
$M$ munie d'une $2$-forme $\omega$ ferm\'{e}e non d\'{e}g\'{e}n\'{e}r\'{e}e. A l'aide de
l'isomorphisme de fibr\'{e}s vectoriels $\omega^\flat:TM\ra T^*M:X\mapsto
\omega(X,~)$ on d\'{e}finit le champ de bivecteurs
$P:=\omega\big((\omega^\flat)^{-1},(\omega^\flat)^{-1}\big)$ qui est une
structure de Poisson.\\
En outre, pour une alg\`{e}bre de Lie r\'{e}elle $(\mathfrak{g},[~,~])$ de dimension
finie, son espace dual $\mathfrak{g}^*$ est une vari\'{e}t\'{e} de Poisson:
$P_\mathfrak{g}$ au point $\alpha\in\mathfrak{g}^*$ est d\'{e}finie par
\[
   P_{\mathfrak{g}~\alpha}:=\alpha([~,~].
\]

Soient $(M,P)$ et $(M',P')$ deux vari\'{e}t\'{e}s de Poisson. Une application
 $\phi:M\ra M'$ de classe $\Cinf$ s'appelle {\em application de Poisson}
 lorsque $P$ et $P'$ sont $\phi$-li\'{e}es, c.-\`{a}-d.
\beq \label{EqDefApplPois}
    \big(T_m\phi\otimes T_m\phi\big)(P_m)=P'_{\phi(m)}
    \mathrm{~~~quel~que~soit~}m\in M,
\end{equation}
d'o\`{u}
\beq \label{EqDefApplPoisAlt}
  P^{\prime\sharp}_{\phi(m)}=T_m\phi ~P_m^{\sharp}~( T_m\phi)^*
    \mathrm{~~~quel~que~soit~}m\in M.
\end{equation}
On voit ais\'{e}ment que l'application
\[
  \phi^*:\CinfK{M'}\ra \CinfK{M}:g\mapsto g\circ \phi
\]
est un homomorphisme d'alg\`{e}bres de Poisson, c.-\`{a}-d.:
\bea
     \phi^*(g_1g_2)     & = & (\phi^*g_1)(\phi^*g_2) \label{EqApPoiCom}\\
     \phi^*\{g_1,g_2\}' & = & \{\phi^*g_1,\phi^*g_2\} \label{EqApPoiLie}
\eea
quelles que soient $g_1,g_2\in \CinfK{M'}$. R\'{e}ciproquement, pour une
application $\korps$-lin\'{e}aire $\Xi:\CinfK{M'}\ra \CinfK{M}$
satisfaisant les deux \'{e}quations (\ref{EqApPoiCom}) et (\ref{EqApPoiLie})
(o\`{u} $\phi^*$ est remplac\'{e} par $\Xi$) il existe toujors une unique
application de Poisson $\phi$ telle que $\Xi=\phi^*$: d'apr\`{e}s l'exercice de
Milnor (voir \cite{KMS93}, p.~301, Cor.~35.9) il existe un unique
$\phi\in\Cinf(M,M')$ telle que $\Xi=\phi^*$. La deuxi\`{e}me \'{e}quation
(\ref{EqApPoiLie}) entra\^{\i}ne que $\phi$ est de Poisson puisque tout champ
de bivecteurs est d\'{e}termin\'{e} par ses valeurs sur $df\otimes dg$,
o\`{u} $f,g\in\CinfK{M'}$.\\
On note que le noyau de $\phi^*$ est un {\em id\'{e}al de Poisson} dans
$\CinfK{M'}$.

Dans le cas $(M',P')=(\mathfrak{g}^*,P_\mathfrak{g})$, une application de
Poisson $J:M\ra\mathfrak{g}^*$ est appel\'{e}e une {\em application moment}.
Une d\'{e}finition \'{e}quivalente est
\beq \label{Eqapplmoment}
     \{ \langle J,\xi\rangle, \langle J,\eta\rangle \} =
                 \langle J, [\xi,\eta] \rangle
                 \mathrm{~quels~que~soient~}\xi,\eta\in\mathfrak{g}.
\end{equation}
 La d\'{e}finition classique d'une application moment par J.-M. Souriau
commence par une action gauche d'un groupe de Lie $G$ sur $M$,
$G\times M\ra M:(g,m)\mapsto gm=:\Phi_g(m)$ telle
que ($i$) l'alg\`{e}bre de Lie de $G$ soit \'{e}gale \`{a} $\mathfrak{g}$, ($ii$)
l'action
pr\'{e}serve la structure de Poisson, c.-\`{a}.-d. toutes les $\Phi_g$ sont
des application de Poisson, ($iii$) il existe une application de classe
$\Cinf$
$J:M\ra\mathfrak{g}^*$ satisfaisant (\ref{Eqapplmoment}), ($iv$)
les champs hamiltoniens $X_{\langle J,\xi\rangle}$ co\"{\i}ncident avec
les g\'{e}n\'{e}rateurs infinit\'{e}simaux $\xi_M(m):=d/dt\big(exp(t\xi)m\big)|_{t=0}$
et ($v$) $J$ est $G$-\'{e}quivariante: $J(gm)=Ad^*(g)\big(m\big)~~\forall~g\in G$.
Ici, les propri\'{e}t\'{e}s ($iv$) et ($v$) impliquent la propri\'{e}t\'{e} ($iii$).\\
Si $M$ est une vari\'{e}t\'{e} symplectique de dimension $2n$, $\mathfrak{g}$
une alg\`{e}bre de Lie r\'{e}elle ab\'{e}lienne de dimension $n$ et si l'ensemble de
tous les points singuliers de $J$ est negligeable (par rapport \`{a} la mesure
de Liouville $\omega^{\wedge n}$) l'application moment est dite un
{\em syst\`{e}me hamiltonien compl\`{e}tement int\'{e}grable} (au sens de Liouville).

La description des applications de Poisson se simplifie beaucoup dans le
cas symplectique (voir aussi \cite[p.76/77]{CHW98}):

\bprop\label{PMorphPoissonSymp}
 Soient $(M,\omega)$ et $(M',\omega')$ deux vari\'{e}t\'{e}s symplectiques et
 $\phi:M\ra M'$ une application de Poisson.
 \ben
  \item $\phi$ est une submersion et son image est une sous-vari\'{e}t\'{e}
    ouverte de $M'$.
  \item Le sous-fibr\'{e} $F:=\mathrm{Ker}~T\phi$ de $TM$ est symplectique
    (et int\'{e}grable) et son compl\'{e}ment symplectique $E:=F^\omega$ est
    symplectique et int\'{e}grable.
  \item Soit $f':M'\ra\real$ une fonction de classe $\Cinf$. Alors
    le champs hamiltonian $X_{f'\circ\phi}$ a ses valeurs dans
    $E$, et $X_{f'\circ\phi}$ et le champs hamiltonien $X'_{f'}$ sont
    $\phi$-li\'{e}s, c.-\`{a}-d. $T\phi~X_{f'\circ\phi}=X'_{f'}\circ \phi$.
  \item La restriction de $\omega$ \`{a} $E$ est \'{e}gale \`{a} la restriction
      de $\phi^*\omega'$ \`{a} $E$.
 \een
\eprop
\bbew
  Soient $P$ et $P'$ les structures de Poisson associ\'{e}es aux formes
  symplectiques $\omega$ et $\omega'$.\\
 1. Soit $m'\in M'$ tel qu'il existe $m\in M$ avec $\phi(m)=m'$ et soit
  $v\in T_{m'}M'$. Puisque $P'$ est nond\'{e}g\'{e}n\'{e}r\'{e}e il existe une
  $1$-forme $\alpha\in {T_{m'}M'}^*$ telle que
  $P^{\prime\sharp}_{m'}\alpha =v$. En utilisant l'equation
  (\ref{EqDefApplPoisAlt}) on voit que
  $T_m\phi~P_m^{\sharp} ~(T_m\phi)^*\alpha=v$
  ce qui
  montre que $T\phi$ est surjective, donc $\phi$ est une submersion.
   Il est bien connu que chaque
  submersion est une application ouverte: en particulier, l'image de
  $\phi$ est un ouvert (et donc automatiquement une sous-vari\'{e}t\'{e}) de
  $M'$.\\
 2. Puisque $\phi$ est une submersion, le noyau $F_m:=\mathrm{Ker}~T_m\phi$ est
 de dimension constante quel que soit $m\in M$ et d\'{e}finit donc un sous-fibr\'{e}
 de $TM$ qui est
 automatiquement int\'{e}grable, les sous-vari\'{e}t\'{e}s int\'{e}grales \'{e}tant les
 composantes
 connexes des images r\'{e}ciproques
 $\phi^{-1}(m')$ pour tout $m'$ dans l'image de $\phi$. Consid\'{e}rons
 pour chaque $m\in M$ le sous-espace
 \[
  E_m:=\{P_m^{\sharp}~(T_m\phi)^*\alpha~|~m\in M,\alpha\in
    {T_{\phi(m)}M'}^*\}.
 \]
  D'apr\`{e}s le premier \'{e}nonc\'{e}, la restriction de $T_m\phi$
  \`{a} $E_m$ est une surjection sur l'espace tangent $T_{\phi(m)}M'$. D'un
  autre c\^{o}t\'{e}, puisque $E_m$ est param\'{e}tr\'{e} par ${T_{\phi(m)}M'}^*$ il s'ensuit
  que $E_m$ et $T_{\phi(m)}M'$ ont la m\^{e}me dimension, et $T_m\phi$
  restreinte \`{a} $E_m$ est un isomorphisme. Par cons\'{e}quent, l'intersection
  de $E_m$ et $F_m$ s'annule et $T_mM=E_m\oplus F_m$ quel que soit $m\in
  M$. Il suffit de montrer que $E_m$ et $F_m$ sont orthogonaux par rapport
  \`{a} $\omega$ pour que $E_m$ et $F_m$ soient de sous-espaces symplectiques
  de $T_mM$: soient $w\in F_m$ et $\alpha\in T_{\phi(m)}M^*$. Alors
  \[
      \omega_m\big(P_m^{\sharp}~(T_m\phi)^*\alpha,~w\big)
      =\omega_m^\flat\big(P_m^{\sharp}~(T_m\phi)^*\alpha\big)(w)
      =\big((T_m\phi)^*\alpha\big)(w)
         = \alpha\big(T_m\phi~w\big) = 0
  \]
  puisque $T_m\phi~w=0$. Le fait que le sous-fibr\'{e} $E$ soit int\'{e}grable
  sera d\'{e}montr\'{e} apr\`{e}s la d\'{e}monstration du prochain \'{e}nonc\'{e}.\\
  3. Soit $V$ un champ de vecteurs \`{a} valeurs dans $F$. Alors:
  \[
     \omega(X_{f'\circ\phi},V)  =  d(f'\circ\phi)\big(V\big)
                               =  df'\circ\phi~(T\phi~V) = 0,
  \]
  et il s'ensuit que $X_{f'\circ\phi}$ est \`{a} valeurs dans $E$. Soit
  $g':M'\ra\real$ une autre fonction de classe $\Cinf$. Alors
  \beas
    dg'\circ\phi~\big(T\phi~X_{f'\circ\phi}\big) & = &
            d(g'\circ\phi)\big(X_{f'\circ\phi}\big)  \\
            & = & \omega\big(X_{g'\circ\phi},X_{f'\circ\phi}\big) \\
            & = & \{g'\circ\phi,f'\circ\phi\}
                  \mathrm{~~par~d\acute{e}finition~du~crochet
                                                       ~de~Poisson}\\
            & = & \{g',f'\}'\circ\phi \mathrm{~~parce~que~}\phi
                                    \mathrm{~est~un~morphisme~de~Poisson}
                                    \\
            & = & \omega'(X'_{g'},X'_{f'})\circ\phi \\
            & = & dg'\circ\phi~\big(X'_{f'}\circ\phi\big)
  \eeas
  ce qui montre que les champs hamiltoniens $X_{f'\circ\phi}$ et $X'_{f'}$
  sont $\phi$-li\'{e}s. Pour l'int\'{e}gra\-bi\-li\-t\'{e} du fibr\'{e} $E$ il suffit de faire
  un calcul local: gr\^{a}ce \`{a} ce qu'on vient de d\'{e}montrer, les champs
  hamiltoniens $X_{f'\circ\phi}$ forment une base locale des sections de
  classe $\Cinf$ de $E$. Alors, puisque
  \[
      [X_{f'\circ\phi},X_{g'\circ\phi}] = -X_{\{f'\circ\phi,g'\circ\phi\}}
                                      = -X_{\{f',g'\}'\circ\phi}
  \]
  on voit l'int\'{e}grabilit\'{e} du fibr\'{e} $E$.\\
  4. Il suffit de v\'{e}rifier l'\'{e}quation \'{e}nonc\'{e}e sur deux champs hamiltoniens
   $X_{f'\circ\phi}$ et $X_{g'\circ\phi}$. Alors:
   \beas \lefteqn{
     \omega\big(X_{f'\circ\phi},X_{g'\circ\phi}\big)
     =\{f'\circ\phi,g'\circ\phi\}
       =\{f',g'\}'\circ\phi=
       \omega'\circ\phi\big(X'_{f'}\circ\phi,X'_{g'}\circ\phi\big)}
             ~~~~~~~~~~~~~    \\
     & & =\omega'\circ\phi\big(T\phi~X_{f'\circ\phi},
                ~T\phi~X_{g'\circ\phi}\big)
       =\phi^*\omega'\big(X_{f'\circ\phi},X_{g'\circ\phi}\big)
   \eeas
   ce qui montre le dernier \'{e}nonc\'{e}.
\ebew

\subsection{Sous-vari\'{e}t\'{e}s}
   \label{SubSecSousvar}

Soit $M$ une vari\'{e}t\'{e} diff\'{e}rentiable de dimension $n$ et $i:C\ra M$
une sous-vari\'{e}t\'{e} ferm\'{e}e de $M$ de dimension $n-k$.
Soit
\beq \label{EqIdAn}
  \mathcal{I}_C:=\mathcal{I}:=\{f\in \CinfK{M}~|~f(c)=0~\forall c\in C\}
\end{equation}
{\em l'id\'{e}al annulateur de} $C$. Soient $TC$ le fibr\'{e} tangent de $C$, $TM|_C$
 (resp. $T^*M|_C$)
 la restriction du fibr\'{e} tangent (resp. du fibr\'{e} cotangent) de $M$ \`{a} $C$ et
 $TC^{\mathrm{ann}}$ le sous-fibr\'{e} de $T^*M|_C$ consistant en toutes les
 $1$-formes in $T^*M|_C$ qui s'annulent sur $TC$.

Par la suite, on aura plusieurs fois besoin du lemme suivant --qui est
sans doute bien connu, voir par exemple \cite[p.127]{BJ73}--
\`{a} l'aide duquel on peut toujours trouver des
voisinages ouverts {\em convexes fibre-par-fibre} de la section nulle d'un
fibr\'{e} vectoriel dans des voisinages
ouverts d\'{e}j\`{a} donn\'{e}s:
\blem \label{LOuvertBeauSectionNulle}
 Soit $N$ une vari\'{e}t\'{e} diff\'{e}rentiable et $\tau_i:E_i\ra N$, $1\leq i\leq
 l$, des fibr\'{e}s vectoriels r\'{e}els sur $N$. Pour tout $i$ soit
 $g^{(i)}\in\Ginf(N,S^2E_i^*)$ une m\'{e}trique d\'{e}finie positive.
 En outre, soit $\tau:K\ra N$ le fibr\'{e}
 vectoriel $K=E_1\oplus\cdots\oplus E_l$ sur $N$ et soit $W$ un
 ouvert de $K$ avec $N\subset W\subset K$.
 Alors il existe $l$ fonctions $f^{(1)},\ldots,f^{(l)}\in\CinfR{N}$
 \`{a} valeurs strictement positives telles que
 pour l'ouvert
 \beas
    \lefteqn{\tilde{W}:=
       \{e=e_1+\ldots+e_l\in K~|~ } \\
  & &        ~~~~~~~~~~  g^{(1)}(e_1,e_1)<f^{(1)}\big(\tau(e)\big),\ldots,
                        g^{(l)}(e_l,e_l)<f^{(l)}\big(\tau(e)\big)\}
 \eeas
 il vient $N\subset \tilde{W}\subset W$. En outre, $\tilde{W}$ est
 convexe fibre-par-fibre, c.-\`{a}-d. que chaque intersection
 $\tilde{W}\cap K_n$ est une partie convexe de la fibre $K_n$ sur $n\in
 N$.
\elem
\bbew
 Soit $(U_\alpha)_{\alpha\in\mathfrak{S}}$ une famille de domaines ouverts
 de cartes
 sur $N$ telle que pour chaque entier positif $1\leq i\leq l$ la famille
 $(U_\alpha,\Psi^{(i)}_\alpha)_{\alpha\in\mathfrak{S}}$ d\'{e}finisse une
 famille de trivialisations locales du fibr\'{e} $E_i$, c.-\`{a}-d. chaque
 $\Psi_\alpha^{(i)}$ est un diff\'{e}omor\-phis\-me de l'ouvert
 $\tau_i^{-1}(U_\alpha)$ sur $U_\alpha\times \real^{n_i}$ avec
 $pr_1\circ \Psi_\alpha^{(i)}=\tau_i$. Evidemment, la trivialisation locale
 $\Psi_\alpha:\tau^{-1}(U_\alpha)\ra
 U_\alpha\times \real^{n_1}\times\cdots\times \real^{n_l}$ du fibr\'{e} $K$
 est form\'{e} de tous les $\Psi_\alpha^{(i)}$. Puisque $N$ est localement
 compacte, on peut supposer que tout $U_\alpha$ est \`{a} adh\'{e}rence
 compacte. Gr\^{a}ce \`{a} la paracompacit\'{e} de $N$ il existe un r\'{e}tr\'{e}cissement
 $(V_\alpha)_{\alpha\in\mathfrak{S}}$ de la famille
 $(U_\alpha)_{\alpha\in\mathfrak{S}}$ (c.-\`{a}-d. on a $U_\alpha\supset
 \overline{V_\alpha}\supset V_\alpha$ pour les ouverts $V_\alpha$).\\
 On d\'{e}montre d'abord le cas $l=1$:
 L'ouvert $W\cap \tau^{-1}(U_\alpha)$ contient une r\'{e}union d'ouverts du
 type $\Psi_\alpha^{-1}\big(U_{\alpha\gamma}\times
 B^{(1)}_{r_{1\gamma}}(0)\big)$ o\`{u}
 tout $U_{\alpha\gamma}$ est un ouvert de $N$ appartenant \`{a} $U_\alpha$ tel
 que $\cup_{\gamma} U_{\alpha\gamma}=U_\alpha$ et
 $B^{(1)}_{r_{1\gamma}}(0)$ d\'{e}signe la boule ouverte de centre $0$ et de
 rayon $r_{1\gamma}>0$ dans $\real^{n_1}$ par rapport au produit scalaire
 canonique $\langle~,~\rangle$ dans $\real^{n_1}$.
 Puisque les $U_{\alpha\gamma}$ recouvrent la partie compacte
 $\overline{V_\alpha}$ on peut choisir un nombre fini des $U_{\alpha\gamma}$
 qui recouvrent $\overline{V_\alpha}$, et si $r_\alpha>0$ est le rayon
 minimal
 de tous les $r_{1\gamma}$ apparaissant dans ce choix, il vient que
 l'ouvert $W\cap \tau^{-1}(V_\alpha)$ contient un ouvert
 $\tilde{W}''_\alpha:=
 \Psi_{\alpha}^{-1}\big(V_\alpha\times B^{(1)}_{r_\alpha}(0)\big)$.
 Puisque la partie $\Psi_{\alpha}^{-1}\big(\overline{V_\alpha}\times
 \overline{B^{(1)}}_1(0)\big)$ est compacte l'application $e\mapsto
 g^{(1)}(e,e)$ y est minor\'{e}e par un nombre r\'{e}el $s_\alpha^2>0$, donc
 $g^{(1)}(e,e)\geq
 s_\alpha^2\langle pr_2\Psi_{\alpha}(e),pr_2\Psi_{\alpha}(e)\rangle$.
 Il s'ensuit que si $e\in\tau^{-1}(V_\alpha)$ tel que
 $g^{(1)}(e,e)<s_\alpha^2r_\alpha^2$ alors $e\in \tilde{W}''_\alpha$.
 Soit $\tilde{W}'_\alpha$ l'ouvert $\{e\in\tau^{-1}(V_\alpha)~|~
 g^{(1)}(e,e)<s_\alpha^2r_\alpha^2\}$. Par cons\'{e}quent,
 $V_\alpha\subset \tilde{W}'_\alpha\subset \tilde{W}''_\alpha\subset W\cap
 \tau^{-1}(V_\alpha)$.
 Evidemment, l'ouvert
 $\tilde{W}':=\cup_{\alpha\in\mathfrak{S}}\tilde{W}'_\alpha$ contient $N$
 et fait partie de $W$.  Soit $(f_\alpha)_{\alpha\in\mathfrak{S}}$ une
 partition de
 l'unit\'{e} subordonn\'{e}e \`{a} la famille $(V_\alpha)_{\alpha\in\mathfrak{S}}$.
 On d\'{e}finit
 $f^{(1)}:=\sum_{\alpha\in\mathfrak{S}}f_\alpha s_\alpha^2r_\alpha^2$.
 Evidemment,
 $f^{(1)}$ est de classe $\Cinf$ et \`{a} valeurs strictement
 positives. Soit
 $e=e_1\in K$ et soit $\hat{g}^{(1)}(e,e)<f^{(1)}\big(\tau(e)\big)$
 Puisque la famille des supports des $f_\alpha$
 est localement finie il existe un entier positif $n_0$ tel que
 $f_\alpha\big(\tau(e)\big)=0$ quel que soit $\alpha\in\mathfrak{S}\setminus
 \{\alpha_1,\ldots,\alpha_{n_0}\}$. Alors
 \[
    0<\sum_{j=1}^{n_0}f_{\alpha_j}\big(\tau(e)\big)\big(
            s_{\alpha_j}^2r_{\alpha_j}^2-g^{(1)}(e,e)\big),
 \]
 et il existe donc au moins un $\alpha_j$ tel que
 $g^{(1)}(e,e)<s_{\alpha_j}^2r_{\alpha_j}^2$. Il s'ensuit que
 $e\in \tilde{W}'_{\alpha_j}$, alors
 \[
    N=\cup_{\alpha\in\mathfrak{S}}V_\alpha \subset
         \tilde{W}:=\{e\in K~|~\hat{g}^{(1)}(e,e)<f^{(1)}\big(\tau(e)\big)\}
         \subset \cup_{\alpha\in\mathfrak{S}}\tilde{W}'_\alpha \subset W,
 \]
 et le cas $l=1$ est montr\'{e}.\\
 Pour le cas g\'{e}n\'{e}ral, on consid\`{e}re d'abord le fibr\'{e} vectoriel $K$ muni de
 la m\'{e}trique $g:(e=e_1+\cdots+e_l,e'=e'_1+\cdots+e'_l)
 \mapsto g^{(1)}(e_1,e'_1)+\cdots+g^{(l)}(e_l,e'_l)$ qui est \'{e}videmment
 d\'{e}finie positive. D'apr\`{e}s
 le cas $l=1$ il existe une fonction $f\in\CinfR{N}$ \`{a} valeurs strictement
 positives tel que l'ouvert $\{e\in K~|~g(e,e)<f\big(\tau(e)\big)\}\subset
 W$. Soit $f^{(i)}:=\frac{1}{l}f$ quel que soit $1\leq i\leq l$. Evidemment
 \beas
    N\subset \tilde{W}
        & := &
        \{e\in K~|~\hat{g}^{(1)}(e_1,e_1)<f^{(1)}\big(\tau(e)\big),\ldots,
                     \hat{g}^{(l)}(e_1,e_1)<f^{(l)}\big(\tau(e)\big)\} \\
        & \subset & \{e\in K~|~g(e,e)<f\big(\tau(e)\big)\}
                     \subset W.
 \eeas
\ebew

 On rappelle le th\'{e}or\`{e}me d'un {\em voisinage tubulaire} de $C$ dans M:
 il existe un sous-fibr\'{e} $E$ de
 $TM|_C$ tel que $TM|_C=TC\oplus E$, un voisinage
 ouvert $V$ de la section nulle de $E$ qu'on peut supposer convexe
 fibre-par-fibre d'apr\`{e}s le lemme \ref{LOuvertBeauSectionNulle},
 un voisinage ouvert $U$ de
 $C$ dans $M$ et un
 diff\'{e}omorphisme $V\mapsto U$ dont la restriction \`{a} la section nulle soit
 l'application identique sur $C$ (voir \cite[p.108-110]{Lan95},
 \cite{BHW00}). On note $\tau:U\ra C$ la submersion surjective induite par la
 projection du fibr\'{e} $E$. En outre, pour $u,u'\in U$ et $t,t'\in\real$ on
 continue \`{a} utiliser la notation $tu+t'u'$ pour la combinaison lin\'{e}aire
 induite de celle dans $E$.
 Soient $U':=M\setminus U$ et $(\psi_U, \psi_{U'})$ une
partition de l'unit\'{e} subordonn\'{e}e au recouvrement ouvert $(U,U')$ de $M$.

 Soit $\xi\in \Ginf(U,TU)$ le champs de vecteurs
 induit par le champ d'Euler dans $TE$ (dont le flot est la multiplication
 avec $\exp(t)$). Soit $\tilde{E}$ le sous-fibr\'{e} vertical de $TU$, c.-\`{a}-d.
 $\tilde{E}:=\mathrm{Ker} T\tau$. Le {\em complexe de Koszul}
 $\big(K:=\bigoplus_{k\in \entier}K_k,\partial\big)$ est d\'{e}fini
 par $K_k:=0~\forall k\leq -1$,
 $K_0:=\CinfK{U}=\Ginf(U,\Lambda^0 \tilde{E}^*)$
 et $K_k:=\Ginf(U,\Lambda^k \tilde{E}^*)$ o\`{u} $\partial_k:K_k\ra K_{k-1}:
 F\mapsto i(\xi)F$. La suite
 \beq \label{EqKos}
   \CinfK{C}\stackrel{i^*}{\longleftarrow}
     \CinfK{U} \stackrel{\partial_1}{\longleftarrow}
      \Ginf(U,\Lambda^1 \tilde{E}^*)
       \stackrel{\partial_2}{\longleftarrow}\cdots
       \stackrel{\partial_k}{\longleftarrow}
       \Ginf(U,\Lambda^k \tilde{E}^*)
        \stackrel{\partial_{k+1}}{\longleftarrow}\cdots
 \end{equation}
 est exacte. Ceci se voit de la fa\c{c}on suivante: pour $\phi\in
 \CinfK{C}$ on d\'{e}finit la prolongation $\prol_U(\phi):=\tau^*\phi$.
 Soit $d_v$ la diff\'{e}rentielle de Cartan verticale (le long des fibres
 de la projection $\tau$) sur $K$. Soit $\Phi_t(u):=tu$ pour $t\in\real$
 suffisamment petit.
 Pour
 $F\in K_k,k\geq 0$ on d\'{e}finit l'homotopie $h_k:K_k\ra K_{k+1}$
 par
 \beq \label{EqHomtopKoszul}
      h_k(F)_u:=\int_{0}^1 \big(\Phi_t^*(d_v F)\big)_{u}dt.
 \end{equation}
 Comme dans la d\'{e}monstration du lemme de Poincar\'{e} (voir par exemple
 \cite[p.67]{KMS93}) on d\'{e}duit sans peine
 que $i^*\prol_U$ est \'{e}gal l'application identique sur $\CinfK{C}$, que
 $\prol i^* + \partial_1 h_0=\id_0$ et $h_{k-1}
 \partial_{k}+\partial_{k+1}h_k=\id_{k}$ o\`{u} $\id_k$ d\'{e}signe l'application
 identique sur $K_k,~k\geq 0$.

Soit
\[
   \mathcal{I}^2:=\{g_1g'_1+\cdots+g_lg'_l~|~l\in\nat,l\geq 1,
    g_1,\ldots,g_l,g'_1,\ldots,g'_l\in\mathcal{I}\}.
\]
Le lemme suivant semble bien connu:
\blem \label{LemIIcarre}
 \ben
  \item La restriction $i^*:\CinfC{M}\ra \CinfK{C}:f\mapsto f\circ i$ est
    un homomorphisme surjectif d'alg\`{e}bres associatives commutatives \`{a}
    noyau $\mathcal{I}$.
  \item $\mathcal{I}$ est un $\CinfK{M}$-module d'un nombre fini de
         g\'{e}n\'{e}rateurs.
  \item Le $\CinfK{C}$-module $\mathcal{I}/\mathcal{I}^2$ (le
    {\em module conormal}) est isomorphe \`{a} l'espace de sections
    $\Ginf(C,TC^{\mathrm{ann}})$. L'isomorphisme est induit par
    l'application $g\mapsto (c\mapsto dg_c)$.
  \item Le $\CinfK{C}$-module
  $\mathrm{Hom}_{\Cinf(C,\korps)}\big(
     \mathcal{I}/\mathcal{I}^2,\CinfK{C}\big)$
    est isomorphe \`{a} l'espace $\Ginf(C,TM|_C/TC)$.
 \een
\elem
\bbew
 1. Il reste \`{a} montrer que la restriction soit surjective: soit $\phi\in
 \CinfK{C}$. On d\'{e}finit une prolongation globale
 $\prol:\CinfK{C}\ra\CinfK{M}$ par $\prol(\phi):=\psi_U(\prol_U(\phi))$
 et on a visiblement $\prol(\phi)(c)=\phi(c)$ quel que soit $c\in C$.\\
 2. Soit $g\in\mathcal{I}$. Alors $g=\psi_U^2(g|_U)+(1-\psi_U^2)g$.
 D'apr\`{e}s le th\'{e}or\`{e}me de Serre et Swan (voir par exemple \cite[p.155]{BJ73}),
 $K_1$ est un $\CinfK{U}$-module
 projectif d'un nombre fini de g\'{e}n\'{e}rateurs $F_1,\ldots,F_N$. De plus,
 $\mathcal{I}^U:=\{g\in\CinfK{U}~|~g(c)=0~\forall c\in C \}$ est \'{e}gal
 \`{a} l'image de l'op\'{e}rateur $\partial_1$ dans la suite exacte (\ref{EqKos}).
 Alors il existent $g_1,\ldots,g_N\in\CinfK{U}$ telles que
 \[
   g=\sum_{j=1}^N~\psi_U~\partial_1(F_j)~\psi_U~g_j~~+~~(1-\psi_U^2)g
 \]
 Puisque $\psi_U~\partial_1(F_1),\ldots,\psi_U~\partial_1(F_N),
 1-\psi_U^2$ appartiennent \`{a} $\mathcal{I}$ et $\psi_Ug_1,\ldots,
 \psi_Ug_N$ appartiennent \`{a} $\CinfK{M}$, l'id\'{e}al annulateur est engendr\'{e}
 par au plus $N+1$ g\'{e}n\'{e}rateurs.\\
 3. Il est clair que $dg_c\in T_cC^{\mathrm{ann}}$ et $d(g_1g_2)_c=0$ quels
 que
 soient $g,g_1,g_2\in \mathcal{I}$. Soit $F\in
 \Ginf(C,TC^{\mathrm{ann}})$. Le fibr\'{e} $TC^{\mathrm{ann}}$ est isomorphe au
 fibr\'{e} $E^*$. Soit $\tilde{F}\in \Ginf(C,\Lambda^1 \tilde{E}^*)$ une
 prolongation verticale de $F$ (par exemple d\'{e}finie pour
 $u,e\in\tau^{-1}(c)\cong E_c$ et $c\in C$ de mani\`{e}re
 $\tilde{F}_u\big(\frac{d}{ds}(u+se)|_{s=0}\big):=F_c(e)$). On d\'{e}finit
 $g:=\psi_U i_\xi(\tilde{F})$. Visiblement $g\in\mathcal{I}$. De plus,
 $\xi_c=0$, $\psi_U(c)=1$ et $\exp(t)c=c$, alors
 \[
   dg_c(e):=d(\psi_U i_\xi(\tilde{F}))_c(e)
          = (L_\xi\tilde{F})_c(e)
          = \frac{d}{dt}F_{c}(\exp(t)e)|_{t=0}
          = F_c(e).
 \]
 Finalement, soit $g\in\mathcal{I}$ tel que $dg_c=0$ quel que soit $c\in
 C$. La restriction $g|_U$ appartient \`{a} $\mathcal{I}^U$, donc, d'apr\`{e}s la
 suite exacte (\ref{EqKos}) il existe $B\in K_1$ tel que $g|_U=i_\xi(B)$.
 Ce $B$ au point $c\in C$ se calcule par $B_c=dg_c$ d'apr\`{e}s ce qui
 pr\'{e}c\`{e}de. Alors $B(X)$ appartient visiblement \`{a} $\mathcal{I}^U$ quel que
 soit le champ de vecteurs $X\in\Ginf(U,TU)$. D'apr\`{e}s (\ref{EqKos}) il
 existe un \'{e}l\'{e}ment $G(X,~)\in K_1$ tel que $B\big(X\big)=G(X,\xi)$. Puisque
 $G$ est une application $\CinfK{U}$-bilin\'{e}aire de $\Ginf(U,TU)\times
 \Ginf(U,\tilde{E})$
 dans $\mathcal{I}^U$ il s'ensuit que $G$ est une section de classe
 $\Cinf$ du fibr\'{e} $TU^*\otimes \tilde{E}^*$. Visiblement
 $g|_U=B(\xi)=G(\xi,\xi)$. Avec les g\'{e}n\'{e}rateurs $F_1,\ldots,F_N$
 du module projectif $K_1$ mentionn\'{e}s ci-dessus on peut r\'{e}crire
 $G(\xi,\xi)$ par une combinaison lin\'{e}aire finie
 $\sum_{i,j}g_{ij}F_i(\xi)F_j(\xi)$ o\`{u} $g_{ij}\in\CinfK{U}$. Il vient
 \[
   g=\sum_{i,j}\psi_Ug_{ij}~\psi_UF_i(\xi)~\psi_UF_j(\xi)~~+~~(1-\psi_U^3)g
 \]
 et $g$ est un \'{e}l\'{e}ment de $\mathcal{I}^2$.\\
 4. Evident.
\ebew

\subsection{Sous-vari\'{e}tes et applications co\"{\i}sotropes}
          \label{SubSecCoiso}

Dans le paragraphe \ref{SubSecPois} on a vu que le noyau de $\phi^*$ pour
une application de Poisson $\phi:M\ra M'$ \'{e}tait un id\'{e}al de Poisson de
$\CinfK{M'}$. Les applications co\"{\i}sotropes qu'on va pr\'{e}senter ici
seront des g\'{e}n\'{e}ralisations des applications de Poisson quand on laisse
tomber la structure de Poisson sur la vari\'{e}t\'{e} $M$:
\bdefi \label{DefAplCoi}
 Soient $(M,P)$ une vari\'{e}t\'{e} de Poisson et $C$ une vari\'{e}t\'{e} quelconque
 et $\phi:C\ra M$ une application de classe $\Cinf$.
 \ben
  \item $\phi$ s'appelle {\em co\"{\i}sotrope} ssi
  \[
   \mathcal{I}:=\{g\in \CinfK{M}~|~\phi^*g=0\}=\mathrm{Ker}\phi^*
  \]
  est une sous-alg\`{e}bre de Poisson de $\CinfK{M}$.
  \item En particulier, si $\phi$ est l'injection canonique d'une
  sous-vari\'{e}t\'{e} ferm\'{e}e $C$ de $M$, alors $C$ s'appelle {\em sous-vari\'{e}t\'{e}
  co\"{\i}sotrope} quand $\phi$ est co\"{\i}sotrope.
 \een
\edefi
A l'aide du lemme \ref{LemIIcarre} on montre la caract\'{e}risation g\'{e}om\'{e}trique
des sous-vari\'{e}t\'{e}s co\"{\i}sotropes (voir \cite{Wei88}):
\blem
 Une sous-vari\'{e}t\'{e} ferm\'{e}e $C$ d'une vari\'{e}t\'{e} de Poisson $(M,P)$ est
 co\"{\i}sotrope si et seulement si
 \beq \label{EqCoisoGeo}
   P_c(F_c,G_c)=0~~\mathrm{quels~que~soient~~}c\in C;
           ~F_c,G_c\in T_cC^{\mathrm{ann}}
 \end{equation}
\elem
\bbew
 Soit $i$ une application co\"{\i}sotrope d'apr\`{e}s la d\'{e}finition
 \ref{DefAplCoi}. Soit $c\in C$ et $F_c,G_c\in T_cC^{\mathrm{ann}}$. D'apr\`{e}s
 le troisi\`{e}me \'{e}nonc\'{e} du lemme \ref{LemIIcarre} il
 existent $g,g'$ dans l'id\'{e}al annulateur $\mathcal{I}$ de $C$ tels que
 $dg_c=F_c$ et $dg'_c=G_c$. Alors $P_c(F_c,G_c)=\{g,g'\}(c)=0$ parce
 que le crochet de Poisson $\{g,g'\}$ appartient \`{a} $\mathcal{I}$.\\
 R\'{e}ciproquement, soit la condition (\ref{EqCoisoGeo}) satisfaite, et
 soient $g,g'\in\mathcal{I}$. Alors pour tout $c\in C$ le crochet de
 Poisson $\{g,g'\}_c=P_c(dg_c,dg'_c)$ s'annule car
 $dg_c,dg'_c\in T_cC^{\mathrm{ann}}$. Alors $\{g,g'\}\in\mathcal{I}$.
\ebew

Pour une vari\'{e}t\'{e} symplectique $(M,\omega)$ et un sous-espace $E$ d'un
espace tangent $T_mM$ il y a la notion du sous-espace $\omega$-orthogonal
\beq
  E^\omega:=\{w\in T_mM~|~\omega_m(v,w)=0~\forall v\in E\},
\end{equation}
et l'on en d\'{e}duit ais\'{e}ment qu'une sous-vari\'{e}t\'{e} ferm\'{e}e $C$ d'une
vari\'{e}t\'{e}
symplectique est co\"{\i}sotrope si et seulement si
 \beq \label{EqCoisoSymp}
      T_cC^\omega\subset T_cC\mathrm{~quel~que~soit~}c\in C.
\end{equation}
Consid\'{e}rons le {\em sous-fibr\'{e} caract\'{e}ristique}
\beq
   TC^\omega:=\bigcup_{c\in C}T_cC^\omega
\end{equation}
du fibr\'{e} tangent de $C$. La proposition suivante est bien connue:
\bprop \label{PIntTComega}
 Soit $C$ une sous-vari\'{e}t\'{e} co\"{\i}sotrope d'une vari\'{e}t\'{e} symplectique
 $(M,\omega)$. Alors le fibr\'{e} caract\'{e}ristique $TC^\omega$ de $C$ est un
 sous-fibr\'{e} int\'{e}gra\-ble de $TC$.
\eprop
Le feuilletage $\mathcal{F}$ r\'{e}gulier de $C$ associ\'{e} \`{a}
$TC^\omega$ \`{a} l'aide du th\'{e}or\`{e}me de Frobenius sera tr\'{e}s important plus tard.

 On rappelle qu'une sous-vari\'{e}t\'{e} co\"{\i}sotrope est dite
{\em lagrangienne} ssi $T_cC^\omega=T_cC\mathrm{~quel~que~soit~}c\in C$.

On voit facilement que toute application $C\ra M$ dont l'image est dense
(par exemple surjective) est co\"{\i}sotrope. En outre, la compos\'{e}e
$\psi\circ\phi$ d'une application co\"{\i}sotrope $\phi:C\ra M$ et d'une
application de Poisson $\psi:M\ra M'$ est toujours co\"{\i}sotrope.
Le singleton $\{m\}$ est une sous-vari\'{e}t\'{e} co\"{\i}sotrope ssi la structure
de Poisson $P$ s'annule, $P_m=0$. De plus, toute sous-vari\'{e}t\'{e} ferm\'{e}e de
codimension $1$ est co\"{\i}sotrope vu que le fibr\'{e} $TC^{\mathrm{ann}}$ est
de dimension $1$. Le lien entre les applications de
Poisson et les sous-vari\'{e}t\'{e}s co\"{\i}sotropes fut \'{e}tabli par Weinstein:
\bsat[A.Weinstein] \label{TheoGrapheWei}
 Soit $\Phi:(M,P)\ra (M',P')$ une application de Poisson. Alors son
  graphe
  \beq
     C:=\{\big(\Phi(m),m\big)\in M'\times M~|~m\in M\}
  \end{equation}
   est une
  sous-vari\'{e}t\'{e} co\"{\i}sotrope de la vari\'{e}t\'{e} de Poisson
  $(M'\times M, P'_{(1)}-P_{(2)})$
\esat
Voir \cite{Wei88} pour une d\'{e}monstration. En particulier, ceci est vrai
pour une application moment $M\ra\mathfrak{g}^*$. Au cas o\`{u} $(M,\omega)$
est symplectique on a une version symplectifi\'{e}e:\\
Soit $G$ un groupe de Lie r\'{e}el, soit $(\mathfrak{g},[~,~])$ son alg\`{e}bre
de Lie et $\mathfrak{g}^*$ l'espace dual de $\mathfrak{g}$. Il est bien
connu que la multiplication gauche $G\times G\ra G:(g,g')\mapsto
gg'=:L_gg'$ induit un isomorphisme de fibr\'{e}s vectoriels $\Xi:G\times
\mathfrak{g}^*\ra T^*G:(g,\beta)\mapsto \beta \circ T_gL_{g^{-1}}$ et
$G\times \mathfrak{g}^*\times \mathfrak{g} \times \mathfrak{g}^*:
(g,\beta,\xi,\gamma)\ra
\frac{d}{dt}\big(\Xi(g\exp(t\xi),\beta+t\gamma)\big)|_{t=0}$. La forme
symplectique canonique $\omega_G$ sur $T^*G$ est de la forme
\[
     \omega_{G(g,\beta)}\big((g,\beta,\xi_1,\gamma_1)
                          ,(g,\beta,\xi_2,\gamma_2)\big)
         = \langle\gamma_2,\xi_1\rangle-\langle\gamma_1,\xi_2\rangle
             +\langle\beta,[\xi_1,\xi_2]\rangle.
\]
Ainsi on peut rendre le graphe d'une application moment co\"{\i}sotrope
dans une vari\'{e}t\'{e} symplectique:
\bprop
 Soit $(M,\omega)$ une vari\'{e}t\'{e} symplectique, $(\mathfrak{g},[~,~])$ une
 alg\`{e}bre de Lie r\'{e}elle de dimension finie et $J:M\ra \mathfrak{g}^*$ une
 application moment. Alors la sous-vari\'{e}t\'{e}
 \[
    C:=\{(g,\beta,m)\in G\times \mathfrak{g}^*\times M~|~\beta=J(m)\}
 \]
 est une sous-vari\'{e}t\'{e} co\"{\i}sotrope de la vari\'{e}t\'{e} symplectique
 $(T^*G\times M,\omega_{G(1)}+\omega_{(2)})$.
\eprop
La d\'{e}monstration est imm\'{e}diate.

Les sous-vari\'{e}t\'{e}s co\"{\i}sotropes des vari\'{e}t\'{e}s symplectiques admettent
toujours une description locale tr\`{e}s simple:
\bprop\label{PCarteDarboux Adaptee}
 Soit $C$ une sous-vari\'{e}t\'{e} co\"{\i}sotrope ferm\'{e}e de dimension $2n-k$
 d'une vari\'{e}t\'{e} symplectique $(M,\omega)$ de dimension $2n$.\\
  Alors autour de tout point $c\in C$ il existe une carte de Darboux
 \[
 \big(U,(q^1,\ldots,q^{n-k},x^1,\ldots,x^k,p_1,\ldots,p_{n-k},
             y_1,\ldots,y_k)\big)=:\big(U,(q,x,p,y)\big)
 \]
 telle que $\omega|_U=\sum_{i=1}^{n-k}dq^i\wedge dp_i +\sum_{j=1}^k
 dx^j\wedge dy_j$ et $C\cap U=\{u\in U~|~y_1(u)=0,\ldots,y_k(u)=0\}$.
 De plus, la carte $\big(U\cap C,
 (q,p,x)\big)$ de $C$ est
 distingu\'{e}e dans le sens que les restrictions des coordonn\'{e}es
 $(q,p)$ \`{a} $U\cap C$ sont constantes sur
 les feuilles locales de $\mathcal{F}$, tandis que les coordonn\'{e}es
 $x$ donnent un param\'{e}trage des feuilles locales.
\eprop

La d\'{e}monstration de cette proposition se d\'{e}duit de celle d'une
caract\'{e}risa\-tion plus globale des sous-vari\'{e}t\'{e}s co\"{\i}sotropes, due \`{a}
Weinstein et Gotay (voir \cite{Wei71} et \cite{Got82}):
\bsat[A.Weinstein/M.Gotay]\label{TWeinGotay}
 Soit $i:C\ra M$ une sous-vari\'{e}t\'{e} co\"{\i}sotrope d'une vari\'{e}t\'{e}
 symplectique $(M,\omega)$, $E=TC^\omega$ le sous-fibr\'{e}
 caract\'{e}ristique de $TC$ et $\tau:E^*\ra C$ le fibr\'{e} dual de $E$.
 \ben
  \item Soit $F\subset TM|_C$ un sous-fibr\'{e} tel que $TM_C=TC\oplus F$ et
     soit $j:F^{\mathrm{ann}}=E^*\ra T^*C$ l'homomorphisme de fibr\'{e}s
     vectoriels qui en r\'{e}sulte. Soit $\omega_0$ la $2$-forme symplectique
     canonique de $T^*C$.
    Alors il existe un ouvert $V\subset E^*$ qui contient $C\subset
    E^*$ et qui est telle que la $2$ forme ferm\'{e}e
    \beq \label{EqFormeWeinsteinGotay}
       \omega_{WG}:=\tau^*(i^*\omega)+j^*\omega_0
    \end{equation}
    soit une forme symplectique sur $V$.
  \item Il existe un voisinage ouvert $V'\subset V$ de $C$, un voisinage
   ouvert $U$ de $C$ dans $M$ et un diff\'{e}omorphisme $\phi:V'\ra U$ (dont
   la restriction \`{a} $C$ induit l'application identique $C\ra C$) tels que
   \[
       \phi^*\omega =\omega_{WG}|_{V'}.
   \]
   En particulier, si $C$ est lagrangienne (i.e. $i^*\omega=0$ et $E=TC$) le
   voisinage $U$ est symplectomorphe \`{a} un voisinage ouvert de la section
   nulle du fibr\'{e} cotangent $T^*C$ de $C$.
 \een
\esat
On rappelle qu'une vari\'{e}t\'{e} diff\'{e}rentiable
$C$ munie d'une $2$-forme ferm\'{e}e $\varpi$ est dite {\em pr\'{e}symplectique}
lorsque le noyau de $\varpi$,
$E:=\cup_{c\in C}\{v\in T_cC~|~\varpi_c(v,w)=0~\forall w\in T_cC\}$, est un
sous-fibr\'{e} vectoriel ({\em le fibr\'{e} caract\'{e}ristique de $\varpi$}) du fibr\'{e}
tangent $TC$. $E$ est automatiquement int\'{e}grable. Donc toute sous-vari\'{e}t\'{e}
co\"{\i}sotrope ferm\'{e}e $C$ d'une vari\'{e}t\'{e} symplectique est toujours
pr\'{e}symplectique avec $\varpi:=i^*\omega$, et le
th\'{e}or\`{e}me \ref{TWeinGotay} montre que un voisinage ouvert d'une sous-vari\'{e}t\'{e}
co\"{\i}sotrope d'une vari\'{e}t\'{e} symplectique est enti\`{e}rement d\'{e}termin\'{e}
par la sous-vari\'{e}t\'{e} pr\'{e}symplectique $(C,i^*\omega)$.

\noindent La caract\'{e}risation suivante des champs de vecteurs verticaux
s'av\`{e}rera utile plus tard:
\bprop \label{PChampsVertIdeal}
 Soit $(M,\omega)$ une vari\'{e}t\'{e} symplectique et $i:C\ra M$ une
 sous-vari\'{e}t\'{e} co\"{\i}sotrope ferm\'{e}e. Alors l'espace de tous les champs
 de vecteurs verticaux est isomorphe \`{a} l'espace de toutes les
 restrictions \`{a} $C$ des champs hamiltoniens $X_g$ des \'{e}l\'{e}ments
 $g\in\mathcal{I}$.
\eprop
\bbew
 L'application $\omega^\flat:TM\ra T^*M$ se restreint \`{a} un
 isomorphisme du fibr\'{e} ${TC^{ann}}^*\ra TM|_C/TC\ra E^*=TC^{\omega*}$
 d\'{e}fini par
 $[v]\mapsto \big(w\mapsto \omega(v,w)\big)$ o\`{u} $v\in T_cM$ et $w\in
 T_cC^\omega$. D'apr\`{e}s le lemme \ref{LemIIcarre}, (3), les restrictions
 de $dg$ \`{a} $C$ de $g\in\mathcal{I}$ donnent toutes les sections de classe
 $\Cinf$ de $TC^{ann}$. Puisque $X_g$ est d\'{e}fini par
 $(\omega^\flat)^{-1}(dg)$ la proposition est claire.
\ebew

Revenant au cas d'une vari\'{e}t\'{e} co\"{\i}sotrope $C$ d'une vari\'{e}t\'{e} de
Poisson $(M,P)$, on rappelle le complexe suivant: puisque $\mathcal{I}$
est \`{a} la fois une alg\`{e}bre de Lie (par rapport au crochet de Poisson)
et un id\'{e}al (par rapport \`{a} la multiplication associatie commutative)
on voit facilement que le $\Cinf(C,\korps)$-module conormal
$\mathcal{I}/\mathcal{I}^2$
(qui est isomorphe \`{a} $\Ginf(C,TC^{\mathrm{ann}})$ d'apr\`{e}s le lemme
\ref{LemIIcarre}) est une alg\`{e}bre de Lie: soient $g_1,g_2\in\mathcal{I}$,
le crochet est d\'{e}fini par
\beq
      \{dg_1|_C,dg_2|_C\}:=(d\{g_1,g_2\})|_C.
\end{equation}
Ensuite, l'alg\`{e}bre de Lie $\mathcal{I}/\mathcal{I}^2$ agit sur
$\Cinf(C,\korps)\cong \Cinf(M,\korps)/\mathcal{I}$ en tant que d\'{e}rivations:
soient $f\in \Cinf(M,\korps)$ et $g\in \mathcal{I}$, alors
\beq
   (dg|_C).i^*f:=i^*\{g,f\}=-i^*X_gf
\end{equation}
Alors le lemme suivant est \'{e}vident:
\blem\label{LCohdP}
  L'espace de sections $\Ginf\big(C,\Lambda(TM|_C/TC)\big)$ --qui est
  isomorphe \`{a}
  $\mathrm{Hom}_{\Cinf(C,\korps)}\big(\Lambda_{\Cinf(C,\korps)}
  \mathcal{I}/\mathcal{I}^2,\Cinf(C,\korps)\big)$-- est un cocomplexe
  \`{a} l'aide de la diff\'{e}rentielle $d_P$ suivante: soit $k\in\nat$,
  $A\in \Ginf\big(C,\Lambda^k(TM|_C/TC)\big)$ (qui est isomorphe \`{a}
  $\Ginf\big(C,\Lambda^k(TC^{\mathrm{ann}})^*\big)$) et
  $g_0,g_1,\ldots,g_k\in \mathcal{I}$, alors
  \bea
    \lefteqn{(d_PA)(dg_0|_C,\ldots,dg_k|_C):=} \nonumber \\
      & & \sum_{i=0}^k(-1)^i~(dg_i|_C).\big(A(dg_0|_C,\ldots,
             \widehat{dg_i|_C},\ldots,dg_k|_C\big) \nonumber \\
      & & +~\sum_{0\leq i<j\leq k}(-1)^{i+j}
               A\big(\{dg_i|_C,dg_j|_C\},dg_0|_C,\ldots,
               \widehat{dg_i|_C},\ldots,\widehat{dg_j|_C},\ldots,
               dg_k|_C\big) \nonumber \\
      & &
  \eea
  o\`{u} $\hat{~}$ d\'{e}signe l'omission de l'argument.
\elem
La formule pour $d_P$ est la fomule usuelle pour la cohomologie de
Chevalley-Eilenberg de l'alg\'{e}bre de Lie $\mathcal{I}/\mathcal{I}^2$
\`{a} valeurs dans le module $\Cinf(C,\korps)$. On va noter les groupes de
cohomologie assoic\'{e}s au complexe pr\'{e}c\'{e}dent par $H_P^k(C,\korps)$: on
en parlera de la {\em cohomologie de Poisson transversale} ou de
la {\em cohomologie BRST}. Ce nom provient de la th\'{e}orie BRST de physique
pour d\'{e}crire des `syst\`{e}mes avec des contraintes', voir \cite{BHW00} et
des r\'{e}f\'{e}rences y mentionn\'{e}es, par exemple \cite{HT92}, \cite{Dub87} et
\cite{FHST89}.

\subsection{Vari\'{e}t\'{e}s feuillet\'{e}es}
          \label{SubSecVarFeuil}

Dans le paragraphe pr\'{e}c\'{e}dent, on a vu qu'une sous-vari\'{e}t\'{e} co\"{\i}sotrope
$C$ d'une vari\'{e}t\'{e} symplectique $(M,\omega)$ est munie d'un feuilletage
$\mathcal{F}$ d\'{e}fini par le fibr\'{e} caract\'{e}ristique $E:=TC^\omega\subset
TC$ (voir la proposition \ref{PIntTComega}). Dans ce paragraphe on va
rappeler quelques d\'{e}tails dont on aura besoin plus tard. On s'appuie
largement sur les livres \cite{KT75}, \cite{Mol88} et surtout \cite{Ton97}
qui contient
une bibliographie tr\`{e}s vaste sur ce sujet.

En g\'{e}n\'{e}ral, soit $C$ une vari\'{e}t\'{e} diff\'{e}rentiable et $E\subset TC$ un
sous-fibr\'{e} int\'{e}grable. Parfois, $E$ est dit le {\em sous-fibr\'{e} vertical}.
Les sections de classe $\Cinf$ de $E$ sont dits les
{\em champs de vecteurs verticaux}, et l'int\'{e}grabilit\'{e} du fibr\'{e} $E$
est \'{e}quivalent \`{a} dire que $\Ginf(C,E)$ est une sous-alg\'{e}bre de Lie de
l'alg\`{e}bre de Lie $\mathcal{X}(C)$ de tous les champs de vecteurs sur $C$.
 On appelle $Q:=TC/E$ le
{\em fibr\'{e} normal du feuilletage} et note $v\mapsto \overline{v}$ la
projection canonique $TC\ra Q$. L'espace de sections $\Ginf(C,Q)$
est isomorphe au $\Cinf(C,\korps)$-module quotient
$\mathcal{X}(C)/\Ginf(C,E)$.

On d\'{e}signe par $\mathcal{F}$ le feuilletage associ\'{e} \`{a} $E$ gr\^{a}ce au
th\'{e}or\`{e}me de Frobenius (voir par exemple \cite[p.24-29]{KMS93} pour une
d\'{e}monstration), c.-\`{a}-d. $\mathcal{F}$ est une partition de $C$
en sous-vari\'{e}t\'{e}s immerg\'{e}es de $C$ telle que pour tout $c\in C$ l'espace
tangent
en $c$ de la feuille passant par $c$ est \'{e}gal \`{a} la fibre $E_c$.
Une carte
$(U,(\xi^1,\ldots,\xi^{n-k},x^1,\ldots,x^k)\big)=:\big(U,(x,\xi)\big)$
de $C$ est dite {\em distingu\'{e}e} lorsque $\langle d\xi^i, V\rangle=0$
quels que soient $1\leq i\leq n-k$ et $V\in \Ginf(C,E)$. Les \'{e}quations
$\xi^1=c^1,\ldots,\xi^{n-k}=c^{n-k}$ avec $c^1,\ldots,c^{n-k}\in\real$
d\'{e}finissent un feuilletage local de l'ouvert $U$ dont les feuilles sont
dites des {\em plaques}. Pour tout $c\in U$ la plaque contenant $c$ fait
partie de la feuille contenant $c$.
Soit $C/\mathcal{F}$ l'espace topologique des feuilles (muni de la topologie
quotient) et $\pi:C\ra C/\mathcal{F}$ la projection canonique.

On rappelle la {\em connexion de Bott correspondant au sous-fibr\'{e} int\'{e}grable
$E$}: soient $V$ un champ de vecteurs
vertical, $X$ un champ de vecteurs quelconque sur $C$
et $\overline{X}$
la section correspondante dans $\Ginf(C,Q)\cong\mathcal{X}(C)/\Ginf(C,E)$,
alors
\beq\label{EqConnexionBott}
  \nabla^{\mathrm{Bott}}_V\overline{X}:=\overline{[V,X]}
\end{equation}
On voit facilement que $\nabla^{\mathrm{Bott}}$ est bien d\'{e}finie,
et satisfait aux axiomes d'une connexion (partielle) le long du fibr\'{e} $E$,
c.-\`{a}-d. la $\korps$-bilin\'{e}arit\'{e} $\Ginf(C,E)\times \Ginf(C,Q)\ra \Ginf(C,Q)$
et pour toutes les fonctions $f\in\Cinf(C,\korps)$ les \'{e}quations
$\nabla^{\mathrm{Bott}}_{fV}\overline{X}
=f\nabla^{\mathrm{Bott}}_V\overline{X}$ et
$\nabla^{\mathrm{Bott}}_V(f\overline{X})
=f\nabla^{\mathrm{Bott}}_V\overline{X}+ (Vf)\overline{X}$. De plus,
puisque $\mathcal{X}(C)$ et $\Ginf(C,E)$ sont des modules de l'alg\`{e}bre de
Lie $\Ginf(C,E)$, il en est de m\^{e}me avec
$\Ginf(C,Q)\cong\mathcal{X}(C)/\Ginf(C,E)$, ce qui veut dire --en termes
g\'{e}om\'{e}triques-- que la connexion de Bott est toujours plate
\beq \label{EConnexionBottEstPlate}
 \nabla^{\mathrm{Bott}}_V\nabla^{\mathrm{Bott}}_W\overline{X}
   -\nabla^{\mathrm{Bott}}_W\nabla^{\mathrm{Bott}}_V\overline{X}
    -\nabla^{\mathrm{Bott}}_{[V,W]}\overline{X}=0.
\end{equation}
Par cons\'{e}quent, les espaces $\Ginf(C,Q)$ et $\Ginf(C,Q^*)$
sont des modules de l'alg\`{e}bre
de Lie $\Ginf(C,E)$ par rapport \`{a} la connexion de Bott. Plus
g\'{e}n\'{e}rale\-ment, soit $\mathsf{F}Q$ une somme directe de produits tensoriels
des fibr\'{e}s
$S^NQ$, $S^{N'}Q^*$, $\Lambda^k Q$ et $\Lambda^l Q^*$
o\`{u} $N,N',k,l$ sont des entiers positifs. Il s'ensuit
directement que l'espace de sections
$\Ginf(C,\mathsf{F}Q\otimes \Lambda E^*)$
est muni de la structure d'un cocomplexe dont la diff\'{e}rentielle $d_v$
est donn\'{e}e par la formule usuelle de Chevalley-Eilenberg: soient
$A\in\Ginf(C,\mathsf{F}Q\otimes \Lambda^k E^*)$ et
$V_0,V_1,\ldots,V_k\in\Ginf(C,E)$, alors
\bea
 \lefteqn{(d_vA)(V_0,\ldots,V_k):=} \nonumber \\
  & & \sum_{i=0}^k(-1)^i~\nabla^{\mathrm{Bott}}_{V_i}
         \big(A(V_0\ldots,\widehat{V_i},\ldots,V_k)\big) \nonumber \\
  & &        + \sum_{0\leq i<j\leq k}(-1)^{i+j}
             A\big([V_i,V_j],V_0,\ldots,\widehat{V_i},\ldots,
             \widehat{V_j},\ldots,V_k\big) \label{EqDefdverticale}
\eea

La cohomologie de ce cocomplexe sera importante plus tard pour quelques
classes caract\'{e}ristiques associ\'{e}es au feuilletage $\mathcal{F}$, notamment
pour la classe d'Atiyah-Molino, voir \'{e}galement le paragraphe
\ref{SubSecClassedAtiyahMolino}.

En particulier, au cas o\`{u} $\mathsf{F}Q$ est le fibr\'{e} trivial \`{a} fibre
$\korps$, la connexion de Bott se r\'{e}duit \`{a} la d\'{e}riv\'{e}e de Lie, et on obtient
{\em l'espace des formes diff\'{e}rentielles verticales}
$\Omega_v(C):=\oplus_{k=0}^{\dim C}\Ginf(C,\Lambda^k E^*)$, et
la cohomologie du cocomplexe
$\big(\Omega_v(C),d_v\big)$ est dite la {\em cohomologie verticale de
$C$}. On dit parfois la {\em cohomologie longitudinale} ou la {\em
cohomologie le long des feuilles}. En choisissant une m\'{e}trique
riemannienne sur $C$ on obtient toujours un {\em sous-fibr\'{e} horizontal}
$F$ de $TC$, c.-\`{a}-d. un sous-fibr\'{e} de $TC$ tel que $TC=F\oplus E$.
Soit $\Omega(C):=\oplus_{k=0}^{\dim C}\Ginf(C,\Lambda^k TC^*)$ l'espace
de toutes les formes diff\'{e}rentielles sur $C$, muni avec la diff\'{e}rentielle
de Cartan, $d$.
On voit ais\'{e}ment que l'application de restriction
\beq \label{EqDefRestVert}
        p_v: \Omega(C)\ra \Omega_v(C)
\end{equation}
(o\`{u} tous les arguments d'une forme diff\'{e}rentielle dans $\Omega(C)$ sont
restreints \`{a} $E$) est surjective et de plus un morphisme de cocomplexes.

 Le noyau de la restriction (\ref{EqDefRestVert}) est appel\'{e}
{\em l'espace
des formes diff\'{e}ren\-tielles relatives}, $\Omega_{rel}(C)$. Il est stable
par la diff\'{e}rentielle de Cartan, donc on a la suite exacte de cocomplexes
\beq \label{EqDefVertSuite}
    \{0\}\ra \big(\Omega_{rel}(C),d\big)
     \stackrel{i_{rel}}{\ra}     \big(\Omega(C),d\big)
     \stackrel{p_v}{\ra} \big(\Omega_v(C),d_v\big)\ra \{0\}.
\end{equation}
La cohomologie de $\big(\Omega_{rel}(C),d\big)$ s'appelle  {\em la
cohomologie relative} de $C$. Cette suite exacte de cocomplexes induit
une suite exacte longue en cohomologie.

Un cas particulier d'une forme relative est une forme diff\'{e}rentielle
$\alpha$ dite
{\em basique}: ceci est le cas lorsque $i_V\alpha=0$ et $L_V\alpha=0$
quels que soient
les champs de vecteurs verticaux $V$.
L'espace $\Omega_{bas}(C)$ de toutes les
formes basiques est un sous-complex du complex de de Rham
$\big(\Omega(C),d\big)$, et sa cohomologie associ\'{e}e est dite la
{\em cohomologie basique} de $C$. Si $C/\mathcal{F}$ est munie d'une
structure diff\'{e}rentiable telle que la projection $\pi$ soit une submersion
surjective, les formes diff\'{e}rentielles basiques co\"{\i}ncident avec
les formes diff\'{e}rentielles retir\'{e}es de la vari\'{e}t\'{e} quotient
$C/\mathcal{F}$. Dans ce cas la cohomologie basique est isomorphe \`{a} la
cohomologie de de Rham de $C/\mathcal{F}$.

Si $C$ est une sous-vari\'{e}t\'{e} co\"{\i}sotrope d'une vari\'{e}t\'{e} symplectique
$(M,\omega)$, alors la $2$-forme $i^*\omega$ est toujours basique.
Pour plus de renseignements concernant les suites spectrales qu'on peut
\'{e}galement d\'{e}finir dans ce domaine, voir \cite{KT75}.

\subsection{R\'{e}duction symplectique}
          \label{SubSecRed1}

Dans ce sous-paragraphe on ne traite que les vari\'{e}t\'{e}s symplectiques
(pour des g\'{e}n\'{e}ralisations voir par exemple \cite{MR86}):

Soit $i:C\ra M$ une sous-vari\'{e}t\'{e} co\"{\i}sotrope ferm\'{e}e d'une vari\'{e}t\'{e}
symplectique
$(M,\omega)$, soit $E=TC^\omega\subset TC$ le sous-fibr\'{e} caract\'{e}ristique
de $TC$, et soit $\mathcal{F}$ le feuilletage  de $C$ correspondant
au fibr\'{e} int\'{e}grable $E$ (\`{a} l'aide du th\'{e}or\`{e}me classique de
Frobenius, voir par exemple \cite{KMS93}, p.~28, Thm.3.25).
Supposons que l'espace des feuilles
\beq
    M_{red}:=C/\mathcal{F}
\end{equation}
est muni
d'une structure diff\'{e}rentiables compatible avec la topologie quotient
telle que la projection canonique $\pi:C\ra M_{red}$ soit une submersion
surjective. Alors on a le th\'{e}or\`{e}me classique suivant:

\bsat\label{RedSymp}
  Avec les hypoth\`{e}ses mentionn\'{e}es ci-dessus, l'espace quotient
  $M_{red}$ est muni d'une structure symplectique canonique,
  $\omega_{red}$, d\'{e}finie par
  \[
      i^*\omega =:\pi^*\omega_{red}.
  \]
  La vari\'{e}t\'{e} symplectique $(M_{red},\omega_{red})$ s'appelle
  la {\em vari\'{e}t\'{e} symplectique r\'{e}duite}.
\esat
Voir par exemple \cite{AM85}, p. 416, Thm. 5.3.23, pour une d\'{e}monstration.

Un cas particulier important s'obtient par une application moment
$J:M\ra \mathfrak{g}^*$ pour laquelle $0$ est une valeur r\'{e}guli\`{e}re
dont l'image r\'{e}ciproque $C:=J^{-1}(0)$ n'est pas vide. Dans ce cas-l\`{a},
$C$ est une sous-vari\'{e}t\'{e} co\"{\i}sotrope, et la vari\'{e}t\'{e} r\'{e}duite
(au cas o\`{u} elle existe) s'obtient en tant qu'espace quotient du groupe
de Lie $G$ (\`{a} alg\`{e}bre de Lie $\mathfrak{g}$) agissant de fa\c{c}on
libre et propre sur $C$. Cette construction importante et extr\^{e}mement
utile est appel\'{e}e la {\em r\'{e}duction de Marsden-Weinstein} \cite{MW74}.
Par exemple l'espace projectif complexe s'obtient en tant que vari\'{e}t\'{e}
symplectique r\'{e}duite de $M=\real^{2n+2},\omega=
\sum_{k=1}^{n+1}dq^k\wedge dp_k$ \`{a} l'aide de l'application moment
$J(q,p):=\frac{\sum_{k=1}^{n+1}(q_k^2+p_k^2)}{2}-\frac{1}{2}$ pour
l'action du groupe $U(1)$ sur $\complex^{n+1}=\real^{2n+2}$.

\noindent Le lemme suivant est bien connu (voir par exemple \cite{BW95})
et permet de
r\'{e}duire quelques questions concernant des sous-vari\'{e}t\'{e}s co\"{\i}sotropes
au cas des sous-vari\'{e}t\'{e}s lagrangiennes:
\blem\label{LCLagdansMMred}
 Soit $i:C\ra M$ une sous-vari\'{e}t\'{e} co\"{\i}sotrope ferm\'{e}e
 d'une vari\'{e}t\'{e} symplectique
 $(M,\omega)$, et on suppose que la vari\'{e}t\'{e} symplectique r\'{e}duite
 $(M_{red},\omega_{red})$ existe o\`{u} $\pi:M\ra M_{red}$ est la projection
 canonique. On consid\`{e}re la vari\'{e}t\'{e} symplectique $(M\times
 M_{red},\omega_{(1)}-\omega_{red~(2)})$.\\
 Alors $i$ et $\pi$ sont des applications co\"{\i}sotropes, et
 l'application suivante
 \beq
    l_C:C\ra M\times M_{red}:c\mapsto \big(i(c),\pi(c)\big)
 \end{equation}
 est un plongement de $C$ sur une sous-vari\'{e}t\'{e} lagrangienne de
 $(M\times M_{red},\hat{\omega}:=\omega_{(1)}-\omega_{red~(2)})$.
\elem
\bbew
 L'application $i$ est un plongement, donc une immersion injective propre.
 Alors $l_C$ reste une immersion injective propre, donc un plongement.
 Puisque $\dim (M_{red})=\dim(C)-\mathrm{codim}(C)=2\dim(C)-\dim(M)$ on a
 $\dim(M\times M_{red})=2\dim(C)$. Alors il suffit de montrer que $l_C(C)$
 soit
 isotrope dans $(M\times M_{red},\hat{\omega})$, i.e.
 $\hat{\omega}_{l_C(c)}\big(w,w'\big)=0$ quels que soient $w,w'$ dans
 l'espace tangent de $l_C(C)$ en $l_C(c)$: il y a $v,v'\in T_cC$ tels que
 $w=T_cl_C v$ et $w'=T_cl_C v'$. Alors
 \beas
    \hat{\omega}_{l_C(c)}\big(w,w'\big)
         & = & \omega_{i(c)}\big(T_ci~ v, T_ci~v'\big)
        -\omega_{red~\pi(c)}\big(T_c\pi ~v,T_c\pi ~v'\big)\\
        & = & (i^*\omega-\pi^*\omega)_c\big(v,v'\big)=0
 \eeas
\ebew

On obtient une caract\'{e}risation tr\`{e}s importante de
l'espace $\CinfK{M_{red}}$ de la fa\c{c}on suivante:
Soit
\beq\label{EqIdLieI}
  \mathcal{N}(\mathcal{I}):=\{f\in\CinfK{M}~|~\{f,g\}\in\mathcal{I}~\forall
  g\in\mathcal{I}\}.
\end{equation}
{\em l'id\'{e}alisateur (de Lie) de $\mathcal{I}$}. On note que
$\mathcal{N}(\mathcal{I})$ est exactement l'espace de fonctions dans
$\CinfK{M}$ dont la restriction \`{a} $C$ est constante sur les feuilles du
feuilletage $\mathcal{F}$ induit par le fibr\'{e} caract\'{e}ristique $E$.
On a la
\bprop \label{PIdLieImodI}
 $\mathcal{B}$ est une sous-alg\`{e}bre de Poisson de $\CinfK{M}$, l'id\'{e}al
 annulateur $\mathcal{I}$ est un id\'{e}al de Poisson de $\mathcal{B}$, et
 la restriction $i^*$ \`{a} $C$ et la projection $\pi$ induit l'isomorphisme
 d'alg\`{e}bres de Poisson suivant:
 \beq
  (\CinfK{M_{red}},\{~,~\}_{\omega_{red}})\cong
           \mathcal{N}(\mathcal{I})/\mathcal{I}
 \end{equation}
 Plus pr\'{e}cis\'{e}ment, soient $f_1,f_2\in \mathcal{N}(\mathcal{I})$, donc il
 existent des uniques fonction $\tilde{f}_1,\tilde{f}_2\in
 \CinfK{M_{\mathrm{red}}}$ telles que $i^*f_1=\pi^*\tilde{f}_1$ et
 $i^*f_2=\pi^*\tilde{f}_2$. Alors
 \beq
   i^*\{f_1,f_2\}=\pi^*\{\tilde{f}_1,\tilde{f}_2\}_{\mathrm{red}}.
 \end{equation}
\eprop
Voir \cite[p.417-418]{AM85} et \cite[p.443]{Kim92} pour la d\'{e}monstration.

\subsection{Op\'{e}rateurs multidiff\'{e}rentiels le long des applications}
 \label{SubSecOpMultidiffAppl}

Dans ce paragraphe on rappelle quelques notions autour des op\'{e}rateurs
multidiff\'{e}rentiels `le long des applications'.

Comme d'habitude, on va utiliser l'abbr\'{e}viation suivante pour les d\'{e}riv\'{e}es
partielles dans $\real^n$: soit $I=(p_1,\ldots,p_n)\in \nat^n$ un multi-indice,
alors on va \'{e}crire $|I|:=p_1+\ldots + p_n$ et
\[
    \frac{\partial}{\partial x^I}
     :=\partial_{x^{I}}:=\frac{\partial^{|I|}}{\partial {x^1}^{p_1}\cdots
                                            \partial {x^n}^{p_n}}
\]
o\`{u} $\partial^0/\partial {x^k}^0$ est \'{e}quivalent \`{a} dire que les d\'{e}riv\'{e}es
partielles par rapport \`{a} $x^k$ ne figurent pas dans l'expression
ci-dessus.

Soit $k$ un entier strictement positif, $\tau=\tau_0:E=E_0\ra M=M_0$,
$\tau_1:E_1\ra M_1,\ldots,\tau_k:E_k\ra M_k$ des fibr\'{e} vectoriels sur
les vari\'{e}t\'{e}s
diff\'{e}rentiables $M,M_1,\ldots,M_k$ et $\phi_1:M\ra
M_1$,$\ldots$,$\phi_k:M\ra M_k$ des applications de classe $\Cinf$.
Pour la d\'{e}finition des op\'{e}rateurs multidiff\'{e}rentiels dans cette situation
on a besoin d'une description locale: soient
$\big(U,(x^1,\ldots,x^n)\big)$,
$\big(U_1,({x_1}^1,\ldots,{x_1}^{n_1})\big)$,
$\ldots$,$\big(U_k,({x_k}^1,\ldots,{x_k}^{n_k})\big)$ des cartes locales
de $M,M_1$, $\ldots,M_k$, respectivement, telles que $\phi_k(U)\subset U_k$.
On suppose que les fibr\'{e}s $E, E_1,\ldots, E_k$ sont trivialisables sur
$U,U_1,\ldots,U_k$, respectivement, et soient
$(e^{(0)}_1,\ldots,e^{(0)}_{d_0})$,
$(e^{(1)}_1,\ldots,e^{(1)}_{d_1})$,$\ldots$,
$(e^{(k)}_1,\ldots,e^{(k)}_{d_k})$ des bases locales des espaces des
sections locales de $E,E_1,\ldots,E_k$, respectivement. Alors toute
section $\psi_r\in\Ginf(M_r,E_r)$ est localement une combinaison lin\'{e}aire
$\sum_{a_r=1}^{d_r}\psi^{a_r}_re^{(r)}_{a_r}$ avec $\psi^{a_r}_r\in
\CinfK{U_r}$ quels que soient $0\leq r\leq k$ et $1\leq a_r\leq d_r$.

Une application $k$-lin\'{e}aire sur $\korps$
\beq
 \mathsf{D}:\Ginf(M_1,E_1)\times \cdots\times\Ginf(M_k,E_k)\ra\Ginf(M,E)
\end{equation}
 est dite un {\em op\'{e}rateur $k$-diff\'{e}rentiel le
long des applications $\phi_1,\ldots,\phi_k$}
lorsqu'il existe un entier positif $N$
telle que $\mathsf{D}$ est de la forme locale suivante: pour tout syst\`{e}me de
cartes et trivialisations des fibr\'{e}s il existe des fonctions
$\mathsf{D}_{a_1\cdots a_k}^{I_1\cdots I_k~a_0}\in\CinfK{U}$ (quels que soient
les multi-indices $I_1\in\nat^{n_1},\ldots,I_k\in\nat^{n_k}$, les indices
$0\leq r\leq k$ et $1\leq a_r\leq d_r$) telles que pour tout $m\in U$
il vient:
\bea \label{EqDefFormeLocale}
 \mathsf{D}(\psi_1,\ldots,\psi_k)_m
     & := & \sum_{a_0,a_1,\ldots,a_k}
     \sum_{\stackrel{I_1,\ldots,I_k}{|I_1|,\ldots,|I_k|\leq N}}
       \mathsf{D}_{a_1\cdots a_k}^{I_1\cdots I_k~a_0}(m) \nonumber \\
       & & ~~~\frac{\partial \psi_1^{a_1}}{\partial {x_1^{I_1}}}
                \big(\phi_1(m)\big)
               \cdots
              \frac{\partial \psi_k^{a_k}}{\partial {x_k^{I_k}}}
                \big(\phi_k(m)\big)
             e^{(0)}_{a_0}(m) \nonumber \\
       & &
\eea
Globalement, tout op\'{e}rateur $k$-diff\'{e}rentiel le long des applications
$\phi_1,\ldots,\phi_k$ peut \^{e}tre regard\'{e} comme
section de classe $\Cinf$ du fibr\'{e} suivant sur $M$:
\beq \label{EqOpMultiGlobalUndSo}
    \mathrm{Hom}\big(\phi_1^*(J^NE_1)\otimes\cdots\otimes
                     \phi_k^*(J^NE_k),E\big)
\end{equation}
o\`{u} $J^NE_r$ d\'{e}signe le fibr\'{e} vectoriel sur $M$ qui est la $N$i\`{e}me
prolongation jet
du fibr\'{e} $E_r$ sur $M_r$,
et $\phi_r^*J^NE_r$ d\'{e}signe le fibr\'{e} retir\'{e} \`{a} $M$ (voir par exemple
\cite[p.124]{KMS93} pour des d\'{e}finitions). Nous allons \'{e}crire
\beq
  \Dop^k\big(\Ginf(M_1,E_1),\ldots,\Ginf(M_k,E_k);\Ginf(M,E)\big)
\end{equation}
pour l'espace vectoriel
de tous les op\'{e}rateurs $k$-diff\'{e}rentiels le long des applications
$\phi_1,\ldots,\phi_k$. Si $M_1=\cdots=M_k$ et $E_1=\cdots=E_k$ on \'{e}crira
$\Dop^k\big(\Ginf(M_1,E_1);\Ginf(M,E)\big)$. Dans le cas des fibr\'{e}s triviaux \`{a}
fibre type $\korps$ on utilisera la notation $\CinfK{M_r}$ au lieu de
$\Ginf(M_r,E_r)$ ($0\leq r\leq k$).

Soient $U_1\subset M_1$,...,$U_k\subset M_k$ des ouverts, soit
$U:=\phi_1^{-1}(U_1)\cap\cdots\cap\phi_k^{-1}(U_k)\subset M$ et soit
$\mathsf{D}$ un op\'{e}rateur $k$-diff\'{e}rentiel le long des applications
$\phi_1,\ldots,\phi_k$. Soient $\psi'_1\in \Ginf(U_1,E_1|_{U_1})$,...,
$\psi'_k\in\Ginf(U_k,E_k|_{U_k}$ des sections locales. On d\'{e}duit de la
forme locale (\ref{EqDefFormeLocale}) de $\mathsf{D}$ qu'il existe
un unique op\'{e}rateur $k$-diff\'{e}rentiel le long de $\phi_1|_U,\ldots,
\phi_k|_U$ --qu'on va noter par $\mathsf{D}|_{U_1\times\cdots\times
U_k}$-- tel que
\beq\label{EqDefRestrOpMultidiffLeLong}
    \mathsf{D}|_{U_1\times\cdots\times U_k}
    (\psi_1|_{U_1},\ldots,\psi_k|_{U_k})
          = \big(\mathsf{D}(\psi_1,\ldots,\psi_k)\big)|_{U}.
\end{equation}
quels que soient les sections globales $\psi_1\in\Ginf(M_1,E_1)$,...,
$\psi_k\in\Ginf(M_k,E_k)$.
On va parler de $\mathsf{D}|_{U_1\times\cdots\times U_k}$ de la
{\em restriction de $\mathsf{D}$ aux ouverts} $U_1$,...,$U_k$.

Il y a une r\`{e}gle \'{e}vidente de composer plusieurs op\'{e}rateurs
multidiff\'{e}rentiels en respectant la composition des applications
sous-jacentes: on laisse au lecteur le soin de les \'{e}crire explicitement.

\subsection{Star-produits}
          \label{SubSecStarprod}

Soit $(M,P)$ une vari\'{e}t\'{e} de Poisson. D'apr\`{e}s \cite{BFFLS78},
un {\em star-produit} $*$ sur
$M$ est la donn\'{e}e d'une suite d'op\'{e}rateurs
bidiff\'{e}rentiels\footnote{Dans ce travail on ne consid\`{e}re que
des star-produits {\em diff\'{e}rentiels} et non pas les star-produits
locaux ou continus.}
$(\mathsf{C}_r)_{r\in\nat}$ o\`{u} $\mathsf{C}_r\in\Dop^2(\CinfK{M};\CinfK{M})$
telle que
pour toutes fonctions $f,g,h\in\CinfK{M}$ on a les conditions suivantes:
\bea
     \mathsf{C}_0(f,g) & = & fg, \\
     \mathsf{C}_1(f,g)-\mathsf{C}_1(g,f) & = & 2P(df,dg), \\
     \mathsf{C}_r(1,f)=\mathsf{C}_r(f,1) & = & 0~~~\forall r\geq 1, \\
     \sum_{a=0}^r\mathsf{C}_a\big(\mathsf{C}_{r-a}(f,g),h\big) & = &
          \sum_{a=0}^r \mathsf{C}_a\big(f,\mathsf{C}_{r-a}(g,h)\big)~~~
            \forall r\in\nat.
\eea
La multiplication $*$ est d\'{e}finie par la s\'{e}rie formelle
\beq
     * := \sum_{r=0}\nu^r \mathsf{C}_r
\end{equation}
vue comme application $\korps[[\nu]]$-bilin\'{e}aire de
$\CinfK{M}[[\nu]]\times \CinfK{M}[[\nu]]$ dans $\CinfK{M}[[\nu]]$. Avec les
conditions pr\'{e}c\'{e}dentes, on voit directement que le $\korps[[\nu]]$-module
$(\CinfK{M}[[\nu]],*)$ est muni de la structure d'une alg\`{e}bre
associative sur $\korps[[\nu]]$. Cette alg\`{e}bre est une d\'{e}formation
associative formelle dans le sens de Gerstenhaber, \cite{Ger63}, \cite{Ger64}.
Le
param\`{e}tre formel $\nu$ correspond \`{a} $\frac{i\hbar}{2}$ o\`{u} $i:=\sqrt{-1}$
et $\hbar$ est la constante de Planck. Pour une partie ouverte $U$ de $M$
on va noter $*|_U$ la s\'{e}rie formelle des restrictions \`{a} $U$ des op\'{e}rateurs
bidiff\'{e}rentiels $\mathsf{C}_s$. Il est \'{e}vident que $*|_U$ est un
star-produit sur la vari\'{e}t\'{e} $(U,P|_U)$.

Sur une vari\'{e}t\'{e} de Poisson arbitraire $(M,P)$ il y a toujours des
star-produits: ceci a \'{e}t\'{e} d\'{e}montr\'{e} par DeWilde et Lecomte pour le cas
symplectique en 1983 \cite{DL83} et plus tard par Fedosov \cite{Fed85}
(voir \'{e}galement le paragraphe \ref{SubSecQuantifMorphVarSymp} pour une
description),
et finalement pour le cas g\'{e}n\'{e}ral en 1997 par Kontsevitch \cite{Kon97b}
comme cas particulier de sa conjecture de formalit\'{e}.

Deux star-produits $*$ et $*'$ sur la m\^{e}me vari\'{e}t\'{e} de Poisson $(M,P)$
sont dits
{\em \'{e}quivalents} (notation: $*~\sim~*'$) lorsqu'il y a un isomorphisme
d'alg\`{e}bres associatives
unitaires $S:\CinfK{M}[[\nu]]\ra\CinfK{M}[[\nu]]$ (qui est dite
`transformation d'\'{e}quivalence') dans le sens suivant:
on a
\[
   f*' g = S\big((S^{-1}f)*(S^{-1}g)\big)~~~
   \forall ~f,g\in\CinfK{M}[[\nu]]
\]
avec
\[
    S= id +\sum_{r=1}^\infty\nu^r S_r
\]
o\`{u} chaque $S_r:\CinfK{M}\ra \CinfK{M}$ est un op\'{e}rateur diff\'{e}rentiel qui
s'annule sur les constantes. Dans ce cas on \'{e}crit $*'=S*$.

Dans le cas d'une vari\'{e}t\'{e} symplectique $(M,\omega)$, on peut associer \`{a}
tout star-produit $*$ sa
{\em classe de Deligne} $[*]$, c.-\`{a}-d. une s\'{e}rie de Laurent formelle \`{a}
coefficients
dans le deuxi\`{e}me groupe de cohomologie de de Rham, $H^2_{dR}(M,\korps)$
d\'{e}finie
comme suit (voir \cite{Del95}, \cite{NT95a},\cite{NT95b}, \cite{BCG97},
et surtout le rapport \cite{GR99}): soit
$\big(U_\alpha,(q,p)\big)_{\alpha\in\mathfrak{S}}$ un
recouvrement ouvert de $M$ par des domaines de cartes de Darboux tel
que toutes les intersections finies sont contractiles. Alors
le champ d'Euler local
$\xi_\alpha:=\sum_{k=1}^n p_k\frac{\partial}{\partial p_k}$ est conform\'{e}ment
symplectique ($L_{\xi_\alpha}\omega=\omega$)
et il existe une s\'{e}rie d'op\'{e}rateurs diff\'{e}rentiels locaux $D'_\alpha$
(dont le terme d'ordre $0$ s'annule) telle que
\beq
    D_\alpha:= \nu \frac{\partial}{\partial \nu} + \xi_\alpha +
                   D'_\alpha
\end{equation}
est une d\'{e}rivation locale de $(\CinfK{U_\alpha}[[\nu]],*|_{U_\alpha})$. Sur
une intersection nonvide $U_\alpha\cap U_\beta$ la diff\'{e}rence
$D_\alpha-D_\beta$ est de la forme
$\frac{1}{\nu}ad_*(d_{\alpha\beta})$ o\`{u} $d_{\alpha\beta}\in
\CinfK{U_\alpha\cap U_\beta}[[\nu]]$ et $ad_*(f)\big(g\big):=f*g-g*f$.
Par cons\'{e}quent, la somme cyclique
$d_{\gamma\beta\alpha}:=d_{\alpha\gamma}+d_{\gamma\beta}+
d_{\beta\alpha}$ est un \'{e}l\'{e}ment de $\korps[[\nu]]$ et d\'{e}finit
un $2$-cocycle de \v{C}ech dont la classe dans
$H^2_{dR}(M,\korps)[[\nu]]$ est d\'{e}not\'{e}e par $d(*)$. La classe de Deligne $[*]$
est d\'{e}finie par
\beq \label{EqThClasSym}
 [*]=\frac{[\omega]}{\nu}+\sum_{r=0}\nu^r [*]_r,~~~~~
 [*]_0:=-[(\mathsf{C}^-_2)^\#],~~~\frac{\partial}{\partial \nu}[*](\nu)
  =\frac{1}{\nu^2}d(*)(\nu)
\end{equation}
o\`{u} la $2$-forme $(\mathsf{C}^-_2)^\#$ est d\'{e}finie par
$2(\mathsf{C}^-_2)^\#(X_f,X_g):=\mathsf{C}_2(f,g)-\mathsf{C}_2(g,f)$. Dans
les travaux cit\'{e}s ci-dessus
il a \'{e}t\'{e} montr\'{e} que deux star-produits sont \'{e}quivalents si et seulement si
leurs classes de Deligne sont \'{e}gales, et que toute s\'{e}rie formelle \`{a}
coefficients dans $H^2_{dR}(M,\korps)$ (avec le terme $\frac{[\omega]}{\nu}$
ajout\'{e}) est la classe de Deligne d'un star-produit. En outre, Neumaier
a montr\'{e} que la s\'{e}rie formelle de classes de de Rham figurant dans la
construction de Fedosov d'un star-produit
(avec le terme $\frac{[\omega]}{\nu}$ ajout\'{e})
co\"{\i}ncide avec la classe de Deligne du star-produit, voir
\cite{Neu02}. De plus, Bonneau donne une construction pour reconna\^{\i}tre les
classes de Fedosov dans les op\'{e}rateurs bidiff\'{e}rentiels du star-produit,
voir \cite{Bon98}. On remarque que nos conventions $\nu\hat{=}\frac{i\hbar}{2}$
sont telles que $*(\nu)=*_{\scriptstyle \mathrm{Neumaier}}(2\nu)
=*_{\scriptstyle \mathrm{Gutt-Rawnsley}}(-2\nu)$ et
$[*](\nu)=2[*]_{\scriptstyle \mathrm{Neumaier}}(2\nu)=
2[*]_{\scriptstyle \mathrm{Gutt-Rawnsley}}(-2\nu)$, voir \cite{Neu02},
\cite{GR99}.\\
La proposition suivante est bien connue et tr\`{e}s utile si deux star-produits
symplectiques co\"{\i}ncident jusqu'\`{a} un certain ordre:
\bprop\label{PDergrosseIrrtummitderexaktenEinsform}
 Soit $(M,\omega)$ une vari\'{e}t\'{e} symplectique. Soient
 $*=\sum_{r=0}^\infty \nu^r \mathsf{C}_r$ et
 $*'=\sum_{r=0}^\infty \nu^r \mathsf{C}'_r$ deux star-produits sur $M$ qui
 co\"{\i}ncident jusqu'\`{a} l'ordre $k\geq 1$, c.-\`{a}-d.
 $\mathsf{C}'_r=\mathsf{C}_r$
 quel que soit $0\leq r\leq k$. Alors:
 \ben
  \item $[*]_r=[*']_r$ quel que soit $-1\leq r\leq k-2$.
  \item il existent un op\'{e}rateur
 diff\'{e}rentiel $B_{k+1}\in\Dop^1\big(\CinfK{M};\CinfK{M}\big)$, une
 $2$-forme ferm\'{e}e $\zeta$ et un nombre rationnel non nul $s$ tels que
 \bea
   [*']_{k-1}-[*]_{k-1} & = & s[\zeta] \\
   \mathsf{C}'_{k+1}(f_1,f_2)-\mathsf{C}_{k+1}(f_1,f_2)
                        & = & (\mathsf{b}B_{k+1})(f_1,f_2)
                       +\zeta(X_{f_1},X_{f_2}) \nonumber \\
                       & & \label{EqOrdrekplusungeneral}
 \eea
 o\`{u} l'on rappelle {\em l'op\'{e}rateur
 cobord de Hochschild} $\mathsf{b}$ pour un op\'{e}rateur diff\'{e}rentiel
 $U:\CinfK{M}\ra\CinfK{M}$:
 \beq \label{EqCobordHochschild}
  (\mathsf{b}U)(f_1,f_2)~~:=~~f_1~U(f_2)~-~U(f_1f_2)~+~U(f_1)~f_2
 \end{equation}
 quels que soient $f_1,f_2\in\CinfK{M}$.\\
 En particulier, si $*$ et $*'$ sont \'{e}quivalents, $\zeta$ est une
 $2$-forme exacte (mais pas forc\'{e}ment nulle!)\footnote{Je dois cette remarque
 \`{a} Nikolai Neumaier. De plus, j'ai \'{e}t\'{e} soulag\'{e}
 d'apprendre que je n'\'{e}tais pas le seul d'avoir oubli\'{e} de tenir en compte
 cette $2$-forme exacte quand $*'\sim *$.}.
 \een
\eprop
\bbew
 L'associativit\'{e} de $*$ et de $*'$ \`{a} l'ordre $k+1$ montre que
 $\mathsf{C}'_{k+1}-\mathsf{C}_{k+1}$ est un $2$-cocycle de Hochschild,
 et l'identit\'{e} de l'associativit\'{e} antisym\'{e}tris\'{e}e \`{a} l'ordre $k+2$ donne
 la deuxi\`{e}me \'{e}quation (\ref{EqOrdrekplusungeneral}), voir par exemple
 \cite{GR99}, proposition 2.13. On calcule la classe relative de Deligne
 $t(*,*')$ de $*$ et $*'$ (voir par exemple \cite{GR99}, paragraphe 4),
 qui est \'{e}gale \`{a} $[*']-[*]$: soit $(U_\alpha)_{\alpha\in\mathfrak{S}}$
 un recouvrement ouvert de $M$ par des cartes d'un bon atlas. Puisque
 les $U_\alpha$ et toutes leurs intersections finies sont contractiles,
 il existe une transformation d'\'{e}quivalence $T_\alpha$ sur chaque
 $U_\alpha$ telle que $T_\alpha(*|_{U\alpha})=*'|_{U_\alpha}$. On peut
 choisir chaque $T_\alpha=id + \nu^k X_\alpha +\nu^{k+1}B_\alpha$
 o\`{u} $X_\alpha$ est un champ de vecteurs (donc $\mathsf{b}X_\alpha=0$)
 et $B_\alpha$ est une s\'{e}rie formelle d'op\'{e}rateurs diff\'{e}rentiels
 (voir la d\'{e}monstration de Proposition 3.1 de \cite{GR99}). Il s'ensuit
  que l'automorphisme $T_\beta^{-1}T_\alpha$ de $*|_{U_\alpha\cap
  U_\beta}$ est de la forme $\exp\big(ad_*(t_{\beta\alpha})\big)$
  avec $t_{\beta\alpha}=\nu^{k-1}t_{\beta\alpha}^{(k-1)}
  +\nu^{k}\rho_{\beta\alpha}$ avec $t_{\beta\alpha}^{(k-1)}
  \in\CinfK{U_\alpha\cap U_\beta}$ et $\rho_{\beta\alpha}\in
  \CinfK{U_\alpha\cap U_\beta}[[\nu]]$. Donc
  \[
     2\{t_{\beta\alpha}^{(k-1)},f\}=(X_\alpha-X_\beta)f,
  \]
  alors
  \[
     dt_{\beta\alpha}^{(k-1)}=\frac{1}{2}i_{(X_\beta-X_\alpha)}\omega.
  \]
  Soit $(\theta_\alpha)_{\alpha\in\mathfrak{S}}$ une partition de l'unit\'{e}
  subordonn\'{e}e \`{a} $(U_\alpha)_{\alpha\in\mathfrak{S}}$. Au cocycle de
  \v{C}ech $t_{\gamma\beta\alpha}=t_{\beta\alpha}+t_{\alpha\gamma}
  +t_{\gamma\beta}$ d\'{e}fini sur $U_\alpha\cap U_\beta \cap U_\gamma$
  on peut associer la $2$-forme ferm\'{e}e
  $\sum_{\gamma\beta\alpha}t_{\gamma\beta\alpha}\theta_\gamma
  d\theta_\beta\wedge d\theta_\alpha$ qui s'annule pour $0\leq r \leq k-2$
  et dont le coefficient d'ordre $k-1$ est cohomologue \`{a}
  \[
     -\sum_{\alpha\beta}\theta_\beta dt_{\beta\alpha}^{(k-1)}\wedge
                             d\theta_\alpha
                             \sim
                             \frac{1}{2}\sum_{\alpha}i_{X_\alpha}\omega\wedge
                             d\theta_\alpha.
  \]
  Alors les coefficients de la classe relative s'annule pour $-1\leq r\leq
  k-2$.
  D'un autre c\^{o}t\'{e}, puisque $\mathsf{C}_1=\mathsf{C}'_1=\mathsf{b}B_1+P$
  (o\`{u} $B_1$ est un op\'{e}rateur diff\'{e}rentiel et $P$ est la structure de
  Poisson associ\'{e}e \`{a} $\omega$) pour tout $\alpha$
  il existe un op\'{e}rateur diff\'{e}rentiel $\tilde{B}_\alpha$ avec
  \beq \label{EqCprimekplusunmoinsCkplusun}
   \big(\mathsf{C}'_{k+1}-\mathsf{C}_{k+1}\big)|_{U_\alpha} =
          -\mathsf{b}\tilde{B}_\alpha + [X_\alpha,P]
  \end{equation}
  avec
  \[
    [X_\alpha,P](f,g) := X_\alpha\{f,g\}-\{X_\alpha f, g\} -\{f,X_\alpha g\}
                      = -\big(di_{X_\alpha}\omega\big)(X_f,X_g).
  \]
  Quand on multiplie les termes dans l'\'{e}qn (\ref{EqCprimekplusunmoinsCkplusun})
  par $\theta_\alpha$ et regarde la somme sur $\alpha$ on obtient le
  r\'{e}sultat.
\ebew

Dans le cas d'une vari\'{e}t\'{e} de Poisson, Kontsevitch a class\'{e} les classes
d'\'{e}qui\-va\-lence des star-produits par des classes de diff\'{e}omorphie formelle
des structures de Poisson, voir \cite{Kon97b}.

Soient $(M,P)$ et $(M',P')$ deux vari\'{e}t\'{e}s de Poisson munies des
star-pro\-duits $*=\sum_{r=0}^\infty\nu^r\mathsf{C}_r$ et
$*'=\sum_{r=0}^\infty\nu^r\mathsf{C}'_r$, respectivement. Pour deux entiers
positifs
$s,t$ on d\'{e}finit dans une carte $(U\times U',x^1,\ldots,x^n,y^1,\ldots,
 y^{n'})$ l'op\'{e}rateur bidiff\'{e}rentiel suivant (o\`{u} l'entier $N$ est le
 maximum des entiers $N_1$ (correspondant \`{a} $\mathsf{C}_s$) et $N_2$
 (correspondant \`{a} $\mathsf{C}'_t$)):
\bea
\lefteqn{\big(\mathsf{C}_s\otimes \mathsf{C}'_t\big)(F,G)(m,m')  := }
               \nonumber \\
    & & \sum_{\stackrel{I_1,I_2,J_1,J_2}{|I_1|,|I_2|,|J_1|,|J_2|\leq N}}
      \mathsf{C}_s^{I_1I_2}(m){C'}_t^{J_1J_2}(m')
      \big(\partial_{x^{I_1}}\partial_{y^{J_1}}F\big)(m,m')
      \big(\partial_{x^{I_2}}\partial_{y^{J_2}}G\big)(m,m')\nonumber \\
      & &
\eea
quelles que soient $F,G\in \Cinf(M\times M',\korps)$. La d\'{e}finition ne
d\'{e}pend pas des cartes choisies. On pose
\beq
      *\otimes *'
      :=\sum_{r=0}^\infty\nu^r\sum_{s+t=r}\mathsf{C}_s\otimes \mathsf{C}'_t
\end{equation}
ce qui d\'{e}finit \'{e}videmment un star-produit sur $(M\times M',
P_{(1)}+P'_{(2)})$ qu'on va appeler le {\em produit tensoriel de $*$ et
$*'$}.
On a la
\bprop \label{PClasseDeligneProduit Tensoriel}
 Soient $(M,\omega)$ et $(M',\omega')$ deux vari\'{e}t\'{e}s symplectiques. Soit
 $*$ un star-produit sur la vari\'{e}t\'{e} symplectique $(M\times
 M',\omega_{(1)}+\omega'_{(2)})$.
 Soient $pr_1$ et $pr_2$ les projections canoniques $M\times M'\ra M$
 et $M\times M'\ra M'$.
 Alors:
 \ben
  \item Soient $\alpha\in H^2_{dR}(M,\korps)[[\nu]]$ et
     $\alpha'\in H^2_{dR}(M',\korps)[[\nu]]$ telles que
     $[*]=pr_1^*\big(\frac{[\omega]}{\nu}+\alpha\big)
     +pr_2^*\big(\frac{[\omega']}{\nu}+\alpha'\big)$. Alors il existe un
     star-produit
     $*_1$ sur $(M,\omega)$ avec
     $[*_1]=\frac{[\omega]}{\nu}+\alpha$ et un star-produit $*_2$ sur
     $(M,\omega')$ avec $[*_2]=\frac{[\omega']}{\nu}+\alpha'$ tels que
     $*~\sim~*_1\otimes *_2$.
  \item Quels que soient les star-produits $*_1$ sur $(M,\omega)$ et
   $*_2$ sur $(M',\omega')$:
   \beq\label{EqFormuleClassesDeligneTensoriel}
      [*_1\otimes *_2]=pr_1^*[*_1]+pr_2^*[*_2].
   \end{equation}
 \een
\eprop
\bbew
 1. Puisque tout star-produit
 est \'{e}quivalent \`{a} un star-produit de Fedosov on peut supposer que
 $*$ soit de Fedosov et est construit sur $M\times M'$ \`{a} l'aide de deux
 connexions symplectiques $\nabla$ (dans $TM$) et $\nabla'$ (dans $TM'$)
 et de deux s\'{e}ries de $2$-formes ferm\'{e}es dont l'une est d\'{e}finie sur $M$
 et repr\'{e}sente $\alpha$ et l'autre est d\'{e}finie sur $M'$ et repr\'{e}sente
 $\alpha'$ d'apr\`{e}s le travail de Neumaier \cite{Neu02}. On voit facilement
 que la construction de Fedosov sur $M\times
 M'$ `factorise' enti\`{e}rement et donne deux star-produits, $*_1$ sur $M$ et
 $*_2$ sur $M'$, tels que $*=*_1\otimes *_2$ avec
 $[*_1]=\frac{[\omega]}{\nu}+\alpha$ et
 $[*_2]=\frac{[\omega']}{\nu}+\alpha'$.\\
 2.
 Pour toute carte de Darboux
 $\big(U_\alpha,(q,p)\big)$ de $M$ (resp. $\big(V_{\alpha'},(q',p')\big)$ de
 $M'$) soient
 $D_\alpha:=\nu\frac{\partial}{\partial \nu} + \xi_\alpha +
 E_\alpha$ (resp.
 $D_{\alpha'}:=\nu\frac{\partial}{\partial \nu} + \xi_{\alpha'} +
 E_{\alpha'}$) des $\nu$-Euler d\'{e}rivations locales pour $*_1$ (resp.
 $*_2$). Alors on voit facilement que
 \[
   D_{(\alpha,\alpha')}:=\nu\frac{\partial}{\partial \nu}
                        + \xi_\alpha + \xi_{\alpha'}
                        + E_\alpha + E_{\alpha'}
 \]
 est une $\nu$-Euler d\'{e}rivation locale pour $*_1\otimes *_2$ dans la
 carte de Darboux produit
 $\big(U_\alpha\times V_{\alpha'},((q,p),(q',p'))\big)$. Par cons\'{e}quent,
 la construction
 de la classe de Deligne de $*_1\otimes *_2$ \`{a} partir des $\nu$-Euler
 d\'{e}rivations locales donn\'{e}e par exemple dans \cite{GR99} `factorise'
 \'{e}galement (ainsi que la classe de la partie antisym\'{e}trique du terme d'ordre
 $2$ du star-produit $*_1\otimes *_2$ car on peut supposer apr\`{e}s une
 \'{e}ventuelle transformation d'\'{e}quivalence que le terme d'ordre $1$ de
 $*_1$ et de $*_2$ soit antisym\'{e}trique)
 pour donner $[*_1\otimes *_2]=pr_1^*[*_1]+pr_2^*[*_2]$.
\ebew

En outre, pour un star-produit $*$ sur une vari\'{e}t\'{e} de Poisson $(M,P)$ donn\'{e}e,
on d\'{e}finit le {\em star-produit oppos\'{e} de $*$}, $*^{\mathrm{opp}}$ par
\beq
   f*^{\mathrm{opp}}g:=g*f~~~~\mathrm{quelles~que~soient~}f,g\in
   \CinfK{M}[[\nu]]
\end{equation}
qui est \'{e}videmment un star-produit pour la vari\'{e}t\'{e} de Poisson $(M,-P)$.
Dans le cas d'une vari\'{e}t\'{e} symplectique $(M,\omega)$ il vient
que
\beq\label{EqClasseDeligneOpp}
    [*^{\mathrm{opp}}]=-[*],
\end{equation}
voir \cite{Neu01}, \cite{Neu02}.

Outre les exemples bien connus des star-produits sur $\real^{2n}$ muni
de la forme symplectique canonique (voir par exemple \cite{BFFLS78}),
le star-produit $*_G$ de Gutt sur l'espace
dual $(\mathfrak{g}^*,P_\mathfrak{g})$ d'une alg\`{e}bre de Lie r\'{e}elle de
dimension finie
$\mathfrak{g}$ est tr\`{e}s important (voir \'{e}galement \cite{Gut83}):
Soit
$\tilde{\mathfrak{g}}:=\mathfrak{g}\otimes_{\real}\korps[[\nu]]$ muni du
crochet
$2\nu[~,~]$. Alors $\tilde{\mathfrak{g}}$ est un
$\korps[[\nu]]$-module libre. Soit $U\tilde{\mathfrak{g}}$ l'alg\`{e}bre
enveloppante de $\tilde{\mathfrak{g}}$ (voir par exemple \cite[p.266]{CE56}).
Le
th\`{e}or\`{e}me de Poincar\'{e}-Birkhoff-Witt (voir par exemple \cite[p.271]{CE56})
entra\^{\i}ne que
$U\tilde{\mathfrak{g}}$ est un $\korps[[\nu]]$-module libre isomorphe
\`{a} l'alg\`{e}bre sym\'{e}trique
$S\tilde{\mathfrak{g}}\cong (S\mathfrak{g})\otimes_{\real}\korps[[\nu]]$
\`{a} l'aide de l'application de sym\'{e}trisation $S\tilde{\mathfrak{g}}\ra
U\tilde{\mathfrak{g}}$. De cette mani\`{e}re, $S\tilde{\mathfrak{g}}$ est
munie d'une multiplication $\korps[[\nu]]$-bilin\'{e}aire $*_G$, isomorphe \`{a} celle
de l'alg\`{e}bre enveloppante, qui se prolonge sur son compl\'{e}t\'{e}
par rapport \`{a} la topologie $\nu$-adique,
$S\big(\mathfrak{g}\otimes_\real \korps\big)[[\nu]]$. Ce dernier espace
s'identifie avec l'espace $\mathcal{P}(\mathfrak{g}^*)[[\nu]]$
de toutes les s\'{e}ries formelles \`{a} coefficients dans
les fonctions polyn\^{o}miales sur $\mathfrak{g}^*$ \`{a} valeurs dans $\korps$.
S.Gutt a montr\'{e} dans \cite{Gut83} que $*_G$ est un star-produit sur
$\mathfrak{g}^*$ (i.e. $*_G$ se prolonge \`{a}
$\CinfK{\mathfrak{g}^*}[[\nu]]$), et on peut le d\'{e}finir \`{a} l'aide des
fonctions exponentielles $e_\xi:\mathfrak{g}^*\ra \korps$ donn\'{e}es pour
tous les \'{e}l\'{e}ments $\xi$ de $\mathfrak{g}\otimes_\real \korps$ par
$e_\xi(\alpha):=e^{\alpha(\xi)}$:
\beq\label{EqStarProduitGutt}
   e_\xi *_G e_\eta:=e_{H(\xi,\eta)}
\end{equation}
o\`{u} $H$ d\'{e}signe la s\'{e}rie de Baker-Campbell-Hausdorff de
$\tilde{\mathfrak{g}}$, voir par exemple \cite[p.40]{KMS93}. $*_G$
s'obtient \'{e}galement comme la restriction d'un certain star-produit bi-invariant
sur la vari\'{e}t\'{e} symplectique $T^*G$ aux fonctions invariantes \`{a} gauche,
voir \cite{Gut83}, \cite{BNW98}. Les deux identit\'{e}s suivantes sont utiles:
\bea
  \xi*_G \eta -\eta *_G \xi & = & 2\nu [\xi,\eta] ~~~~\forall \xi,\eta\in
                                              \tilde{\mathfrak{g}},
                     \label{EGuttCommutateur}\\
  \xi*_G\cdots *_G \xi ~(k~\mathrm{facteurs}~) & = & \xi^k~~~~
                                                \forall \xi \in
                                              \tilde{\mathfrak{g}}.
                     \label{EGuttPuissance}
\eea

\section{Homomorphismes, re\-pr\'{e}\-sen\-ta\-tions et r\'{e}duc\-tion}
  \label{SecHomoRepRed}

\subsection{Homomorphismes de star-produits et applications de Poisson}
          \label{SubSecMor1}

\bdefi \label{DDefHomStarGen}
 Soient $(M,P)$ et $(M',P')$ deux vari\'{e}t\'{e}s de Poisson munies des
 star-produits $*$ et $*'$, respectivement.\\
 Une application $\korps[[\nu]]$-lin\'{e}aire
$\Phi:\CinfK{M}[[\nu]]\ra\CinfK{M'}[[\nu]]$ est appel\'{e}e
 {\em homomorphisme de star-produits} ssi $\Phi$ est un homorphisme
 d'alg\`{e}bres associatives unitaires sur $\korps[[\nu]]$:
\[
   \Phi(f*g)=\big(\Phi(f)\big)*'\big(\Phi(g)\big)
   ~~~\mathrm{et}~~~\Phi(1)=1
\]
 quels que soient $f,g\in \CinfK{M}[[\nu]]$.
\edefi

\noindent Le lien avec les applications de Poisson est contenu dans le
lemme suivant:

\blem \label{LHomoEgalPoissonOrdreZero}
 Soit $\Phi=\sum_{r=0}^\infty \nu^r \Phi_r:\CinfK{M}[[\nu]]
 \ra\CinfK{M'}[[\nu]]$ un homomorphisme de
 star-produits.\\
 Alors il existe une application de Poisson $\phi:(M',P')\ra (M,P)$
 telle que $\Phi_0(f)=\phi^*f:=f\circ\phi$ quelle que soit
 $f\in\CinfK{M}$.
\elem
\bbew
  Soient $f,g\in\CinfK{M}$.
  La propri\'{e}t\'{e} d'homomorphisme de $\Phi$ s'\'{e}crit \`{a} l'ordre
  $0$ de $\nu$:
  \[
    \Phi_0(fg)=\big(\Phi_0(f)\big)\big(\Phi_0(g)\big).
  \]
  Alors $\Phi_0$ est un homomorphisme d'alg\`{e}bres commutatives associatives
  unitaires $\CinfK{M}\ra \CinfK{M'}$. D'apr\`{e}s l'exercice
  de Milnor (voir \cite{KMS93}, p. 301, Cor. 35.9) il existe une application
  de classe
  $\Cinf$ $\phi:M'\ra M$ telle que  $\Phi_0(f)=\phi^*f$. Ensuite,
  le commutateur
  de la propri\'{e}t\'{e} d'homomorphismes \`{a} l'ordre $1$ de $\nu$
  s'\'{e}crit
  \[
    \Phi_0\{f,g\}_P= \{\Phi_0(f),\Phi_0(g)\}_{P'}
  \]
  d'o\`{u} le fait que $\phi$ est une application de Poisson.
\ebew

Puisque $\Phi_0$ est donn\'{e}e par une application $\phi$ de classe $\Cinf$ on
peut restreindre la d\'{e}finition \ref{DDefHomStarGen} de la mani\`{e}re
suivante:
\bdefi \label{DDefHomStarGenDiff}
 Un homomorphisme de star-produits $\Phi=\phi^*+\sum_{r=1}^\infty \nu^r
 \Phi_r:
 \CinfK{M}[[\nu]]\ra\CinfK{M'}[[\nu]]$ est dit
 {\em diff\'{e}rentiel} lorsque toutes les applications $\Phi_r$ sont des
 op\'{e}rateurs
 diff\'{e}rentiels le long de $\phi$, i.e. ils apartiennent \`{a} l'espace
 $\Dop^1(\CinfK{M};\CinfK{M'})$.
\edefi
Une transformation d'\'{e}quivalence $S$ d'un star-produit $*$ \`{a} un
star-produit $*'$ est \'{e}videmment un homomorphisme de star-produit
diff\'{e}rentiel qui d\'{e}forme l'application identique.

Il me semble important de savoir si tout homomorphisme de star-produits
est diff\'{e}rentiel.

La question r\'{e}ciproque de savoir quand une application de Poisson donn\'{e}e
se d\'{e}forme dans un homomorphisme de star-produits est sans doute
int\'{e}res\-sante:

\bprob[Quantification des applications de Poisson] \label{prob1}
 Quelles \\
 sont les conditions sur une application de Poisson
 $\phi:(M',P')\ra (M,P)$ pour qu'elle soit
 {\em (diff\'{e}rentiellement) quantifiable}, i.e. pour
 qu'il existent des star-produits
 $*'$ et $*$ sur les vari\'{e}t\'{e}s de Poisson $(M',P')$ et
 $(M,P)$, respectivement, et des applications lin\'{e}aires
 $\Phi_1,\Phi_2,\ldots:\CinfK{M}\ra \CinfK{M'}$
 tels que
 \[
    \Phi:= \phi^* + \sum_{r=1}^\infty\nu^r\Phi_r
 \]
 soit un homomorphisme de star-produits (diff\'{e}rentiel)?
\eprob

Un premier cas particulier tr\`{e}s simple, mais important de cette situation
est le cas d'une partie
ouverte $U$ de la vari\'{e}t\'{e} de Poisson $(M,P)$: l'injection canonique
$i_U:U\ra M$ est visiblement une application de Poisson $(U,P|_U)\ra
(M,P)$. On voit facilement que la restriction
$i_U^*:\CinfK{M}\ra\CinfK{U}$ d\'{e}finit directement un
homomorphisme de star-produits diff\'{e}rentiels de l'alg\`{e}\-bre
$\big(\CinfK{M}[[\nu]],*\big)$ \`{a} $\big(\CinfK{U}[[\nu]],*|_U\big)$.
\bprop \label{PQuantRestLocal}
 Soit $*$ un star-produit sur la vari\'{e}t\'{e} de Poisson $(M,P)$ et $*_U$
 un star-produit sur la partie ouverte $U$ de $M$ munie de la restriction
 de $P$ \`{a} $U$. Alors $i_U:U\ra M$ est diff\'{e}rentiellement
 quantifiable si et seulement si
 la restriction de $*$ \`{a} $U$, $*|_U$, est \'{e}quivalente \`{a} $*_U$ sur $U$.
\eprop
\bbew
 Si $i_U$ est diff\'{e}rentiellement quantifiable, il existe alors une s\'{e}rie
 d'op\'{e}rateurs diff\'{e}rentiels $S_U=id_U+\sum_{r=1}^\infty\nu^r S_r$
 avec $S_r\in \Dop^1\big(\CinfK{U},\CinfK{U}\big)$ s'annulant sur les
 constantes telle que  $\Phi:=S_U\circ (i_U)^*$ est un homomorphisme
 de star-produits. Il vient
 \[
   \big(S_U(i_U^*f)\big)*_U  \big(S_U(i_U^*g)\big)
     = \big(S_U(i_U^*(f*g))\big)
     = S_U\big((i_U^*f)*|_U (i_U^*g)\big).
 \]
 Par cons\'{e}quent, $*_U=S_U(*|_U)$: en fait, il suffit de v\'{e}rifier cette \'{e}quation
 localement et --en utilisant une
 partition de l'unit\'{e}-- pour toutes les fonctions $\phi,\psi\in\CinfK{U}$
 et pour tout point $x\in U$ il existe un voisinage ouvert $V$ de $x$
 (dont l'adh\'{e}rence fait partie de $U$) et des fonctions $f,g\in\Cinf{M}$
 tels que $f(y)=\phi(y)$ et $g(y)=\psi(Y)$ quel que soit $y\in V$.\\
 La r\'{e}ciproque est \'{e}vidente: s'il existe une transformation d'\'{e}quivalence
 $S_U$ avec $*_U=S_U(*|_U)$, alors $\Phi:=S_U\circ (i_U)^*$ est un
 homomorphisme de star-produits.
\ebew

Un deuxi\`{e}me cas particulier du probl\`{e}me \ref{prob1} est l'action \`{a} gauche
$\phi:G\times M\ra M:(g,m)\mapsto \phi_g(m)$ d'un groupe de Lie $G$ sur une
vari\'{e}t\'{e} de Poisson
$(M,P)$ en tant qu'isomorphismes de Poisson, i.e. $\phi_g^*P=P$ quel que soit
$g\in G$: le probl\`{e}me
des {\em star-produits $G$-invariants $*$} est la question de savoir si les
applications de Poisson $\Phi_g$ se quantifient de telle fa\c{c}on qu'on
obtienne une repr\'{e}sentation de $G$ en tant qu'automorphismes de l'alg\`{e}bre
d\'{e}form\'{e}e $\big(\CinfK{M}[[\nu]],*\big)$.

Un troisi\`{e}me cas particulier important est donn\'{e} par les applications
moments $J:M\ra \mathfrak{g}^*$ (voir eqn (\ref{Eqapplmoment}) du paragraphe
\ref{SubSecPois}). Le probl\`{e}me de la quantification de cette application
de Poisson est li\'{e} \`{a} la situation suivante: \\
Une s\'{e}rie formelle de fonctions $\mathbb{J}\in
\mathfrak{g}^*\otimes\CinfK{M}[[\nu]]$ satisfaisant
\beq \label{DAplMomQ}
  \langle \mathbb{J},\xi\rangle*\langle \mathbb{J},\eta\rangle
  -\langle \mathbb{J},\eta\rangle*\langle \mathbb{J},\xi\rangle
            =   2\nu \langle \mathbb{J},[\xi,\eta]\rangle.
\end{equation}
quels que soient $\xi,\eta\in\mathfrak{g}$
s'appelle une {\em application moment quantique} d'apr\`{e}s Xu \cite{Xu98}.
Pour un star-produit $*$ donn\'{e}, on voit facilement que les obstructions
r\'{e}currentes \`{a} l'existence d'une application moment quantique se trouvent dans
le deuxi\`{e}me groupe de cohomologie de Chevalley-Eilenberg
$H^2_{CE}\big(\mathfrak{g},\CinfK{M}\big)$ de l'alg\`{e}bre de Lie
$\mathfrak{g}$ qui agit sur $\CinfK{M}$ via
$\xi\mapsto -X_{\langle J,\xi\rangle}$.
Si $\mathbb{J}=J=\mathbb{J}_0$ le star-produit
$*$ s'appelle {\em $\mathfrak{g}$-covariant} d'apr\`{e}s Arnal, Cortet, Molin et
Pinczon \cite{ACMP83}.
Plus particuli\`{e}rement, un star-produit $*$ s'appelle {\em fortement
$\mathfrak{g}$-invariant} d'apr\`{e}s ces auteurs si
\beq
 \langle \mathbb{J},\xi\rangle*f
  -f*\langle \mathbb{J},\xi\rangle
            =   2\nu \{ \langle \mathbb{J},\xi\rangle,f\}
\end{equation}
quelle que soit la fonction $f\in\Cinf(M,\korps)[[\nu]$. On a le
crit\`{e}re suffisant suivant pour l'existence de ces star-produits:
\bsat[Fedosov,1996]
 Soit $(M,\omega)$ une vari\'{e}t\'{e} symplectique, $J:M\ra\mathfrak{g}^*$
une application
 moment et $\nabla$ une connection dans le fibr\'{e} tangent telle que
 \[
   0=\big(L_{X_{\langle J,\xi\rangle}}\nabla\big)_XY
         :=[X_{\langle J,\xi\rangle},\nabla_XY]
                -\nabla_{[X_{\langle J,\xi\rangle},X]}Y-
                    \nabla_X[X_{\langle J,\xi\rangle},Y].
 \]
 quel que soit $\xi$ dans l'alg\`{e}bre de Lie $\mathfrak{g}$.
 Alors il existe un star-produit fortement $\mathfrak{g}$-invariant $*$
 avec $\mathbb{J}=J$.
\esat
Voir \cite[Sec.5.8]{Fed96} pour la d\'{e}monstration.
Par exemple, si les flots des champs de vecteurs $X_{\langle J,\xi\rangle}$
d\'{e}finissent
l'action d'un groupe de Lie compacte ou plus g\'{e}n\'{e}ralement une action
propre d'un groupe de Lie, il r\'{e}sulte d'un th\'{e}or\`{e}me classique de R.Palais
(voir \cite[p.316]{Pal61}) que ces champs
de vecteurs pr\'{e}servent une m\'{e}trique riemannienne sur $M$, alors
sa connection Levi-Civita, et le th\'{e}or\`{e}me de Fedosov est applicable.
\bprop \label{PMorphApplMomentQ}
 Soit $J:M\ra\mathfrak{g}^*$ une application moment.
 \ben
  \item
    Si $J$ --vue comme application de Poisson--
    admet une quantification $\mathbf{J}$ de l'alg\`{e}bre
    $\CinfK{\mathfrak{g}^*}[[\nu]]$ munie du star-produit de Gutt $*_G$,
    voir l'\'{e}q. (\ref{EqStarProduitGutt}),
    dans l'alg\`{e}bre d\'{e}form\'{e}e $\big(\CinfK{M}[[\nu]],*\big)$,
    alors l'application
    $\xi\mapsto\langle \mathbb{J},\xi\rangle:=\mathbf{J}(\xi)$ d\'{e}finit une
    application moment quantique.
  \item
    Si $J$ est le terme d'ordre $0$ d'une application moment quantique
    $\mathbb{J}$, alors il existe un morphisme d'alg\`{e}bres associatives
    unitaires sur $\korps[[\nu]]$, $\mathbf{J}$, de l'alg\`{e}bre
    $\big(\mathcal{P}(\mathfrak {g}^*)[[\nu]],*_G\big)$ (les s\'{e}ries formelles
    \`{a} coefficients dans l'espace des fonctions polyn\^{o}miales sur
    $\mathfrak{g}^*$) dans l'alg\`{e}bre d\'{e}form\'{e}e
    $\big(\CinfK{M}[[\nu]],*\big)$.
 \een
\eprop
\bbew
 1. Soient $\xi,\eta\in\mathfrak{g}$. Alors $\xi *_G\eta-\eta *_G \xi
 = 2\nu [\xi,\eta]$ (voir l'\'{e}q. (\ref{EGuttCommutateur})), donc $\mathbb{J}$
 est une application moment quantique.\\
 2. Une application moment quantique est un morphisme d'alg\`{e}bres de Lie
  $\tilde{\mathfrak{g}}\ra \CinfK{M}[[\nu]]$ o\`{u} $\tilde{\mathfrak{g}}$
  est l'espace $\mathfrak{g}\otimes_\real \korps[[\nu]]$ muni du crochet
  $2\nu[~,~]$ et $\CinfK{M}[[\nu]]$ est vu comme alg\`{e}bre de Lie muni du
  commutateur du star-produit $*$. Ceci induit un morphisme d'alg\`{e}bres
  associatives de l'alg\`{e}bre enveloppante $U\tilde{\mathfrak{g}}$ dans
  l'alg\`{e}bre associative $\big(\CinfK{M}[[\nu]],*\big)$ qui se prolonge sur
  le compl\'{e}t\'{e} (par rapport \`{a} la topologie $\nu$-adique)
  $\big(\mathcal{P}(\mathfrak {g}^*)[[\nu]],*_G\big)$ de
  $U\tilde{\mathfrak{g}}$, voir aussi le paragraphe \ref{SubSecStarprod}.
\ebew
\brem
 Je ne sais pas si on peut en g\'{e}n\'{e}ral prolonger le morphisme
d'alg\`{e}bres
associatives de l'\'{e}nonc\'{e} $2.$ \`{a} $\CinfK{\mathfrak{g}^*}[[\nu]]$ sans
imposer des conditions additionnelles \`{a} l'application moment quantique.
\erem

Puisque les syst\`{e}mes hamiltoniens int\'{e}grables constituent une
sous-classe de la classe des applications moment,
on peut appeler un syst\`{e}me hamiltonien
$(M,\omega,H)$ un {\em syst\`{e}me int\'{e}grable quantique} lorsqu'il y a
une application $\mathbb{J}=\sum_{r=0}^\infty\nu^r\mathbb{J}_r~\in
\real^n\otimes\CinfK{M}[[\nu]]$ telle que
$\big(M,\omega,H,\mathbb{J}_0=:J=:(J_1,\ldots,J_n)\big)$ soit un syst\`{e}me
int\'{e}grable classique et
\beq
   \mathbb{J}_k*\mathbb{J}_l-\mathbb{J}_l*\mathbb{J}_k=0~~~\forall
         1\leq k,l\leq n
\end{equation}
quels que soient $1\leq k,l\leq n$. La plupart des syst\`{e}mes
hamiltoniens int\'{e}grables connus sont \'{e}galement int\'{e}grables quantiques
(voir par exemple \cite{BW99}).

Un quatri\`{e}me cas particulier est donn\'{e} par une application $\phi$ de classe
$\Cinf$ d'une vari\'{e}t\'{e} diff\'{e}rentiable $M$ dans une vari\'{e}t\'{e} de Poisson
$(M',P')$ munie d'un star-produit $*'=\sum_{r=0}\nu^r\mathsf{C}'_r$:
si $\phi$ est une subimmersion,
i.e. un diff\'{e}omorphisme local, on peut d\'{e}finir la {\em structure de Poisson
retir\'{e}e} $P:=\phi^*P'$ sur $M$ par
\[
   P_m:=(T_m\phi)^{-1}\otimes (T_m\phi)^{-1}\big(P_{\phi(m)}\big)
\]
(quel que soit $m\in M$),
et puisque les fibr\'{e}s en jets $J^N(M',\real)_0$ et leur fibr\'{e}s duaux
$\mathrm{Hom}\big(J^N(M',\real)_0,\korps\big)$ (voir par exemple
\cite[p.117-119]{KMS93} pour la d\'{e}finition)
sont des fibr\'{e}s vectoriels naturels (voir par exemple \cite[p.138]{KMS93})
on peut \'{e}galement construire {\em le star-produit retir\'{e} $\phi^*(*')$} en
retirant tous les op\'{e}rateurs bidiff\'{e}rentiels $\mathsf{C}'_r$ qui sont des
sections de classe $\Cinf$ dans
$\Ginf\big(\mathrm{Hom}\big(J^{N_r}(M',\real)_0,\korps\big)\otimes
            \mathrm{Hom}\big(J^{N_r}(M',\real)_0,\korps\big)\big)$
(pour un certain entier positif $N_r$) par
\beq \label{EStarProduitRetire}
    (\phi^*\mathsf{C}'_r)_m:=(J^{N_r}\phi)_m^{-1}\otimes (J^{N_r}\phi)_m^{-1}
                               \big(\mathsf{C}'_{r~\phi(m)}\big)
    ~~\mathrm{et}~~\phi^*(*'):=\sum_{r=0}^\infty\nu^r\phi^*\mathsf{C}'_r
\end{equation}
o\`{u} $m\in M$ et $(J^{N_r}\phi)_m$ d\'{e}signe le jet d'ordre $N_r$ de $\phi$ en $m$
(qui est inversible car $T_m\phi$ l'est). Il s'ensuit que
\[
    (\phi^*\mathsf{C}'_r)(\phi^*f',\phi^*g')
          =\phi^*\big(\mathsf{C_r}'(f',g')\big)~~~\forall~f',g'\in
                     \CinfK{M'}[[\nu]],
\]
ce qui permet \'{e}galement de d\'{e}finir $\phi^*\mathsf{C}_r$ dans le cas o\`{u}
$\phi$ est un diff\'{e}omor\-phis\-me global, et il vient que le star-produit
retir\'{e} $*:=\phi^*(*')$ est un star-produit de la vari\'{e}t\'{e} de Poisson
$(M,P)$ et que $\Phi:=\phi^*$ est une quantification de l'application
de Poisson $\phi$.

\subsection{Repr\'{e}sentations des star-produits et applications co\"{\i}sotropes}
          \label{SubSecRep1}

Soit $(M,P)$ une vari\'{e}t\'{e} de Poisson, $C$ une vari\'{e}t\'{e} diff\'{e}rentiable et
$\zeta:E\ra C$ un fibr\'{e} vectoriel (sur $\korps$) sur $C$.
On consid\`{e}re l'espace  de tous les
op\'{e}rateurs diff\'{e}rentiels $\Dop^1\big(\Ginf(C,E);\Ginf(C,E)\big)$ sur
l'espace des sections $\Ginf(C,E)$. Ceci est une alg\`{e}bre
associative unitaire
munie de la multiplication usuelle d'op\'{e}rateurs diff\'{e}rentiels. Il en est
de m\^{e}me avec l'espace de toutes les s\'{e}ries formelles
$\Dop^1\big(\Ginf(C,E);\Ginf(C,E)\big)[[\nu]]$.
Soit $*$ un star-produit sur $M$.

\bdefi \label{DefRep1}
  Une {\em re\-pr\'{e}\-sen\-ta\-tion du star-produit $*$ dans $E$} est un
  homomorphisme d'alg\`{e}bres associatives unitaires
  \[
  \rho:\CinfK{M}[[\nu]]\ra \Dop^1\big(\Ginf(C,E);\Ginf(C,E)\big)[[\nu]]
  \]
  sur $\korps[[\nu]]$. Dans ce cas-l\`{a}, l'espace $\Ginf(C,E)[[\nu]]$
  est appel\'{e} un {\em module (\`{a} gauche) de} $*$.
\edefi
Dans le cas particulier tr\`{e}s important $E=C\times \korps$, on parlera
\'{e}galement d'une {\em repr\'{e}sentation de $*$ sur $C$}.
Il est clair qu'on obtient des {\em modules \`{a} droite} en rempla\c{c}ant
$*$ par le star-produit oppos\'{e} $*^\mathrm{opp}$ dans la d\'{e}finition
pr\'{e}c\'{e}dente. En outre, dans le cas de la donn\'{e}e de deux vari\'{e}t\'{e}s de
Poisson $(M_1,P_1)$ et $(M_2,P_2)$
munies des star-produits $*_1$ et $*_2$, respectivement, et d'une
repr\'{e}sentation du star-produit $*_1\otimes *_2^\mathrm{opp}$
(sur la vari\'{e}t\'{e}
de Poisson $(M_1\times M_2,P_{1~(1)}-P_{2~(2)})$) on parlera de
$\Ginf(C,E)[[\nu]]$ d'un {\em $*_1$-$*_2$-bimodule}. Dans ce cas, on
obtient \'{e}videmment par restriction une repr\'{e}sentation $\rho_1$ de $*_1$
dans $E$ et une repr\'{e}sentation $\rho_2$ de $*_2^\mathrm{opp}$ dans $E$
telles que
\beq \label{EVraiBimoduleQuant}
    \rho_1(f_1)\rho_2(f_2)=\rho_2(f_2)\rho_1(f_1)
     ~~~\forall~f_1\in\CinfK{M_1}[[\nu]]~~\forall~f_2\in\CinfK{M_2}[[\nu]].
\end{equation}
donc $\Ginf(C,E)[[\nu]]$ est un bimodule dans le sens alg\'{e}brique
par rapport aux anneaux
$\big(\CinfK{M_1}[[\nu]],*_1\big)$ et
$\big(\CinfK{M_2}[[\nu]],*_2\big)$.
\bexem
On obtient un exemple simple
en choisissant $M=C$ et $E=M\times\korps$ et en d\'{e}finissant
$\rho(f):g\mapsto f*g$.
\eexem

\noindent Pour les modules, il y a la fa\c{c}on usuelle d'isomorphie:
\bdefi \label{DefModEquiv}
 Soient $\rho$ une re\-pr\'{e}\-sen\-ta\-tion de $*$ sur $E$ et $\rho'$ une
 re\-pr\'{e}\-sen\-ta\-tion de $*$ sur $E'$. Les deux modules $\Ginf(C,E)[[\nu]]$
 et $\Ginf(C,E')[[\nu]]$ sont dits {\em isomorphes} lorsqu'il existe
 une s\'{e}rie d'op\'{e}rateurs diff\'{e}rentiels inversible $T$ dans
 $\mathbf{D}^1\big(\Ginf(C,E),\Ginf(C,E')\big)[[\nu]]$ telle que
 \[
      T\rho(f)=\rho'(f)T~~~\forall~f\in\CinfK{M}[[\nu]]
 \]
\edefi
On peut combiner les repr\'{e}sentations avec les morphismes de star-produits
de la mani\`{e}re suivante
\footnote{Cette id\'{e}e s'\'{e}tait form\'{e}e dans une des nombreuses conversations
fructueuses avec Stefan Waldmann}:

\bdefi \label{DefRepEquiv}
 Soient $(M,P)$ et $(M',P')$ deux vari\'{e}t\'{e}s de Poisson munies de deux
 star-produits $*$ et $*'$, respectivement. Soit $\phi:M\ra M'$ une application
 de Poisson et $\Phi$ une quantification de $\phi$ par rapport \`{a} $*$ et
 $*'$. Soient $\zeta:E\ra C$
 et $\zeta':E'\ra C'$ deux fibr\'{e}s vectoriels sur les vari\'{e}t\'{e}s
 diff\'{e}rentiables $C$ et $C'$, respectivement. {\em Une repr\'{e}sentation $\rho$
 de $*$ dans $E$ est dite li\'{e}e \`{a} une repr\'{e}sentation $\rho'$ de $*'$
 dans $E'$ par rapport \`{a} $\phi$}
 lorsqu'il existe une application $\korps[[\nu]]$-lin\'{e}aire
 $T:\Ginf(C',E')[[\nu]]\ra \Ginf(C,E)[[\nu]]$ telle que
 \beq\label{EqDeuxRepLiees}
    \rho\big(\Phi(f')\big)\big(T(\psi')\big)
          =
          T\big(\rho'(f')(\psi')\big)
 \end{equation}
 quels que soient $f'\in\CinfK{M'}[[\nu]]$ et $\psi'\in\Ginf(C',E')[[\nu]]$.
 Dans le cas particulier $(M',P')=(M,P)$, $\phi=id$ et $C'=C$
 on dit que les {\em couples $(*,\rho)$ et $(*',\rho')$ sont \'{e}quivalents}
 lorsque $\Phi$ est une transformation d'\'{e}quivalence
 et l'application $T$ est donn\'{e}e par une s\'{e}rie formelle
 $\sum_{r=0}\nu^r T_r$ o\`{u} $T_r$ appartient \`{a} l'espace d'op\'{e}rateurs
 diff\'{e}rentiels
 $\Dop^1\big(\Ginf(C,E');\Ginf(C,E)\big)$ quel que soit l'entier
 positif $r$.
% \bea
%     f*'g & = & S\big((S^{-1}f)*(S^{-1}g\big)
%            ~~~\mathrm{quels~que~soient~}f,g\in\CinfK{M}[[\nu]].\nonumber \\
%          &   &     \\
%    \rho'(f) & = &
%    T\rho(S^{-1}f)T^{-1}~~~\mathrm{quel~que~soit~}f\in\CinfK{M}[[\nu]].
% \eea
\edefi
%En particulier, les deux star-produits sont \'{e}quivalents si les couples
%$(*,\rho)$ et $(*',\rho')$ sont \'{e}quivalents.

Le lien entre les re\-pr\'{e}\-sen\-ta\-tions de star-produits et les
sous-vari\'{e}t\'{e}s
co\"{\i}so\-tropes est contenu dans la proposition suivante:

\bprop \label{PRepCois}
 Soit $(M,P)$ une vari\'{e}t\'{e} de Poisson munie d'un star-produit $*$, $C$
 une vari\'{e}t\'{e} diff\'{e}rentiable et
 \[
  \rho=\sum_{r=0}^\infty\nu^r\rho_r:
  \CinfK{M}[[\nu]]\ra \Dop^1(\CinfK{C};\CinfK{C})[[\nu]]
 \]
 une re\-pr\'{e}\-sen\-ta\-tion de star-produits. \\
 Alors il existe une application de classe $\Cinf$, $i:C\ra M$
 telle que $\rho_0(f)(\psi)=(i^*f)\psi:=(f\circ i)\psi$ quelles
 que soient $f\in\CinfK{M}$ et $\psi\in\CinfK{C}$.
 L'application $i$ est co\"{\i}sotrope.
 En particulier, au cas o\`{u} $i$ est un plongement sur une sous-vari\'{e}t\'{e} ferm\'{e}e
 $i(C)$ de $M$, alors $i(C)$ est une sous-vari\'{e}t\'{e} co\"{\i}sotrope de $(M,P)$.
\eprop
\bbew
 La propri\'{e}t\'{e} de re\-pr\'{e}\-sen\-ta\-tion s'\'{e}crit \`{a} l'ordre $0$:
 $\rho_0(f)\rho_0(g)=\rho_0(fg)$. Donc $\rho_0$ est un homomorphisme
 de l'alg\`{e}bre associative commutative $\Cinf(M,\korps)$ dans
 l'alg\`{e}bre associative $\Dop^1(\CinfK{C};\CinfK{C})$.
 Puisque $\rho_0(1)=$ l'application
 identique, alors $\rho_0$ envoie des fonctions qui ne s'annulent nulle
 part sur des op\'{e}rateurs diff\'{e}rentiels inversibles. Si l'on regarde la
 forme locale (\ref{EqDefFormeLocale}) d'un op\'{e}rateur
 diff\'{e}rentiel inversible dans des coordonn\'{e}es locales, on voit rapidement
 qu'il ne peut
 contenir aucune puissance strictement positive d'une d\'{e}riv\'{e}e partielle.
 Alors, un tel op\'{e}rateur diff\'{e}rentiel est de la forme $\psi\mapsto
 \chi\psi$ o\`{u} $\chi\in\Cinf(C,\korps)$. Soit $f\in \Cinf(M,\real)$.
 Alors la fonction $1+f^2$ est un \'{e}l\'{e}ment inversible dans l'alg\`{e}bre
 $\CinfK{M}$. Par cons\'{e}quent, il existe une fonction
 $\chi\in\Cinf(C,\korps)$ telle que
 \[
    \psi+\rho_0(f)^2(\psi)=\rho_0(1+f^2)(\psi)=\chi\psi
 \]
 quelle que soit $\psi\in\Cinf(C,\korps)$. Il s'ensuit qu'il existe
 une fonction $\chi'\in\Cinf(C,\korps)$ telle que $\rho(f)(\psi)
 =\chi'\psi$. Pour une fonction \`{a} valeurs complexes on arrive \`{a} la m\^{e}me
 conclusion tout en s\'{e}parant en partie r\'{e}elle et partie imaginaire.
 Alors il existe un homomorphisme d'alg\`{e}bres associatives commutatives
 $\tilde{\rho}_0:\Cinf(M,\korps)\ra \Cinf(C,\korps)$ tel que
 $\rho_0(f)\psi=\tilde{\rho}_0(f)~\psi$. D'apr\`{e}s l'exercice de Milnor
 (voir par exemple \cite{KMS93}, p.~301, Corollary 35.10) il existe une
 application
 de classe $\Cinf$ $i:C\ra M$ telle que $\tilde{\rho}_0(f)=f\circ i$.\\
 Consid\'{e}rons deux fonctions $f,g\in\CinfK{M}$ telles que $i^*f=0=i^*g$.
 Le commutateur de l'identit\'{e} de re\-pr\'{e}\-sen\-ta\-tion \`{a} l'ordre $1$
 s'\'{e}crit
 \[
   [\rho_1(f),\rho_0(g)]-[\rho_1(g),\rho_0(f)]=\rho_0(\{f,g\})
 \]
 Puisque $\tilde{\rho}_0(f)=i^*f=0=i^*g=\tilde{\rho}_0(g)$ il s'ensuit que
 $\{f,g\}\circ i=0$, alors $i$ est une application co\"{\i}sotrope.
\ebew

Dans ce cas particulier on peut restreindre la d\'{e}finition \ref{DefRep1}
de la fa\c{c}on suivante:
\bdefi\label{DRepDiff}
 Une representation de star-produits
 \[
  \rho=\sum_{r=0}^\infty\nu^r\rho_r:
   \CinfK{M}[[\nu]]\ra \Dop^1\big(\CinfK{C};\CinfK{C}\big)[[\nu]]
 \]
 est dite {\em diff\'{e}rentielle} lorsque tout $\rho_r$ est un op\'{e}rateur
 bidiff\'{e}rentiel le long de $i:C\ra M$ et $\id:C\ra C$ o\`{u} $\rho_0(f)(\psi)
 =(i^*f)\psi$, c.-\`{a}-d. tout $\rho_r$ appartient \`{a}
 $\Dop^2(\CinfK{M},\CinfK{C};\CinfK{C})$. On parlera du module
 $\Cinf{C}[[\nu]]$ comme {\em module diff\'{e}rentiel de $*$}, et au cas
 particulier o\`{u} $*=*_1\otimes *_2^{\mathrm{opp}}$ d'un
 {\em $*_1$-$*_2$-bimodule diff\'{e}rentiel}.
\edefi
Dans ce cas particulier il me semble int\'{e}ressant de savoir si toute
re\-pr\'{e}\-sen\-ta\-tion de star-produits est diff\'{e}rentielle.

\brem
 Le cas d'un fibr\'{e} vectoriel $E$ autre que $C\times \korps$
 peut \^{e}tre `non diff\'{e}rentiel': il suffit de prendre deux repr\'{e}sentations
 diff\'{e}rentielles $\rho^{(1)}$ et $\rho^{(2)}$ de $*$ dans $C\times \korps$
 dont les op\'{e}rateurs bidiff\'{e}rentiels sont le long de $i_1:C\ra M$ et
 le long de $i_2:C\ra M$ o\`{u} les applications co\"{\i}sotropes
 $i_1$ et $i_2$ sont distinctes. Ainsi la somme
 directe $\CinfK{C}[[\nu]]\oplus \CinfK{C}[[\nu]]$ est toujours un $*$-module
 (avec $\rho(f)(\psi_1+\psi_2):=\rho^{(1)}(f)\psi_1+\rho^{(2)}\psi_2$)
 qui est \'{e}gal \`{a} l'espace des sections $\Ginf(C,C\times \korps^2)[[\nu]]$,
 mais les coefficients de $\rho$ ne sont pas les op\'{e}rateurs diff\'{e}rentiels
 le long d'une seule application.
 Cet exemple montre que la notion d'une `repr\'{e}sentation diff\'{e}rentielle'
 pour les fibr\'{e}s vectoriels g\'{e}n\'{e}raux sera plus compliqu\'{e}e.
\erem

Encore une fois, la question r\'{e}ciproque de savoir quand l'injection
canonique $i$ d'une sous-vari\'{e}t\'{e} co\"{\i}sotrope ferm\'{e}e $C$
d'une vari\'{e}t\'{e} de Poisson donn\'{e}e
se d\'{e}forme dans une re\-pr\'{e}\-sen\-ta\-tion de star-produits me semble
aussi int\'{e}res\-sante:

\bprob[Quantification des sous-vari\'{e}t\'{e}s co\"{\i}sotropes] \label{prob3}
 Quelles sont les conditions sur l'injection canonique
 $i:C\ra (M,P)$  d'une sous-vari\'{e}t\'{e} co\"{\i}sotrope ferm\'{e}e $C$
d'une vari\'{e}t\'{e} de Poisson $(M,P)$ pour qu'il existent un star-produit
$*$ sur la vari\'{e}t\'{e} de Poisson $(M,P)$ et des applications lin\'{e}aires
 $\rho_1,\rho_2,\ldots:\Cinf(M,\korps)\ra \Dop^1(\CinfK{C};\CinfK{C})$
 tels que l'application suivante
 \[
    (f,\psi)\mapsto
    \rho(f)(\psi):= (i^*f)\psi + \sum_{r=1}^\infty\nu^r\rho_r(f)(\psi)
 \]
 (quels que soient $f\in\CinfK{M}[[\nu]], \psi\in
 \CinfK{C}[[\nu]]$)
 soit une re\-pr\'{e}\-sen\-ta\-tion (diff\'{e}rentielle) de star-produits?
\eprob
On arrive \`{a} la d\'{e}finition suivante:
\bdefi \label{DRepresentabilite}
 Soit $i:C\ra M$ une sous-vari\'{e}t\'{e} co\"{\i}sotrope ferm\'{e}e d'une vari\'{e}t\'{e}
 de Poisson $(M,P)$. Un star-produit $*$ sur $M$ est dit {\em
 repr\'{e}sentable sur $C$} lorsqu'il existe une repr\'{e}sentation diff\'{e}rentielle
 de $*$ sur $C$.
\edefi
On peut donner une `version diff\'{e}rentielle' de la liaison des
repr\'{e}sentations, voir \ref{DRepDiff}:
\bdefi\label{DRepDiffDifferentielle}
 Soient $(M,P)$ et $(M',P')$ deux vari\'{e}t\'{e}s de Poisson et $C,C'$ deux
 vari\'{e}t\'{e}s diff\'{e}rentiables. Soit $\phi:M\ra M'$ une application de de
 Poisson, soient $i:C\ra M$ et $i':C'\ra M'$ des applications
 co\"{\i}sotropes et soit $\varphi:C\ra C'$ une application de classe
 $\Cinf$ telles que le diagramme suivant soit commutatif:
 \beq\label{EqDiagrammeLiaison}
   \begin{array}{rcl}
     M & \stackrel{\phi}{\ra} & M' \\
     i\uparrow &              & \uparrow i' \\
     C & \stackrel{\varphi}{\ra} & C'
   \end{array}
 \end{equation}
 Soient $*$ et $*'$ des star-produits sur $M$ et $M'$, respectivement,
 soit $\Phi$ une quantification diff\'{e}rentielle de $\phi$
 (par rapport \`{a} $*$ et $*'$) et soient $\rho$ et $\rho'$ des
 repr\'{e}sentations diff\'{e}rentielles de $*$ sur $C$ (le long de $i$) et de
 $*'$ sur $C'$ (le long de $i'$), respectivement. \\
 {\em La repr\'{e}sentation $\rho$ est dite li\'{e}e \`{a} la repr\'{e}sentation $\rho'$
 par rapport \`{a} $\phi$ et $\varphi$} lorsqu'il existe une s\'{e}rie
 d'op\'{e}rateurs diff\'{e}rentiels le long de $\varphi$,
 $T=\varphi^*+\sum_{r=1}\nu^r T_r$ (c.-\`{a}-d. chaque $T_r$ appartient \`{a}
 $\Dop^1\big(\CinfK{C'},\CinfK{C}\big)$) telle que
 \beq\label{EqDeuxRepDiffLiees}
    \rho\big(\Phi(f')\big)\big(T(\psi')\big)
          =
          T\big(\rho'(f')(\psi')\big)
 \end{equation}
 quels que soient $f'\in\CinfK{M'}[[\nu]]$ et $\psi'\in\CinfK{C}[[\nu]]$.
\edefi
La proposition suivante est \'{e}vidente:
\bprop \label{PLiaisonEtRepresentabilite}
 Avec les hypoth\`{e}ses g\'{e}om\'{e}triques de la d\'{e}finition
 \ref{DRepDiffDifferentielle}, on suppose que $\varphi$ soit
 l'application identique de
 $C'=C$. Soient $*$ et $*'$ des star-produits sur $M$ et $M'$, respectivement,
 et soit $\Phi$ une quantification diff\'{e}rentielle de $\phi$
 (par rapport \`{a} $*$ et $*'$).
 \ben
  \item Soit $\rho$ une repr\'{e}sentation diff\'{e}rentielle de $*$ sur $C$
   le long de $i$. Alors la formule
 \[
    \rho'(f')\chi:=\rho\big(\Phi(f')\big)\chi~~~\forall~
    f'\in\CinfK{M'}[[\nu]]~~\forall~ \chi\in\CinfK{C}[[\nu]]
 \]
 d\'{e}finit une repr\'{e}sentation diff\'{e}rentielle de $*$ sur $C$ le long de $i$.
  \item Soit l'application de Poisson $\phi$ un diff\'{e}omorphisme de $M$
   sur un ouvert $U'$ de $M'$ (qui contient forc\'{e}ment $i'(C)$).
   Soit $\rho'$ une repr\'{e}sentation
   diff\'{e}rentielle de $*'$ sur $C$ le long de $i'$. Soit $*'|_U$
   la restriction de $*$ \`{a} $U$ et
   $\rho'|_{U'\times C}$ la restriction de la s\'{e}rie d'op\'{e}rateurs
   bidiff\'{e}rentiels $\rho'$
   le long de $i',id_C$ \`{a} $U'\times C$, voir l'\'{e}qn
   (\ref{EqDefRestrOpMultidiffLeLong}). Alors la restriction \`{a} $U'$ de la
   s\'{e}rie d'op\'{e}rateurs
   diff\'{e}rentiels le long de $\phi$, $\Phi|_{U'}$, admet une application
   r\'{e}ciproque
   $\Phi|_{U'}^{-1}$ de l'alg\`{e}bre associative $\big(\CinfK{M}[[\nu]],*\big)$
   dans $\big(\CinfK{U'}[[\nu]],*'|_{U'}\big)$ qui est
   une quantification diff\'{e}rentielle de $\phi^{-1}|_{U'}$. En outre,
   la formule
 \[
    \rho(f)\chi:=\rho'|_{U'\times C}\big(\Phi^{-1}(f)\big)\chi~~~\forall~
    f\in\CinfK{M}[[\nu]]~~\forall~ \chi\in\CinfK{C}[[\nu]]
 \]
 d\'{e}finit une repr\'{e}sentation diff\'{e}rentielle de $*$ sur $C$ le long de $i$.
 \een
\eprop
\bbew
 1. Ceci est \'{e}vident car $\phi\circ i =i'$.\\
 2. $\Phi|_{U'}$ est bien d\'{e}finie et inversible car son terme d'ordre
   $r=0$, $\phi^*$, l'est apr\`{e}s avoir \'{e}t\'{e} restreint \`{a} $U'$. Il est \'{e}vident
   que $\rho'|_{U'\times C}$ est une repr\'{e}sentation diff\'{e}rentielle
   du star-produit $*'|_{U'}$ sur $C$ (car $i'(C)\subset U'$. La formule
   est \'{e}vidente.
\ebew
Une cons\'{e}quence directe de cette proposition est le fait que la
repr\'{e}sentabilit\'{e} d'un star-produit ne d\'{e}pend que d'un voisinage ouvert
de la sous-vari\'{e}t\'{e} co\"{\i}sotrope; plus pr\'{e}cis\'{e}ment:
\bcor\label{CRepresentabiliteLocale}
 Soit $i:C\ra M$ une sous-vari\'{e}t\'{e} co\"{\i}sotrope ferm\'{e}e de la
 vari\'{e}t\'{e} de
 Poisson $(M,P)$. Soit $*$ un star-produit sur $M$.
 \ben
  \item Soit $U$ un ouvert de $M$ avec $C\subset U$ et soit $*_U$ un
  star-produit sur $U$ qui est repr\'{e}sentable sur $C$. Si la restriction de
  $*$ \`{a} $U$ est \'{e}quivalent \`{a} $*_U$, alors $*$ est repr\'{e}sentable sur $U$.
  \item Soit $*$ repr\'{e}sentable sur $C$, soit $U$ un ouvert de $M$ avec
    $C\subset U$ et soit $*_U$ un star-produit sur $U$ qui est \'{e}quivalent
    \`{a} la restriction $*|_U$ de $*$ \`{a} $U$. Alors $*_U$ est repr\'{e}sentable
    sur $C$.
 \een
\ecor

Un exemple important et bien connu pour des re\-pr\'{e}\-sen\-ta\-tions des
star-produits est le {\em calcul symbolique} (voir par exemple \cite{DL83a},
 \cite{BNW98},
\cite{BNW99} et \cite{BNPW03} dans un context de la quantification
par d\'{e}formation): soit $Q$ une vari\'{e}t\'{e} diff\'{e}rentiable et $\alpha$
une $1$-forme. Le rel\`{e}vement vertical de $\alpha$ est un champ de
vecteurs sur le fibr\'{e} cotangent $T^*Q$ de $Q$ dont le flot consiste
en translations le long des fibres de $T^*Q$ dans la direction de
$\alpha$. Par it\'{e}ration et en observant que deux translations commutent
on obtient ainsi un homomorphisme d'alg\`{e}bres associatives
$\mathsf{F}:\Ginf(Q,ST^*Q)\ra
\mathbf{D}^1\big(\CinfK{T^*Q},\CinfK{T^*Q}\big)$. On fixe une connexion
sans torsion $\nabla$ dans le fibr\'{e} tangent de $Q$ et on note
$\mathsf{D}:\Ginf(Q,S^\bullet T^*Q)\ra\Ginf(Q,S^{\bullet + 1}T^*Q)$ la
d\'{e}riv\'{e}e covariante
sym\'{e}tris\'{e}e (qui est l'analogue de la diff\'{e}rentielle de Cartan pour
les champs de tenseurs sym\'{e}triques). Soit $A\in\nu\Ginf(Q,T^*Q)[[\nu]]$.
Alors pour $f\in\CinfK{M}[[\nu]]$ et $\phi\in\CinfK{C}[[\nu]]$
la formule
\beq \label{EqRepRhoAConnGen}
    \rho^A_0(f)\phi
    :=i^*\big(\mathsf{F}\big(e^{-2\nu(D+\frac{1}{2\nu}A)}\phi\big)(f)\big)
\end{equation}
(o\`{u} $i:Q\ra T^*Q$ est la section nulle, donc un plongement sur une
sous-vari\'{e}t\'{e} lagrangienne)
d\'{e}finit une re\-pr\'{e}\-sen\-ta\-tion diff\'{e}rentielle
d'un star-produit sur $T^*Q$ (qui est
\'{e}quivalent au star-produit $\star_0$ dit {\em standard} au cas $A=0$) sur $Q$.
On peut associer \`{a} $A$ une s\'{e}rie
d'op\'{e}rateurs diff\'{e}rentiels $\mathcal{A}$ (voir \cite{BNPW03}, p.12,
eqn (3.11)) telle que
\beq \label{EqDefRhoAStand}
   \rho^A_0(f)\phi=\rho_0(\mathcal{A}^{-1}f)\phi,
\end{equation}
o\`{u} $\rho_0$ est une re\-pr\'{e}\-sen\-ta\-tion pour $\star_0$, et $\mathcal{A}$
est un
automorphisme de $\star_0$ ssi $dA=0$. Soit $B$ une s\'{e}rie de $2$-formes
ferm\'{e}es
dans $\nu\Ginf(Q,\Lambda^2T^*Q)[[\nu]]$. En \'{e}crivant localement
$B|_U=dA$ on calcule que la construction ci-dessus permet de
construire un star-produit $\star_0^B$ dont la classe de Deligne
est \'{e}gale \`{a} $\frac{[B]}{2\nu}=[\star_0^B]$ (voir \cite{BNPW03}, Theorem 4.6).
Dans ce cadre on trouve \'{e}galement d'autres re\-pr\'{e}\-sen\-ta\-tions
diff\'{e}rentielles param\'{e}tr\'{e}es
par un param\`{e}tre `d'ordre' comme la re\-pr\'{e}\-sen\-ta\-tion de Schr\"{o}dinger,
une interpr\'{e}tation de l'approximation
WKB et des re\-pr\'{e}\-sen\-ta\-tions diff\'{e}rentielles
dans un fibr\'{e} en droites complexes
holomorphe pour des monop\^{o}les magn\'{e}tiques (voir \cite{BNPW03}).

Au cas o\`{u} $Q=\real^k$ et $A=0$ l'\'{e}quation (\ref{EqRepRhoAConnGen})
se simplifie :
\beq \label{EqRhoStandardRk}
  \big(\rho_0(f)\phi\big)(q) = \sum_{r=0}^\infty\frac{(-2\nu)^r}{r!}
         \sum_{i_1,\ldots,i_r=1}^k
         \frac{\partial^r f}{\partial p_{i_1}\cdots\partial p_{i_r}}(q,0)
         \frac{\partial^r \phi}{\partial q^{i_1}\cdots\partial
         q^{i_r}}(q),
\end{equation}
et le star-produit $\star_0$ est donn\'{e} par
\beq \label{EqProdStarStandardRk}
    f\star_0 g = \sum_{r=0}^\infty\frac{(-2\nu)^r}{r!}
         \sum_{i_1,\ldots,i_r=1}^k
         \frac{\partial^r f}{\partial p_{i_1}\cdots\partial p_{i_r}}
         \frac{\partial^r g}{\partial q^{i_1}\cdots\partial
         q^{i_r}}.
\end{equation}
Cette formule bine connue est un cas particuier de la formule exponentielle
de M.Gerstenhaber, voir \cite[p.13]{Ger68} et \cite{DM93} pour une d\'{e}monstration
\'{e}l\'{e}gante.\\
Revenant au cas d'une sous-vari\'{e}t\'{e} co\"{\i}sotrope ferm\'{e}e $C$ d'une
vari\'{e}t\'{e} symplectique, on peut utiliser les deux formules pr\'{e}c\'{e}dentes
et l'existence d'une carte de Darboux adapt\'{e}e (voir la proposition
\ref{PCarteDarboux Adaptee}) pour arriver
au suivant
\bcor \label{CRepresentCarteSympl}
 Avec les notations mentionn\'{e}s ci-dessus, soit $*$ un star-produit
 sur $M$. Alors autour de tout point $c\in C$ il existe un ouvert $U$
 de $M$ et une re\-pr\'{e}\-sen\-ta\-tion diff\'{e}rentielle
 $\rho_U$ de $*|_U$ sur $C\cap U$.
\ecor
\bbew
 On choisit une carte de Darboux adapt\'{e}e $\big(U,(q^1,\ldots,q^m,p_1,\ldots,$
 $p_m,
 x^1,\ldots,x^k,y_1,\ldots,y_k)\big)$  avec un domaine contractile
 $U$. Alors $*|_U$ est \'{e}quivalent \`{a} tout star-produit local sur $U$,
 par exemple \`{a} $\star_0$ o\`{u} on remplace $q$ par $(q,x)$ et $p$ par
 $(p,y)$, et il est facile \`{a} calculer que la formule suivante d\'{e}finit
 une re\-pr\'{e}\-sen\-ta\-tion $\hat{\rho}_U$ de $\star_0$ sur $C\cap U$ qui est
 donn\'{e} par $y=0$: soit $f\in\CinfK{U}[[\nu]]$ et
 $\phi\in\CinfK{C}[[\nu]]$,
 \bea
  \lefteqn{\big(\hat{\rho}_U(f)\phi\big)(q,p,x) = \sum_{r,s=0}^\infty
        \frac{(-2\nu)^{r+s}}{r!s!} \sum_{i_1,\ldots,i_r=1}^k
           \sum_{a_1,\ldots,a_s=1}^m}  \nonumber \\
        & & \frac{\partial^{r+s} f}{\partial p_{i_1}\cdots\partial p_{i_r}
                                 \partial y_{a_1}\cdots\partial y_{a_s}}
                                 (q,p,x,0)
         \frac{\partial^{r+s} \phi}{\partial q^{i_1}\cdots\partial
         q^{i_r}\partial x^{a_1}\cdots\partial x^{a_s}}(q,p,x). \nonumber \\
        & &
 \eea
 ce qui montre le corollaire puisque $*|_U$ est \'{e}quivalent \`{a} $\star_0$.
\ebew
\brem
 Supposons qu'on ait une application de Poisson $\phi$ du fibr\'{e}
cotangent $T^*C$ d'une vari\'{e}t\'{e} diff\'{e}rentiable $C$
dans une vari\'{e}t\'{e} de Poisson $(M,P)$ qui admette une quantification $\Phi$
par rapport \`{a} un star-produit $*$ sur $M$ et au star-produit $*_0$
(relatif \`{a} une connexion sans torsion sur $C$) sur $T^*C$. Evidemment
la compos\'{e}e $f\mapsto \rho_0\big(\Phi(f)\big)$ d\'{e}finit une repr\'{e}sentation
diff\'{e}rentielle de $*$ sur $C$. Mais je ne sais pas si la r\'{e}ciproque
est forc\'{e}ment vraie, c'est-\`{a}-dire si une repr\'{e}sentation diff\'{e}rentielle
$\rho$ sur $C$ d'un star-produit $*$ sur $M$ d\'{e}finit automatiquement un
morphisme de star-produits $\Phi$ par rapport \`{a} $*$ et un star-produit
\'{e}quivalent \`{a} un $*_0$ sur $T^*C$: on rencontre des probl\`{e}mes de
convergence.
\erem

Une classe de re\-pr\'{e}\-sen\-ta\-tions particuli\`{e}res plus ou moins li\'{e}e au
contexte du calcul formel s'obtient par un proc\'{e}d\'{e}
analogue \`{a} celui qu'on utilise pour les re\-pr\'{e}\-sen\-ta\-tions de
Gel'fand, Naimark et Segal (GNS) des alg\`{e}bres
stellaires, voir \cite{BW98}, \cite{Bor99v}, \cite{Wal00},
\cite{waldmann:2002a}.

\subsection{R\'{e}duction des star-produits}
  \label{SubSecReducStarProduit}

Dans ce paragraphe on suppose le suivant:
Soit $i:C\ra M$ une sous-vari\'{e}t\'{e} co\"{\i}sotrope ferm\'{e}e de la vari\'{e}t\'{e}
symplectique $(M,\omega)$ telles que l'espace symplectique r\'{e}duit
$(M_{\mathrm{red}},\omega_{\mathrm{red}}$ est une vari\'{e}t\'{e} symplectique
pour laquelle la projection canonique $\pi:M\ra M_{\mathrm{red}}$ soit une
submersion. Soit $\mathcal{I}$ l'id\'{e}al annulateur de $C$ et
$\mathcal{N}(\mathcal{I})$ sont id\'{e}alisateur, voir (\ref{EqIdLieI}).

Il n'y a pas d'unique d\'{e}finition de la quantification de la r\'{e}duction
symplectique. Il y a au moins les deux variantes suivantes:

Dans la premi\`{e}re, on essaie de d\'{e}former l'isomorphisme d'alg\`{e}bres de
Poisson entre $\big(\CinfK{M_{\mathrm{red}}},\{~,~\}_{\mathrm{red}}\big)$
et l'alg\`{e}bre quotient $\mathcal{N}(\mathcal{I})/\mathcal{I}$, voir
la proposition \ref{PIdLieImodI}:

\bdefi \label{DStarProjet}
Avec les notations donn\'{e}es ci-dessus, soit $*$ un star-produit sur $M$.
\ben
 \item $*$ est dit {\em projetable} lorsque
 $(\mathcal{N}(\mathcal{I})[[\nu]],*)$ et $(\mathcal{I}[[\nu]],*)$ sont des
 sous-alg\`{e}bres de $\big(\CinfK{M}[[\nu]],*\big)$ telles que
 $\mathcal{I}[[\nu]]$ soit un id\'{e}al bilat\`{e}re dans l'alg\`{e}bre
 $(\mathcal{N}(\mathcal{I})[[\nu]],*)$.
 \item $*$ est dit {\em r\'{e}ductible} lorsqu'il est \'{e}quivalent \`{a} un
 star-produit projetable.
\een
\edefi
On arrive ainsi \`{a} une g\'{e}n\'{e}ralisation de la
formule (100) de \cite[p.138]{BHW00}:
\bprop\label{PStarProjet}
 Avec les notations donn\'{e}es ci-dessus, soit $*$ un star-produit projetable
 sur $M$. Alors l'espace $\CinfK{M_{\mathrm{red}}}[[\nu]]$ est muni d'un
 star-produit canonique $*_{\mathrm{red}}$ d\'{e}fini de la mani\`{e}re suivante:
 soient $f_1,f_2\in \mathcal{N}(\mathcal{I})[[\nu]]$, donc il
 existent des uniques fonctions $\tilde{f}_1,\tilde{f}_2\in
 \CinfK{M_{\mathrm{red}}}[[\nu]]$ telles que $i^*f_1=\pi^*\tilde{f}_1$ et
 $i^*f_2=\pi^*\tilde{f}_2$. Alors
 \beq \label{EqDefStarRed}
   i^*(f_1*f_2)=:\pi^*(\tilde{f}_1 *_{\mathrm{red}} \tilde{f}_2).
 \end{equation}
 On appelle $*_{\mathrm{red}}$ {\em le star-produit r\'{e}duit de $*$} et
 $\big(\CinfK{M_{\mathrm{red}}}[[\nu]],*_{\mathrm{red}}\big)$ {\em l'alg\`{e}bre
 r\'{e}duite}.
\eprop
\bbew
 Gr\^{a}ce au fait que $\mathcal{I}[[\nu]]$ est un id\'{e}al bilat\`{e}re dans
 la sous-alg\`{e}bre $\mathcal{N}(\mathcal{I})[[\nu]]$, l'alg\`{e}bre quotient
 $\mathcal{N}(\mathcal{I})[[\nu]]/\mathcal{I}[[\nu]]$ est bien d\'{e}finie
 et est isomorphe --en tant que $\korps[[\nu]]$-module-- \`{a}
 $\CinfK{M_{\mathrm{red}}}[[\nu]]$. Ceci justifie l'\'{e}quation
 (\ref{EqDefStarRed}) dont le membre droit est bien d\'{e}fini. Une
 \'{e}tude locale de cette \'{e}quation montre que la multiplication
 $*_{\mathrm{red}}$ associative est une s\'{e}rie formelle d'op\'{e}rateurs
 bidiff\'{e}rentiels sur $M_{\mathrm{red}}$ et que
 $1*_{\mathrm{red}}\tilde{f}_2=\tilde{f}_2=1*_{\mathrm{red}}\tilde{f}_2$.
 Par cons\'{e}quent, $*_{\mathrm{red}}$ est une star-produit sur
 $M_{\mathrm{red}}$.
\ebew
\brem
 Bien qu'un star-produit projetable induise un unique star-produit r\'{e}duit
 ceci n'est plus vrai pour un star-produit r\'{e}ductible: un tel star-produit
 peut \^{e}tre \'{e}quivalent \`{a} plusieurs star-produits projetables qui induisent
 des star-produits r\'{e}duits deux-\`{a}-deux non \'{e}quivalents, voir l'exemple
 de la r\'{e}duction pour l'espace projectif complexe en paragraphe
 \ref{SubSecEspProj}.
\erem

\brem
Dans la d\'{e}finition \ref{DStarProjet} on ne suppose a priori
pas que $\mathcal{I}[[\nu]]$ soit un id\'{e}al \`{a} gauche: donc elle peut
\^{e}tre vue comme un affaiblissement de la condition que $*$ admet
une {\em r\'{e}duction quantique coh\'{e}rente} donn\'{e}e dans le contexte
des r\'{e}ductions de Marsden-Weinstein dans \cite{BHW00}, p.135, Def.~31.
On verra pourtant en paragraphe \ref{SubSecISousAlg} que dans le cas
symplectique le fait que $\mathcal{I}[[\nu]]$ est une sous-alg\`{e}bre
implique que $\mathcal{I}[[\nu]]$ est ou bien un id\'{e}al \`{a} gauche ou bien
un id\'{e}al \`{a} droite (si $C$ est connexe), voir le th\'{e}or\`{e}me \ref{TIAlgId}.
\erem

Pour la deuxi\`{e}me d\'{e}finition, qui m'a \'{e}t\'{e} communiqu\'{e}e par Alan Weinstein,
on essaie de `quantifier' le fait que
l'espace des fonctions $\CinfK{C}$ est un bimodule par rapport aux
anneaux $\CinfK{M}$ et $\CinfK{M_\mathrm{red}}$: soient $f\in\CinfK{M}$,
$\phi\in\CinfK{C}$ et $\tilde{h}\in\CinfK{M_\mathrm{red}}$, alors
l'application
\beq
%    \CinfK{M}\times\CinfK{C}\times\CinfK{M_\mathrm{red}}\ra \CinfK{C}:
      (f,\phi,\tilde{h})\mapsto (i^*f)\phi(\pi^*\tilde{h})
             \label{EqBimodClasRed}
\end{equation}
d\'{e}finit la structure d'un $\CinfK{M}$-$\CinfK{M_\mathrm{red}}$-bimodule
sur $\CinfK{C}$.
Ceci provoque la suivante
\bdefi \label{DStarStarredBimod}
 Avec les notations donn\'{e}es ci-dessus, soient $*$ un star-produit sur
 $(M,\omega)$
 et $*_r$ un star-produit sur $(M_\mathrm{red},\omega_\mathrm{red})$.
 On va dire que {\em $\CinfK{C}[[\nu]]$ est un $*$-$*_r$-bimodule}
 (diff\'{e}rentiel) lorsqu'il
 existe une re\-pr\'{e}\-sen\-ta\-tion (diff\'{e}\-ren\-tielle) $\rho$ de $*$ dans
 $C\times\korps$ qui d\'{e}forme $i^*$
 et une re\-pr\'{e}\-sen\-ta\-tion (diff\'{e}rentielle) $\tilde{\rho}$ de
 $*_r^\mathrm{opp}$ dans $C\times \korps$ qui d\'{e}forme $\pi^*$ telles que
 \beq\label{EqRhoRhoComm}
     \rho(f)\tilde{\rho}(\tilde{h})\big(\phi\big)
        =\tilde{\rho}(\tilde{h})\rho(f)\big(\phi\big)
 \end{equation}
 quels que soient $f\in \Cinf{M}[[\nu]],
 ~\tilde{h}\in\Cinf{M_\mathrm{red}}[[\nu]],~\phi\in\Cinf{C}[[\nu]]$.
\edefi
En particulier, $\rho$ est une representation diff\'{e}rentielle de $*$ sur
$C$, alors, comme on verra en paragraphe \ref{SubSecRepDeformIdealAnnul},
l'id\'{e}al \`{a} gauche $\{f\in\CinfK{M}[[\nu]]~|~\rho(f)1=0\}$ sera une
d\'{e}formation de l'id\'{e}al annulateur $\mathcal{I}$. Ceci correspond bien au
`coisotropic creed' de Weinstein et Lu, voir \cite{Lu93}.

Au cas o\`{u} l'espace r\'{e}duit n'existe plus en tant que vari\'{e}t\'{e}
diff\'{e}rentiable, on peut introduire une version locale d'un star-produit
projetable: soit $U\subset M$ un ouvert. On d\'{e}finit
\beas
  \mathcal{I}_U & := & \{g\in\CinfK{U}~|~g(c)=0~\forall c\in C\cap U\} \\
  \mathcal{N}(\mathcal{I})_U & := &
                       \{ f\in \CinfK{U}~|~\{f,g\}\in\mathcal{I}_U~\forall
                              g\in\mathcal{I}_U \}
\eeas
et la version locale de la d\'{e}finition \ref{DStarProjet} est la suivante
\bdefi\label{DStarProjetLocal}
 Avec les notations donn\'{e}es ci-dessus, soit $*$ un star-produit sur $M$.
\ben
 \item $*$ est dit {\em localement projetable} lorsque quel que soit
 l'ouvert $U$ de $M$ les espaces
 $(\mathcal{N}(\mathcal{I})_U[[\nu]],*|_U)$ et $(\mathcal{I}_U[[\nu]],*|_U)$
 sont des
 sous-alg\`{e}bres de $\big(\CinfK{U}[[\nu]],*|_U\big)$ telles que
 $\mathcal{I}_U[[\nu]]$ soit un id\'{e}al bilat\`{e}re dans l'alg\`{e}bre
 $(\mathcal{N}(\mathcal{I})_U[[\nu]],*)$.
 \item $*$ est dit {\em localement r\'{e}ductible} lorsqu'il est \'{e}quivalent \`{a} un
 star-produit localement projetable.
\een
\edefi
Il ne faut pas oublier que le star-produit $*$ est globalement donn\'{e}.
De la m\^{e}me fa\c{c}on que pour la proposition \ref{PStarProjet}, on montre
l'existence des star-produits r\'{e}duits $*_{U\mathrm{red}}$
--cette fois locaux-- sur les quotients locaux
$\big(\mathcal{N}(\mathcal{I})_U/\mathcal{I}_U\big)[[\nu]]$.

\section{Premiers r\'{e}sultats}
 \label{SecPremiersResultats}

\subsection{Quantification des applications de Poisson comme
       probl\`{e}me de re\-pr\'{e}\-sen\-ta\-tion}
          \label{SubSecMorRep}

Le lien entre la quantification des application de Poisson et la
quantification des sous-vari\'{e}t\'{e}s co\"{\i}sotropes se trouve dans la
suivante

\bprop \label{PDefoMorphDefoBimod}
 Soit $\phi:(M',P')\ra (M,P)$ une application de Poisson entre deux
 vari\'{e}t\'{e}s de Poisson. Soit $C:=M'$ et $i:C:=M'\ra M\times M'$ le plongement
 canonique de $M'$ dans le graphe de $\phi$, c.-\`{a}-d.
 $i(p'):=(\phi(p'),p')\in M\times M'$ quel que soit $p'\in M'$.
 Soient $*$ et $*'$ des star-produits sur $(M,P)$ et $(M',P')$,
 respectivement.\\
 Alors les \'{e}nonc\'{e}s suivants sont \'{e}quivalents:
 \ben
  \item Il existe un homomorphisme de star-produits diff\'{e}rentiel
   \[
    \Phi=\sum_{r=0}\nu^r\Phi_r:
     \big(\Cinf(M,\korps)[[\nu]],*\big)\ra
        \big(\Cinf(M',\korps)[[\nu]],*'\big)\mathrm{~avec~}
           \Phi_0=\phi^*.
   \]
   \item Il existe la structure d'un
   $*$-$*^{\prime \mathrm{opp}}$-bimodule sur $C$ qui d\'{e}forme $i$, c.-\`{a}-d.
   une re\-pr\'{e}\-sen\-ta\-tion de star-produits diff\'{e}rentielle
   $\rho=\sum_{r=0}\nu^r\rho_r$ de
 $\big(\Cinf(M\times M',\korps)[[\nu]],*\otimes *^{\prime\mathrm{opp}}\big)$
 dans $\CinfK{C}[[\nu]]$ telle que $\rho_0=i^*$.
 \een
\eprop
\bbew
 D'abord, $C$ est
 une sous-vari\'{e}t\'{e} co\"{\i}sotrope de $(M\times M',P_{(1)}-P'_{(2)}$
 d'apr\`{e}s le th\'{e}or\`{e}me de Weinstein \ref{TheoGrapheWei}.

 \noindent ``$\Longrightarrow$'' Pour $f\in \CinfK{M}$ et
  $g_1,g_2\in \CinfK{M'}$ on d\'{e}finit
  \[
    \rho(f\otimes g_1)(g_2):=\Phi(f)*'g_2*'g_1.
  \]
  Puisque $\Phi$ est diff\'{e}rentielle et $*'$ est bidiff\'{e}rentiel, on voit
  sans peine que cette d\'{e}finition s'\'{e}tend de $f\otimes g_1$ \`{a} toutes les
  fonctions $F\in \CinfK{M\times M'}[[\nu]]$ \`{a} une application
  bidiff\'{e}rentielle le long de $i$ et $\mathrm{id_C}$. Il est \'{e}vident que ceci
  nous donne une re\-pr\'{e}\-sen\-ta\-tion de star-produits diff\'{e}rentielle
  sur $M'=C$, et
  \[
     \rho_0(f\otimes g_1)(g_2)=\Phi_0(f)g_1 g_2=i^*(f\otimes g_1)g_2.
  \]
 ``$\Longleftarrow$'':
 Soient $f,f_1,f_2\in\CinfK{M}$ et $g,g_1,g_2\in \CinfK{M'}$.
 On consid\`{e}re les applications
 $r$ et $l$: $\Cinf(M',\korps)[[\nu]]\ra\Cinf(M',\korps)[[\nu]]$ donn\'{e}es par
 $r(g):=\rho(1\otimes g)(1)$ et
 $l(f):=\rho(f\otimes 1)(1)$.
 On a $r_0(g)=\rho_0(1\otimes g)(1)=i^*(1\otimes g)=g$, alors $r_0$ est
 l'application identique ce qui entra\^{\i}ne que $r$ est inversible. On d\'{e}finit
 \[
     \Phi:=r^{-1}\circ l:
     \Cinf(M,\korps)[[\nu]]\ra\Cinf(M',\korps)[[\nu]]
 \]
 Ensuite,
 $l_0(f)=\rho_0(f\otimes 1)(1)=i^*(f\otimes 1)=\phi^*f$, par cons\'{e}quent
 $\Phi_0=(r^{-1}\circ l)_0=l_0=\phi^*$. On a
 \beas
   \rho(1\otimes g_1)r(g_2)=\rho(1\otimes g_1)\rho(1\otimes g_2)(1)
               & = &
        \rho\big(1\otimes (g_1*^{\prime \mathrm {opp}} g_2)\big)(1) \\
               & = & r(g_1*^{\prime \mathrm {opp}} g_2) = r(g_2*'g_1).
 \eeas
 d'o\`{u}
 \[
      g_2*'g_1=\big(r^{-1}\circ \rho(1\otimes g_1) \circ r\big)(g_2)
 \]
 En outre,
 \beas
   \big(r^{-1}\circ \rho(f\otimes 1)\circ r\big)(g)
     & = & r^{-1}\big(\rho(f\otimes 1)\rho(1\otimes g)(1) \big)
       =   r^{-1}\big(\rho(f\otimes g)(1)\big) \\
     & = & r^{-1}\big(\rho(1\otimes g)\rho(f\otimes 1)(1)\big) \\
     & = & \big(r^{-1}\circ \rho(1\otimes g)\circ r\big)
                   \big(r^{-1}(l(f))\big)\\
     & = & \Phi(f) *' g
 \eeas
 d'apr\`{e}s l'\'{e}quation pr\'{e}c\'{e}dente. Par cons\'{e}quent
 \beas
   \Phi(f_1*f_2)*'g & = &
          \big(r^{-1}\circ \rho((f_1*f_2)\otimes 1)\circ r\big)(g) \\
          & = & \big(r^{-1}\circ \rho(f_1\otimes 1)\circ r\big)\left(
            \big(r^{-1}\circ \rho(f_2\otimes 1)\circ r\big)(g)\right)\\
          & = & \Phi(f_1)*'\big(\Phi(f_2)*'g\big)
                 = \big(\Phi(f_1)*'\Phi(f_2)\big)*'g
 \eeas
 ce qui montre en particulier $(g=1)$ que $\Phi$ est un homomorphisme de
 star-produits.
\ebew
\brem
 La cat\'{e}gorie des $\korps[[\nu]]$-alg\`{e}bres associatives
est parfois
remplac\'{e}e par une cat\'{e}gorie `plus vaste'
qui a les m\^{e}mes objets, mais dont les morphismes
de $\mathcal{A}$ dans $\mathcal{B}$ sont donn\'{e}s par des classes
d'isomorphismes
de $\mathcal{A}$-$\mathcal{B}$-bimodules qui soient projectifs et d'un
nombre fini de g\'{e}n\'{e}rateurs en tant que $\mathcal{B}$-modules,
voir par exemple \cite[p.20]{Lod92}. La proposition pr\'{e}c\'{e}dente montre que la
d\'{e}formation d'un morphisme de Poisson dans la cat\'{e}gorie plus vaste,
c.-\`{a}-d. en tant que bi-module, reste dans la cat\'{e}gorie usuelle o\`{u} les
morphismes sont des morphismes d'alg\`{e}bres associatives.
\erem

%{\tt Hier k\"{o}nnte die Sache mit Impulsabbildungen und $T^*G$ und so kommen}

\subsection{Repr\'{e}sentations, la d\'{e}formation de l'id\'{e}al annulateur et
  star-produits adapt\'{e}s}
  \label{SubSecRepDeformIdealAnnul}

Soit $(M,P)$ une vari\'{e}t\'{e} de Poisson, $i:C\ra M$ une sous-vari\'{e}t\'{e}
co\"{\i}sotrope ferm\'{e}e et $\mathcal{I}$ l'id\'{e}al annulateur de $C$,
voir \ref{EqIdAn}. Soit $*$ un star-produit sur $M$ et
$\rho$ une re\-pr\'{e}\-sen\-ta\-tion diff\'{e}rentielle de $*$ dans $C\times \korps$.
On d\'{e}finit
\beq \label{EqDefIdrho}
  \mathcal{I}_\rho:=\{f\in \CinfK{M}[[\nu]]~|~\rho(f)1=0\}.
\end{equation}
Il est \'{e}vident que $\mathcal{I}_\rho$ est un id\'{e}al \`{a} gauche de
$\big(\CinfK{M}[[\nu]],*\big)$. Pour le cas d'une re\-pr\'{e}\-sen\-ta\-tion GNS
formelle (voir \cite{BW98}) o\`{u} l'espace pr\'{e}hilbertien est engendr\'{e} par
un vecteur `vacuum' (par exemple pour la re\-pr\'{e}\-sen\-ta\-tion de Schr\"{o}dinger
sur un espace des configurations compact) l'id\'{e}al $\mathcal{I}_\rho$
co\"{\i}ncide avec {\em l'id\'{e}al de Gel'fand}.
On peut se demander quand cet id\'{e}al
co\"{\i}ncide avec l'espace $\mathcal{I}[[\nu]]$:
\bprop\label{PKerRho1congId}
 Soient $(M,P)$, $C$, $*$ et $\rho$ comme ci-dessus. Alors il existe
 une s\'{e}rie d'op\'{e}rateurs diff\'{e}rentiels $S=id+\sum_{r=1}^\infty \nu^r S_r$
 (o\`{u} $S_r$ appartient \`{a} $\Dop^1\big(\CinfK{M};\CinfK{M}\big)$ et s'annule
 sur les
 constantes quel que soit l'entier
 positif $r$) telle que pour le star-produit $*':=S(*)$ et
 la re\-pr\'{e}\-sen\-ta\-tion diff\'{e}rentielle
 $\rho'(f):=\rho(S^{-1}f)$ (o\`{u} $f\in \CinfK{M}[[\nu]]$) on a
 \ben
  \item Pour tous $f,g\in \CinfK{M}[[\nu]]$:
     \[
         \rho'(f)\big(i^*g\big) = i^*(f*'g).
     \]
  \item  L'id\'{e}al \`{a} gauche
    $\mathcal{I}_{\rho'}$ de l'alg\`{e}bre
    $\big(\CinfK{M}[[\nu]],*'\big)$ est \'{e}gal \`{a} l'espace
    $\mathcal{I}[[\nu]]$.
  \item L'id\'{e}al \`{a} gauche $\mathcal{I}_\rho$ de l'alg\`{e}bre
    $\big(\CinfK{M}[[\nu]],*\big)$ est \'{e}gal \`{a}
   l'espace $S^{-1}\big(\mathcal{I}[[\nu]]\big)$.
 \een
\eprop
\bbew
 Puisque la repr\'{e}sentation $\rho$ est diff\'{e}rentielle, alors l'application
 $D:\CinfK{M}[[\nu]]\ra \CinfK{C}[[\nu]]$ d\'{e}finie par
 \[
    Df := \rho(f)1
 \]
 quel que soit $f\in \CinfK{M}[[\nu]]$ est une s\'{e}rie d'op\'{e}rateurs
 diff\'{e}rentiels $D=\sum_{r=0}\nu^r D_r$ o\`{u} $D_r$ appartient
 \`{a} $\Dop^1\big(\CinfK{M};\CinfK{C}\big)$ (le long de $i$) quel que soit
 l'entier positif $r$. Evidemment, $D_0=i^*$ et $D_r$ s'annule sur les
 constantes
 quel que soit $r\geq 1$. A l'aide d'un voisinage
 tubulaire $U$ de $C$ on va construire vers la fin de cette d\'{e}monstration
 une s\'{e}rie d'op\'{e}rateurs
 diff\'{e}rentiels $S=id +\sum_{r=1}^\infty\nu^r S_r$ o\`{u}
 tout $S_r$ appartient \`{a} $\Dop^1\big(\CinfK{M};\CinfK{M}\big)$ et s'annule sur
 les constantes, telle que
 \[
     D_rf= i^*S_rf~~~\forall f\in \CinfK{M}~\forall r\in \nat.
 \]
 Alors $S$ est inversible et
 \[
    \rho(f)1=i^*(Sf)~~~\mathrm{quel~que~soit}~f\in\CinfK{M}[[\nu]].
 \]
 Alors $\rho'$ est \'{e}videmment une re\-pr\'{e}\-sen\-ta\-tion de $*'$ dans
 $C\times\korps$. Puisque $\rho'(f)1=i^*f$ il s'ensuit que
 $\mathcal{I}_{\rho'}$ est \'{e}gal au noyau de $i^*$, alors \`{a}
 $\mathcal{I}[[\nu]]$. De plus, on d\'{e}duit de $\rho(f)1=i^*(Sf)$ que
 $\mathcal{I}_{\rho}=S^{-1}\big(\mathcal{I}[[\nu]]\big)$. En outre
 \beas
  \rho'(f)\big(i^*g\big) = \rho(S^{-1}f)i^*\big(S(S^{-1}g)\big)
              & = & \rho(S^{-1}f)\rho(S^{-1}g)1 =
                    \rho\big(S^{-1}(f*'g)\big)1 \\
              & = &
                   i^*\big(S(S^{-1}(f*'g))\big) = i^*(f*'g),
 \eeas
 ce qui montre l'\'{e}nonc\'{e} 1.\\
 Pour construire $S$ il faut prolonger les valeurs de $D$:
 d'abord il est clair
 que $D_0=i^*=i^*~id=i^*~S_0$. Soit $\tau:U\ra C$ un
 voisinage tubulaire de $C$ dans $M$. Soit $(U_\alpha)_{\alpha\in
 \mathbf{A}}$ un recouvrement ouvert de $C$ tel que
 tout $\tau^{-1}(U_\alpha)$ soit le domaine d'une carte
 $x_\alpha=(x_\alpha^1,\ldots,x_\alpha^n)$ de
 $M$. Soit $(\chi_\alpha)_{\alpha\in\mathbf{A}}$ une partition de l'unit\'{e}
 subordonn\'{e}e \`{a} $(U_\alpha)_{\alpha\in \mathbf{A}}$.
 Pour tout $r\in\nat,r\geq 1$
 il existe un entier positif $N_r$ tel que l'op\'{e}rateur diff\'{e}rentiel $D_r$
 est de la forme
 locale suivante pour tout $f\in \CinfK{M}$ et $c\in U_\alpha$:
 \[
      D_r(f)|_{U_\alpha}(c) = \sum_{|I|\leq N_r}D^I_{\alpha r}(c)
                      \frac{\partial^{|I|}f}{\partial x_\alpha^I}(c)
 \]
 pour des multi-indices $I\in N^n$ et certaines fonctions $D^I_{\alpha r}$
 appartenant \`{a} $\CinfK{U_\alpha}$. On d\'{e}finit pour tout $u\in U$
 \[
      S_{r,U}(f)(u):= \sum_{\alpha\in\mathbf{A}}
                            \chi_\alpha\big(\tau(u)\big)
                             \sum_{|I|\leq N_r}
                            D^I_{\alpha r}\big(\tau(u)\big)
                                  \frac{\partial^{|I|}f}{\partial
                                  x_\alpha^I}(u).
 \]
 Visiblement, tout $S_{r,U}$ est un op\'{e}rateur diff\'{e}rentiel bien d\'{e}fini
 appartenant \`{a}
 $\Dop^1\big(\CinfK{M};\CinfK{U}\big)$ tel que $S_{r,U}(f)(c)=D_r(f)(c)$ quel
 que soit $c\in C$ et $f\in\CinfK{M}$. Soit maintenant $V$ un ouvert de
 $M$ tel que $C\subset V\subset \bar{V} \subset U$ (puisque $M$ est un
 espace topologique normal) et soit $W$ l'ouvert $M\setminus \bar{V}$
 (ici $\bar{V}$ etc. d\'{e}signe l'adh\'{e}rence de $V$ dans $M$).
 Alors $M=U\cup W$ et $C\cap \bar{W}=\emptyset$. Soit $(\psi_U,1-\psi_U)$
 une partition de l'unit\'{e} subordonn\'{e}e au recouvrement ouvert $(U,W)$ de
 $M$. Finalement, les op\'{e}rateurs diff\'{e}rentiels
 \[
     S_r(f):=\psi_U~S_{r,U}(f)~~~~~~\forall r\in\nat, r\geq 1
 \]
 appartiennent \`{a} l'espace d'op\'{e}rateurs diff\'{e}rentiels
 $\Dop^1\big(\CinfK{M};\CinfK{M}\big)$ et satisfont \`{a}
 l'\'{e}quation $S_r(f)(c)=D_r(f)(c)$ quel que soient $c\in C$ et
 $f\in\CinfK{M}$.
 \ebew

\noindent On voit qu'il existe toujours un couple \'{e}quivalent $(*',\rho')$
tel que l'id\'{e}al annulateur est directement un id\'{e}al \`{a} gauche pour le
star-produit \'{e}quivalent $*'$.

R\'{e}ciproquement, soit $*$ un star-produit sur $M$ pour lequel l'id\'{e}al
annulateur $\mathcal{I}[[\nu]]$ soit un id\'{e}al \`{a} gauche. Il s'ensuit
que l'\'{e}quation
\beq
   \rho(f)\big(i^*g\big):= i^*(f*g)~~~\mathrm{quels~que~soient~}
                          f,g\in\CinfK{M}[[\nu]]
\end{equation}
 est une re\-pr\'{e}\-sen\-ta\-tion diff\'{e}rentielle de
$*$ dans $C\times\korps$ bien d\'{e}finie parce que la restriction $i^*$ est
surjective et $f*g$ appartient \`{a} $\mathcal{I}[[\nu]]$ si $g$ est un
\'{e}l\'{e}ment de $\mathcal{I}[[\nu]]$. \\
Les consid\'{e}rations pr\'{e}c\'{e}dentes montrent le
\bsat\label{TIdAnIdGauche}
 Soit $(M,P)$ une vari\'{e}t\'{e} de Poisson, $i:C\ra M$ une sous-vari\'{e}t\'{e}
 co\"{\i}sotrope et $*=\sum_{r=0}^\infty\nu^r \mathsf{C}_r$ un star-produit
 sur $M$.\\
 Alors $*$ est repr\'{e}sentable sur $C$, i.e.
 il existe une re\-pr\'{e}\-sen\-ta\-tion diff\'{e}rentielle
 de $*$ dans $C\times\korps$ si et
 seulement s'il existe un star produit
 $*'=\sum_{r=0}^\infty\nu^r \mathsf{C}'_r$
 \'{e}quivalent \`{a} $*$ tel que tous les op\'{e}rateurs bidiff\'{e}rentiels $\mathsf{C}'_r$
 sont tels que
 \beq \label{EqIIdealGauche}
       \mathsf{C}'_r(f,g)\in\mathcal{I}~~\mathrm{quels~que~soient~}
            f\in\CinfK{M}, g\in\mathcal{I}.
 \end{equation}
\esat
Cette proposition provoque la suivante
\bdefi \label{DDefStarAdapte}
 Soit $(M,P)$ une vari\'{e}t\'{e} de Poisson et $i:C\ra M$ une sous-vari\'{e}t\'{e}
 co\"{\i}sotrope ferm\'{e}e de $M$.\\
 Un star-produit $*$ sur $M$ est dit
 {\em adapt\'{e} \`{a} la sous-vari\'{e}t\'{e} co\"{\i}sotrope $C$} lorsque
 la condition (\ref{EqIIdealGauche}) est satisfaite, c.-\`{a}-d.
 $\mathcal{I}[[\nu]]$ est un id\'{e}al \`{a} gauche de $*$ o\`{u}
 $\mathcal{I}$ est l'id\'{e}al annulateur de $C$. On appellera {\em sa
 repr\'{e}senta\-tion canonique} $\rho$ l'application
 \[
     \rho(f)i^*g:=i^*(f*g)
 \]
 quels que soient $f,g\in\CinfK{M}[[\nu]]$.
\edefi
L'exemple clef pour un star-produit adapt\'{e} est le star-produit du type
standard $*_0$ pour un fibr\'{e} cotangent, voir l'\'{e}qn
(\ref{EqRepRhoAConnGen}) pour sa repr\'{e}sentation canonique au cas o\`{u} $A=0$.
Un cas particulier plus \'{e}lementaire est le star-produit
(\ref{EqProdStarStandardRk}) sur $\real^{2n}$.

La notion du commutant s'av\`{e}rera tr\`{e}s importante plus tard:
\bdefi \label{DefCommutant}
 Soit $(M,P)$ une vari\'{e}t\'{e} de Poisson et $i:C\ra M$ une sous-vari\'{e}t\'{e}
 co\"{\i}sotrope ferm\'{e}e de $M$. Soit $*$ un star-produit repr\'{e}sentable
 sur $M$ et $\rho$ une repr\'{e}sentation diff\'{e}rentielle de $*$ sur $C$.\\
 Le {\em commutant de $*$ par rapport \`{a} $\rho$} est d\'{e}fini par
 \bea
  \mathrm{Hom}_{(*,\rho)}\big(C,C\big) & :=&
                          \{D:\CinfK{C}[[\nu]]\ra \CinfK{C}[[\nu]]~|~
                                                            \nonumber\\
    & & ~~~~ D~\mathrm{~est~}\korps[[\nu]]-\mathrm{lin\acute{e}aire~et~}
                      \nonumber \\
    & & ~~~~~~~D\rho(f)=\rho(f)D~\forall~f\in
    \CinfK{M}[[\nu]]\}.\nonumber\\
    & &
 \eea
\edefi
Pour un star-produit adapt\'{e} $*$ on d\'{e}finit {\em l'id\'{e}alisateur de
$\mathcal{I}[[\nu]]$ par rapport \`{a} $*$} par
\beq \label{EqDefIdIStar}
  \mathcal{N}_*(\mathcal{I}):=\{h\in \CinfK{M}[[\nu]]
           ~|~g*h\in\mathcal{I}[[\nu]]~\forall~g\in\mathcal{I}[[\nu]]\}.
\end{equation}
Il est \'{e}vident que $\mathcal{N}_*(\mathcal{I})$ est une sous-alg\`{e}bre
de $\big(\CinfK{M}[[\nu]],*\big)$ et que $\mathcal{I}[[\nu]]$ est un id\'{e}al
bilat\`{e}re de $\mathcal{N}_*(\mathcal{I})$. On a
\bprop \label{PCommutantNdeI}
 Soit $(M,P)$ une vari\'{e}t\'{e} de Poisson, $i:C\ra M$ une sous-vari\'{e}t\'{e}
 co\"{\i}sotrope ferm\'{e}e de $M$ et soit $*$ un star-produit adapt\'{e}
 sur $M$.\\
 \ben
  \item Le commutant de $*$ par rapport \`{a} sa repr\'{e}sentation canonique
 $\rho$ est donn\'{e} par
 \[
    \mathrm{Hom}_{(*,\rho)}(C,C)=
      \{ i^*f\mapsto i^*(f*h)~|~h\in \mathcal{N}_*(\mathcal{I})\}
        \cong
        \big(\mathcal{N}_*(\mathcal{I})/\mathcal{I}[[\nu]]\big)^{\mathrm{opp}}.
 \]
 \item En outre, si $*$ est projetable, alors $\mathcal{N}_*(\mathcal{I})
 =\mathcal{N}(\mathcal{I})[[\nu]]$, et le commutant est donc anti-isomorphe
 \`{a} l'alg\`{e}bre r\'{e}duite.
 \een
\eprop
\bbew
 1. Soit $D\in \mathrm{Hom}_{(*,\rho)}(C,C)$.
   Soit
   $f'\in\CinfK{M}[[\nu]]$ tel que $i^*f'=D1$. Il vient
   \[
         i^*(f*f')=\rho(f)i^*f'=\rho(f)\big(D1\big)=D\rho(f)1=D(i^*f),
   \]
   donc $D$ est d\'{e}termin\'{e} par sa valeur sur $1$.
   En particulier, pour $f=g\in\mathcal{I}[[\nu]]$ la fonction $f'$ est
   telle que $g*f'\in\mathcal{I}[[\nu]]$, donc $f'\in
   \mathcal{N}_*(\mathcal{I})$ et $D(i^*f)=i^*(f*f')$. R\'{e}ciproquement,
   soit $h\in \mathcal{N}_*(\mathcal{I})$. Alors l'application
   \[
      D_h(i^*f):=i^*(f*h)
   \]
    est bien d\'{e}finie car $i^*(g*h)=0$ quel que soit
   $g\in \mathcal{I}[[\nu]]$ par d\'{e}finition de
   $\mathcal{N}_*(\mathcal{I})$. De plus
   \[
      D_h\rho(f)i^*\hat{f}= D_h(i^*(f*\hat{f}) = i^*(f*\hat{f}*h)
         = \rho(f)i^*(\hat{f}*h) = \rho(f) D_h i^*\hat{f}
   \]
   et $D_h\in \mathrm{Hom}_{(*,\rho)}(C,C)$. On voit facilement que $D_h=0$
   si $h\in\mathcal{I}[[\nu]]$ et que $D_h=0$ implique $0=D_h(1)=i^*h$,
   donc $h\in\mathcal{I}[[\nu]]$
   \[
      D_{h*h'}=D_{h'}D_h~~~~~\forall h,h'\in \mathcal{N}_*(\mathcal{I}).
   \]
   Ceci montre que $D:h\mapsto D_h$ d\'{e}finit un antihomomorphisme
   surjective de $\mathcal{N}_*(\mathcal{I})$ sur l'alg\`{e}bre associative
   $\mathrm{Hom}_{(*,\rho)}(C,C)$. Le noyau de $D$ \'{e}tant
   $\mathcal{I}[[\nu]]$, l'isomorphisme annonc\'{e} est montr\'{e}.

   2. Puisque $*$ est projetable, $\mathcal{I}[[\nu]]$ est un id\'{e}al \`{a} gauche
   de $\mathcal{N}(\mathcal{I})[[\nu]]$. Alors il est \'{e}vident que
   $\mathcal{N}(\mathcal{I})[[\nu]]\subset \mathcal{N}_*(\mathcal{I})$.\\
   R\'{e}ciproquement, soit
   $f=\sum_{r=0}^\infty\nu^rf_r$ un \'{e}l\'{e}ment de $\mathcal{N}_*(\mathcal{I})$
   (avec $f_r\in\CinfK{M}$). Puisque pour tout $f'\in \CinfK{M}[[\nu]]$ il
   vient
   $f'*g\in\mathcal{I}[[\nu]]$ (car $\mathcal{I}[[\nu]]$ est un id\'{e}al \`{a}
   gauche) il suffit de consid\'{e}rer l'\'{e}quation
   \[
          0=i^*(f*g-g*f)~~~~\forall g\in\mathcal{I}.
   \]
   A l'ordre $r=1$ on obtient $0=2i^*\{f_0,g\}~\forall g\in\mathcal{I}$,
   donc $f_0\in\mathcal{N}(\mathcal{I})$. Supposons par r\'{e}currence
   que $f_0,\ldots,f_r\in \mathcal{N}(\mathcal{I})$. Donc
   $f^{(r)}:=f_0+\nu f_1+\cdots+\nu^r f_r\in\mathcal{N}(\mathcal{I})[[\nu]]$,
   et puisque $\mathcal{I}[[\nu]]$ est un id\'{e}al bilat\`{e}re dans ce dernier
   espace il vient que $0=i^*(f^{(r)}*g-g*f^{(r)})$, donc
   $h^{(r)}:=f'-f^{(r)}$ est un multiple de $\nu^{r+1}$, et on a
   $0=i^*(h^{(r)}*g-g*h^{(r)})$. Cette condition \`{a} l'ordre $r+2$ donne
   $0=2i^*\{f_{r+1},g\}~\forall g\in\mathcal{I}$, alors $f_{r+1}$ est
   un \'{e}l\'{e}ment de $\mathcal{N}(\mathcal{I})$. Il s'ensuit que $f\in
    \mathcal{N}(\mathcal{I})[[\nu]]$, alors
    $\mathcal{N}_*(\mathcal{I})\subset \mathcal{N}(\mathcal{I})[[\nu]]$.
\ebew

Soit $\Dop_I\big(\CinfK{M};\CinfK{M}\big)$
{\em l'espace de tous les op\'{e}rateurs
diff\'{e}rentiels adapt\'{e}s \`{a} $C$}, c.-\`{a}-d. ceux qui pr\'{e}servent $\mathcal{I}$:
\bea
  \lefteqn{\Dop_I(\CinfK{M};\CinfK{M}) :=}\nonumber \\
      & & ~~~~~~~~~\{D\in\Dop^1(\CinfK{M};\CinfK{M})~|~
                 Dg\in\mathcal{I}~\forall~g\in\mathcal{I}\}.
\eea
Puisque $\CinfK{C}\cong \CinfK{M}/\mathcal{I}$ il existe une application
naturelle
\bea
   \Dop_I(\CinfK{M};\CinfK{M})& \ra & \Dop^1(\CinfK{C};\CinfK{C}):\nonumber
   \\
     D & \mapsto & \big(i^*f\mapsto i^*(Df)\big). \label{EqDopIMDopC}
\eea
Le sous-espace d'op\'{e}rateurs diff\'{e}rentiels
$\Dop_I\big(\CinfK{M},\CinfK{M}\big)$
servira d'espace naturel pour les
transformations d'\'{e}quivalence entre deux star-pro\-duits adapt\'{e}s \`{a} $C$: d'un
c\^{o}t\'{e}, il est \'{e}vident que toute s\'{e}rie
$S=\mathrm{id}+\sum_{r=1}^\infty \nu^rS_r$
 telle que tous les op\'{e}rateurs diff\'{e}rentiels $S_r$ soient adapt\'{e}s \`{a} $C$
 transforme tout star-produit adapt\'{e} \`{a} $C$ dans un star-produit adapt\'{e}
 \`{a} $C$. Dans ce cas on a la proposition suivante:
\bprop \label{PCouplesEquiSiTrafoAdap}
 Soit $i:C\ra M$ une sous-vari\'{e}t\'{e} co\"{\i}sotrope de la vari\'{e}t\'{e} de Poisson
 $(M,P)$. Soient $*$ et $*'$ deux star-produits adapt\'{e}s \`{a} $C$, $\rho$ et
 $\rho'$ leurs repr\'{e}sentations canoniques, et soit
 $S$ une transformation d'\'{e}quivalence adapt\'{e}e, i.e. $*'=S(*)$ et
 $S\big(\mathcal{I}[[\nu]]\big)=\mathcal{I}[[\nu]]$. \\
 Alors les couples $(*,\rho)$ et $(*',\rho')$ sont \'{e}quivalents.
\eprop
\bbew
 Soit $T$ la s\'{e}rie d'op\'{e}rateurs diff\'{e}rentiels de $C$ induite par $S$,
 i.e. $Ti^*g=i^*(Sg)$ quel que soit $g\in\CinfK{M}[[\nu]]$.
 Soient $f,g\in\CinfK{M}[[\nu]]$. Alors
 \beas
  \rho'(f)i^*g & = & i^*(f*'g) = i^*\big(S\big((S^{-1}f)*(S^{-1}g)\big)\big)
                         \\
               & = & Ti^*\big((S^{-1}f)*(S^{-1}g)\big)
                     = T\rho(S^{-1}f)i^*(S^{-1}g) \\
               & = & T\rho(S^{-1}f)T^{-1}i^*g,
 \eeas
 donc $\rho'(f)=T\rho(S^{-1}f)T^{-1}$, et les couples $(*,\rho)$ et
 $(*',\rho')$ sont \'{e}quivalents.
\ebew

D'un autre c\^{o}t\'{e}, il n'est pas vrai en g\'{e}n\'{e}ral que deux star-produits
adapt\'{e}s \'{e}quivalents sont \'{e}quivalents par rapport \`{a} une transformation
d'\'{e}quivalence adapt\'{e}e comme on peut
 voir dans des exemples, voir le paragraphe \ref{SecExemples}.
 Mais le cas particulier suivant est important:
\bsat \label{Tdaserstemalfalsch}
 Soit $(M,\omega)$ une vari\'{e}t\'{e} symplectique et $i:C\ra M$ une sous-vari\'{e}t\'{e}
 co\"{\i}sotrope ferm\'{e}e de $M$. Soient $*$ et $*'$ deux star-produits
 adapt\'{e}s \`{a} $C$. On suppose que le premier groupe de cohomologie verticale
 de $C$ s'annule:
 \[
    H^1_v(C,\korps)= \{0\}.
 \]
 Si $*$ et $*'$ sont \'{e}quivalents, alors il existe une
 transformation d'\'{e}quivalence $S=\mathrm{id}+\sum_{r=1}^\infty \nu^rS_r$
 telle que tous les op\'{e}rateurs diff\'{e}rentiels $S_r$ soient adapt\'{e}s \`{a} $C$.
\esat
\bbew
 Soit $*=\sum_{r=0}^\infty\nu^r \mathsf{C}_r$ et
 $*'=\sum_{r=0}^\infty\nu^r \mathsf{C}'_r$.
 Le cas $r=0$ \'{e}tant \'{e}vident, supposons qu'on ait d\'{e}j\`{a} trouv\'{e} $r+1$
 op\'{e}rateurs diff\'{e}rentiels adapt\'{e}s $T_0:=0,T_1,\ldots,T_r$ ($r\geq 0$)
 tels que $T^{(r)}:=(\mathrm{id}+\nu^rT_r)\cdots(\mathrm{id}+\nu T_1)
 \mathrm{id}$ envoie $*$ sur un star-produit adapt\'{e}
 $T^{(r)}(*)=:*^{(r)}=:\sum_{s=0}^\infty\nu^s \mathsf{C}^{(r)}_s$ qui
 co\"{\i}ncide avec $*'$ jusqu'\`{a} l'ordre $r$. Puisque $*$ et $*'$ sont
 \'{e}quivalents, $*'$ et $*^{(r)}$ le sont aussi.
 %On rappelle {\em l'op\'{e}rateur
 %cobord de Hochschild} $\mathsf{b}$ pour un op\'{e}rateur diff\'{e}rentiel
 %$U:\CinfK{M}\ra\CinfK{M}$:
 %\beq
 % (\mathsf{b}U)(f_1,f_2)~~:=~~f_1~U(f_2)~-~U(f_1f_2)~+~U(f_1)~f_2
 %\end{equation}
 %quels que soient $f_1,f_2\in\CinfK{M}$.
 D'apr\`{e}s la th\'{e}orie g\'{e}n\'{e}rale
 des classes de Deligne (voir la proposition
 \ref{PDergrosseIrrtummitderexaktenEinsform} et \cite{GR99}),
 il existe un
 op\'{e}rateur diff\'{e}rentiel $B:\CinfK{M}\ra\CinfK{M}$ et une $1$-forme
 $\gamma$ sur $M$ tels que
 \beas
    \mathsf{C}'_{r+1}(f_1,f_2) & = &\mathsf{C}^{(r)}_{r+1}(f_1,f_2)
           -(\mathsf{b}B)(f_1,f_2) +(d\gamma)(X_{f_1},X_{f_2})
                 \\
                 & & ~~~~~~~~~~~~~~~~~~~~~~\forall~f_1,f_2\in\CinfK{M}
 \eeas
 o\`{u} $\mathsf{b}$ d\'{e}signe l'op\'{e}rateur cobord de Hochschild, voir
 \'{e}qn (\ref{EqCobordHochschild}).
 Puisque les star-produits $*$ et $*^{(r)}$ sont adapt\'{e}s, alors la
 restriction $i^*$ appliqu\'{e}e \`{a} l'\'{e}quation pr\'{e}c\'{e}\-dente (pour $f_1=g_1$ et
 $f_2=g_2$ avec $g_1,g_2\in\mathcal{I}$) entra\^{\i}ne que la somme
 $-i^*B(g_1g_2)+i^*(d\gamma)(X_{g_1},X_{g_2})$ s'annule. En s\'{e}parant
 la partie antisymetrique, il vient
 \[
     d_v(p_v i^*\gamma)=0
 \]
 (voir les \'{e}quations (\ref{EqDefdverticale}) et (\ref{EqDefVertSuite})
 pour la notation),
 car les champs hamiltoniens $X_g$ de $g\in\mathcal{I}$
 induisent tous les champs verticaux, voir la proposition
 \ref{PChampsVertIdeal}. Puisque $H^1_v(C,\korps)=\{0\}$ il existe une
 fonction $f\in\CinfK(M)$ telle que
 \[
      \gamma_c\big(V)=(df)_c\big(V)~~~~~~\forall ~c\in C
                            ~~\forall~V\in T_cC^\omega.
 \]
 La restriction de la $1$-forme modifi\'{e}e $\gamma':=\gamma-df$ \`{a} $C$
 s'annule donc sur tous les vecteurs verticaux. Soit $Y$ l'unique champ
 de vecteurs tel que $\omega^\flat(Y)=-\gamma'$ et soit $g\in\mathcal{I}$.
 Il vient
 \[
   i^*Y(g)=i^*\big(dg(Y)\big)=i^*\big(\omega(X_g,Y)\big)
                 = i^*\big(\gamma'(X_g)\big) = 0
 \]
 alors $Y(g)\in\mathcal{I}$, et $Y$ est un champ de vecteurs adapt\'{e} (en
 tant qu'op\'{e}rateur diff\'{e}rentiel d'ordre $1$). En modifiant l'op\'{e}rateur
 adapt\'{e} $T_r$ par l'op\'{e}rateur adapt\'{e} $T'_r:=T_r+Y$, les coefficients
 $\mathsf{C}^{(r)}_k$ du star-produit $*^{(r)}$ ne changent pas
 pour $0\leq k\leq r$ (car $Y$ est une d\'{e}rivation de la multiplication
 point-par-point, donc $\mathsf{b}Y=0$), mais \`{a} l'ordre $r+1$ le star-produit
 modifi\'{e} adapt\'{e} $\hat{*}^{(r)}$ est \'{e}gal \`{a}
 \[
   \hat{\mathsf{C}}^{(r)}_{r+1}(f_1,f_2)=
    \mathsf{C}^{(r)}_{r+1}(f_1,f_2)+Y\big(\mathsf{C}_1(f_1,f_2)\big)
                       -\mathsf{C}_1\big(Y(f_1),f_2\big)
                       -\mathsf{C}_1\big(f_1,Y(f_2)\big)
 \]
 Puisque il est connu que
 $\mathsf{C}_1(f_1,f_2)=\{f_1,f_2\}+(\mathsf{b}B_1)(f_1,f_2)$
 pour un op\'{e}rateur diff\'{e}rentiel $B_1$ il s'ensuit que
 \[
  \mathsf{C}^{(r)}_{r+1}(f_1,f_2)=
    \hat{\mathsf{C}}^{(r)}_{r+1}(f_1,f_2)
      -(d\gamma')(X_{f_1},X_{f_2})-(\mathsf{b}[Y,B_1])(f_1,f_2),
 \]
 donc
  \beas
    \mathsf{C}'_{r+1}(f_1,f_2) & = & \hat{\mathsf{C}}^{(r)}_{r+1}(f_1,f_2)
           -\big(\mathsf{b}(B+[Y,B_1])\big)(f_1,f_2)
                 \\
                 & & ~~~~~~~~~~~~~~~~~~~~~~\forall~f_1,f_2\in\CinfK{M}.
 \eeas
 En restreignant cette \'{e}quation \`{a} $C$ et en observant que $*'$ et
 $\hat{*}^{r}$ sont adapt\'{e}s, on obtient les \'{e}quations suivantes pour
 l'op\'{e}rateur diff\'{e}rentiel $B':=B+[Y,B_1]$
 \beq \label{EqistarBfgg}
   \begin{array}{cccl}
     i^*B'(fg) & = & i^*f~i^*B(g) &\forall f\in\CinfK{M}~\forall
     g\in\mathcal{I} \\
     i^*B'(g_1g_2) & = & 0 & \forall g_1,g_2\in\mathcal{I}.
   \end{array}
 \end{equation}
 Il s'ensuit que la restriction de $i^*B'$ \`{a} $\mathcal{I}$ d\'{e}finit un
 morphisme $\tilde{B}$ de $\CinfK{C}$-modules du module conormal
 $\mathcal{I}/\mathcal{I}^2\cong\Ginf(C,TC^\mathrm{ann})$ \`{a}
 $\CinfK{C}$, par cons\'{e}quent une section du fibr\'{e} vectoriel $TM|_C/TC$.
 En choisissant un sous-fibr\'{e} $F$ de $TM|_C$ compl\'{e}mentaire \`{a} $TC$ et
 en utilisant un voisinage tubulaire de $C$ on trouve un champ de vecteurs
 $X$ sur $M$ tel que $X_c~\mathrm{mod}~T_cC~=~\tilde{B}_c$ quel que soit
 $c\in C$. Par cons\'{e}quent, pour l'op\'{e}rateur diff\'{e}rentiel
 \[
    T_{r+1}:=B' - L_X
 \]
 (o\`{u} $L_X$ d\'{e}signe la d\'{e}riv\'{e}e de Lie par rapport au champ de vecteurs $X$)
 il vient
 \[
    i^*T_{r+1}(g)=i^*B'(g)-i^*L_X(g)=0~~~\forall g\in\mathcal{I},
 \]
 donc $T_{r+1}$ est adapt\'{e} \`{a} $C$. Puisque $L_X$ est une d\'{e}rivation
 de $\CinfK{M}$ on a
 \beas
    \mathsf{C}'_{r+1}(f_1,f_2) & = & \hat{\mathsf{C}}^{(r)}_{r+1}(f_1,f_2)
                                  -(\mathsf{b}T_{r+1})(f_1,f_2) \\
                  & &~~~~~~~~~~~~~~~~~~~\forall~f_1,f_2\in\CinfK{M}
 \eeas
 et le star-produit adapt\'{e}
 $*^{(r+1)}:=(\mathrm{id}+\nu^{r+1}T_{r+1})(\hat{*}^{(r)})$ co\"{\i}ncide avec
 $*'$ jusqu'\`{a} l'ordre $r+1$. Ceci montre que le produit infini
 $S:=\lim_{r\to\infty}(\mathrm{id}+\nu^rT'_r)\cdots (\mathrm{id}+\nu T'_1)$
 d\'{e}finit une transformation d'\'{e}quivalence adapt\'{e}e de $*$ \`{a} $*'$.
\ebew
\bcor \label{Cdaserstemalfalsch}
 Soit $(M,\omega)$ une vari\'{e}t\'{e} symplectique et $i:C\ra M$ une sous-vari\'{e}t\'{e}
 co\"{\i}sotrope ferm\'{e}e de $M$. Soient $*$ et $*'$ deux star-produits
 sur $M$ et $\rho$ (resp. $\rho'$) une repr\'{e}sentation diff\'{e}rentielle de $*$
 (resp. de $*'$) sur $C$. On suppose que le premier groupe de cohomologie
 verticale s'annule, $H^1_v(C,\korps)= \{0\}$.\\
 Si $*$ et $*'$ sont \'{e}quivalents, alors les couples $(*,\rho)$ et
 $(*',\rho')$ sont \'{e}galement \'{e}quivalents.
\ecor
\bbew
 D'apr\`{e}s la proposition \ref{PKerRho1congId} il y a deux transformations
 d'\'{e}quiva\-lence $S$ et $S'$ tels que les star-produits $\hat{*}:=S(*)$ et
 $\hat{*}':=S'(*')$ soient adapt\'{e}s et le couple
 $(*,\rho)$ (resp. $(*',\rho')$) est \'{e}quivalent au couple
 $(\hat{*},\hat{\rho})$ (resp. $(\hat{*}',\hat{\rho})'$) o\`{u}
 $\hat{\rho}$ (resp. $\hat{\rho}'$) est la rep\'{e}sentation canonique
 de $\hat{*}$ (resp. de $\hat{*}'$). Puisque $*$ et $*'$ sont \'{e}quivalents,
 alors $\hat{*}$ et $\hat{*}'$ le sont, et d'apr\`{e}s le th\'{e}or\`{e}me
 \ref{Tdaserstemalfalsch} il existe une transformation d'\'{e}quivalence
 adapt\'{e}e $U$ telle que $\hat{*}'=U(\hat{*})$. D'apr\`{e}s la proposition
 \ref{PCouplesEquiSiTrafoAdap} les couples $(\hat{*},\hat{\rho})$ et
 $(\hat{*}',\hat{\rho})'$
 sont \'{e}quivalents. Par cons\'{e}quent, les couples $(*,\rho)$ et
 $(*',\rho')$ sont \'{e}galement \'{e}quivalents.
\ebew

\brem
L'analogue multidiff\'{e}rentiel de $\Dop_\mathcal{I}$, l'espace
$\mathfrak{G}_\mathcal{I}^r$, d\'{e}fini par le sous-espace de tous les
op\'{e}rateurs multidiff\'{e}rentiels (de rang $r\in\nat$) $\mathsf{D}$ tels que
pour tous $f_1,\ldots,f_{r-1}\in\CinfK{M}$ et $g\in\mathcal{I}$
la fonction $\mathsf{D}(f_1,\ldots,f_{r-1},g)\in\mathcal{I}$ s'av\`{e}re
tr\`{e}s important pour les consid\'{e}rations de formalit\'{e} dans \cite{BGHHW03p}.
\erem

\subsection{Repr\'{e}sentations sur les sous-vari\'{e}t\'{e}s lagrangiennes}
  \label{SubSecRepSousvarLag}

Soit $(M,\omega)$ une vari\'{e}t\'{e} symplectique et $L\subset M$ une
sous-vari\'{e}t\'{e} lagrangienne ferm\'{e}e. D'apr\`{e}s un th\'{e}or\`{e}me de Weinstein
 (voir \cite{Wei71} ou \cite{AM85}, p.411, thm 5.3.18) il existe un
 voisinage ouvert $U\supset L$ dans $M$,
 un voisinage ouvert $L\subset V \subset T^*L$ de la section nulle
 du fibr\'{e} cotangent de $L$ et un diff\'{e}omorphisme symplectique
 $\phi:U\ra V$ dont la restriction \`{a} $L$ donne l'identification usuelle
 de $L$ avec la section nulle $L\ra T^*L$.

Alors, pour \'{e}tudier les re\-pr\'{e}\-sen\-ta\-tions des alg\`{e}bres d\'{e}form\'{e}es
sur les
sous-vari\'{e}t\'{e}s lagrangiennes, on regarde d'abord le cas $M=T^*L$ \`{a} l'aide
du calcul symbolique esquiss\'{e} dans paragraphe \ref{SubSecRep1}.

\bsat \label{TRepLagTStarQ}
 Soit $V$ un voisinage ouvert fibre-par-fibre convexe de $L$ dans la vari\'{e}t\'{e}
 symplectique $T^*L$.
 Soit $*$ un star-produit sur $V$.
 \ben
  \item Il existe une re\-pr\'{e}\-sen\-ta\-tion diff\'{e}rentielle
  de $*$ sur $L$ si et seulement si $[*]=0$.
  \item Soit $\rho$ une re\-pr\'{e}\-sen\-ta\-tion diff\'{e}rentielle
   de $*$ sur $L$ et $*'$ un autre
    star-produit sur $V$ avec une representation diff\'{e}rentielle
    $\rho'$ sur $L$.
    Alors les couples $(*,\rho)$ et $(*',\rho')$ sont \'{e}quivalents.
  \item On suppose que $*$ soit repr\'{e}sentable sur $L$. Alors
   l'ensemble de toutes les classes d'isomorphisme de modules de $*$
     sur $L$ est en bijection avec
     \[
       \nu H^{1\prime}_{dR}(L,\korps)\oplus \nu^2 H^1_{dR}(L,\korps)[[\nu]],
    \]
    o\`{u} on a \'{e}crit $H^{1\prime}_{dR}(L,\korps):=H^1_{dR}(L,\real)$ au cas
    o\`{u} $\korps=\real$
    et
    $H^{1\prime}_{dR}(L,\korps):=
    H^1_{dR}(L,\complex)/2\pi iH^1_{dR}(L,\mathbb{Z})$ au cas o\`{u}
    $\korps=\complex$.
 \een
\esat
\bbew
 1. Soit $B\in \nu\Omega^2(L)[[\nu]]$ une s\'{e}rie formelle de $2$-formes
 ferm\'{e}es telle que $\frac{B}{2\nu}$ repr\'{e}sente la classe de Deligne $[*]$ de
 $*$: ceci est possible
 car $V$ se retracte sur $L$ et $[\omega]=0$ pour la forme symplectique
 canonique $\omega$ de $V\subset T^*L$. Soit $\star_0^B$ le star-produit
 sur $T^*L$ d\'{e}fini dans \cite{BNPW03}, paragraphe 4, eqn (4.2). Puisque
 sa classe de Deligne est \'{e}gale \`{a} $\frac{[B]}{2\nu}=[*]$
 (voir \cite{BNPW03}, Theorem 4.6) alors la restriction de
 $\star_0^B$ \`{a} $V$ est \'{e}quivalente \`{a} $*$. Soit $\tilde{S}$ une
 transformation d'\'{e}quivalence telle que $\tilde{S}*=\star_0^B$.\\
 ``$\Longleftarrow$'': If $[*]=0$, alors $*$ est \'{e}quivalent au star-produit
 $\star_0$ (o\`{u} $B=0$) qui admet toujours une re\-pr\'{e}\-sen\-ta\-tion sur $L$
 gr\^{a}ce au calcul symbolique (voir le paragraphe \ref{SubSecRep1} et
 \cite{BNW98}, \cite{BNW99} et
 \cite{BNPW03}, \'{e}quation (2.1)). Donc il existe une re\-pr\'{e}\-sen\-ta\-tion
 diff\'{e}rentielle de $*$ sur $L$.\\
 ``$\Longrightarrow$'': R\'{e}ciproquement, supposons qu'une re\-pr\'{e}\-sen\-ta\-tion
 diff\'{e}rentielle $\rho$
 existe pour le star-produit $*$.
 Par cons\'{e}quent, $\CinfK{V}[[\nu]]\ni
 f\mapsto\tilde{\rho}(f):=\rho(\tilde{S}^{-1}f)$ est une
 re\-pr\'{e}\-sen\-ta\-tion
 de la restriction de $\star_0^B$ \`{a} $V$ (qu'on va noter d\'{e}sormais par
 $\star_0^B$) dans $L$. Soit $\pi:V\ra L$ la restriction de la projection
 du fibr\'{e} $T^*L\ra L$ \`{a} $V$. On consid\`{e}re l'op\'{e}rateur diff\'{e}rentiel
 \[
      D:\CinfK{L}[[\nu]]\ni \phi\mapsto \tilde{\rho}(\pi^*\phi)1\in
      \CinfK{L}[[\nu]].
 \]
 Gr\^{a}ce au fait que $\pi\circ i$ est l'application identique de $L$ il
 s'ensuit que $D\phi = \phi +$ termes d'ordre sup\'{e}rieur en $\nu$, donc
 $D$ est une s\'{e}rie formelle d'op\'{e}rateurs diff\'{e}rentiels inversible.
 Puisque
 \[
    (\pi^*\phi)\star_0^B(\pi^*\psi)=\pi^*(\phi\psi)~~~\forall
 \phi,\psi\in\CinfK{L}[[\nu]]
 \]
 (voir \cite{BNPW03}, Theorem 4.1 ii)) il
 vient $D(\phi\psi)=\tilde{\rho}(\pi^*\phi)D\psi$. En d\'{e}finissant la
 repr\'{e}senta\-tion
 $\CinfK{V}[[\nu]]\ni
 f\mapsto\hat{\rho}(f):=D^{-1}\tilde{\rho}(f)D$ de $\star_0^B$ dans $L$,
 on voit que les couples $(*,\rho)$
 et $(\star_0^B,\hat{\rho})$ sont \'{e}quivalents et que
 \beq \label{EHatRhoPhi}
     \hat{\rho}(\pi^*\phi)\psi=\phi\psi ~~~\forall
 \phi,\psi\in\CinfK{L}[[\nu]].
 \end{equation}
 Pour un champ de vecteurs $X$ sur $L$ soit $J(X)\in\CinfK{T^*L}$ la
  fonction $\zeta\mapsto \langle\zeta,X_{\pi(\zeta)}\rangle$. On utilisera
  les m\^{e}mes symboles pour sa restriction \`{a} $V$. Soit $Y$
  un autre champ de vecteurs sur $L$ et $\phi\in\CinfK{L}[[\nu]]$. A
  l'aide des \'{e}quations (4.2), (4.1), (3.11) et (2.1) de \cite{BNPW03}
  on calcule:
  \bea
       \pi^*\phi~\star_0^B J(X) & = & \pi^*\phi~J(X) = J(\phi X)
                                            \label{EphistarBX}\\
       J(X)\star_0^B \pi^*\phi & = & \pi^*\phi ~J(X) -2\nu \pi^*(X\phi)
                                            \label{EXstarBphi} \\
       J(X)\star_0^B J(Y)-J(Y)\star_0^B J(X) & = &
               -2\nu\bigg(J([X,Y])+\pi^*\big(B(X,Y)\big)\bigg)
                             \label{EXstarBcommYphi} \\
  \eea
  Par cons\'{e}quent, pour tout $\psi\in\CinfK{L}[[\nu]]$ il vient de
  (\ref{EHatRhoPhi}), (\ref{EphistarBX}) et (\ref{EXstarBphi}):
  \beas
     -2\nu(X\phi)\psi & = &\hat{\rho}\big(-2\nu\pi^*(X\phi)\big)\psi
        =\hat{\rho}\big(J(X)\star_0^B \pi^*\phi
                     -~\pi^*\phi~\star_0^B J(X)\big) \psi \\
            & = & \hat{\rho}\big(J(X)\big)\hat{\rho}\big(\pi^*\phi\big)\psi
             -\hat{\rho}\big(\pi^*\phi\big)\hat{\rho}\big(J(X)\big)\psi \\
            & = & \hat{\rho}\big(J(X)\big)(\phi\psi)
                  -\phi\big(\hat{\rho}\big(J(X)\big)\psi\big)
  \eeas
  En particulier, pour $\psi=1$ on obtient
  \beq \label{EHatRhoX}
   \hat{\rho}\big(J(X)\big)\phi = -2\nu X\phi
                                -A(X)\phi
  \end{equation}
  o\`{u} l'application $A:\Ginf(L,TL)\ra \CinfK{L}[[\nu]]:
        X\mapsto \hat{\rho}\big(J(X)\big)1$ est une s\'{e}rie d'$1$-formes
        sur $L$, car gr\^{a}ce \`{a} (\ref{EphistarBX}) on a
  \beas
   A(\phi X) & = & \hat{\rho}\big(\pi^*\phi J(X)\big)1
            =\hat{\rho}\big(\pi^*\phi ~\star_0^B J(X)\big)1
            =\hat{\rho}\big(\pi^*\phi\big)\hat{\rho}\big(J(X)\big)1
            =\phi\hat{\rho}\big(J(X)\big)1 \\
             & = & \phi A(X)
  \eeas
  Finalement, les \'{e}quations (\ref{EXstarBcommYphi}) et (\ref{EHatRhoX})
  entra\^{\i}nent
  \beas
   -2\nu\hat{\rho}\big(J([X,Y]\big)\psi-2\nu B(X,Y)\psi
    & = & \hat{\rho}\big(J(X)\big)\hat{\rho}\big(J(Y)\big)\psi
      -\hat{\rho}\big(J(Y)\big)\hat{\rho}\big(J(X)\big)\psi \\
    & = & 4\nu^2[X,Y]\psi \\
    &   & +2\nu\bigg(X\big(A(Y)\big)-Y\big(A(X)\big)-A([X,Y])\bigg)\psi\\
    &   & +2\nu A([X,Y])\psi \\
    & = & -2\nu\hat{\rho}\big(J([X,Y])\big)\psi +2\nu (dA)(X,Y)\psi,
  \eeas
  donc $B=dA$, alors $[*]=\frac{1}{2\nu}[B]=0$.\\
  2. D'apr\`{e}s la premi\`{e}re partie, on a $[*]=0$, et $*$ est \'{e}quivalent
  \`{a} $\star_0$. Dans la d\'{e}monstration de 1) on vu des cas particuliers de
  la re\-pr\'{e}\-sen\-ta\-tion
  $\hat{\rho}$ de $\star_0$: quels que soient $\phi,\psi\in\CinfK{L}[[\nu]]$
  et $X\in\Gamma(TL)[[\nu]]$ il vient $\hat{\rho}(\pi^*\phi)\psi=\phi\psi$
  et $\hat{\rho}\big(J(X)\big)\psi=-2\nu\big(X+\frac{1}{2\nu}A(X)\big)\psi$.
  Par cons\'{e}quent, $\hat{\rho}$ con\"{\i}ncide avec $\rho_0^A$ (voir
  l'\'{e}quation (\ref{EqDefRhoAStand})) sur les sous-espaces
  $\pi^*\CinfK{L}[[\nu]]$ et $J\big(\Ginf(L,TL)\big)[[\nu]]$ o\`{u} $dA=0$.
  Puisque ces
  deux espaces engendrent localement l'espace de toutes les fonctions
  qui sont polyn\^{o}miales le long des fibres de $T^*L$, et puisque le fait
  que $\rho$ et $*$ sont diff\'{e}rentielles entra\^{\i}ne que toute
  re\-pr\'{e}\-sen\-ta\-tion est d\'{e}termin\'{e}e par ses valeurs sur le sous-espace des
  fonctions
  polyn\^{o}miales, il vient que $\hat{\rho}$ et $\rho_0^A$ co\"{\i}ncident.
  Donc les couples $(*,\rho)$ et $(\star_0,\rho_0^A)$ sont \'{e}quivalents,
  et puisque $(*,\rho)$ \'{e}tait arbitrairement choisi, alors tous les
  couples sont \'{e}quivalents.\\
  3. Gr\^{a}ce \`{a} 2) il reste \`{a} v\'{e}rifier lesquels des modules $\rho_0^A$ et
  $\rho_0^{A'}$ pour $dA=0$ et $dA'=0$ sont isomorphes. L'\'{e}nonc\'{e} pour
  $\korps=\complex$ est le contenu du th\'{e}or\`{e}me 7.3 de \cite{BNPW03},
  c.-\`{a}-d. la classe de de Rham de $A-A'$ est n\'{e}cessairement un \'{e}l\'{e}ment de
  $\nu 2\pi H^1_{dR}(L,\entier)$ parce qu'on peut
  avoir de l'holonomie int\'{e}grale. Pour $\korps=\real$ ce ph\'{e}nom\`{e}ne
  n'existe pas, d'o\`{u} le fait que la classe de $A-A'$ doit s'annuler.
\ebew

\noindent Le th\'{e}or\`{e}me de Weinstein mentionn\'{e} ci-dessus permet de d\'{e}duire le
suivant
\bcor \label{CorRepLagArbi}
Soit $*$ un star-produit sur la vari\'{e}t\'{e} symplectique $(M,\omega)$ et
 $i:L\ra M$
 une sous-vari\'{e}t\'{e} lagrangienne de $M$.
 \ben
  \item Il existe une re\-pr\'{e}\-sen\-ta\-tion diff\'{e}rentielle de $*$ sur $L$
   si et seulement si $i^*[*]=0$.
  \item Soit $\rho$ une re\-pr\'{e}\-sen\-ta\-tion diff\'{e}rentielle
  de $*$ sur $L$ et $*'$ un autre
    star-produit sur $M$ avec une representation diff\'{e}rentielle
    $\rho'$ sur $L$.
    Alors il existe un voisinage ouvert $U$ de $L$ tel que
    les couples $(*|_U,\rho|_{U\times L})$ et $(*'|_U,\rho'|_{U\times L})$
    sont \'{e}quivalents.
  \item Soit $\rho$ une re\-pr\'{e}\-sen\-ta\-tion diff\'{e}rentielle de $*$ sur $L$.
    Alors l'ensemble de toutes les classes d'isomorphisme de modules de $*$
     sur $L$ est en bijection avec
     \[
       \nu H^{1\prime}_{dR}(L,\korps)\oplus \nu^2 H^1_{dR}(L,\korps)[[\nu]],
    \]
    voir le th\'{e}or\`{e}me \ref{TRepLagTStarQ} pour la notation.
 \een
\ecor
\bbew
 Le th\'{e}or\`{e}me de Weinstein \cite{Wei71} mentionn\'{e} au d\'{e}but de ce paragraphe
 nous donne un ouvert $U$ de $M$ (qui est convexe fibre-par-fibre et contient
 $L$) et un symplectomorphisme $\phi:U\ra V\subset T^*L$ o\`{u} $V$ est un ouvert
 de $T^*L$ qui contient $L$ comme section nulle. On voit d'abord que
 $*$ est repr\'{e}sentable sur $L$ ssi $*|_U$ est repr\'{e}sentable sur $L$
 d'apr\`{e}s le corollaire \ref{CRepresentabiliteLocale}. Pour les vari\'{e}t\'{e}s $U$
 munie du star-produit $*|_U$ et $V$ munie du star-produit retir\'{e}
 $\phi^{-1~*}(*|_U)$ et l'application $\phi$ les
 hypoth\`{e}ses de la proposition \ref{PLiaisonEtRepresentabilite} sont
 satisfaites. Par cons\'{e}quent, $*|_U$ est repr\'{e}sentable sur $L$ ssi
 $\phi^{-1~*}(*|_U)$ est repr\'{e}sentable sur $L$. Puisque la cohomologie
 de de Rham de $U$ et de $V$ co\"{\i}ncident avec celle de $L$ o\`{u}
 l'injection $i$ de $L$ dans $U$ ou $V$ induit un isomorphisme, il
 s'ensuit que $\phi^{-1~*}(*|_U)$ est repr\'{e}sentable sur $L$ ssi
 la classe de Deligne de $[\phi^{-1~*}(*|_U)]$ s'annule d'apr\`{e}s le th\'{e}or\`{e}me
 \ref{TRepLagTStarQ}, donc ssi
 \[
    0=[\phi^{-1~*}(*|_U)]=[*|_U]=i^*[*].
 \]
 A l'aide d'un raisonnement analogue on trouve que les classes
 d'isomorphisme de $*$-modules sur $L$ sont en bijection avec
 celles de $\phi^{-1~*}(*|_U)$-modules sur $L$, alors le reste du
 corollaire se d\'{e}duit du th\'{e}or\`{e}me \ref{TRepLagTStarQ}.
\ebew

\subsection{D\'{e}formation de l'id\'{e}al annulateur comme sous-\-al\-g\`{e}bre}
  \label{SubSecISousAlg}

Soit $i:C\ra M$ une sous-vari\'{e}t\'{e} co\"{\i}sotrope connexe et ferm\'{e}e
de la vari\'{e}t\'{e}
symplectique $(M,\omega)$, soit $\mathcal{I}$ l'id\'{e}al annulateur de $C$
et soit $*$ un star-produit sur $M$. Le but de ce paragraphe est la
d\'{e}monstration du th\'{e}or\`{e}me suivant:

\bsat \label{TIAlgId}
 Avec les notations donn\'{e}es ci-dessus, soit $*$ tel que l'espace
 $\mathcal{I}[[\nu]]$ est une sous-alg\'{e}bre de l'alg\`{e}bre
 $\big(\CinfK{M}[[\nu]],*\big)$. Alors $\mathcal{I}[[\nu]]$ est ou bien
 un id\'{e}al \`{a} gauche ou bien un id\'{e}al \`{a} droite de
 $\big(\CinfK{M}[[\nu]],*\big)$.
\esat
\bbew
 Soit $*=\sum_{r=0}^\infty \nu^r \mathsf{C}_r$. D'apr\`{e}s l'hypoth\`{e}se on a
 $\mathsf{C}_r(g_1,g_2)\in\mathcal{I}$ quels que soient $g_1,g_2\in\mathcal{I}$
 et $r\in\nat$
 donc
 \beq \label{EqIAlgId1}
    i^*\mathsf{C}_r(g_1,g_2)=0~~~~~~~~~~~~~~~~
    \forall ~g_1,g_2\in\mathcal{I}~~\forall~r\in\nat.
 \end{equation}
 On va proc\'{e}der par r\'{e}currence: soient $f\in\CinfK{M}$ et $g,g_1,g_2\in
 \mathcal{I}$.\\
 0. Pour $r=0$, il est clair que
 $\mathsf{C}_0(f,g)=fg=\mathsf{C}_0(g,f)\in\mathcal{I}$.\\
 1. Le traitement du cas $r=1$, ce qui est le plus difficile, sera remis
    \`{a} la fin de cette d\'{e}monstration.\\
 2. Supposons d'abord que $\mathcal{I}[[\nu]]$ soit un id\'{e}al \`{a} gauche
    jusqu'\`{a} l'ordre $r\geq 1$, c-\`{a}-d.
    \beq \label{EqIAlgId2}
    i^*\mathsf{C}_s(f,g)=0~~~~~~~~~~~~~~
    \forall~f\in\CinfK{M}~~\forall g\in\mathcal{I}
       ~~\forall~0\leq s\leq r.
    \end{equation}
    L'associativit\'{e} \`{a} l'ordre $r+1$ de $*$, suivie de la restriction $i^*$,
    donne
    \beas
     0 & = & i^*\mathsf{C}_{r+1}(fg_1,g_2)-i^*\mathsf{C}_{r+1}(f,g_1g_1)
              +i^*\big(\mathsf{C}_{r+1}(f,g_1)g_2\big)
                -i^*\big(f\mathsf{C}_{r+1}(g_1,g_2)\big) \\
       &   & +\sum_{s=1}^r i^* \Big(
               \mathsf{C}_s\big(\mathsf{C}_{r+1-s}(f,g_1),g_2\big)
                 -\mathsf{C}_s\big(f,\mathsf{C}_{r+1-s}(g_1,g_2)\big)
                       \Big),
     \eeas
     d'o\`{u}, gr\^{a}ce \`{a} (\ref{EqIAlgId1}) et (\ref{EqIAlgId2}),
     \beq \label{EqIAlgId3}
       i^*\mathsf{C}_{r+1}(f,g_1g_1)=0 ~~~~~~~~~~~
                \forall~f\in\CinfK{M}~~\forall g_1,g_2\in\mathcal{I}
    \end{equation}
    Ensuite, l'associativit\'{e} \`{a} l'ordre $r+2$ de $*$, suivie de la restriction
    $i^*$, donne
    \beas
     0 & = & i^*\mathsf{C}_{r+2}(fg_1,g_2)-i^*\mathsf{C}_{r+2}(f,g_1g_1)
              +i^*\big(\mathsf{C}_{r+2}(f,g_1)g_2\big)
                -i^*\big(f\mathsf{C}_{r+2}(g_1,g_2)\big) \\
       &   & +i^*\mathsf{C}_{r+1}\big(\mathsf{C}_1(f,g_1),g_2\big)
             -i^*\mathsf{C}_{r+1}\big(f,\mathsf{C}_1(g_1,g_1)\big) \\
       &   &       +i^*\mathsf{C}_1\big(\mathsf{C}_{r+1}(f,g_1),g_2\big)
                -i^*\mathsf{C}_1\big(f,\mathsf{C}_{r+1}(g_1,g_2)\big) \\
       &   & +\sum_{s=2}^r i^* \Big(
               \mathsf{C}_s\big(\mathsf{C}_{r+2-s}(f,g_1),g_2\big)
                 -\mathsf{C}_s\big(f,\mathsf{C}_{r+2-s}(g_1,g_2)\big)
                       \Big)                 \\
       & = & - i^*\mathsf{C}_{r+2}(f,g_1g_1)
               -i^*\mathsf{C}_{r+1}\big(f,\mathsf{C}_1(g_1,g_2)\big)
               ~~~~
                \forall~f\in\CinfK{M}~~\forall g_1,g_2\in\mathcal{I}
    \eeas
    gr\^{a}ce \`{a} $(\ref{EqIAlgId1})$ et $(\ref{EqIAlgId2})$.
    En retranchant de cette \'{e}quation la m\^{e}me
    \'{e}quation avec $g_1$ et $g_2$ \'{e}chang\'{e}s, on obtient
    \beq \label{EqIAlgId4}
         i^*\mathsf{C}_{r+1}\big(f,\{g_1,g_2\}\big)=0
            ~~~~~~~~\forall~f\in\CinfK{M}~~\forall g_1,g_2\in\mathcal{I}.
    \end{equation}
    Puisque $\mathcal{I}/\mathcal{I}^2\cong \Ginf(C,TC^{\mathrm{ann}})\cong
    \Ginf(C,E)$ (via $[g]\mapsto -X_g|_C$) en tant qu'alg\`{e}bre de Lie,
    les \'{e}quations $(\ref{EqIAlgId1})$ et  $(\ref{EqIAlgId3})$ nous disent que
    $(f,g)\mapsto i^*\mathsf{C}_{r+1}(f,g)$ induit
    un op\'{e}rateur diff\'{e}rentiel $D_{r+1}$ de $\CinfK{C}\times \Ginf(C,E) \ra
    \CinfK{C}$. De plus, l'\'{e}quation (\ref{EqIAlgId4})
    entra\^{\i}ne que
    \[
         D_{r+1}(f,[V,W])=0~~~~~~~\forall ~V,W\in\Ginf(C,E)~~\forall
         ~f\in\CinfK{C}.
    \]
    Puisque $D_{r+1}$ en tant qu'op\'{e}rateur bidiff\'{e}rentiel est local, il
    est d\'{e}termin\'{e} par ses restrictions \`{a} toutes les cartes. Soit
    $\big(U,(\xi^1,\ldots,\xi^{2m},x^1,\ldots,x^k)\big)=:\big(U,(\xi,x)\big)$
    une carte distingu\'{e}e de la vari\'{e}t\'{e} feuillet\'{e}e $C$
    (dont les coordonn\'{e}es $\xi$ sont constantes sur les
    feuilles locales et les coordonn\'{e}es $x$ donnent un param\'{e}trage des
    feuilles locales). On peut en outre supposer que l'image de la carte
    soit un domaine convexe de $\real^{2m+k}$ qui contient l'origine.
    Soit $W=\sum_{i=1}^kW^i\frac{\partial}{\partial
    x^i}$ un champ vertical local, c.-\`{a}-d. un \'{e}l\'{e}ment de $\Ginf(U,E)$.
    On d\'{e}finit
    \[
        \hat{W}(\xi,x):=\sum_{i=1}^k\hat{W}^i(x,\xi)
                            \frac{\partial}{\partial x^i}
                      :=\sum_{i=1}^k\int_0^{x^1}W^i(\xi,t,x^2,\ldots,x^{2m})
                            \frac{\partial}{\partial x^i}.
    \]
    On calcule sans peine que
    \beq \label{EqVertCrochetVV}
       \left[\frac{\partial}{\partial x^1},\hat{W}\right]=W
    \end{equation}
    et la restriction de $D_{r+1}$ \`{a} $U$ s'annule sur $W=[V,\hat{W}]$
    car $V:=\frac{\partial}{\partial x^1}$ est \'{e}galement vertical. Par
    cons\'{e}quent, $D_{r+1}$ s'annule et donc
    \[
      i^*\mathsf{C}_{r+1}(f,g)=D_{r+1}(i^*f,-X_g|_C)=0,
    \]
    ce qui montre la r\'{e}currence \`{a} condition que le cas $r=1$ soit \'{e}tabli.\\
    Le cas o\`{u} l'\'{e}quation $(\ref{EqIAlgId2})$ est remplac\'{e}e par les
    conditions analogues
    pour un id\'{e}al \`{a} droite est trait\'{e} de mani\`{e}re enti\`{e}rement analogue
    \`{a} la fa\c{c}on pr\'{e}c\'{e}dente.\\
    3. Le cas $r=1$: dans la suite, soient $f,f_1,f_2\in\CinfK{M}$ et
    $g,g_1,g_2\in\mathcal{I}$.
    Il est connu (voir par exemple \cite{GR99}, Prop. 2.23)
    que l'op\'{e}rateur bidiff\'{e}rentiel $\mathsf{C}_1$
    est de la forme g\'{e}n\'{e}rale
    \beq \label{EqDefDecC1}
      \mathsf{C}_1(f_1,f_2)=\{f_1,f_2\}+ (\mathsf{b} B_1)(f_1,f_2)
                  := \{f_1,f_2\} + f_1 (B_1f_2) -B_1(f_1f_2) + (Bf_1)f_2.
    \end{equation}
    pour un op\'{e}rateur diff\'{e}rentiel $B_1:\CinfK{M}\ra\CinfK{M}$
    quelles que soient $f_1,f_2\in\CinfK{M}$.
    La partie sym\'{e}trique de (\ref{EqIAlgId1}) donne
    $i^*\big(B_1(g_1g_2)\big)=0$ quels que soient $g_1,g_2\in\mathcal{I}$,
    donc $i^*B_1$ induit un op\'{e}rateur diff\'{e}rentiel $B:\Ginf(C,E)\ra \CinfK{C}$
    par $B(-X_g|_C):=i^*B_1(g)$. Il s'av\`{e}rera utile de d\'{e}finir l'op\'{e}rateur
    bidiff\'{e}rentiel suivant: quels que soient $V\in\Ginf(C,E)$ et $\phi\in
    \CinfK{C}$
    \beq \label{EqDefDVphi}
       D_V\phi:=B(\phi V)-\phi B(V)
    \end{equation}
    Puisque pour tous $f\in \CinfK{M}$ et $g_1,g_2,g\in\mathcal{I}$ il vient
    \[
     i^*(\mathsf{b} B_1)(f,g_1g_2)=i^*\Big(fB_1(g_1g_2)\Big)
                            -i^*\Big(B_1\big((fg_1)g_2\big)\Big)
                                +i^*(B_1f)i^*(g_1g_2)=0
    \]
    (car $fg_1\in\mathcal{I}$ et $i^*B(\mathcal{I}^2)=\{0\}$), et puisque
    $i^*\big((\mathsf{b} B_1)(f,g)\big)$ ne d\'{e}pend que de $i^*f$ (car
    $i^*(\mathsf{b} B_1)(\mathcal{I},\mathcal{I})=\{0\}$) il vient
    \beq \label{EqDefDVphi2}
       D_V\phi:=-i^*(\mathsf{b} B_1)(f,g)~~~\mathrm{si}~~~i^*f=\phi~
                                     \mathrm{et}~V=-X_g|_C.
    \end{equation}
    En notant par $V$ \'{e}galement la d\'{e}riv\'{e}e de Lie par rapport au champ de
    vecteurs $V$, on a
    \beq \label{EqDefVphi2}
         V\phi = i^*\{g,f\} ~~~\mathrm{si}~~~i^*f=\phi~
                                     \mathrm{et}~V=-X_g|_C,
    \end{equation}
    donc on peut \'{e}crire
    \bea
       i^*\big(\mathsf{C}_1(f,g)\big) & = & -(V+D_V)\phi,
                                       \label{EqAlternative1}\\
       i^*\big(\mathsf{C}_1(g,f)\big) & = & (V-D_V)\phi.
                                         \label{EqAlternative2}
    \eea
    En regardant l'associativit\'{e} de $*$ \`{a} l'ordre $2$ pour
    $(f,g_1,g_2)$ on obtient
    \[
    0 = -i^*\mathsf{C}_2(f,g_1g_2)
                    + i^*\mathsf{C}_1\big(\mathsf{C}_1(f,g_1),g_2\big)
                         - i^*\mathsf{C}_1\big(f,\mathsf{C}_1(g_1,g_2)\big),
    \]
    alors, en retranchant de cette \'{e}quation celle avec $g_1$ et $g_2$
    \'{e}chang\'{e}es,
    \[
      2i^*\mathsf{C}_1(f,\{g_1,g_2\})
            =i^*\mathsf{C}_1\big(\mathsf{C}_1(f,g_1),g_2\big)
            -i^*\mathsf{C}_1\big(\mathsf{C}_1(f,g_2),g_1\big).
    \]
    En d\'{e}composant $\mathsf{C}_1$ \`{a} l'aide de (\ref{EqDefDecC1}) et
    en utilisant
    (\ref{EqDefDVphi2}) et (\ref{EqDefVphi2}) nous sommes men\'{e}s
    \`{a} l'\'{e}quation d'op\'{e}rateurs diff\'{e}rentiels suivante
    \beq \label{EqDVDWprelim}
      [D_V,D_W]-[V,W]+[V,D_W]+[D_V,W]-2D_{[V,W]}=0
    \end{equation}
    quels que soient les champs de vecteurs $V,W\in\Ginf(C,E)$.
    On regarde encore une fois l'associativit\'{e} de $*$ \`{a} l'ordre $2$,
    mais pour $(g_1,f,g_2)$ cette fois:
    \beas
     0 & = & i^*\mathsf{C}_2(g_1f,g_2)-i^*\mathsf{C}_2(g_1,fg_2)
        +i^*\big(\mathsf{C}_2(g_1,f)g_2\big)
        -i^*\big(g_1\mathsf{C}_2(f,g_2)\big) \\
       & &
          +i^*\mathsf{C}_1\big(\mathsf{C}_1(g_1,f),g_2\big)
          -i^*\mathsf{C}_1\big(g_1,\mathsf{C}_1(f,g_2)\big),
    \eeas
    alors, gr\^{a}ce \`{a} (\ref{EqIAlgId1}) il vient
    \[
       0 = i^*\mathsf{C}_1\big(\mathsf{C}_1(g_1,f),g_2\big)
           -i^*\mathsf{C}_1\big(g_1,\mathsf{C}_1(f,g_2)\big).
    \]
    A l'aide de (\ref{EqDefDecC1}), de (\ref{EqDefDVphi2}) et de
    (\ref{EqDefVphi2}) on arrive \`{a} l'\'{e}quation d'op\'{e}rateurs diff\'{e}rentiels
    suivante:
    \[
       [D_V,D_W]-[V,W]=[V,D_W]+[W,D_V]
    \]
    quels que soient les champs de vecteurs $V,W\in\Ginf(C,E)$. Puisque le
    membre gauche de cette \'{e}quation est antisym\'{e}trique par rapport \`{a}
    l'\'{e}change de $V$ et $W$, tandis que le membre droit est sym\'{e}trique
    par rapport \`{a} cet \'{e}change, il vient que les deux membres doivent
    s'annuler s\'{e}par\'{e}ment. Avec l'\'{e}quation (\ref{EqDVDWprelim}) on obtient
    les trois \'{e}quations suivantes
    \bea
     ~[V,D_W] & = & D_{[V,W]}, \label{EqVDWgDuVW} \\
     ~[V,D_W] & = & [D_V,W],  \label{EqVDWgDVW}\\
     ~[D_V,D_W] & = & [V,W],   \label{EqDVDWgVW}
    \eea
    quels que soient les champs de vecteurs $V,W\in\Ginf(C,E)$. La
    derni\`{e}re \'{e}quation (\ref{EqDVDWgVW}) montre en particulier que
    $D_V$ ne s'annule pas en g\'{e}n\'{e}ral.\\
    On \'{e}tudiera les trois derni\`{e}res \'{e}quations \`{a} l'aide d'un calcul
    symbolique dans une carte distingu\'{e}e $\big(U,(\xi,x)\big)$ de $C$
    mentionn\'{e}e ci-dessus:
    dans la suite, soient $(\eta,y)$ et $(\zeta,z)$ dans
    $\real^{2m}\times \real^k$. Soit $e_{(\eta,y)}$ la fonction
    exponentielle
    \[
     e_{(\eta,y)}(\xi,x):=\exp\left(\sum_{a=1}^{2m}\eta_a\xi^a+
                               \sum_{i=1}^ky_a x^a\right),
    \]
    et on d\'{e}finit le symbole de l'op\'{e}rateur diff\'{e}rentiel $B$ par
    \[
      B_i(\eta,y):=B\big(e_{(\eta,y)}\frac{\partial}{\partial x^i}\big)
                     e_{(-\eta,-y)}~~~~~~\forall~~ 1\leq i \leq k.
    \]
    Il est \'{e}vident que la restriction de $B$ \`{a} $U$ est d\'{e}termin\'{e} par son
    symbole et que $B_i$ est une application polyn\^{o}miale dans les
    variables $(\eta,y)$. Il s'ensuit que
    \beq \label{EqSymbDVPhi}
      D_{e_{(\eta,y)}\frac{\partial}{\partial x^i}}\big(
              e_{(\zeta,z)}\big)
              =
      \big(B_i(\eta+\zeta,y+z)-B_i(\eta,y)\big)
          e_{(\eta+\zeta,y+z)}.
    \end{equation}
    Puisque $[\frac{\partial}{\partial x^i},e_{(\eta,y)}
    \frac{\partial}{\partial x^j}]=y_i e_{(\eta,y)}
    \frac{\partial}{\partial x^j}$, l'\'{e}quation (\ref{EqVDWgDuVW}) entra\^{\i}ne
     pour le choix $V=\frac{\partial}{\partial x^i}$ et
     $W=e_{(\eta,y)}\frac{\partial}{\partial x^j}$ le suivant:
     \beas
        \lefteqn{
           y_i\big(B_j(\eta+\zeta,y+z)-B_j(\eta,y)\big)
             e_{(\eta+\zeta,y+z)} ~~=}~~~~~~~~~~~~~~~~~~~~~~~~~~~~ \\
          D_{[\frac{\partial}{\partial x^i},
                    e_{(\eta,y)}\frac{\partial}{\partial x^j}]}
                    e_{(\zeta,z)}
          & \stackrel{(\ref{EqVDWgDuVW})}{=} &
            \left[\frac{\partial}{\partial x^i},
               D_{e_{(\eta,y)}\frac{\partial}{\partial x^j}} \right]
               e_{(\zeta,z)}\\
          & = &
          (y_i+\frac{\partial}{\partial x^i})
            \big(B_j(\eta+\zeta,y+z)-B_j(\eta,y)\big)e_{(\eta+\zeta,y+z)}
     \eeas
     d'o\`{u}
    \[
      \frac{\partial\big(B_j(\eta+\zeta,y+z)-B_j(\eta,y)\big)}{\partial x^i}
         = 0.
    \]
    Par cons\'{e}quent, on a la d\'{e}composition $B_j(\eta,y)=\hat{B}_j(\eta,y)
    +B_j(0,0)$ o\`{u} $\hat{B}_j(\eta,y)$ ne d\'{e}pend plus de
    $x=(x^1,\ldots,x^k)$ et s'annule pour $(\eta,y)=(0,0)$. En outre, dans
    l'\'{e}quation pour le symbole de $D_V$, (\ref{EqSymbDVPhi}), on peut
    remplacer $B_i$ par $\hat{B}_i$. Ensuite, il vient
     \beas
        \lefteqn{
           y_i\big(\hat{B}_j(\eta+\zeta,y+z)-\hat{B}_j(\eta,y)\big)
             e_{(\eta+\zeta,y+z)} ~~=}~~~~~~~~~~~~~~~~~~~~~~~~~~~~ \\
          ~\left[\frac{\partial}{\partial x^i},
               D_{e_{(\eta,y)}\frac{\partial}{\partial x^j}} \right]
               e_{(\zeta,z)}
          & \stackrel{(\ref{EqVDWgDVW})}{=} &
          \left[D_{\frac{\partial}{\partial x^i}},
                e_{(\eta,y)}\frac{\partial}{\partial x^j} \right]
                e_{(\zeta,z)}\\
          & = &
            z_j\big(\hat{B}_i(\eta+\zeta,y+z)
                  -\hat{B}_i(\zeta,z)\big)e_{(\eta+\zeta,y+z)}
     \eeas
     d'o\`{u}
    \beq \label{EqWichtigerAusgangspunkt}
     y_i\big(\hat{B}_j(\eta+\zeta,y+z)-\hat{B}_j(\eta,y)\big)
       = z_j\big(\hat{B}_i(\eta+\zeta,y+z)
                  -\hat{B}_i(\zeta,z)\big).
    \end{equation}
    En particulier, le choix $y=0$ entra\^{\i}ne que $\hat{B_i}$ ne d\'{e}pend pas
    de $\eta$, alors il existe une unique application
    $E_i(y):=\hat{B}_i(0,y)$ polyn\^{o}miale en $y$. Dans l'\'{e}quation qui
    en r\'{e}sulte de (\ref{EqWichtigerAusgangspunkt}),
    \beq \label{EqWichtigerAusgangspunkt2}
      y_i\big(E_j(y+z)-E_j(y)\big)=z_j\big(E_i(y+z)-E_i(z)\big),
    \end{equation}
    la d\'{e}riv\'{e}e par rapport \`{a} $y_l$ en $y=0$ nous donne
    \beq \label{EqWichtigerAusgangspunkt3}
         \delta^l_iE_j(z)=z_j\frac{\partial E_i}{\partial z_l}(z)
    \end{equation}
    quels que soient $1\leq i,j,l \leq k$. La somme sur $i=l$ de $1$
    jusqu'\`{a} $k$ fait appara\^{\i}tre le champ d'Euler
    $\sum_{i=1}^kz_k\frac{\partial}{\partial z_k}$ dans le membre droit et
    montre que chaque $E_i$ est un polyn\^{o}me lin\'{e}aire en $z$,
    et l'\'{e}quation (\ref{EqWichtigerAusgangspunkt3}) pour le cas
    $i\neq l$ entra\^{\i}ne que $E_i$ ne d\'{e}pend que de la variable $z_i$.
    Par cons\'{e}quent, il existe une unique fonction $f_i$ (ne d\'{e}pendant
    que des variables $\xi$) telle que $E_i(z)=f_iz_i$. Finalement,
    si on remplace $E_i(y)$ par $f_iy_i$ dans l'\'{e}quation
    (\ref{EqWichtigerAusgangspunkt2}) on voit que $f_i=f_j=:f$ quels que
    soient $1\leq i,j\leq k$. Il s'ensuit
    \[
        D_V|_U=f V|_U.
    \]
    La troisi\`{e}me \'{e}quation (\ref{EqDVDWgVW}) montre que
    \[
       ~[V,W]|_U=[D_V,D_W]|_U=[fV,fW]|_U=f^2[V,W]|_U
    \]
    quels que soient les champs de vecteurs $V,W\in\Ginf(C,E)$, car
    $Vf=0$ ($f$ n'est une fonction que des variables $\xi$). Puisque cet
    espace de champs de vecteurs est localement donn\'{e} par des sommes
    finies des crochets de Lie $[V,W]$ (comme on a montr\'{e} pr\'{e}c\'{e}demment,
    voir (\ref{EqVertCrochetVV})) il vient que $f^2=1$, et puisque $C$
    est connexe, il
    ne restent que les possibilit\'{e}s suivantes
    \[
       \forall V\in \Ginf(C,E): D_V=V~~~~~\mathrm{ou}~~~~~
       \forall V\in \Ginf(C,E): D_V=-V.
    \]
    Dans le premier cas, $I[[\nu]]$ est un id\'{e}al \`{a} droite jusqu'\`{a} l'ordre $1$
    d'apr\`{e}s l'\'{e}quation (\ref{EqAlternative2}), et dans le deuxi\`{e}me cas
    un id\'{e}al \`{a} gauche jusqu'\`{a} l'ordre $1$, voir (\ref{EqAlternative1}).\\
    Ceci finit la d\'{e}monstration.
\ebew

\subsection{R\'{e}duction (et bimodules) des star-produits quand l'espace
            r\'{e}duit existe}
  \label{SubSecEspRedReduction}

\bsat \label{TRepUndAllesEspRedExist}
 Soit $(M,\omega)$ une vari\'{e}t\'{e} symplectique munie d'un star-produit $*$.
 On note $[*]$ sa classe de Deligne.
 Soit $i:C\ra M$ une sous-vari\'{e}t\'{e} co\"{\i}sotrope ferm\'{e}e et connexe
 de $M$ telle que la vari\'{e}t\'{e} symplectique r\'{e}duite
 $\pi:C\ra(M_{\mathrm{red}},\omega_{\mathrm{red}})$
 (voir th\'{e}or\`{e}me \ref{RedSymp}) existe. Alors les \'{e}nonc\'{e}s suivants sont
 \'{e}quivalents:
 \begin{description}
  \item[i.]
   Il existe
   une s\'{e}rie formelle $\beta$ \`{a} coefficients dans le deuxi\`{e}me groupe de
   cohomologie de de Rham de $M_{\mathrm{red}}$ telle que
   \[
     i^*[*]=\pi^*\left(\frac{[\omega_{\mathrm{red}}]}{\nu}+\beta\right).
   \]
  \item[ii.] Le star-produit $*$ est r\'{e}ductible.
  \item[iii.] Il existe un star-produit $*_r$ sur $M_{\mathrm{red}}$ et
    la structure d'un $*$-$*_r$-bimodule diff\'{e}rentiel (ou d'un
    $*_r$-$*$-bimodule diff\'{e}rentiel) sur $\CinfK{C}[[\nu]]$
    (voir la d\'{e}finition \ref{DStarStarredBimod}).
 \end{description}
 Si une de ces conditions est satisfaite, alors
 \ben
  \item  La classe de Deligne de $[*_r]$ satisfait \`{a} la condition
  suivante:
   \[
        i^*[*]=\pi^*[*_r]
   \]
  \item Soit $*_{\mathrm{red}}$ un star-produit sur $M_{\mathrm{red}}$ qui
    est \'{e}quivalent au star-produit $*_r$ et soit $\rho$ la repr\'{e}sentation
    diff\'{e}rentielle de $*$ sur $C$ provenant de la structure du
    $*$-$*_r$-bimodule (resp $*_r$-$*$-bimodule) diff\'{e}rentiel.
    Alors l'alg\`{e}bre r\'{e}duite est antiisomorphe au commutant
   du module \`{a} gauche (resp. du module \`{a} droite) de
   $(\CinfK{M}[[\nu]],*)$:
   \[
    \big(\CinfK{M_{\mathrm{red}}}[[\nu]],
                 *_{\mathrm{red}}^{\mathrm{opp}}\big)
       \cong
               \mathrm{Hom}_{(*,\rho)}\big(C,C \big).
   \]
  \item L'ensemble de toutes les classes d'isomorphismes de
   $*$-$*_r$--bimodules diff\'{e}rentiels sur
  $\CinfK{C}[[\nu]]$ est en bijection avec l'espace
   \[
      \nu H^{1\prime}_{dR}(C,\korps)\oplus \nu^2 H^1_{dR}(C,\korps)[[\nu]],
   \]
   voir le th\'{e}or\`{e}me \ref{TRepLagTStarQ} pour la notation.
 \een
\esat
\bbew
 ``{\bf i.} $\Longleftrightarrow$ {\bf iii.}'': il est d'abord
  \'{e}vident que l'existence des re\-pr\'{e}\-sen\-ta\-tions diff\'{e}rentielles
  $\rho$ et $\tilde{\rho}$ satisfaisant la propri\'{e}t\'{e} (\ref{EqRhoRhoComm})
  est \'{e}quivalente \`{a} l'existence $*$-$*_r$-bimodule diff\'{e}rentiel, i.e.
  d'une seule re\-pr\'{e}\-sen\-ta\-tion
  diff\'{e}rentielle pour l'alg\`{e}bre $\CinfK{M\times M_{\mathrm{red}}}[[\nu]]$
  munie
  du star-produit $*\otimes *_r^{\mathrm{opp}}$. L'\'{e}nonc\'{e} analogue pour
  les $*_r$-$*$-bimodules diff\'{e}rentiels est \'{e}vident. \\
 On consid\`{e}re la vari\'{e}t\'{e} symplectique $(M\times M_{\mathrm{red}},
 \omega_{(1)}-\omega_{\mathrm{red}(2)})$. On rappelle que
  \[
    j:C\ra M\times M_{\mathrm{red}}:c\mapsto \big(i(c),\pi(c)\big)
  \]
  est un plongement de $C$ comme sous-vari\'{e}t\'{e} co\"{\i}sotrope
  dans $M\times M_{\mathrm{red}}$, voir le lemme \ref{LCLagdansMMred}.
   Soit $*_r$ un star-produit arbitraire
  sur $M_{\mathrm{red}}$, donc on peut consid\'{e}rer le star-produit
  $*\otimes *_r^{\mathrm{opp}}$ sur $M\times M_{\mathrm{red}}$.
  Sa classe de Deligne vaut (voir la proposition
  \ref{PClasseDeligneProduit Tensoriel})
  \[
      [*\otimes *_r^{\mathrm{opp}}]=pr_1^*[*]-pr_2^*[*_r]
  \]
  o\`{u} $pr_1:M\times M_{\mathrm{red}}\ra M$ et $pr_2:M\times M_{\mathrm{red}}
  \ra M_{\mathrm{red}}$ d\'{e}signent les projections canoniques. Puisque
  $pr_1\circ j=i$ et $pr_2\circ j=\pi$ il vient
  \[
    j^*[*\otimes *_r^{\mathrm{opp}}]=
        j^*pr_1^*[*]-j^*pr_2^*[*_r]=i^*[*]-\pi^*[*_r]
  \]
  D'apr\`{e}s le corollaire \ref{CorRepLagArbi} il existe une
  re\-pr\'{e}\-sen\-ta\-tion de
  $*\otimes *_r^{\mathrm{opp}}$ sur $C$ si et seulement si
  $j^*[*\otimes *_r^{\mathrm{opp}}]=0$ d'o\`{u} l'\'{e}nonc\'{e} 1.
  D'apr\`{e}s les th\'{e}or\`{e}mes d'existence et de classification des star-produits
  symplectiques il existe toujours un star-produit $*_r$ sur
  $M_{\mathrm{red}}$
  dont la classe de Deligne est \'{e}gale \`{a} $\frac{[\omega_{\mathrm{red}}]}{\nu}$
  plus une classe $\beta\in H^2_{dR}(M_{\mathrm{red}},\korps)[[\nu]]$
  arbitraire. Ceci montre l'\'{e}quivalence avec la condition {\bf i.}.
  L'ensemble des classes d'isomorphismes de bimodules diff\'{e}rentiels
  correspondent exactement
  \`{a} celui des classes d'isomorphismes de modules diff\'{e}rentiels pour
  $*\otimes *_r^{\mathrm{opp}}$. Ce dernier est en bijection avec
  l'espace du premier groupe de la cohomologie de de Rham annonc\'{e}, d'apr\`{e}s
   le corollaire \ref{CorRepLagArbi}. Ceci montre l'\'{e}nonc\'{e} 3.\\
   Le cas du $*_r$-$*$-bimodule diff\'{e}rentiel se d\'{e}duit de fa\c{c}on
   enti\`{e}rement analogue.

  ``{\bf ii.} $\Longrightarrow$ {\bf iii.}'': Soit $*$ un star-produit
  r\'{e}ductible sur $M$. Par d\'{e}finition, il existe un star-produit
  projetable $\hat{*}$ qui est \'{e}quivalent \`{a} $*$. L'espace
  $\mathcal{I}[[\nu]]$ (o\`{u} $\mathcal{I}$ est l'id\'{e}al annulateur de $C$)
  est donc une sous-alg\`{e}bre de
  $\big(\CinfK{M}[[\nu]],\hat{*}\big)$. Puisque $C$ est connexe,
  il vient que $\mathcal{I}[[\nu]]$ est automatiquement ou bien un id\'{e}al
  \`{a} gauche ou bien un id\'{e}al \`{a} droite d'apr\`{e}s le th\'{e}or\`{e}me \ref{TIAlgId}.
   On suppose d'abord que $\mathcal{I}[[\nu]]$ soit un id\'{e}al \`{a} gauche. Soient
   $f\in \CinfK{M}[[\nu]]$ et $\phi\in\CinfK{C}[[\nu]]$. Soit
   $f'\in \CinfK{M}[[\nu]]$ telle que $i^*f'=\phi$. On d\'{e}finit
   \[
      \hat{\rho}(f)\phi:=i^*(f\hat{*}f').
   \]
   Puisque $\mathcal{I}[[\nu]]$ est un id\'{e}al \`{a} gauche cette d\'{e}finition
   ne d\'{e}pend pas du repr\'{e}sentant $f'$ modulo $\mathcal{I}[[\nu]]$ de $\phi$
   choisi. L'associativit\'{e} et la nature bidiff\'{e}rentielle
   de $\hat{*}$ et le fait que $f\hat{*}f'$ est un repr\'{e}sentant modulo
   $\mathcal{I}[[\nu]]$ de $\hat{\rho}(f)\phi$
   montrent que $\hat{\rho}$ est une re\-pr\'{e}\-sen\-ta\-tion diff\'{e}rentielle
   de $\hat{*}$
   sur $C$ qui d\'{e}forme visiblement $i^*$. En utilisant la transformation
   d'\'{e}quivalence entre $*$ et $\hat{*}$ on obtient une re\-pr\'{e}\-sen\-ta\-tion
   diff\'{e}rentielle $\rho$ de $*$ qui d\'{e}forme $i^*$.
   De plus, soit $\mathcal{N}(\mathcal{I})$ l'id\'{e}alisateur de
   $\mathcal{I}$, alors soit $\tilde{h}\in
   \CinfK{M_{\mathrm{red}}}[[\nu]]$ et $h\in
   \mathcal{N}(\mathcal{I})[[\nu]]$ tel que $i^*h=\pi^*\tilde{h}$. On
   d\'{e}finit
   \[
      \tilde{\rho}(\tilde{h})\phi := i^*(f'\hat{*}h).
   \]
   Puisque $\mathcal{I}[[\nu]]$ est un id\'{e}al bilat\`{e}re dans l'espace
   $\mathcal{N}(\mathcal{I})[[\nu]]$ qui est une sous-alg\`{e}bre de
   $\big(\CinfK{M}[[\nu]],\hat{*}\big)$ (par d\'{e}finition de la
   projetabilit\'{e} de $\hat{*}$)
   $\tilde{\rho}$ est une application bien d\'{e}finie. L'associativit\'{e} et la
   nature bidiff\'{e}rentielle de
   $\hat{*}$ et le fait que $f'\hat{*}h$ est un repr\'{e}sentant modulo
   $\mathcal{I}[[\nu]]$ de $\tilde{\rho}(\tilde{h})\phi$
   montrent que $\tilde{\rho}$ est une re\-pr\'{e}\-sen\-ta\-tion diff\'{e}rentielle
   du star-produit r\'{e}duit
   $*_{\mathrm{red}}^{\mathrm{opp}}$ sur $C$ qui d\'{e}forme $\pi^*$. On a
   \[
      \hat{\rho}(f)\tilde{\rho}(\tilde{h})\phi
     = i^*(f\hat{*}f'\hat{*}h) = \tilde{\rho}(\tilde{h})\hat{\rho}(f)\phi,
   \]
   ce qui montre l'\'{e}nonc\'{e} {\bf iii.}\\
   Le cas o\`{u} $\mathcal{I}[[\nu]]$ est un id\'{e}al \`{a} droite d\'{e}finit un module
   \`{a} droite de $(\CinfK{M}[[\nu]],\hat{*})$, et on d\'{e}duit la structure
   d'un $*_r$-$*$-bimodule de mani\`{e}re enti\`{e}rement analogue.

   ``{\bf ii.} $\Longleftarrow$ {\bf iii.}'': Soient $\rho$ et
   $\tilde{\rho}$ les re\-pr\'{e}\-sen\-ta\-tions diff\'{e}rentielles
   de $*$ et $*_r^{\mathrm{opp}}$, respectivement, qui proviennent de la
   structure du $*$-$*_r$-bimodule diff\'{e}rentiel sur
   $\CinfK{C}[[\nu]]$. Alors l'application
   \[
      T:\CinfK{M_{\mathrm{red}}}[[\nu]]\ra \CinfK{C}[[\nu]]:
          \tilde{h}\ra \tilde{\rho}(\tilde{h})1
   \]
   est une s\'{e}rie d'op\'{e}rateurs diff\'{e}rentiels le long de la projection $\pi$.
   D'apr\`{e}s le lemme 2.2 de \cite[p.11]{KMS93}, il existe un atlas
   $\big(U_\alpha,(\xi_\alpha,x_\alpha)\big)_{\alpha\in\mathfrak{S}}$ de $C$,
   un atlas $\big(U'_\beta,\xi'_\beta\big)_{\beta\in\mathfrak{S}'}$ de
   $M_{\mathrm{red}}$ et une application $p:\mathfrak{S}\ra\mathfrak{S}'$
   tels que $\pi(U_\alpha)\subset U'_{p(\alpha)}$ et
   $\xi_{p(\alpha)}\circ \pi= \xi_{\alpha}$. D'apr\`{e}s
   la structure locale (\ref{EqDefFormeLocale}) de $T$ chaque coefficient
   de $T=\sum_{r=0}^\infty \nu^r T_r$ est de la forme locale suivante
   \[
       T_r\big(\tilde{h}\big)(u) =
         \sum_{|I|\leq N_r} T_{r\alpha}^I(u)
              \frac{\partial^{|I|}f}{\partial
              \xi^I_{p(\alpha)}}\big(\pi(u)\big)
        =\sum_{|I|\leq N_r} T_{r\alpha}^I(u)
              \frac{\partial^{|I|}(\pi^*\tilde{h})}{\partial
              \xi^I_{\alpha}}(u)
   \]
   quels que soient $u\in U_\alpha$ et
   $\tilde{h}\in\CinfK{M_{\mathrm{red}}}$.
   Soit $(\chi_\alpha)_{\alpha\in\mathfrak{S}}$
   une partition de l'unit\'{e}
   subordonn\'{e}e \`{a} $\big(U_\alpha\big)_{\alpha\in\mathfrak{S}}$. Alors il
   s'ensuit que les op\'{e}rateurs diff\'{e}rentiels
   \[
      D_r:=\sum_{\alpha\in\mathfrak{S}}\chi_\alpha
             \sum_{|I|\leq N_r} T_{r\alpha}^I
              \frac{\partial^{|I|}}{\partial
              \xi^I_{\alpha}}
   \]
   dans $\Dop^1\big(\CinfK{C};\CinfK{C}\big)$
   sont tels que $D:=\sum_{r=0}^\infty \nu^r D_r$ satsifait \`{a} l'\'{e}quation
   \[
      T\tilde{h}=D\pi^*\tilde{h}.
   \]
   Puisque $T_0=\pi^*$ et $T(1)=1$ il s'ensuit que
   $D$ est \'{e}gale \`{a} l'application identique \`{a} l'ordre $\nu^0$ et est donc
   inversible avec des termes d'ordre sup\'{e}rieur qui s'annulent sur les
   constantes. Soit $f\in\CinfK{M}[[\nu]]$. On d\'{e}finit
   $\rho'(f):=D^{-1}\rho(f)D$ et $\tilde{\rho}'(\tilde{h}):=
   D^{-1}\tilde{\rho}(\tilde{h})D$. Il est clair que le couple
   $(\rho',\tilde{\rho}')$ d\'{e}finit \'{e}galement la structure d'un
   $*$-$*_r$--bimodule sur $C$. De plus,
   \beq\label{EqRhoprimePistar}
       \tilde{\rho}'(\tilde{h})1=D^{-1}\tilde{\rho}(\tilde{h})D1
                  = D^{-1}\tilde{\rho}(\tilde{h})1
                  = D^{-1}D\pi^*\tilde{h}=\pi^*\tilde{h}.
   \end{equation}
   D'apr\`{e}s la proposition \ref{PKerRho1congId} il existe une
   transformation d'\'{e}quivalence $S$ de $*$ telle que pour le star-produit
   $*':=S(*)$ et la re\-pr\'{e}\-sen\-ta\-tion $\rho''(f):=\rho'(S^{-1}f)$ on a
   pour tout $f'\in\CinfK{M}[[\nu]]$:
   \beq\label{EqRhoprimeprimeIstar}
       \rho''(f)i^*f'=i^*(f*'f').
   \end{equation}
   De la m\^{e}me proposition \ref{PKerRho1congId} on d\'{e}duit que
   $\mathcal{I}[[\nu]]$ est un id\'{e}al \`{a} gauche pour le star-produit $*'$.
   Soient $h_1,h_2\in\mathcal{N}(\mathcal{I})[[\nu]]$. Alors il existent des
   uniques fonctions $\tilde{h}_1,\tilde{h}_2\in
   \CinfK{M_{\mathrm{red}}}[[\nu]]$ telles que $i^*h_1=\pi^*\tilde{h_1}$
   et $i^*h_2=\pi^*\tilde{h_2}$. Par cons\'{e}quent
   \beas
    i^*(h_1*'h_2) & \stackrel{(\ref{EqRhoprimeprimeIstar})}{=}  &
                   \rho''(h_1*'h_2)1 = \rho''(h_1)\rho''(h_2)1
                          = \rho''(h_1)i^*h_2 \\
                  & = & \rho''(h_1)\pi^*\tilde{h}_2
                    \stackrel{(\ref{EqRhoprimePistar})}{=}
                        \rho''(h_1)\tilde{\rho}'(\tilde{h}_2)1
                        =\tilde{\rho}'(\tilde{h}_2)\rho''(h_1)1 \\
                  & \stackrel{(\ref{EqRhoprimeprimeIstar})}{=} &
                        \tilde{\rho}'(\tilde{h}_2)i^*h_1
                        = \tilde{\rho}'(\tilde{h}_2)\pi^*\tilde{h}_1
                       \stackrel{(\ref{EqRhoprimePistar})}{=}
                  \tilde{\rho}'(\tilde{h}_2)\tilde{\rho}'(\tilde{h}_1)1 \\
                  & = & \tilde{\rho}'(\tilde{h}_1 *_r \tilde{h}_2)1
                        = \pi^*(\tilde{h}_1 *_r \tilde{h}_2).
   \eeas
   Cette \'{e}quation montre en m\^{e}me temps que $\mathcal{N}(\mathcal{I})[[\nu]]$
   est une sous-alg\`{e}bre pour le star-produit $*'$, que
   $\mathcal{I}[[\nu]]$ est un id\'{e}al bilat\`{e}re dans
   $\mathcal{N}(\mathcal{I})[[\nu]]$ (car $\tilde{h}_1=0$ pour $h_1$ dans
   $\mathcal{I}[[\nu]]$) et que $*_r$ est le star-produit r\'{e}duit du
   star-produit $*'$ qui s'av\`{e}re projetable. Par d\'{e}finition, $*$ est
   r\'{e}ductible.\\
   On arrive au m\^{e}me r\'{e}sultat pour un $*_r$-$*$-bimodule diff\'{e}rentiel.

   Finalement, l'\'{e}nonc\'{e} 2. se d\'{e}duit directement du fait que le star-produit
   $*'$ du raisonnement pr\'{e}c\'{e}dent est projetable et que
   $\CinfK{M_{\mathrm{red}}}[[\nu]],*_{\mathrm{red}})$ est isomorphe \`{a}
   $\mathcal{N}(\mathcal{I})[[\nu]]/\mathcal{I}[[\nu]]$. Puisque
   $\mathcal{N}(\mathcal{I})[[\nu]]$ est \'{e}gal \`{a}
   $\mathcal{N}_{*'}(\mathcal{I})[[\nu]]$ d'apr\`{e}s l'\'{e}nonc\'{e} 2. de
   la proposition \ref{PCommutantNdeI} l'isomorphisme avec le commutant
   est \'{e}vident.
%   {\tt na, was jetzt kommt, kann gek\"{u}rzt werden}
%
%  D'un autre c\^{o}t\'{e}, pour tout
%   $\tilde{h}\in \CinfK{M_{\mathrm{red}}}[[\nu]]$ il existe $h\in
%    \mathcal{N}(\mathcal{I})[[\nu]]$ tel que $i^*h=\pi^*\tilde{h}$
%    et la prescription
%    \[
%            \tilde{h}\mapsto \big(D_{\tilde{h}}:i^*f\mapsto i^*(f*h)\big)
%    \]
%    est une application bien d\'{e}finie de $\CinfK{M_{\mathrm{red}}}[[\nu]]$
%    dans l'espace des homomomorphismes de $*$-modules. On calcule sans
%    peine que $D_{\tilde{h_1}*_r\tilde{h_2}} =
%    D_{\tilde{h}_2}D_{\tilde{h}_1}$ quels que soient
%    $\tilde{h}_1,\tilde{h}_2\in\CinfK{M_{\mathrm{red}}}[[\nu]]$. Donc
%    on a montr\'{e} l'isomorphisme souhait\'{e}.\\
    On obtient le m\^{e}me r\'{e}sultat en supposant que $\mathcal{I}[[\nu]]$
    soit un id\'{e}al \`{a} droite.
\ebew

\brem
 Il est possible qu'un star-produit r\'{e}ductible puisse
\^{e}tre \'{e}qui\-val\-ent \`{a} plusieurs star-produits projetables distincts dont
les star-produits r\'{e}duits sont deux-\`{a}-deux non-equivalents, voir
l'exemple de la r\'{e}duction de l'espace projectif complexe en paragraphe
\ref{SubSecEspProj}.
Dans ce cas, les structures de bimodules et m\^{e}mes les structures de modules
(\`{a} gauche ou \`{a} droite) qu'on construit \`{a} l'aide des star-produits projetables
sont non-isomorphes et ont des commutants non-isomorphes.\\
D'un autre c\^{o}t\'{e}, si $H^1_v(C,\korps)=\{0\}$ (ce qui n'est pas le cas pour
la fibration de Hopf $S^{2n+1}\ra \complex P(n)$), alors deux star-produits
repr\'{e}sentables \'{e}quivalents ont des repr\'{e}sentations \'{e}quivalentes
(voir le corollaire \ref{Cdaserstemalfalsch}), et la structure du commutant
ne d\'{e}pend que du star-produit $*$. Un calcul rapide montre que l'application
$p^*$ induit une application {\em injective} de
$H^2_{dR}(M_{\mathrm{red}},\korps)$ dans $H^2_{dR}(C,\korps)$ si
$H^1_v(C,\korps)=\{0\}$, donc la classe de Deligne $[*_r]$ est uniquement
d\'{e}termin\'{e}e par $[*]$ via la condition $p^*[*_r]=i^*[*]$.
\erem

\brem
 En \'{e}crivant $R:=(\CinfK{M}[[\nu]],*)$ et
$B:=\CinfK{C}[[\nu]]$ ce th\'{e}or\`{e}me nous permet d'identifier
 $R_{\mathrm{red}}:=
\big(\CinfK{M_{\mathrm{red}}}[[\nu]],*_r^{\mathrm{opp}}\big)$,
i.e. l'alg\`{e}bre r\'{e}duite,
avec l'anneau de tous les $R$-homomorphismes $B\ra B$, alors
$R_{\mathrm{red}}\cong \mathrm{Hom}_{R}(B,B)$. Le triplet
$(R,R_{\mathrm{red}},B)$ est donc le commencement d'un {\em contexte de
Morita}, voir par exemple \cite[p.485]{Lam99}. Il reste \`{a} calculer
le module $R_{\mathrm{red}}$-$R$-bimodule $\mathrm{Hom}_{R}(B,R)$:
soit $\Psi:\CinfK{C}[[\nu]]\ra\CinfK{M}[[\nu]]$ un morphisme de
$R$-modules. Puisque $i^*$ est surjective, l'application $\Psi$ est
d\'{e}termin\'{e}e par sa valeur $h$ en $1\in B$, \`{a} savoir $\Psi(\phi)=
\Psi(i^*f)=\Psi\big(\rho(f)1\big)=f*h$. En particulier, pour tout $f=g\in
\mathcal{I}$ il vient $g*h=0$ car $i^*g=0$. Si l'id\'{e}al annulateur
ne s'annule pas, alors \`{a} l'ordre $r=0$ ceci entra\^{\i}ne
$gh_0=0$, donc $h_0=0$ parce dans une carte on a $y_ih_0=0$ quels que
soient $1\leq i\leq k$ pour les coordonn\'{e}es transversales $y$. Par
r\'{e}currence il s'ensuit $h=0$, donc
\[
  \mathrm{Hom}_{R}(B,R)=\{0\},\mathrm{~~si~}\mathcal{I}\neq\{0\}
\]
et le contexte de Morita est donn\'{e} par le quadruplet
$(R,R_{\mathrm{red}},B,\{0\})$. Pour le cas $\mathcal{I}=\{0\}$ il
vient $C=M$ si $M$ est connexe, et on obtient $(R,R^{\mathrm{opp}},R,R)$
comme contexte de Morita.
\erem

\section{Obstructions et feuilletages}
   \label{SecObstrucFeuilletage}

\subsection{Obstructions r\'{e}currentes pour l'existence des
                re\-pr\'{e}\-sen\-tations}
     \label{SubSecObstRecurExisRep}

Soit $(M,P)$ une vari\'{e}t\'{e} de Poisson, $i:C\ra M$ une sous-vari\'{e}t\'{e}
co\"{\i}sotrope ferm\'{e}e et $*=\sum_{r=0}^\infty \nu^r \mathsf{C}_r$ un
star-produit
sur $M$. Soit $\mathcal{I}$ l'id\'{e}al annulateur de $C$. D'apr\`{e}s le th\'{e}or\`{e}me
\ref{TIdAnIdGauche} il existe une re\-pr\'{e}\-sen\-ta\-tion diff\'{e}rentielle de
$*$ dans
 $C\times\korps$ si et seulement si $*$ est \'{e}quivalent \`{a} un star-produit
 $*'$ qui est adapt\'{e} \`{a} $C$, c.-\`{a}-d. pour
lequel $\mathcal{I}[[\nu]]$ est un id\'{e}al \`{a} gauche. Dans ce paragraphe,
on va \'{e}tudier les obstructions r\'{e}currentes pour la construction ordre par
ordre de $*'$.

On d\'{e}finit la suite d'applications bilin\'{e}aires
$(B_r)_{r\in\nat}:\CinfK{M}\times \mathcal{I}\ra \CinfK{C}$ par
\beq
      B_r(f,g):= i^*\mathsf{C}_r(f,g)~~\mathrm{quels~que~soient~}f\in\CinfK{M},
      g\in\mathcal{I}
\end{equation}
et $B^-_{r+1}:\mathcal{I}\times \mathcal{I}\ra \CinfK{C}$ par
\beq
  B^-_{r+1}(g_1,g_2):=B_{r+1}(g_1,g_2)-B_{r+1}(g_2,g_1)
                         ~~\mathrm{quels~que~soient~}g_1,g_2\in\mathcal{I}
\end{equation}
Puisque $\mathcal{I}$ est un id\'{e}al de l'alg\`{e}bre associative commutative
$\CinfK{M}$ il est clair que $B_0=0$.

On suppose d\'{e}sormais que le star-produit $*$ soit {\em adapt\'{e} \`{a} $C$ jusqu'\`{a}
l'ordre $r$}, c.-\`{a}-d. que $*$ soit tel que
\beq
  B_0=0,\cdots,B_r=0
\end{equation}
pour un entier positif $r$ donn\'{e}.

\noindent L'associativit\'{e} de $*$ entra\^{\i}ne le lemme suivant:
\blem\label{LObstructionsRecursives}
  Pour tous $f,f_1,f_2\in\CinfK{M}$ et $g,g_1,g_2,g_3\in\mathcal{I}$
  on a les \'{e}quations suivantes:
  \bea
   (i^*f_1) B_{r+1}(f_2,g) -B_{r+1}(f_1f_2,g)+B_{r+1}(f_1,f_2g)
        &   =   &     0,   \label{EqLOb1} \\
   B^-_{r+1}(fg_1,g_2)=B^-_{r+1}(g_1,fg_2) & = & (i^*f)B^-_{r+1}(g_1,g_2),
                           \nonumber \\
        & & \label{EqLOb2} \\
   B^-_{r+1}(g_1g_2,g_3)
        &   =   &     0 \label{EqLOb3}
  \eea
  et
  \bea
   -X_{g_1}\big(B^-_{r+1}(g_2,g_3)\big)
   -X_{g_2}\big(B^-_{r+1}(g_3,g_1)\big)
   -X_{g_3}\big(B^-_{r+1}(g_1,g_2)\big)
          &   &  \nonumber \\
                -B^-_{r+1}(\{g_1,g_2\},g_3)-B^-_{r+1}(\{g_2,g_3\},g_1)
                 -B^-_{r+1}(\{g_3,g_1\},g_2)
        &   =   & 0. \nonumber \\
        & & \label{EqLOb4}
  \eea
  o\`{u} $X_{g_1}$ etc. d\'{e}signe le champ hamiltonien de $g_1$ (qui --en tant
  que champ de vecteurs vertical-- agit sur $\CinfK{C}$.
\elem
\bbew
 Pour $f_1,f_2\in\CinfK{M}$ et $g\in\mathcal{I}$ l'associativit\'{e} de $*$
 \[
     i^*\big((f_1*f_2)*g\big)-i^*\big(f_1*(f_2*g)\big)=0
 \]
 s'\'{e}crit \`{a} l'ordre $r+1$
 \beq \label{EqAssOrdr}
   \sum_{a=0}^{r+1}\left(i^*\mathsf{C}_a\big(\mathsf{C}_{r+1-a}(f_1,f_2),g\big)
        -i^*\mathsf{C}_a\big(f_1,\mathsf{C}_{r+1-a}(f_2,g)\big)\right)
                         = 0.
 \end{equation}
 Par hypoth\`{e}se il vient
 $\mathsf{C}_a(f,g')\in\mathcal{I}$ quel que soit $g'\in\mathcal{I}$ pour tout
 $0\leq a\leq r$. Il s'ensuit que seuls les termes de la somme ci-dessus
 avec $a=0$ et $a=r+1$ ne sont pas forc\'{e}ment nuls, d'o\`{u} l'\'{e}nonc\'{e}
 (\ref{EqLOb1}) car $i^*g=0$.\\
 Le cas particulier $f_1=g_1\in\mathcal{I}$,
 $f_2=f$ et $g=g_2$  du premier \'{e}nonc\'{e} nous donne
 \beq \label{EqBTens}
    B_{r+1}(g_1f,g_2)=B_{r+1}(g_1,fg_2)
 \end{equation}
  dont la partie antisym\'{e}trique entra\^{\i}ne la premi\`{e}re
 \'{e}quation de l'\'{e}nonc\'{e} (\ref{EqLOb2}). D'autre part, si l'on met
 $f_1=f$, $f_2=g_1\in\mathcal{I}$ et $g=g_2$ dans le premier \'{e}nonc\'{e}, on
 obtient
 en antisym\'{e}trisant en $g_1$ et $g_2$
 \beas
   0 & = & (i^*f)B_{r+1}(g_1,g_2)-(i^*f)B_{r+1}(g_2,g_1)
      -B_{r+1}(fg_1,g_2) + B_{r+1}(fg_2,g_1) \\
     &   &
      + B_{r+1}(f,g_1g_2) - B_{r+1}(f,g_2g_1) \\
     & \stackrel{(\ref{EqBTens})}{=} & (i^*f)B^-_{r+1}(g_1,g_2)
                             -B_{r+1}(fg_1,g_2) + B_{r+1}(g_2,fg_1)
 \eeas
 ce qui donne la deuxi\`{e}me \'{e}quation de (\ref{EqLOb2}).\\
 Pour obtenir l'\'{e}nonc\'{e} (\ref{EqLOb3}) on utilise alternativement
 l'antisym\'{e}trie de $B^-_{r+1}$ et la premi\`{e}re \'{e}quation de l'\'{e}nonc\'{e}
 (\ref{EqLOb2}):
 \beas
  B^-_{r+1}(g_1g_2,g_3) & = & B^-_{r+1}(g_1,g_2g_3)
                              = -B^-_{r+1}(g_2g_3,g_1)
                              = -B^-_{r+1}(g_2,g_3g_1) \\
                        & =  & B^-_{r+1}(g_3g_1,g_2)
                              = B^-_{r+1}(g_3,g_1g_2)
                              = -B^-_{r+1}(g_1g_2,g_3)
 \eeas
 d'o\`{u} $B^-_{r+1}(g_1g_2,g_3)=0$.\\
 Si l'on regarde l'\'{e}quation de l'associativit\'{e} (\ref{EqAssOrdr}) \`{a} l'ordre
 $r+2$
 pour $f_1=g_1,f_2=g_2,g=g_3$ o\`{u} $g_1,g_2,g_3\in\mathcal{I}$ et si l'on
 antisym\'{e}trise en $g_1,g_2,g_3$ on obtient --gr\^{a}ce \`{a} la sym\'{e}trie de
 $\mathsf{C}_0$-- l'\'{e}quation suivante o\`{u}
 $C^-_s(f_1,f_2):=\mathsf{C}_s(f_1,f_2)-\mathsf{C}_s(f_2,f_1)$ quel que
 soit l'entier positif $s$:
 \beas
   i^*C^-_1\big(g_1,C^-_{r+1}(g_2,g_3)\big)
   +i^*C^-_1\big(g_2,C^-_{r+1}(g_3,g_1)\big)
          +i^*C^-_1\big(g_3,C^-_{r+1}(g_1,g_2)\big)
          &   &   \\
                -B^-_{r+1}\big(C^-_1(g_1,g_2),g_3\big)
                -B^-_{r+1}\big(C^-_1(g_2,g_3),g_1\big)
                 -B^-_{r+1}\big(C^-_1(g_3,g_1),g_2\big)
        &   =   & 0,
 \eeas
 d'o\`{u} l'\'{e}nonc\'{e} (\ref{EqLOb4}) car $i^*C^-_1(g_1,f)$ est \'{e}gal au crochet de
 Poisson
 $2i^*\{g_1,f\}=-2X_{g_1}(i^*f)$ par d\'{e}finition.
\ebew

\bcor
 \ben
  \item
  Il existe un unique champ de tenseurs $\tilde{B}_{r+1}$ appartenant
 \`{a} $\Ginf\big(C,\Lambda^2(TM|_C/TC)\big)$ tel que
 \beq\label{EqDefTildeBR+1}
    B^-_{r+1}(g_1,g_2)=\tilde{B}_{r+1}(dg_1|_C,dg_2|_C)~~\forall
    g_1,g_2\in\mathcal{I}
 \end{equation}
 qui soit un $2$-cocycle par rapport \`{a} la diff\'{e}rentielle $d_P$ (voir
 le lemme \ref{LCohdP}), c.-\`{a}-d. $d_P\tilde{B}_{r+1}=0$.
 \item Si $P$ est une structure de Poisson venant d'une forme symplectique
 $\omega$ sur $M$: il existe une unique $2$-forme verticale
 $\beta_{r+1}$ appartenant \`{a} $\Ginf(C,\Lambda^2 TC^{\omega*})$ qui est
 ferm\'{e}e et telle que
 \beq\label{EqDefBetaR+1}
    B^-_{r+1}(g_1,g_2)=\beta_{r+1}(X_{g_1}|_C,X_{g_2}|_C)~~\forall
    g_1,g_2\in\mathcal{I}
 \end{equation}
 o\`{u} $X_{g_1}$ etc. d\'{e}signe le champ hamiltonien associ\'{e} \`{a} $g_1$.
 \een
\ecor
\bbew
 Les \'{e}nonc\'{e}s (\ref{EqLOb2}) et (\ref{EqLOb3}) du lemme pr\'{e}c\'{e}dent montrent
 que $B^-_{r+1}$ induit une unique application $\CinfK{C}$-lin\'{e}aire
 \[
     \Lambda^2_{\CinfK{C}}(\mathcal{I}/\mathcal{I}^2)\ra \CinfK{C}.
 \]
 Mais d'apr\`{e}s le lemme \ref{LemIIcarre}, \'{e}nonc\'{e} (3), le $\CinfK{C}$-module
 $\mathcal{I}/\mathcal{I}^2$ est isomorphe \`{a} $\Ginf(C,TC^{ann})$. Puisque
 le fibr\'{e} dual de $TC^{ann}$ est donn\'{e} par le fibr\'{e} quotient $TM|_C/TC$
 l'existence et l'unicit\'{e} du champ de tenseurs $\tilde{B}_{r+1}$ sont
 claires. Ceci montre le premier \'{e}nonc\'{e} de ce corollaire.\\
 Pour le cas symplectique on utilise la premi\`{e}re partie et la proposition
 \ref{PChampsVertIdeal} pour l'existence et l'unicit\'{e} de la $2$-forme
 verticale $\beta_{r+1}$. D'apr\`{e}s la proposition \ref{PChampsVertIdeal} on
 peut repr\'{e}senter trois champs de vecteurs verticaux $V_1,V_2,V_3$ par les
 restrictions \`{a} $C$ des champs hamiltoniens $X_{g_1},X_{g_2}$ et $X_{g_3}$
 de $g_1,g_2,g_3\in\mathcal{I}$. On a (tout en \'{e}crivant $X_{g_1}$ pour
 $X_{g_1}|_C$):
 \beas
   (d_v\beta_{r+1})\big(X_{g_1},X_{g_2},X_{g_3}\big) & = &
       X_{g_1}\left(\beta_{r+1}(X_{g_2},X_{g_3})\right)
       + X_{g_2}\left(\beta_{r+1}(X_{g_3},X_{g_1})\right) \\
    & &       ~~~~~~~~~+ X_{g_3}\left(\beta_{r+1}(X_{g_1},X_{g_2})\right) \\
    & &   -\beta_{r+1}\left([X_{g_1},X_{g_2}],X_{g_3}\right)
       -\beta_{r+1}\left([X_{g_2},X_{g_3}],X_{g_1}\right) \\
    & &  ~~~~~~~~~-\beta_{r+1}\left([X_{g_3},X_{g_1}],X_{g_2}\right).
 \eeas
 Puisque $[X_{g_1},X_{g_2}]=-X_{\{g_1,g_2\}}$ on voit que le membre droit de
 cette
 \'{e}quation co\"{\i}ncide avec $-1$ fois le membre gauche de (\ref{EqLOb4}).
\ebew

\noindent Supposons qu'on fasse maintenant la transformation d'\'{e}quivalence
$S^{(r)}:=\mathrm{id}+\nu^rT_r$ pour le star-produit $*$ de fa\c{c}on que le
star-produit transform\'{e} $*':=S^{(r)}(*)$ soit toujours adapt\'{e} \`{a} $C$ jusqu'\`{a}
l'ordre $r$.
\bprop
Avec les hypoth\`{e}ses faites ci-dessus,
 il existe alors une unique section $\xi_r\in\Ginf(C,TM|_C/TC)$
 telle que $\tilde{B}'_{r+1}=\tilde{B}_{r+1}-d_P\xi_r$. Dans le cas
 symplectique il existe une unique $1$-forme verticale $\gamma_r$ telle
 que la $2$-forme verticale ferm\'{e}e $\beta_{r+1}$ est modifi\'{e}e par
 $d_v\gamma_{r}$
\eprop
\bbew
 La transformation d'\'{e}quivalence ne change pas les op\'{e}rateurs
 bidiff\'{e}rentiels $\mathsf{C}_0,\ldots,\mathsf{C}_{r-1}$. Le terme
 $\mathsf{C}_r$ est modifi\'{e}
 par $\mathsf{C}'_{r}=\mathsf{C}_r-\mathsf{b} T_r$ avec
 $\mathsf{b} T_r(f_1,f_2):=f_1T_r(f_2)-T_r(f_1f_2)+T_r(f_1)f_2$ et $T_r$
 un op\'{e}rateur diff\'{e}rentiel qui s'annule sur $1$.
 Soient $f\in\CinfK{M}$ et $g\in\mathcal{I}$. Il vient
 $0=i^*(\mathsf{b} T_r)(f,g)=i^*fi^*T_r(g)-i^*T_r(fg)$. Le raisonnement
 qui suit l'\'{e}quation (\ref{EqistarBfgg}) montre l'existence d'une
 section $\xi_r$ telle que $\langle \xi, dg|_C\rangle:= i^*T_r(g)$.
 Le terme $\mathsf{C}_{r+1}$ est remplac\'{e} par
 \[
   \mathsf{C}'_{r+1}(f_1,f_2)=\mathsf{C}_{r+1}(f_1,f_2)
        -\big(\mathsf{C}_1(T_rf_1,f_2)
        -T_r\big(\mathsf{C}_1(f_1,f_2)\big)+\mathsf{C}_1(f_1,T_rf_2\big).
 \]
 Si on prend la partie antisym\'{e}trique de cette \'{e}quation pour
 $f_1=g_1,f_2=g_2\in\mathcal{I}$ on arrive \`{a} l'\'{e}quation pour
 $\tilde{B}'_{r+1}$. Le cas symplectique est \'{e}vident.
\ebew

\bprop
 Soit $*$ un star-produit adapt\'{e} $C$ jusqu'\`{a} l'ordre $r$ et supposons
 que
 \[
   i^*\big(\mathsf{C}_{r+1}(g_1,g_2)
    -\mathsf{C}_{r+1}(g_2,g_1)\big)=0\mathrm{~quels~que~
            soient~~}g_1,g_2\in\mathcal{I}.
 \]
 Alors il existe un op\'{e}rateur diff\'{e}rentiel $T_{r+1}$ dans
 $\mathbf{D}^1\big(\CinfK{M},\CinfK{M}\big)$ qui s'annule sur les
 constantes tel que le star-produit $(\mathrm{id}+\nu^{r+1}T_{r+1})(*)$
 soit adapt\'{e} \`{a} $C$ jusqu'\`{a} l'ordre $r+1$.
\eprop
\bbew
 Pour un op\'{e}rateur diff\'{e}rentiel $T_{r+1}$ arbitrairement choisi
 l'op\'{e}\-ra\-teur bidiff\'{e}rentiel $\mathsf{C}_{r+1}$ changera en
 \[
     \mathsf{C}'_{r+1}(f_1,f_2)=\mathsf{C}_{r+1}(f_1,f_2)
         - \big(f_1T_{r+1}(f_2)-T_{r+1}(f_1,f_2)+T_{r+1}(f_1)f_2\big),
 \]
 donc, si on veut que $0\stackrel{!}{=}i^*\mathsf{C}'_{r+1}(f,g)$ quels que
 soient
 $f\in\CinfK{M},g\in\mathcal{I}$ il faut demander que
 \beq\label{EqCeQuonVeutDeTRplusun}
  i^*\mathsf{C}_{r+1}(f,g)\stackrel{!}{=} i^*T_{r+1}(fg)-i^*f i^*T_{r+1}(g)
 \end{equation}
 On fait d'abord une \'{e}tude locale: soit $\big(U,(\eta,y)\big)$ une carte
 de sous-vari\'{e}t\'{e} o\`{u} les coordonn\'{e}s $\eta=(\eta^1,\ldots,\eta^l)$
 param\'{e}trisent l'intersection $U\cap C\neq\emptyset$ et les coordonn\'{e}es
 $y=(y_1,\ldots,y_k)$ d\'{e}finissent $U\cap C$ comme $\{p\in U~|~y=0\}$.
 Soit $g\in\mathcal{I}$. Puisque dans $U$ il vient $g(\eta,0)=0$, on
 a
 \beq\label{EqHomotopieH0locale}
    g(\eta,y)=g(\eta,y)-g(\eta,0)
          =\sum_{i=1}^k\int_0^1\frac{\partial g}{\partial
          y_i}(\eta,ty)y_i~dt =:\sum_{i=1}^kg^i(\eta,y)y_i,
 \end{equation}
 ce qui est la version locale de l'homotopie $h_0$ (\ref{EqHomtopKoszul})
 du complexe de Koszul (\ref{EqHomtopKoszul}). Suppsons qu'il y ait une
 autre d\'{e}composition $g(\eta,y)=\sum_{i=1}^k\hat{g}^i(\eta,y)y_i$, alors
 on a $0=\sum_{i=1}^k\big(g^i(\eta,y)-\hat{g}^i(\eta,y)\big)y_i$, et gr\^{a}ce
 \`{a} l'acyclicit\'{e} du complexe de Koszul il existe des fonctions
 $h^{ij}\in\CinfK{U}$, $1\leq i,j\leq k$ avec $h^{ij}=-h^{ji}$ telles que
 \beq\label{EqDeuxRepreKoszulH}
    g^i(\eta,y)-\hat{g}^i(\eta,y)=\sum_{j=1}^kh^{ij}(\eta,y)y_j.
 \end{equation}
 On d\'{e}finit l'op\'{e}rateur $E_U:\mathcal{I}|_U\ra \CinfK{C\cap U}$ par
 \beq\label{EqDefEU}
    E_U(g):=\sum_{i=1}^ki^*\mathsf{C}_{r+1}(g^i,y_i).
 \end{equation}
 Cette d\'{e}finition ne d\'{e}pend pas des fonctions $g^1,\ldots,g^k$ choisies:
 en fait, en utilisant (\ref{EqDeuxRepreKoszulH}) et le fait que
  $y_i,y_j\in\mathcal{I}_U$ on peut appliquer la symm\'{e}trie
 de la restriction de $i^*\mathsf{C}_{r+1}$ \`{a} $\mathcal{I}\times\mathcal{I}$
 pour conclure
 \beas
  \sum_{i,j=1}^ki^*\mathsf{C}_{r+1}(h^{ij}y_j,y_i) & = &
       \sum_{i,j=1}^ki^*\mathsf{C}_{r+1}(y_i,h^{ij}y_j)
       \stackrel{(\ref{EqBTens})}{=}
           \sum_{i,j=1}^ki^*\mathsf{C}_{r+1}(h^{ij}y_i,y_j) \\
           & = & \sum_{i,j=1}^ki^*\mathsf{C}_{r+1}(h^{ji}y_j,y_i)
              = -\sum_{i,j=1}^ki^*\mathsf{C}_{r+1}(h^{ij}y_j,y_i),
 \eeas
 donc $\sum_{i,j=1}^ki^*\mathsf{C}_{r+1}(h^{ij}y_j,y_i)=0$ et il vient
 \beq \label{EqIndepDecIdlocale}
  \sum_{i=1}^ki^*\mathsf{C}_{r+1}(g^i,y_i)  =
    \sum_{i=1}^ki^*\mathsf{C}_{r+1}(\hat{g}^i,y_i)
 \end{equation}
 On a pour tout $f\in\CinfK{M}$ et $g\in\mathcal{I}$ que
 \beq\label{Eqfgigleichfig}
     \sum_{i=1}^k(fg)^i(\eta,y)y_i=(fg)(\eta,y)
     =f(\eta,y)\sum_{i=1}^kg^i(\eta,y)y_i
     =\sum_{i=1}^k(fg^i)(\eta,y)y_i,
 \end{equation}
 donc
 \bea
   i^*f~E_U(g)-E_U(fg) & = & \sum_{i=1}^k\big(i^*f~i^*\mathsf{C}_{r+1}(g^i,y_i)
                            -i^*\mathsf{C}_{r+1}\big((fg)^i,y_i\big)\big)
                            \nonumber\\
         & \stackrel{(\ref{Eqfgigleichfig}),(\ref{EqIndepDecIdlocale})}{=} &
           \sum_{i=1}^k\big(i^*f~i^*\mathsf{C}_{r+1}(g^i,y_i)
                            -i^*\mathsf{C}_{r+1}\big(f~ g^i,y_i\big)\big)
                            \nonumber\\
         & \stackrel{(\ref{EqLOb1})}{=} &
              -\sum_{i=1}^ki^*\mathsf{C}_{r+1}(f,g^iy_i)=
              -i^*\mathsf{C}_{r+1}(f,g),
 \eea
 donc, en vue de (\ref{EqCeQuonVeutDeTRplusun}) $E_U$ est un bon candidat
 pour la restriction de $i^*T_{r+1}$ \`{a} $U$ et \`{a} $\mathcal{I}$. Puisque
 pour tous multi-indices $I$,$J$ la d\'{e}finition (\ref{EqHomotopieH0locale})
 implique
 \[
    i^*\frac{\partial^{|I|+|J|}g^i}{\partial \eta^I\partial y_J}
      = i^*\frac{1}{|J|+1}
        \frac{\partial^{|I|+|J|+1}g}{\partial \eta^I\partial y_J\partial y_i}
 \]
 on voit --en utilisant la forme locale de $\mathsf{C}_{r+1}$ comme op\'{e}rateur
 bidiff\'{e}rentiel-- qu'il existe un op\'{e}rateur differentiel $T^U_{r+1}$ dans
 $\mathbf{D}^1\big(\CinfK{U},\CinfK{U}\big)$ tel que
 \[
     i^*T^U_{r+1}(g)=E_U(g)~~~\forall~g\in\mathcal{I}|_U.
 \]
 Finalement, on recolle les op\'{e}rateurs $T^U_{r+1}$ \`{a} l'aide d'une
 partition de l'unit\'{e} d'une mani\`{e}re analogue \`{a} celle utilis\'{e}e dans la fin
 de la d\'{e}monstration de la proposition \ref{PKerRho1congId}, tout en
 observant que
 la multiplication avec des fonctions de la partition de l'unit\'{e}
 n'influence pas l'\'{e}quation (\ref{EqCeQuonVeutDeTRplusun}).
\ebew

\noindent On est arriv\'{e} \`{a} la caract\'{e}risation compl\`{e}te des obstructions
r\'{e}currentes pour la construction d'un star-produit adapt\'{e} \`{a} $C$:

\bsat \label{TObstructionsRecurrentes}
 Soit $*=\sum_{r=0}^\infty \nu^r \mathsf{C}_r$ un star-produit sur une
 vari\'{e}t\'{e} de Poisson $(M,P)$ et soit
 $i:C\ra M$ une sous-vari\'{e}t\'{e} co\"{\i}sotrope ferm\'{e}e.\\
 Alors pour tout entier $r\geq 1$ on a le suivant: soit $*$ adapt\'{e} \`{a}
 $C$ jusqu'\`{a} l'ordre $r$. Alors on peut modifier $*$ par une
 transformation d'\'{e}quivalence de la forme
 $(\mathrm{id}+\nu^{r+1}T_{r+1})(\mathrm{id}+\nu^{r}T_{r})$ de telle
 sorte que le star-produit transform\'{e} soit adapt\'{e} jusqu'\`{a} l'ordre $r+1$
 si et seulement si le $2$-cocycle $\tilde{B}_{r+1}$
 (voir (\ref{EqDefTildeBR+1})) est un cobord
 dans la cohomologie d\'{e}finie par $d_P$. Dans le cas d'une vari\'{e}t\'{e}
 symplectique $(M,\omega)$ ceci est le cas si et seulement si
 la $2$-forme verticale ferm\'{e}e $\beta_{r+1}$ (voir (\ref{EqDefBetaR+1}))
 est exacte.
\esat
\bcor \label{CHDeuxNulleRepresentable}
 Soit $*=\sum_{r=0}^\infty \nu^r \mathsf{C}_r$ un star-produit sur une
 vari\'{e}t\'{e} de Poisson $(M,P)$ et soit
 $i:C\ra M$ une sous-vari\'{e}t\'{e} co\"{\i}sotrope ferm\'{e}e.\\
 Si le deuxi\`{e}me groupe de cohomologie BRST de $C$,
 $H_P^2(C,\korps)$, s'annule, alors le star-produit $*$ est repr\'{e}sentable.
 Dans le cas particulier d'une vari\'{e}t\'{e} symplectique, si le deuxi\`{e}me groupe
 de cohomologie de de Rham verticale, $H^2_v(C,\korps)$, s'annule,
 alors $*$ est repr\'{e}sentable.
\ecor
\bbew
 Puisque les obstructions r\'{e}currentes pour rendre l'id\'{e}al annulateur un id\'{e}al
 \`{a} gauche appartiennent \`{a} $H_P^2(C,\korps)$, ce corollaire est \'{e}vident.
\ebew

\brem
Au moins dans le cas symplectique on peut ajouter une
modification de $*$
\`{a} l'ordre $r+1$, \`{a} savoir une $2$-forme ferm\'{e}e $\nu^{r+1}\alpha_{r+1}$
(qui modifie la classe de Deligne de $*$),
voir la construction de DeWilde par exemple dans \cite{GR99}, Theorem 7.1.
Par cons\'{e}quent,
l'espace de toutes les obstructions r\'{e}currentes est donn\'{e} par le quotient
\[
        H^2_v(C,\korps)/p_v\big(i^*H^2_{dR}(M,\korps)\big).
\]
\erem

\subsection{Le cas de codimension $1$}
   \label{SubSecCodimUn}

Les m\'{e}thodes de la section pr\'{e}c\'{e}dente permettent de red\'{e}montrer un cas
particulier qui a \'{e}t\'{e} fait en 1998 par Peter Gl\"{o}{\ss}ner dans sa th\`{e}se
\cite{Glo98}:

\bsat[P.Gl\"{o}{\ss}ner 1998]\label{TGloessnerCodimUn} Soit $*$ un star-produit sur une
 vari\'{e}t\'{e} de Poisson
  $(M,P)$ et $C$ une sous-vari\'{e}t\'{e} co\"{\i}sotrope ferm\'{e}e de codimension $1$
  de $M$.\\
  Alors $*$ est toujours repr\'{e}sentable sur $C$.
\esat
\bbew
 Puisque le fibr\'{e} $TM|_C/TC$ est de dimension $1$ dans ce cas, toutes
 les obstructions r\'{e}currentes --qui se d\'{e}duisent du fibr\'{e}
 $\Lambda^2(TM|_C/TC)$-- sont automatiquement r\'{e}duites \`{a} z\'{e}ro,
 donc on peut construire par r\'{e}currence un star-produit \'{e}quivalent \`{a} $*$
 qui soit adapt\'{e} \`{a} $C$.
\ebew

\noindent Gl\"{o}{\ss}ner a montr\'{e} ce r\'{e}sultat pour une vari\'{e}t\'{e}
symplectique dont la
sous-vari\'{e}t\'{e} co\"{\i}sotrope est donn\'{e} par l'ensemble des z\'{e}ros
d'une fonction $J\in\Cinf(M,\real)$. La d\'{e}monstration de son Lemma 1, par.3,
dans lequel il construit par r\'{e}currence une transformation d'\'{e}quivalence,
reste valable texto dans le cas d'une vari\'{e}t\'{e} de Poisson. Le cas
de $\real^{n-1}\subset \real^n$ a \'{e}galement \'{e}t\'{e} trait\'{e} dans
\cite{CF03p} par des m\'{e}thodes graphiques.

\subsection{Connexions symplectiques adapt\'{e}es \`{a} une sous-vari\'{e}t\'{e}
    co\"{\i}sotrope} \label{SubSecConnSympAdaptCois}

 Soit $C$ une vari\'{e}t\'{e} diff\'{e}rentiable et $E\subset TC$ un sous-fibr\'{e}
 int\'{e}grable. Soit $F\subset TC$ un sous-fibr\'{e} compl\'{e}mentaire \`{a} $E$,
 c.-\`{a}-d. $TC=E\oplus F$. On d\'{e}signe par
 $\mathsf{P}_E\in \Ginf\big(C,\mathrm{Hom}(TC,TC)\big)$
 (resp. $\mathsf{P}_F\in\Ginf\big(C,\mathrm{Hom}(TC,TC)\big)$)
 le champs d'endomorphismes
 du fibr\'{e} tangent qui projette $TC$ sur $E$ (resp. sur $F$) le long de
 $F$ (resp. $E$). La {\em courbure de $F$},
 $\hat{R}^F\in\Ginf(C,\Lambda^2T^*C\otimes TC)$ est d\'{e}finie par
 \beq
      \hat{R}^F(X,Y):=\mathsf{P}_E[\mathsf{P}_FX,\mathsf{P}_FY]
 \end{equation}
 quels que soient les champs de vecteurs $X,Y\in\Ginf(C,TC)$.

\bprop\label{PExistConnAdaptees}
 Avec les notations introduites ci-dessus:
 \ben
  \item Il existe une connexion sans
 torsion $\nabla$ dans le fibr\'{e} tangent $TC$ avec les propri\'{e}t\'{e}s
 suivantes pour tous $X\in\Ginf(C,TC)$, $V,W\in\Ginf(C,E)$ et
 $H,H_1,H_2,H_3\in\Ginf(C,F)$:
 \begin{itemize}
  \item[i)] $\nabla_XV\in\Ginf(C,E)$,
  \item[ii)] $\nabla_VH=\mathsf{P}_F[V,H]\in\Ginf(C,F)$,
  \item[iii)] $\nabla_HV=\mathsf{P}_E[H,V]\in\Ginf(C,E)$,
  \item[iv)] $\mathsf{P}_E\nabla_{H_1}H_2=\frac{1}{2}\hat{R}^F(H_1,H_2)
                \in\Ginf(C,E)$.
 \end{itemize}
 Pour toute connexion sans torsion $\nabla$ satisfaisant les conditions
 ci-dessus son
 tenseur de courbure $R$ a la {\em propri\'{e}t\'{e} de Bott} suivante:
 \beq\label{EqConnExtCourbureNulle}
     R(V,W)H=0~~~~\forall~V,W\in\Ginf(C,E)~~\forall~H\in\Ginf(C,F).
 \end{equation}
 \item
 Dans le cas particulier o\`{u} $E$ est le sous-fibr\'{e} caract\'{e}ristique d'une
 $2$-forme
 pr\'{e}symplectique $\varpi$ sur $C$ la connexion $\nabla$ peut \^{e}tre choisie
 d'une telle mani\`{e}re qu'elle pr\'{e}serve \'{e}galement $\varpi$, c.-\`{a}-d.
 \begin{itemize}
  \item[v)] $\nabla\varpi = 0$.
 \end{itemize}
 \een
\eprop
\bbew
 1. Soit $\nabla'$ une connexion sans torsion arbitraire sur $C$
 (par exemple la connexion Levi-Civita d'une m\'{e}trique riemannienne).
 Pour tous les champs de vecteurs $X,Y\in\Ginf(C,TC)$
 on d\'{e}finit

 \bea
  \nabla_X Y & := & \nabla'_XY
                      -\mathsf{P}_F\big(
                      \nabla'_{\mathsf{P}_EX}(\mathsf{P}_EY)\big)
                      -\mathsf{P}_F\big(
                      \nabla'_{\mathsf{P}_FX}(\mathsf{P}_EY)+
      \nabla'_{\mathsf{P}_FY}(\mathsf{P}_EX)\big) \nonumber \\
              &    &  ~~-\mathsf{P}_E\big(
                \nabla'_{\mathsf{P}_EX}(\mathsf{P}_FY)
                +\nabla'_{\mathsf{P}_EY}(\mathsf{P}_FX)\big)\nonumber
                               \\
              &    &   ~~~~~ -\frac{1}{2}
                         \mathsf{P}_E\big(
                         \nabla'_{\mathsf{P}_FX}(\mathsf{P}_FY)
                         +\nabla'_{\mathsf{P}_FY}(\mathsf{P}_FX)\big).
 \eea
 Il est \'{e}vident que $\nabla_XY-\nabla'_XY$ sont des champs de tenseurs
 sym\'{e}triques par rapport \`{a} l'\'{e}change de $X$ et de $Y$,
 et --gr\^{a}ce au fait que la torsion de $\nabla'$ s'annule et
 \`{a} l'identit\'{e} $\mathsf{P}_E+\mathsf{P}_F=\mathrm{id}$-- on v\'{e}rifie rapidement
 les quatre propri\'{e}t\'{e}s $i)$ - $iv)$
 et la condition (\ref{EqConnExtCourbureNulle}) pour le tenseur de
 courbure.\\
 2. Pour le cas pr\'{e}symplectique, on calcule sans peine ($d\varpi=0$!)
 que pour la connexion
 $\nabla$ construite ci-dessus la d\'{e}riv\'{e}e covariante
 $(\nabla_X\varpi)(Y,Z)$ s'annule automatiquement si un des champs de
 vecteurs $X,Y,Z$
 est vertical. Puisque la restriction de $\varpi$ \`{a} $F\times F$
 est nond\'{e}g\'{e}n\'{e}r\'{e}e, la section $S\in\Ginf(C,F^*\otimes F^*\otimes F)$
 suivante est bien d\'{e}finie (l'astuce de Tondeur, Lichnerowicz, He{\ss}, voir
 par exemple \cite{Hes81}):
 \beq \label{EqAstuceTondeur}
    \varpi\big(S(H_1,H_2),H_3\big):=
    \frac{1}{3}\big(\nabla_{H_1}\varpi\big)(H_2,H_3)
    +\frac{1}{3}\big(\nabla_{H_2}\varpi\big)(H_1,H_3),
 \end{equation}
 et \`{a} l'aide de $d\varpi=0$ il s'ensuit que la connexion
 \[
     \nabla_XY+S(\mathsf{P}_FX,\mathsf{P}_FY)
 \]
 pr\'{e}serve $\varpi$ et satisfait toujours $i)$ - $iv)$.
\ebew

\noindent Pour toute vari\'{e}t\'{e} munie d'un feuilletage r\'{e}gulier (resp.
pr\'{e}symplectique) on va appeler le couple $(F,\nabla)$ qui remplit les
condions $i)$ - $iv)$ (resp. $i)$ - $v)$) une {\em connexion adapt\'{e}e}.
On voit \'{e}galement que toute connexion adapt\'{e}e induit une connexion
$\bar{\nabla}$ dans le fibr\'{e} quotient $Q=TC/E$ par
\beq \label{EqConnAdapteeInduite}
  \overline{\nabla}_X\overline{Y}:=\overline{\nabla_XY}
\end{equation}
qui est une prolongation de la connexion de Bott en vue de l'\'{e}quation
$ii)$ de la proposition \ref{PExistConnAdaptees} et
de l'\'{e}quation (\ref{EqConnexionBott}). L'action de l'alg\`{e}bre de Lie
$\Ginf(C,E)$ sur $\Ginf(C,F)$ par $\nabla_V$ est visiblement isomorphe
\`{a} celle de Bott sur $\Ginf(C,Q)$, et on va en faire usage beaucoup
de fois.

Pour une vari\'{e}t\'{e} pr\'{e}symplectique $(C,\varpi)$ on
consid\`{e}re le fibr\'{e} dual $\tau:E^*\ra C$ au fibr\'{e} caract\'{e}ristique $E$ de
$\omega$.
Soient $c\in C$ et $\varphi,\psi\in E^*$ avec $\tau(\varphi)=\tau(\psi)=c$.
On rappelle le {\em rel\`{e}vement vertical de
$\psi$} comme vecteur tangent en $\varphi$: $\psi^v_\varphi:=
\frac{d}{dt}(\varphi+t\psi)|_{t=0}$. On d\'{e}finit le rel\`{e}vement vertical
des sections dans $\Ginf(C,E^*)$ par celui de leurs valeurs.
Soit $v\in T_cC$ et $\nabla$ une connexion dans le fibr\'{e} $E^*$. On
rappelle la d\'{e}finition
du {\em rel\`{e}vement horizontal de $v$} comme vecteur tangent en $\varphi$:
$v^h_\varphi:=\frac{d}{dt}(\mathsf{Tp}_\gamma(t)\varphi)|_{t=0}$ o\`{u} $\gamma$
est une
courbe lisse dans $C$ telle que $\gamma(0)=c$ et $\frac{d\gamma}{dt}(0)=v$,
et $\mathsf{Tp}_\gamma(t)\varphi$ d\'{e}signe le transport parall\`{e}le de $\varphi$
le long de la courbe $\gamma$ (par rapport \`{a} la connexion $\nabla$).
On d\'{e}finit le rel\`{e}vement horizontal des champs de vecteurs sur $C$
par celui de leurs valeurs. On note les formules
\[
 \begin{array}{ccl}
  ~[X^h,Y^h]_\varphi & = &
    {[X,Y]^h}_\varphi +\langle\varphi,R(X,Y)\mathsf{P}_E~\underline{~~}~
               \rangle^v_\varphi \\
  ~[X^h,\chi^v] & = & (\nabla_X\chi)^v \\
  ~[\psi^v,\chi^v] & = & 0
 \end{array}
\]

\bprop
 Soit $(C,\varpi)$ une vari\'{e}t\'{e} pr\'{e}symplectique et $\tau:E^*\ra C$ le fibr\'{e}
 dual du fibr\'{e} caract\'{e}ristique de $(C,\varpi)$. Soit $\nabla$ une
 connexion adapt\'{e}e dans le fibr\'{e} tangent de $C$. Quels que soient les
 champs de vecteurs $X,Y$ sur $C$ et les sections $\psi,\chi\in\Ginf(C,E^*)$
 on a la formule suivante pour la $2$-forme symplectique de
 Weinstein-Gotay (voir le th\'{e}or\`{e}me \ref{TWeinGotay}) au point $\varphi\in
 E^*_c$
 \bea
     \omega_{WG}(X^h,Y^h)_\varphi & = & \varpi(X,Y)_c
                           +\langle \varphi, \hat{R}^F(X,Y)_c \rangle
                           =:\hat{\varpi}_\varphi(X,Y), \label{EqHatOmegaPhi}
                           \\
     \omega_{WG}(X^h,\chi^v)_\varphi & = &
       \langle \chi, \mathsf{P}_EX\rangle_c, \\
     \omega_{WG}(\psi^v,\chi^v)_\varphi & = & 0.
 \eea
\eprop
\bbew
 Calcul direct \`{a} partir de la formule (\ref{EqFormeWeinsteinGotay}).
\ebew

Soit maintenant $\check{M}$ un ouvert de $E^*$ qui contient la z\'{e}ro-section
$C$
et qui est telle que la forme $\omega_{WG}$ soit nond\'{e}g\'{e}n\'{e}r\'{e}e. On peut et
supposer que $\check{M}$ soit convexe fibre-par-fibre.
%(i.e. pour
%tout
%$\varphi\in \check{M}$ le segment
% $\{t\varphi~|~t\in[0,1]\}\subset \check{M}$).
La non-d\'{e}g\'{e}n\'{e}rescence de $\omega_{WG}$ en $\varphi\in \check{M}$ est
\'{e}quivalente \`{a} celle de
la restriction de $\hat{\varpi}_\varphi$ \`{a} $F\times F$
(voir (\ref{EqHatOmegaPhi})): puisque la restriction de $\varpi$ \`{a}
$F\wedge F$ est non d\'{e}g\'{e}n\'{e}r\'{e}e et la diff\'{e}rence $\hat{\varpi}_\varphi-\varpi$
est lin\'{e}aire en $\varphi$ le voisinage $\check{M}$ de $C$ existe.
Il y a donc un unique champ de bivecteurs
$\varphi\mapsto
\hat{P}_\varphi\in F_{\tau(\varphi)}\wedge F_{\tau(\varphi)}\subset
T_{\tau(\varphi)}C\wedge T_{\tau(\varphi)}C$ tel que
\beq
    \hat{P}^\sharp_\varphi\hat{\varpi}_\varphi=\mathsf{P}_F
        \mathrm{~~~et~~~}
    \hat{\varpi}_\varphi\hat{P}^\sharp_\varphi={\mathsf{P}_F}^*.
\end{equation}
o\`{u} $\hat{P}^\sharp_\varphi$ est consid\'{e}r\'{e}e comme application
$T_{\tau(\varphi)}C^*\ra T_{\tau(\varphi)}C$
dont le noyau est \'{e}gal \`{a} $E_{\tau(\varphi)}$.
On peut en d\'{e}duire le suivant

\bsat \label{TConnSympAutourCoisotrope}
 Soit $(C,\varpi)$ une vari\'{e}t\'{e} pr\'{e}symplectique et $(F,\nabla)$ une con\-nexion
 adapt\'{e}e.
 Alors la con\-nex\-ion $\tilde{\nabla}$ suivante
 dans l'ouvert $\check{M}\subset E^*$ est sans torsion et symplectique dans
 la vari\'{e}t\'{e} symplectique $(\check{M},\omega_{WG})$: soient $X,Y$ deux champs de
 vecteurs sur $C$, $\psi,\chi\in\Ginf(C,E^*)$ et $\varphi\in \check{M}$ avec
 $c:=\tau(\varphi)$, alors
 \bea
    \big(\tilde{\nabla}_{X^h}Y^h\big)_\varphi & = &
                \big(\nabla_XY)\big)^h_\varphi +
               \frac{1}{2}\langle \varphi,
               R(X,Y)\mathsf{P}_E~\underline{~~}~)_c\rangle^v_\varphi
                                            \nonumber \\
          & &      + \langle \varphi,A(X,Y,~\underline{~~}~)_c\rangle^v_\varphi
                +\big(\hat{P}_\varphi^\sharp\langle \varphi,
                   B(X,Y,~\underline{~~}~)_c\rangle\big)^h_\varphi,
                                      \label{EqConnSympHorHor} \\
    \big(\tilde{\nabla}_{X^h}\chi^v\big)_\varphi & = &
                \big(\nabla_X\chi\big)^v_\varphi +
    \frac{1}{2}\big(\hat{P}_\varphi^\sharp\langle\chi_c,
       \hat{R}^F_c(X,~\underline{~~}~)\rangle\big)^h_\varphi,
                                       \label{EqConnSympHorVer} \\
    \big(\tilde{\nabla}_{\chi^v}X^h\big)_\varphi & = &
    \frac{1}{2}\big(\hat{P}_\varphi^\sharp\langle\chi_c,
       \hat{R}^F_c(X,~\underline{~~}~)\rangle\big)^h_\varphi,
                                       \label{EqConnSympVerHor} \\
    \tilde{\nabla}_{\psi^v}\chi^v & = & 0, \label{EqConnSympVerVer}
 \eea
 avec des champs de tenseurs
 $A\in\Ginf(C,TC^*\otimes TC^*\otimes E^*\otimes E)$ et
 $B\in\Ginf(C,TC^*\otimes TC^*\otimes TC^*\otimes E)$ donn\'{e}s par
 ($V\in\Ginf(C,E),Z\in\Ginf(C,TC)$):
 \bea
  A(X,Y,V) & := & \frac{1}{6}\big(R(\mathsf{P}_FX,\mathsf{P}_EY)V
                                 +R(\mathsf{P}_FY,\mathsf{P}_EX)V\big)
                      \nonumber \\
           &    &          +\frac{1}{6}\big(R(\mathsf{P}_EX,V)\mathsf{P}_EY
                                     +R(\mathsf{P}_EY,V)\mathsf{P}_EX\big),
                                           \\
  B(X,Y,Z) & := &    \frac{1}{3}\big(
      (\nabla_{\mathsf{P}_FX}\hat{R}^F)(\mathsf{P}_FY,\mathsf{P}_FZ)
      +(\nabla_{\mathsf{P}_FY}\hat{R}^F)(\mathsf{P}_FX,\mathsf{P}_FZ)\big)
                                                  \nonumber \\
           &    & -\frac{1}{2}\big(R(\mathsf{P}_FX,Z)\mathsf{P}_EY
                                 +R(\mathsf{P}_FY,Z)\mathsf{P}_EX\big)
                                                   \nonumber \\
           &    &  -\frac{1}{6}\big(R(\mathsf{P}_EX,Z)\mathsf{P}_EY
                                 +R(\mathsf{P}_EY,Z)\mathsf{P}_EX\big).
 \eea
\esat
\bbew
 La pr\'{e}scription
 $\big(\hat{\nabla}_{X^h}Y^h\big)_\varphi:=(\nabla_XY)^h_\varphi
 +\frac{1}{2}\langle \varphi,
 R(X,Y)\mathsf{P}_E~\underline{~~}~)_c\rangle^v_\varphi$,
 $\hat{\nabla}_{X^h}\chi^v:=(\nabla_{X}\chi)^v$,
 $\hat{\nabla}_{\chi^v}X^h:=0$ et $\hat{\nabla}_{\psi^v}\chi^v:=0$
 d\'{e}finit une connexion sans torsion dans le fibr\'{e} tangent de $E^*$, voir par
 exemple \cite[p.410]{KMS93}. Les termes qui restent r\'{e}sultent au bout d'un
 long calcul \`{a} l'aide de l'astuce
 de Tondeur, Lichnerowicz et He{\ss} (une formule analogue \`{a}
 (\ref{EqAstuceTondeur})), appliqu\'{e}e \`{a} $\omega_{WG}$ et $\hat{\nabla}$.
\ebew

On rappelle qu'une application affine entre deux espaces vectoriels r\'{e}els
de dimension finie est d\'{e}finie par la condition que sa deuxi\`{e}me d\'{e}riv\'{e}e
s'annule. On rappelle \'{e}galement la notion d'une {\em application affine}
entre deux vari\'{e}t\'{e}s diff\'{e}rentiables $M$ et $M'$ munies des connexions sans
torsion
$\nabla$ et $\nabla'$ dans leurs fibr\'{e}s tangents, respectivement: soit
$\phi:M\ra M'$ de classe
$\Cinf$. On consid\`{e}re le fibr\'{e} retir\'{e} $\phi^*TM'$. Dans ce fibr\'{e} vectoriel
sur $M$ on a la connexion retir\'{e}e $\phi^*\nabla'$. Les deux connexions
$\nabla$
et $\phi^*\nabla'$ induisent une connexion $\hat{\nabla}$ dans le fibr\'{e}
vectoriel $\mathrm{Hom}(TM,\phi^*TM')$ sur $M$. L'application tangente
$T\phi$ est
une section de $\mathrm{Hom}(TM,\phi^*TM')$. Alors $\phi$ est dite affine
lorsque
sa deuxi\`{e}me d\'{e}riv\'{e}e s'annule dans le sens
\beq
      \hat{\nabla}(T\phi) = 0.
\end{equation}
Si $\big(U,x=(x^1,\ldots,x^m)\big)$ est une carte de $M$ et
$\big(U',y=(y^1,\ldots,y^n)\big)$ une carte de $M'$ telle que
$\phi(U)\subset U'$, on obtient la condition locale (o\`{u} $\phi^b:=y^b\circ
\phi$)
\[
    \frac{\partial^2 \phi^a}{\partial x^i\partial x^j}
       +\sum_{b,c=1}^n\Gamma^{\prime a}_{bc}(\phi)
                           \frac{\partial \phi^b}{\partial x^i}
                           \frac{\partial \phi^c}{\partial x^j}
       -\sum_{k=1}^m\frac{\partial \phi^a}{\partial x^k}\Gamma^k_{ij} =0
\]
avec les symboles de Christoffel usuels $\Gamma^k_{ij}$ de $\nabla$ et
$\Gamma^{\prime a}_{bc}$ de
$\nabla'$. Un autre crit\`{e}re util est le
suivant: soient $X,Y\in\Ginf(M,TM)$ et $X',Y'\in\Ginf(M',TM')$ tels que
$(X,X')$ et $(Y,Y')$ sont deux couples $\phi$-li\'{e}s, c.-\`{a}-d.
$T\phi~X=X'\circ\phi$ et $T\phi~Y=Y'\circ \phi$. On calcule sans peine
que $\nabla'_{X'}Y'\circ\phi-T\phi\nabla_XY=(\hat{\nabla}_XT\phi)Y$.
Par cons\'{e}quent, \`{a} condition que pour tout $p\in M$ et pour tout $v,w\in T_pM$
il existe des couples locaux $(X,X')$ et $(Y,Y')$ $\phi$-li\'{e}s tels que
$X(p)=v$ et
$Y(p)=w$, $\phi$ est affine si et seulement si
$(\nabla_XY,\nabla'_{X'}Y')$ est $\phi$-li\'{e}.

Comme exemple on peut prendre toute g\'{e}od\'{e}sique de $(M',\nabla')$: c'est
une application affine
d'un intervalle r\'{e}el ouvert (muni de la connexion usuelle) dans $M'$. Plus
g\'{e}n\'{e}ralement, une sous-vari\'{e}t\'{e} $i:M\ra M'$ de $M'$ telle que $M$ soit munie
d'une connexion $\nabla$ dans le fibr\'{e} tangent est dite {\em totalement
g\'{e}od\'{e}sique} lorsque l'injection $i$ est affine.

Le prochain corollaire g\'{e}n\'{e}ralise un r\'{e}sultat de Xu pour les sous-vari\'{e}t\'{e}s
lagrangiennes (voir \cite{Xu98}):

\bcor
 Soit $i:C\ra M$ une sous-vari\'{e}t\'{e} co\"{\i}sotrope ferm\'{e}e d'une vari\'{e}t\'{e}
 symplectique
 $(M,\omega)$. Soit $\nabla$ une connexion adapt\'{e}e sur la vari\'{e}t\'{e}
 pr\'{e}symplectique $(C,\varpi=i^*\omega)$. Alors il existe une connexion
 symplectique $\tilde{\nabla}$ dans le fibr\'{e} tangent de $M$ telle que $C$
 soit une sous-vari\'{e}t\'{e} totalement g\'{e}od\'{e}sique de $M$.
\ecor
\bbew
 Gr\^{a}ce au th\'{e}or\`{e}me de Weinstein/Gotay \ref{TWeinGotay} il existe un
 voisinage tubulaire $U$ de $C$ qui est symplectomorphe \`{a} un voisinage
 $V$ de la z\'{e}ro-section $C$ dans $E^*$ munie de la forme $\omega_{WG}$.
 Par cons\'{e}quent, le th\'{e}or\`{e}me \ref{TConnSympAutourCoisotrope} nous donne
 une connexion symplectique $\tilde{\nabla}^U$ avec les propri\'{e}t\'{e}s
 souhait\'{e}es dans $U$: en fait, pour tout $X\in\Ginf(C,TC)$ le couple
 $(X,X^h)$ est toujours $i$-li\'{e}, donc l'\'{e}quation (\ref{EqConnSympHorHor})
 montre que $(\nabla_XY, \tilde{\nabla}_{X^h}Y^h)$ est $i$-li\'{e} (on met
 $\varphi=0\in E_{\tau(\varphi)}$). Soit maintenant $U'$ le voisinage ouvert
 correspondant \`{a}
 $V':=\frac{1}{2}V\subset V$. Alors l'adh\'{e}rence $\overline{U'}$ de $U'$
 est contenu dans
 $U$, donc les ouverts $U$ et $W:=M\setminus \overline{U'}$ recouvrent $M$.
 En
 utilisant une partition de l'unit\'{e} $1=\psi_U+\psi_W$ subordonn\'{e}e \`{a}
 $(U,W)$ et une connexion symplectique $\tilde{\nabla}^W$ arbitraire
 dans $W$, on obtient une connexion symplectique globale
 $\tilde{\nabla}:=\psi_U\tilde{\nabla}^U+\psi_W\tilde{\nabla}^W$.
\ebew

\subsection{La classe d'Atiyah-Molino d'une vari\'{e}t\'{e} pr\'{e}symplec\-tique}
 \label{SubSecClassedAtiyahMolino}

\noindent Nous d\'{e}montrons d'abord le lemme suivant:

\blem \label{LSubmersionAffine}
 Soit $\phi:M\ra M'$ une submersion surjective entre deux vari\'{e}t\'{e}s
 diff\'{e}rentiables.
 On note $E:=Ker T\phi$.
 \ben
 \item Soit $\nabla$ une connexion sans torsion dans le
 fibr\'{e} tangent de $M$. On consid\`{e}re les deux \'{e}nonc\'{e}s suivants:
  \ben
   \item Il existe une unique connexion $\nabla'$ sans torsion dans le fibr\'{e}
  tangent de $M'$ telle que $\phi$ soit une application affine.
   \item La connexion $\nabla$ pr\'{e}serve $E$ (i.e. $\nabla_XV\in\Ginf(M,E)$
  quels que soient $X\in\Ginf(M,TM),V\in\Ginf(M,E)$), et le tenseur de
  courbure
  $\overline{R}$  de la connexion
  $\overline{\nabla}$ induite dans le fibr\'{e} $TM/E$ a la propri\'{e}t\'{e}
  suivante:
  \[
     \overline{R}(V,X)\overline{Y}=0~~~~\forall~ X,Y\in\Ginf(M,TM),
                                         \forall~ V\in\Ginf(M,E)
  \]
  \een
  Alors $(a)$ implique toujours $(b)$. Dans le cas o\`{u}
  toutes les fibres de $\phi$ soient connexes, $(b)$ implique $(a)$.
 \item Soit $\nabla'$ une connexion sans torsion dans le fibr\'{e} tangent de
   $M'$. Alors il existe une connexion sans torsion $\nabla$ dans le fibr\'{e}
   tangent de $M$ telle que $\phi$ soit une application affine.
 \een
\elem
\bbew
 Il est clair que $E$ est un sous-fibr\'{e} int\'{e}grable.
 On choisit un sous-fibr\'{e} $F$ de $TM$ compl\'{e}mentaire \`{a} $E$.
 Par
 cons\'{e}quent, pour tout champ de vecteurs $X'$ sur $M'$ la prescription
 $T_m\phi~X^{\prime h}_m:=X'_{\phi(m)}$ d\'{e}finit un unique champ de
 vecteurs $X^{\prime h}$ sur $M$ \`{a} valeurs dans $F$ (un rel\`{e}vement
 horizontal). Par d\'{e}finition, le couple $(X^{\prime h},X')$ est toujours
 $\phi$-li\'{e}. Puisque $(V,0)$ est $\phi$-li\'{e} pour tout champ de vecteurs
 vertical, le crochet de Lie $[X^{\prime h},V]$ est vertical. De plus,
 tout vecteur tangent de $M$ se repr\'{e}sente comme une somme de la valeur
 d'un rel\`{e}vement horizontal et de la valeur d'un champ de vecteurs
 vertical.\\
 1) ``$(a)\Longrightarrow (b)$'': Soit $\phi$ affine par rapport \`{a} une
  connexion
 $\nabla'$ dans le fibr\'{e} tangent de $M'$. Soient $V,W\in\Ginf(M,E)$ et
 $X',Y'\in\Ginf(M',TM')$. Il vient que les couples
 $(\nabla_{X^{\prime h}}W,0)$ et $(\nabla_VW,0)$ sont $\phi$-li\'{e}s, donc
 $\nabla_XW$ est un champ de vecteurs vertical pour tout
 $X\in\Ginf(M,TM)$, et $\nabla$ pr\'{e}serve $E$. De plus, le couple
 $(\nabla_VX^{\prime h},0)$ est \'{e}galement $\phi$-li\'{e}, alors le champ de
 vecteurs $\nabla_VX^{\prime h}$ est aussi vertical. Puisque $\nabla$
 pr\'{e}serve $E$ il s'ensuit que $R(X,Y)V$ est vertical. De plus, les champs
 $\nabla_WY^{\prime h}$, $\nabla_VY^{\prime h}$ et
 $\nabla_{[V,W]}Y^{\prime h}$ sont verticaux, donc le tenseur de
 courbure
 \[
   R(V,W)Y^{\prime h}=\nabla_V\nabla_WY^{\prime h}
                        -\nabla_W\nabla_VY^{\prime h}
                          -\nabla_{[V,W]}Y^{\prime h}
 \]
 est vertical, d'o\`{u} $\overline{R}(V,W)\overline{Y}=0$ quel que soit
 $Y\in\Ginf(M,TM)$. Finalement, puisque le couple
 $(\nabla_{X^{\prime h}}Y^{\prime h},\nabla'_{X'}Y')$ est $\phi$-li\'{e},
 alors le couple $(\nabla_V\nabla_{X^{\prime h}}Y^{\prime h},0)$ l'est,
 donc $\nabla_V\nabla_{X^{\prime h}}Y^{\prime h}$ est vertical et
 le tenseur de courbure
 \[
   R(V,X^{\prime h})Y^{\prime h}=\nabla_V\nabla_{X^{\prime h}}Y^{\prime h}
                        -\nabla_{X^{\prime h}}\nabla_VY^{\prime h}
                          -\nabla_{[V,X^{\prime h}]}Y^{\prime h}
 \]
 est donc vertical, d'o\`{u} $\overline{R}(V,X^{\prime h})\overline{Y}=0$ quel
 que soit $Y\in\Ginf(M,TM)$.\\
 ``$(a)\Longleftarrow (b)$'': Soient $(X,X')$ et $(Y,Y')$ deux couples
 $\phi$-li\'{e}s. En utilisant la d\'{e}composition $TM=F\oplus E$ on voit que
 $X=X^{\prime h}+W_1$ et $Y=Y^{\prime h}+W_2$ pour certains
 $W_1,W_2\in\Ginf(M,E)$. Par hypoth\`{e}se, pour tout champ de vecteurs
 vertical $V$ le champ de vecteurs $\nabla_ZV$ est vertical, donc
 $\nabla_{X^{\prime h}}W_2$ est vertical. Puisque $(X^{\prime h},X')$
 est $\phi$-li\'{e}, alors $([X^{\prime h},W_2],0)$ l'est, alors
 $[X^{\prime h},W_2]=\nabla_{X^{\prime h}}W_2-\nabla_{W_2}X^{\prime h}$
 est vertical, donc $\nabla_{W_2}X^{\prime h}$ est vertical. Par
 cons\'{e}quent $T\phi~\nabla_XY=T\phi~\nabla_{X^{\prime h}}Y^{\prime h}$.
 Soient $m_0,m_1\in M$ dans la m\^{e}me fibre ($\phi(m_0)=\phi(m_1)=:m'$). Puisque
 chaque fibre est connexe, on peut utiliser une cha\^{\i}ne de cartes
 $U_1,\ldots,U_k$ avec $U_i\cap U_{i+1}\neq \emptyset$ et $m_0\in U_1$, $m_1\in
 U_k$, pour construire un champ vertical $W$ \`{a} support compact
 dont le flot (forc\'{e}ment complet) $\Psi_t$ repr\'{e}sente cette courbe,
 i.e. $t\mapsto \Psi_t(m_0)$ avec $\Psi_1(m_0)=m_1$. Il vient
 \beas
 \lefteqn{ \frac{d}{dt}\big(\Psi_t^*(\nabla_{X^{\prime h}}Y^{\prime h}\big)
     =  \Psi_t^*\big([W,\nabla_{X^{\prime h}}Y^{\prime h}]\big)} \\
     & = &
       \Psi_t^*\left(R(W,X^{\prime h})Y^{\prime h}
          +\nabla_{X^{\prime h}}\nabla_W Y^{\prime h}
          +\nabla_{[W,X^{\prime h}]}Y^{\prime h}
          -\nabla_{\nabla_{X^{\prime h}}Y^{\prime h}}W\right) \\
     & =: & \tilde{W}(t)
 \eeas
 Puisque $R(W,X^{\prime h})Y^{\prime h}$ est vertical par hypoth\`{e}se et les
 trois autres termes le sont d'apr\`{e}s ce qui pr\'{e}c\`{e}de, il vient que
 $\tilde{W}(t)$ consiste en champs verticaux. En int\'{e}grant cette \'{e}quation
 de $0$ \`{a} $1$ on arrive \`{a}
 \[
  \big(\nabla_{X^{\prime h}}Y^{\prime h}\big)_{m_1} =
     T_{m_0}\Psi_1 \big(\nabla_{X^{\prime h}}Y^{\prime h}\big)_{m_0}
     +T_{m_0}\Psi_1 \big(\int_0^1~\tilde{W}(s)~ds\big),
 \]
 et cette \'{e}quation montre -- gr\^{a}ce \`{a} $\phi\circ\Psi_t=\phi$ -- pour
 tous $m_0,m_1\in M$ avec $\phi(m_0)=\phi(m_1)=:m'$
 \[
   T_{m_1}\phi~\big(\nabla_{X^{\prime h}}Y^{\prime h}\big)_{m_1} =
    T_{m_0}\phi~\big(\nabla_{X^{\prime h}}Y^{\prime h}\big)_{m_0}
    =:(\nabla'_{X'}Y')_{m'}.
 \]
 On v\'{e}rifie sans peine que $\nabla'$ est une connexion sans torsion
 sur $M'$. Donc le couple $(\nabla_XY,\nabla'_{X'}Y')$ est $\phi$-li\'{e}
 et $\phi$ est affine.\\
 Puisque $T\phi$ est surjective, il est clair que $\nabla'$ est unique.\\
 2) On note que le fibr\'{e} retir\'{e} $\phi^*TM'$ est isomorphe \`{a} $F$, donc
 la connexion retir\'{e}e $\phi^*\nabla'$ d\'{e}finit une connexion $\check{\nabla}$
 dans le fibr\'{e} $F$ sur $M$. Soient $X',Y'$ deux champs de vecteurs sur
 $M'$. Leurs rel\`{e}vements horizontaux $X^{\prime h},Y^{\prime h}$
 s'interpr\`{e}tent comme des sections retir\'{e}es dans $\phi^*TM'\cong F$ et on
 a $\check{\nabla}_{X^{\prime h}}Y^{\prime h}=(\nabla'_{X'}Y')^h$. Il
 vient que $\check{\nabla}_{X^{\prime h}}Y^{\prime h}
 -\check{\nabla}_{Y^{\prime h}}X^{\prime h}
 =\mathsf{P}_F[X^{\prime h},Y^{\prime h}]$, donc pour avoir une connexion
 sans torsion le long de $F$ on peut d\'{e}finir pour deux champs horizontaux
 $H_1,H_2\in\Ginf(C,F)$
 \[
      \nabla_{H_1}H_2:=\check{\nabla}_{H_1}H_2+\frac{1}{2}
                          \hat{R}^F(H_1,H_2).
 \]
 De plus, pour $V,W\in\Ginf(C,E)$ on d\'{e}finit
 \[
    \nabla_{H_1}W:=\mathsf{P}_E[H_1,W] ~~~\mathrm{et}~~~
         \nabla_{W}H_1:=\mathsf{P}_F[W,H_1],
 \]
 et on compl\`{e}te $\nabla$ sur $V$ et $W$ par la partie verticale-verticale
 d'une connexion adapt\'{e}e arbitraire (voir par exemple la proposition
 \ref{PExistConnAdaptees}). Soient $(X,X')$ et $(Y,Y')$ deux couples
 $\phi$-li\'{e}s. Donc --comme ci-dessus-- $X=X^{\prime h}+V$ et
 $Y=Y^{\prime h}+W$ pour certains $V,W\in\Ginf(C,E)$. Il vient que
 $\nabla_VY^{\prime h}=\mathsf{P}_F[V,Y^{\prime h}]=0$, alors
 \[
   T\phi\nabla_XY=T\phi\nabla_{X^{\prime h}}Y^{\prime h}
                    +T\phi\nabla_{X^{\prime h}}W +T\phi\nabla_VW
                 =\nabla'_{X'}Y'\circ \phi,
 \]
 donc $\phi$ est affine.
\ebew

\noindent Dans la situation du lemme pr\'{e}c\'{e}dent, on va dire que la connexion
$\nabla$ soit {\em projetable} (par rapport \`{a} $\phi$) si
$\phi$ est affine,
et dans ce cas on appelle $\nabla'$ la {\em connexion induite
par $\nabla$ et $\phi$}.

Soit $(C,\mathcal{F})$ une vari\'{e}t\'{e} diff\'{e}rentiable munie d'un feuilletage
r\'{e}gulier. Soit $E\subset TC$ son fibr\'{e} caract\'{e}ristique. On regarde une
carte distingu\'{e}e du feuilletage $\big(U,(\xi,x)\big)$ o\`{u}
$(\xi,x)=(\xi^1,\ldots,\xi^m,x^1,\ldots,x^k)\big)$ (dont
 les coordonn\'{e}es
$\xi$ sont transverses \`{a} et les coordonn\'{e}es $x$ sont le long des plaques).
Ces derni\`{e}res d\'{e}finissent une submersion locale $\pi_U$ de $U$ \`{a}
{\em l'espace quotient local} $U_{\mathrm{red}}$ qui est de la forme
$\pi_U(\xi,x):=\xi$. Dans le cas o\`{u} $(C,\varpi)$ est une vari\'{e}t\'{e}
pr\'{e}symplectique on note que la prescription
\beq
    \pi_U^*\varpi_{U_{\mathrm{red}}}:=\varpi|_{U}
\end{equation}
d\'{e}finit une forme symplectique $\varpi_{U_{\mathrm{red}}}$ sur l'espace
quotient local.

\bdefi[d'\`{a}pr\`{e}s P.Molino, 1971] Une connexion sans torsion $\nabla$ dans le
  fibr\'{e} tangent de $C$ est dite {\em localement projetable} lorsque pour toute
  carte distingu\'{e}e du feuilletage $\big(U,(\xi,x)\big)$
  la restriction de $\nabla$ \`{a} $U$
  est projetable par rapport \`{a} $\pi_U$.
\edefi
Voir les travaux de P.Molino \cite{Mol71a} et \cite{Mol71b} pour
une d\'{e}finition plus g\'{e}n\'{e}rale.
Soit maintenant $F$ un sous-fibr\'{e} de $TC$ compl\'{e}mentaire \`{a} $E$ et
$\nabla$ une connexion adapt\'{e}e dans le fibr\'{e} tangent de $C$
(qui existe toujours gr\^{a}ce \`{a} la proposition \ref{PExistConnAdaptees}). Pour le
tenseur de courbure $R$ de $\nabla$ on  a montr\'{e} que $R(V,W)W'$ est vertical
et que $R(V,W)H=0$ (la propri\'{e}t\'{e} de Bott, voir aussi l'\'{e}q.
(\ref{EConnexionBottEstPlate})) quels que soient les champs de
vecteurs
verticaux $V,W,W'$ et horizontal $H$. Donc sur le fibr\'{e} transversal
$Q=TC/E$ avec la
connexion induite $\overline{\nabla}$ la section suivante
\beq
  (V,\overline{X},\overline{Y})\mapsto r_{AM}(V)(\overline{X},\overline{Y})
                 := \overline{R}(V,X)\overline{Y}
\end{equation}
dans $\Ginf(C,\Lambda^1E^*\otimes S^2Q^*\otimes Q)$ est bien d\'{e}finie:
la sym\'{e}trie
en $(\overline{X},\overline{Y})$ se d\'{e}duit de la premi\`{e}re identit\'{e} de Bianchi
pour $R$ et du fait que $R(X,Y)V$ est vertical. Le th\'{e}or\`{e}me suivant
est d\^{u} \`{a} P.Molino, \cite{Mol71a}:
\bsat[d'\`{a}pr\`{e}s P.Molino, 1971] \label{TMolino}
 Avec les hypoth\`{e}ses faites ci-dessus on a:
 \ben
  \item $r_{AM}$ est un $1$-cocycle de la cohomologie longitudinale \`{a} valeurs
        dans
          $\Ginf(C,S^2Q^*\otimes Q)$.
  \item La classe de cohomologie $c_{AM}(C,E)$ de $r_{AM}$ ne d\'{e}pend pas de
     la connexion
      adapt\'{e}e $(F,\nabla)$ choisie et est donc un invariant de la
      topologie dif\-f\'{e}\-ren\-tielle de la vari\'{e}t\'{e} feuillet\'{e}e.
  \item $c_{AM}(C,E)=0$ si et seulement s'il existe une connexion adapt\'{e}e
      localement projetable sur $C$.
 \een
\esat
\bbew
 1. Il vient
  \beas
   \lefteqn{(d_vr_{AM})\big(V,W\big)\big(H_1,H_2\big) = } \\
     & & \overline{\nabla_V\big(\mathsf{P}_FR(W,H_1)H_2\big)}
       -\overline{\mathsf{P}_FR(W,\nabla_VH_1)H_2}
       -\overline{\mathsf{P}_FR(W,H_1)\nabla_VH_2 }  \\
     & &  -\overline{\nabla_W\big(\mathsf{P}_FR(V,H_1)H_2\big)}
            +\overline{\mathsf{P}_FR(V,\nabla_WH_1)H_2}
            +\overline{\mathsf{P}_FR(V,H_1)\nabla_WH_2}   \\
     & & -\overline{\mathsf{P}_FR([V,W],H_1)H_2} \\
     & = & \overline{\mathsf{P}_F(\nabla_VR)(W,H_1)H_2}
              +\overline{\mathsf{P}_F(\nabla_WR)(H_1,V)H_2} \\
     & = & 0
  \eeas
  gr\^{a}ce \`{a} la deuxi\`{e}me identit\'{e} de Bianchi pour
  le tenseur de courbure $R$ de $\nabla$ et au fait que $\nabla\mathsf{P}_F=0$,
  car $\mathsf{P}_F(\nabla_{H_1}R)(V,W)H_2=0$ (Bott).\\
  2. Soit $(F',\nabla')$ une autre connexion sans torsion adapt\'{e}e.
  La diff\'{e}rence
  $\nabla'-\nabla$ d\'{e}finit un champ de tenseurs
  $S\in\Ginf(C,S^2T^*C\otimes TC)$ pour lequel on a
  $S(V,X)\in\Ginf(C,E)$ quels que soient $V\in\Ginf(C,E), X\in\Ginf(C,TC)$
  parce que $\nabla$ et $\nabla'$ pr\'{e}servent $E$.
  Par cons\'{e}quent, il vient pour tous $X,Y\in\Ginf(C,TC)$ et
  $V\in\Ginf(C,E)$:
  \beas
    \overline{R'(V,X)Y} & = & \overline{(\nabla'_{V}\nabla'_{X}Y
                                    -\nabla'_{X}\nabla'_{V}Y
                                    -\nabla'_{[V,X]}Y)} \\
                   & = & \overline{(\nabla_{V}\nabla'_{X}Y
                                    -\nabla'_{X}\nabla_{V}Y
                                    -\nabla'_{[V,X]}Y)} \\
                   & = & \overline{R(V,X)Y} \\
                   &   &   +\overline{\big(\nabla_V\big(S(X,Y)
                           -S(X,\nabla_VY)-S(\nabla_VX,Y\big)\big)}
  \eeas
  Avec la d\'{e}finition $\overline{S}(\overline{X},\overline{Y}):=
  \overline{S(X,Y)}$ ceci implique
  \beq\label{EqFormeAMChangeConn}
   r'_{AM}(V)(\overline{X},\overline{Y}) =
       r_{AM}(V)(\overline{X},\overline{Y})+
          (d_v\overline{S})(V)(\overline{X},\overline{Y}).
  \end{equation}
  Alors la classe de $r_{AM}$ ne d\'{e}pend pas de la connexion adapt\'{e}e.\\
  3. Soit d'abord $\nabla$ localement projetable. D'apr\`{e}s le lemme
  \ref{LSubmersionAffine} il s'ensuit que $r_{AM}$ s'annule, donc
  $c_{AM}(C,E)=0$.\\
  R\'{e}ciproquement, soit $(F,\nabla)$ une connexion sans
  torsion adapt\'{e}e telle que la classe d'Atiyah-Molino $c_{AM}(C,E)$
  s'annule. Alors la connexion $\nabla$ pr\'{e}serve $E$
  et il existe une section $T\in\Ginf(C,S^2Q^*\otimes Q)$ telle que
  $r_{AM}=d_vT$. Soit $S\in\Ginf(C,S^2T*C\otimes F)$ l'unique champ de
  tenseurs tel que
  $\overline{S(X,Y)}=T(\overline{X},\overline{Y})$ et $S(V,X)=0$ quels que
  soient $X,Y\in\Ginf(C,TC)$, $V\in\Ginf(C,E)$ et $S(H_1,H_2)\in\Ginf(C,F)$
  pour tous $H_1,H_2\in\Ginf(C,F)$. Soit
  $\nabla'$ la connexion sans torsion $\nabla-S$. Alors $(F,\nabla')$ est
  adapt\'{e}e (donc pr\'{e}serve $E$) et $r'_{AM}=r_{AM}-d_vT=0$ d'apr\`{e}s l'\'{e}quation
  (\ref{EqFormeAMChangeConn}), donc $\nabla'$
  est localement projetable d'apr\`{e}s le lemme \ref{LSubmersionAffine}.
\ebew

\noindent On va appeler $r_{AM}$ le {\em $1$-cocycle d'Atiyah-Molino
(relatif \`{a} la connexion adap\-t\'{e}e $(F,\nabla)$)} et $c_{AM}(C,E):=[r_{AM}]$ la
{\em classe d'Atiyah-Molino de la vari\'{e}t\'{e} feuille\-t\'{e}e $(C,\mathcal{F})$}.
En outre, la partie $2)$ du lemme \ref{LSubmersionAffine} montre que la
classe
d'Atiyah-Molino s'annule toujours si $E$ est donn\'{e} par le noyau
de l'application tangente d'une submersion surjective. Le feuilletage
de codimension $1$
de la $3$-sph\`{e}re d\^{u} \`{a} G.Reeb est un exemple d'une vari\'{e}t\'{e}
feuillet\'{e}e dont la classe d'Atiyah-Molino ne s'annule pas, voir
\cite[p.66]{Mol88}. D'un autre c\^{o}t\'{e}, les {\em feuilletages riemanniens}
sont toujours de classe d'Atiyah-Molino nulle: ce sont des feuilletages
qui admettent une m\'{e}trique riemannienne $\gamma$ sur le fibr\'{e} transverse
$Q=TC/E$ telle que $\bar{\nabla}^{\mathrm{Bott}}_V\gamma=0$ quels que
soient les champs de vecteurs verticaux, voir \cite[p.81]{Mol88}. On en
d\'{e}duit le
\bcor \label{CActionsPropresClasseAMNulle}
 Soit $C$ une vari\'{e}t\'{e} diff\'{e}rentiable et $\Phi:G\times C\ra C$ l'action
 propre d'un groupe de Lie connexe dont toutes les orbites ont la m\^{e}me
 dimension. Alors le feuilletage de $C$ en orbites de $G$ est de classe
 d'Atiyah-Molino nulle.
\ecor
\bbew
 D'apr\`{e}s le th\'{e}or\`{e}me classique de R.Palais \cite{Pal61}, il existe une
 m\'{e}trique riemannienne $\mathsf{g}$ sur $C$ qui est invariante par
 l'action de $G$. Soit $E$ le fibr\'{e} tangent aux orbites et $F$ le sous-fibr\'{e}
 de $TC$ orthogonal \`{a} $E$. Soit $\gamma$ la restriction de
 $\mathsf{g}$ \`{a} $F\otimes F$. La condition $L_{\xi_C}g=0$ quel que soit
 $\xi\in\mathfrak{g}$, l'alg\'{e}bre de Lie de $G$, implique directement
 que $\bar{\nabla}^{\mathrm{Bott}}_V\gamma=0$ quels que soient les
 sections $V$ de $\Ginf(C,E)$. Alors le feuilletage est riemannien et
 la classe d'Atiyah-Molino s'anulle.
\ebew
Si $G$ est compact, son action est automatiquement propre.

Pour les vari\'{e}t\'{e}s pr\'{e}symplectiques $(C,\varpi)$ il y a l'analogue suivant du
$1$-cocycle d'Atiyah-Molino: soit $(F,\nabla)$ une connexion adapt\'{e}e
et $r_{AM}$ le $1$-cocycle d'Atiyah-Molino correspondant.
Puisque le noyau de $\varpi$ est donn\'{e} par le fibr\'{e} caract\'{e}ristique $E$,
le fibr\'{e} normal $Q=TC/E$ est muni d'une $2$-forme $\overline{\varpi}$
symplectique d\'{e}finie par $\overline{\varpi}(\overline{X},\overline{Y})
:=\varpi(X,Y)$ quels que soient les champs de vecteurs $X,Y$ sur $C$.
Soient $Z\in\Ginf(C,TC)$ et $V\in\Ginf(C,E)$. On d\'{e}finit
\beq\label{EqCocycleAtiyahMolinoPresymp}
  \rho_{AM}(V)(\overline{X},\overline{Y},\overline{Z})
    := \overline{\varpi}\big(\overline{X},
              r_{AM}(V)(\overline{Y},\overline{Z})\big)
              =
              \varpi\big(X,R(V,Y)Z\big).
\end{equation}
On a le
\bsat\label{TMolinoPresymp}
 Avec les hypoth\`{e}ses faites ci-dessus on a:
 \ben
  \item $\rho_{AM}$ est un $1$-cocycle de la cohomologie longitudinale
     \`{a} valeurs dans $\Ginf(C,S^3Q^*)$.
  \item La classe de cohomologie $\kappa_{AM}(C,E):=[\rho_{AM}]$ de $\rho_{AM}$
    ne d\'{e}pend pas de
     la connexion adapt\'{e}e $(F,\nabla)$ choisie et est donc un invariant
     de la topologie diff\'{e}rentielle de la vari\'{e}t\'{e} pr\'{e}symplectique $(C,\varpi)$.
  \item $\kappa_{AM}(C,E)=0$ si et seulement s'il existe une connexion adapt\'{e}e
      localement projetable sur $(C,\varpi)$. En outre
      $\kappa_{AM}(C,E)=0$ si et seulement $c_{AM}(C,E)=0$. Si c'est le cas,
      la connexion induite
      sur chaque espace r\'{e}duit local
      $(U_{\mathrm{red}},\varpi_{U_{\mathrm{red}}})$ pr\'{e}serve la forme
      symplectique r\'{e}duite $\varpi_{U_{\mathrm{red}}}$.
 \een
\esat
\bbew
1. Puisque $\nabla$ pr\'{e}serve $\varpi$ on a
\[
  \varpi\big(X,R(V,Y)Z\big)=-\varpi\big(R(V,Y)X,Z\big)
=\varpi\big(Z,R(V,Y)X)\big),
\]
donc $\rho_{AM}\in\Ginf(C,\Lambda^1E^*\otimes S^3Q^*)$. L'\'{e}quation
$d_v\rho_{AM}=0$ se d\'{e}duit de $d_vr_{AM}=0$ et du fait que
$\nabla\varpi=0$.\\
2. Analogue \`{a} l'enonc\'{e} 2) du th\'{e}or\`{e}me \ref{TMolino} de Molino.\\
3. La premi\`{e}re phrase se d\'{e}montre comme l'\'{e}nonc\'{e} analogue du th\'{e}or\`{e}me
  \ref{TMolino}. \\
  Soit $\kappa_{AM}(C,E)=0$. Il existe donc une connexion sans torsion $\nabla$
  adapt\'{e}e \`{a} $(C,\varpi)$ qui est localement projetable. Alors
  $c_{AM}(C,E)=[r_{AM}]=0$.
  R\'{e}ciproquement, soit $c_{AM}(C,E)=0$. Alors il existe une connexion
  sans torsion $\nabla'$ adapt\'{e}e localement projetable, alors $r'_{AM}=0$.
  En ajoutant un champ de tenseurs $S$ comme (\ref{EqAstuceTondeur}) on
  peut faire en sorte que $\nabla:=\nabla'+S$ pr\'{e}serve $\varpi$. On a
  $r_{AM}=r'_{AM}+d_v\overline{S}=d_v\overline{S}$. Il vient que
  $\rho'_{AM}=d_v\overline{T}$ avec $T(X,Y,Z):=\varpi\big(X,S(Y,Z)\big)$.
  $\overline{T}$ n'est pas forc\'{e}ment totalement sym\'{e}trique en
  $\overline{X},\overline{Y},\overline{Z}$. Soit $6\mathsf{S}$ la somme
  sur toutes les permutations d'un ensemble de trois \'{e}l\'{e}ments, repr\'{e}sent\'{e}e
  sur $\Ginf(C,Q^*\otimes Q^*\otimes Q^*)$. Par cons\'{e}quent,
  \[
  \rho'_{AM}=\mathsf{S}\rho'_{AM}=\mathsf{S}d_v\overline{T}
         =d_v\mathsf{S}\overline{T}=:d_v\overline{T'}
  \]
  o\`{u} $\overline{T'}\in\Ginf(C,S^3Q^*)$, donc $[\rho'_{AM}]=0$.\\
  Puisque $\pi_U^*\varpi_{U_{\mathrm{red}}}=\varpi|_{U_{\mathrm{red}}}$
  on peut utiliser les rel\`{e}vements horizontaux des champs de vecteurs
  sur $U_{\mathrm{red}}$ \`{a} $U$ pour montrer que la connexion induite
  pr\'{e}serve $\varpi_{U_{\mathrm{red}}}$.
\ebew

\noindent On va appeler $\rho_{AM}$ le {\em $1$-cocycle d'Atiyah-Molino
(relatif \`{a} la connexion adapt\'{e}e $(F,\nabla)$)} et
$\kappa_{AM}(C,E):=[\rho_{AM}]$ la
{\em classe d'Atiyah-Molino de la vari\'{e}t\'{e} pr\'{e}sym\-plec\-tique $(C,\varpi)$}.

\brem
 La g\'{e}n\'{e}ralisation de la
{\em classe de Godbillon-Vey} d'une vari\'{e}t\'{e}
pr\'{e}symplectique $(C,\varpi)$ en tant que vari\'{e}t\'{e} feuillet\'{e}e (voir
par exemple \cite[p.30]{Ton88}) s'annule toujours: en fait, si les feuilles
sont de dimension $k$ et $\dim C=2m+k$, alors la $m$i\`{e}me
puissance de $\varpi$ est une $2m$-forme ferm\'{e}e dont le noyau est
\'{e}gal au sous-fibr\'{e} caract\'{e}ristique $E$ de $\varpi$. Donc la forme de
Godbillon-Vey (qui est de rang $4m+1$) s'annule.
\erem

Pour les vari\'{e}t\'{e}s pr\'{e}symplectiques on a la structure quadratique
suivante sur les classes de cohomologie longitudinale: la $2$-forme
$\varpi$ induit une forme symplectique $\overline{\varpi}$ sur
le fibr\'{e} normal $Q=TC/E$. Il existe donc un unique champ de bivecteurs
$\overline{P}\in\Ginf(C,\Lambda^2Q)$ tel que
\[
   \overline{\varpi}^\flat~\overline{P}^\sharp = \mathrm{id}_{Q*}~~~
   \mathrm{et}~~~
   \overline{P}^\sharp~\overline{\varpi}^\flat = \mathrm{id}_{Q}.
\]
Pour tout entier positif $k$, $\overline{P}$ d\'{e}finit une section
$\overline{P}^{(k)}\in\Ginf(C,S^kQ\otimes S^kQ)$ d\'{e}finie de mani\`{e}re
usuelle par
\beq\label{EqDefProduitSymetriquePoisson}
  \overline{P}^{(k)}\big(\zeta_1^k,\zeta_2^k\big):=
         \big(\overline{P}(\zeta_1,\zeta_2)\big)^k
\end{equation}
quels que soient $\zeta_1,\zeta_2\in\Ginf(C,Q^*)$. Ceci nous donne
un accouplement
\beq
  \overline{P}^{(k)}:\Ginf(C,\Lambda^pE^*\otimes S^kQ^*)\times
         \Ginf(C,\Lambda^qE^*\otimes S^kQ^*)\ra
         \Ginf(C,\Lambda^{p+q}E^*)
\end{equation}
d\'{e}fini de fa\c{c}on usuelle:
soient
$\alpha\in\Ginf(C,\Lambda^pE^*)$,
$\beta\in\Ginf(C,\Lambda^qE^*)$ et
$\tau_1,\tau_2\in\Ginf(C,S^kQ^*)$, alors
\beq
  \overline{P}^{(k)}(\alpha\otimes \tau_1,\beta\otimes \tau_2)
  :=   \overline{P}^{(k)}(\tau_1,\tau_2)\alpha\wedge\beta.
\end{equation}
Puisque $L_V\varpi=0$ quel que soit le champ vertical $V$, alors
$\overline{\nabla}_V\overline{\varpi}=0$ pour la connexion de Bott,
donc $\overline{\nabla}_V\overline{P}^{(k)}=0$, et il vient
pour tous $A\in\Ginf(C,\Lambda^pE^*\otimes S^kQ^*)$ et
$B\in\Ginf(C,\Lambda^qE^*\otimes S^kQ^*)$
\beq
   d_v\big(\overline{P}^{(k)}(A,B)\big)=
   \overline{P}^{(k)}(d_vA,B)+(-1)^p\overline{P}^{(k)}(A,d_vB),
\end{equation}
donc l'accouplement $\overline{P}^{(k)}$ descend \`{a} un accouplement
--aussi not\'{e} $\overline{P}^{(k)}$-- des groupes de cohomologie
\beq
 \overline{P}^{(k)}:H^p_v(C,S^kQ^*)\times H^q_v(C,S^kQ^*)\ra
                    H^{p+q}_v(C).
\end{equation}

\subsection{Les obstructions \`{a} la repr\'{e}sentabilit\'{e}
       jusqu'\`{a} l'ordre $3$}
  \label{SubSecObstrRepOrdreTrois}

Le but de ce paragraphe est d'esquisser la d\'{e}monstration du th\'{e}or\`{e}me
suivant:
\bsat\label{TObstructionsOrdreTrois}
 Soit $(M,\omega)$ une vari\'{e}t\'{e} symplectique et $i:C\ra M$ une
 sous-vari\'{e}t\'{e} co\"{\i}sotrope ferm\'{e}e. Soit $E$ le sous-fibr\'{e} caract\'{e}ristique
 de $\varpi:=i^*\omega$. Soit $*$ un star-produit sur
 $M$ et $[*]=\sum_{r=-1}^\infty\nu^r[*]_r$ sa classe de Deligne.
 Pour que $*$ soit adapt\'{e} \`{a} $C$ jusqu'\`{a} l'ordre $3$ il faut et il
 suffit que les deux conditions suivantes soient satisfaites:
 \ben
  \item Il existe un repr\'{e}sentant $\alpha_0$ de
    $[*]_0$ telle que $i^*\alpha_0$ est une $2$-forme relative, i.e.
    \beq\label{EqObstrClasseDelZero}
        p_v i^*[*]_0=0
    \end{equation}
  \item De plus, $\alpha_0$ a la propri\'{e}t\'{e} suivante:
    \bea
      \lefteqn{\frac{1}{12}\overline{P}^{(3)}
              \big(\kappa_{AM}(C,E),\kappa_{AM}(C,E)\big)} \nonumber \\
        & & ~~~~~~~~~~~+\frac{1}{2}
          \overline{P}^{(1)}([i^*\alpha_0]_{(1,1)},[i^*\alpha_0]_{(1,1)})
        -p_vi^*[*]_1=0 \label{EqObstructionImportante}
    \eea
 \een
 o\`{u} $\kappa_{AM}(C,E)$ d\'{e}signe la classe d'Atiyah-Molino de la sous-vari\'{e}t\'{e}
 pr\'{e}\-sym\-plec\-tique $(C,\varpi:=i^*\omega)$ et $[~~]_{(1,1)}$ d\'{e}signe la
 classe de cohomologie d'une $2$-forme relative vue comme $1$-forme
 logitudinale \`{a} valeurs dans $Q^*$.
\esat

On va d'abord donner la formule pour un star-produit symplectique
g\'{e}n\'{e}ral jusqu'\`{a} l'ordre $3$: soit $P$ la structure de Poisson correspondant
\`{a} $\omega$, et soit $P^{(k)}\in\Ginf(S^kTM\otimes S^kTM)$ d\'{e}fini comme
$\overline{P}^{(k)}$, voir (\ref{EqDefProduitSymetriquePoisson}).
Soit $\tilde{\nabla}$ une connexion sans
torsion symplectique dans le fibr\'{e} tangent de $M$. On rappelle l'op\'{e}rateur
bidiff\'{e}rentiel suivant pour $f,h\in\CinfK{M}$ (voir \cite{Lic82},
\cite{DL88}):
\beq \label{EqDefP2Gamma}
 P^2_{\tilde{\Gamma}}(f,h)  :=
        P^{(2)}(\tilde{\mathsf{D}}^2f,\tilde{\mathsf{D}}^2h)
\end{equation}
avec la d\'{e}riv\'{e}e covariante sym\'{e}tris\'{e}e $\tilde{\mathsf{D}}$
(voir le paragraphe
\ref{SubSecRep1}). Ensuite, pour le champ hamiltonien $X_f$ de
$f\in\CinfK{M}$ on regarde la d\'{e}riv\'{e}e de Lie de la connexion,
\[
    \big(L_{X_f}\tilde{\nabla}\big)_XY:=
      ~ [X_f,\tilde{\nabla}_XY]-\tilde{\nabla}_{[X_f,X]}Y
       -\tilde{\nabla}_X[X_f,Y]~~~\forall~X,Y\in\Ginf(M,TM)
\]
ce qui est un champ de tenseur dans $\Ginf(M,S^2T^*M\otimes TM)$, et
\beq
  \big(L_{X_f}\tilde{\Gamma}\big)(X,Y,Z):=
    \omega\big(\big(L_{X_f}\tilde{\nabla}\big)_XY,Z\big)
             ~~~\forall~X,Y,Z\in\Ginf(M,TM),
\end{equation}
ce qui est un champ de tenseur dans $\Ginf(M,S^3T^*M)$. On rappelle
l'op\'{e}rateur bidiff\'{e}rentiel suivant pour $f,h\in\CinfK{M}$ (voir \cite{Lic82},
\cite{DL88}):
\beq \label{EqDefS3Gamma}
 S^3_{\tilde{\Gamma}}(f,h)  :=  P^{(3)}\big(L_{X_f}\tilde{\Gamma},
                                    L_{X_h}\tilde{\Gamma}\big).
\end{equation}
Ensuite, on rappelle la multiplication diamant $\diamond$, voir
\cite{Bon98}, pour deux $2$-formes $\beta$ et $\beta'$:
pour garder la simplicit\'{e} on les donnera
dans une carte $\big(U,(x^1,\ldots,x^n)\big)$;
soit $P=\frac{1}{2}\sum_{i,j}P^{ij}
\frac{\partial}{\partial x^i}\wedge\frac{\partial}{\partial x^j}$ et
$\beta=\frac{1}{2}\sum_{k,l}\beta_{kl}dx^k\wedge dx^l$:
\[
     (\beta\diamond\beta')_{ij}:=\sum_{r,s}P^{rs}\beta_{ri}\beta'_{sj}.
\]
De plus, l'op\'{e}rateur bidiff\'{e}rentiel $B^3_{\tilde{\Gamma}}[\beta]$
--qui d\'{e}pend de fa\c{c}on lin\'{e}aire
d'une $2$-forme ferm\'{e}e $\beta$-- est donn\'{e} en coordonn\'{e}es
par la formule suivante:
\beas
   B^3_{\tilde{\Gamma}}[\beta](f,h) & := & \sum_{i,j,k,l,r,s}
                \frac{1}{4}\beta_{rs}\big(P^{ir}P^{ks}P^{jl}
                              +P^{jr}P^{ks}P^{il}+P^{ir}P^{ls}P^{jk}\\
          & & ~~~~~~~~~~~~~~~~~+P^{jr}P^{ls}P^{ik}\big)
          ~\tilde{\nabla}^2_{(\frac{\partial}{\partial x^i},
                   \frac{\partial}{\partial x^j})}f~
   \tilde{\nabla}^2_{(\frac{\partial}{\partial x^k},
                  \frac{\partial}{\partial x^l})}h
                                     \\
         &    & +\sum_{i,j,k,u,r,s}\frac{1}{6}
              \big(\tilde{\nabla}_{\frac{\partial}{\partial x^u}}\beta
                           \big)_{rs}\big(P^{ju}P^{ir}P^{ks}
                                         +P^{iu}P^{jr}P^{ks}\big) \\
         & & ~~~~~~~~~~~~~~~\big(
               \tilde{\nabla}^2_{(\frac{\partial}{\partial x^i},
                   \frac{\partial}{\partial x^j})}f~
              \tilde{\nabla}_{\frac{\partial}{\partial x^k}}h
              +\tilde{\nabla}^2_{(\frac{\partial}{\partial x^i},
                   \frac{\partial}{\partial x^j})}h~
               \tilde{\nabla}_{\frac{\partial}{\partial x^k}}f\big).
\eeas
%Finalement, on rappelle l'op\'{e}rateur cobord de Hochschild $\mathsf{b}$
%pour des op\'{e}rateurs diff\'{e}rentiels $A:\CinfK{M}\ra\CinfK{M}$:
%\beq\label{EqOpCobordHochschildb}
%  (\mathsf{b}A)(f,h):=f~(Ah)-A(fh)+(Af)h
%\end{equation}
Le th\'{e}or\`{e}me suivant est plus ou moins bien connu:
\bsat\label{TStarprodOrdreTrois}
 Soit $(M,\omega)$ une vari\'{e}t\'{e} symplectique et $*=\sum_{r=0}^\infty\nu^r\mathsf{C}_r$
 un star-produit sur
 $M$. Soit $[*]$ sa classe de Deligne et
 $\alpha=\sum_{r=0}^\infty\nu^r\alpha_r$ une
 s\'{e}rie formelle de $2$-formes ferm\'{e}es repr\'{e}sentant
 $[*]-\frac{[\omega]}{\nu}$. Soit $\tilde{\nabla}$ une connexion sans
 torsion symplectique. Alors il existe des op\'{e}rateurs diff\'{e}rentiels
 $S_1,S_2,S_3$ dans $\Dop^1\big(\CinfK{M},\CinfK{M}\big)$ s'annulant sur les
 constantes tels que les
 premiers quatre termes $\mathsf{C}_0,\mathsf{C}_1,\mathsf{C}_2$ et
 $\mathsf{C}_3$ de $*$ soient de la forme
 suivante
 (pour tous $f,h\in\CinfK{M}$):
 \beas
  \mathsf{C}_0(f,h) & = & fh \label{EqCNull}  \\
  \mathsf{C}_1(f,h) & = & \{f,h\}+(\mathsf{b}S_1)(f,h) \label{EqCUn} \\
  \mathsf{C}_2(f,h) & = & \frac{1}{2}P^2_{\tilde{\Gamma}}(f,h)
                    -S_1\big((\mathsf{b}S_1)(f,h)\big)+(S_1f)(S_1h)
                    +(\mathsf{b}S_2)(f,h)\nonumber \\
           &   & ~~~-\alpha_0(X_f,X_h) +\{S_1f,h\}+\{f,S_1h\}-S_1\{f,h\}
                              \label{EqCDeux} \\
  \mathsf{C}_3(f,h) & = &
      \frac{1}{6}S^3_{\tilde{\Gamma}}(f,h)-\alpha_1(X_f,X_h)
                                   \nonumber \\
           &   & +\alpha_0\diamond\alpha_0
                 -\alpha_0(X_{S_1f},X_h)-\alpha_0(X_f,X_{S_1h})
                    +S_1\big(\alpha_0(X_f,X_h)\big) \nonumber \\
           &   & +\{S_2f,h\}+\{f,S_2h\}-S_2\{f,h\} \nonumber \\
           &   & -S_1\big(\{S_1f,h\}+\{f,S_1h\}-S_1\{f,h\}\big)
                   +\{S_1f,S_1h\} \nonumber \\
           &   & -B^3_{\tilde{\Gamma}}[\alpha_0](f,h)
                 +P^2_{\tilde{\Gamma}}(S_1f,h)+P^2_{\tilde{\Gamma}}(f,S_1h)
                 -S_1\big(P^2_{\tilde{\Gamma}}(f,h)\big)\nonumber \\
           &   &  +(\mathsf{b}S_3)(f,h)
                     +\big(S_1^2-S_2\big)\big((\mathsf{b}S_1)(f,h)\big)
                     -S_1\big((\mathsf{b}S_2)(f,h)\big) \nonumber \\
           &   &  +(S_1f)(S_2h)+(S_2f)(S_1h)-S_1\big((S_1f)(S_1h)\big)
                      \label{EqETrois}
 \eeas
 o\`{u} $\mathsf{b}$ est l'op\'{e}rateur cobord de Hochschild, voir \'{e}qn
 (\ref{EqCobordHochschild}).
\esat
\bbew
 Puisque tout star-produit est \'{e}quivalent \`{a} un star-produit construit
 par la m\'{e}thode de Fedosov $*_F$ (voir \cite{Fed96} et le paragraphe
 \ref{SecMorphRepAtiyahMolinoNulle}), alors il existe une
 transformation d'\'{e}quivalence $S=\mathrm{id}+\sum_{r=1}^\infty\nu^rS_r$
 telle que $*=S^{-1}(*_F)$. En utilisant le r\'{e}sultat de Neumaier \cite{Neu02}
 reliant
 la classe de Deligne avec la s\'{e}rie de $2$-formes ferm\'{e}es entrant dans
 la construction de Fedosov et le r\'{e}sultat de Bonneau \cite{Bon98} qui
 a montr\'{e} comment ces $2$-formes de Fedosov apparaissent dans les
 op\'{e}rateurs bidiff\'{e}rentiels du star-produit, on obtient $*_F$ qui est
 donn\'{e} par la formule ci-dessus avec $S_1=0=S_2=S_3$, voir \'{e}galement
 \cite{DL88} pour le cas $\alpha_0=0$. Le reste est un
 calcul \'{e}vident.
\ebew

\noindent Soit maintenant $i:C\ra M$ une sous-vari\'{e}t\'{e} co\"{\i}sotrope
ferm\'{e}e de $M$ de codimension $k$ et soit $E$ son fibr\'{e} caract\'{e}ristique.
Soit $\tau:\check{M}\ra C$ un voisinage tubulaire de $C$ dans $M$
que l'on va identifier de mani\`{e}re symplectique avec un voisinage ouvert
--\'{e}galement not\'{e} $\check{M}$-- de la section nulle de $E^*$ gr\^{a}ce au
th\'{e}or\`{e}me de Weinstein/Gotay \ref{TWeinGotay}. Soit $*$ un star-produit
sur $M$. Pour \'{e}tudier les
r\'{e}pr\'{e}sentations de $*$ sur $C$ il suffit de regarder les repr\'{e}sentations
de $*|_{\check{M}}$ sur $C$, ou --en vue de la proposition
 \ref{PKerRho1congId}-- de regarder les star-produits adapt\'{e}s dans
 $\check{M}$. Soit $(F,\nabla)$ une connexion sans torsion adapt\'{e}e dans le
 fibr\'{e} tangent de la vari\'{e}t\'{e} pr\'{e}symplectique $(C,\varpi=i^*\omega)$.
 On va utiliser la connexion symplectique $\tilde{\nabla}$ (voir la th\'{e}or\`{e}me
 \ref{TConnSympAutourCoisotrope}) dans
 $(\check{M},\omega_{WG})$ pour calculer explicitement $S_1$ et $S_2$
 du th\'{e}or\`{e}me pr\'{e}cedent tels que $*$ soit adapt\'{e} \`{a} $C$ jusqu'\`{a} l'ordre $2$.
 Soit $e_1,\ldots,e_k$ une base locale de $E$ et $e^1,\ldots,e^k$ la base
 duale dans $E^*$. Soit \'{e}galement $f_1,\ldots,f_{2m}$ une base locale
 du fibr\'{e} $F$ et $f^1,\ldots,f^{2m}$ la base duale dans $F^*$. Soit
 $R$ le tenseur de courbure de $\nabla$. On note pour
 $V_1,V_2\in\Ginf(C,E)$
 \[
    Ric_v(V_1,V_2):=\sum_{i=1}^k\langle e^i,R(e_i,V_2)V_1\rangle
 \]
 le {\em tenseur de Ricci vertical} de $\nabla$. De plus, il est connu que
 la trace du tenseur de courbure est une $2$-forme exacte (voir par
 exemple \cite{BNW98}, Lemma 16), c.-\`{a}-d. il existe une $1$-forme $\xi$
 sur $C$ telle que pour tous $X,Y,Z\in\Ginf(C,TC)$

 \[
    \mathrm{trace}\big(Z\mapsto R(X,Y)Z\big)=:d\xi(X,Y).
 \]
 Les op\'{e}rateurs diff\'{e}rentiels suivants ne d\'{e}pendent pas des bases choisies:
 soit $f\in\CinfK{\check{M}}$ et $\zeta\in\Ginf(C,T^*C\otimes T^*C)$:
 \bea
   \Delta (f)  & := &  \sum_{i=1}^k ~(e^i)^ve_i^h (f), \\
    \zeta^{E^* E^*}(f)  & := & \sum_{i,j=1}^k
    \zeta(e_i,e_j)~(e^i)^v(e^j)^v(f),
 \eea
 \bea
   \zeta^{E^* E}(f) & := & \sum_{i,j=1}^k\zeta\big(e_i^h,(e^j)^v\big)
                                      ~(e^i)^v e_j^h (f), \\
   \big(\zeta^{E^* F} (f)\big)_\varphi & := & -\sum_{i=1}^k~
      (e^i)^v~
      \big(\hat{P}_\varphi^\#\big(\zeta(\underline{~~},e_i^h)\big)\big)^h
               (f)_\varphi.
 \eea
\bsat \label{TObsRepOrdreDeux}
  Avec les notations introduites ci-dessus:
  \ben
   \item Pour que $*$ soit adapt\'{e} \`{a} $C$ jusqu'\`{a} l'ordre $1$ il faut et il
     suffit
     que l'op\'{e}rateur diff\'{e}rentiel $S_1$ soit de la forme suivante:
     \beq\label{EqDefSUnAdapte}
      S_1 = \Delta + \gamma_1^v +\tilde{S_1}
     \end{equation}
     o\`{u}
     $\gamma_1$ est une section de classe $\Cinf$ arbitraire de $E^*$ et
     $\tilde{S}_1$ est un op\'{e}rateur diff\'{e}rentiel arbitraire de
     $\CinfK{\check{M}}$ dans $\CinfK{\check{M}}$ qui s'annule sur les
     constantes et soit adapt\'{e} \`{a} $C$.
   \item Pour que $*$ soit adapt\'{e} \`{a} $C$ jusqu'\`{a} l'ordre $2$ il faut et il
   suffit que
    \ben
      \item $S_1$ soit de la forme (\ref{EqDefSUnAdapte}),
      \item il existe un repr\'{e}sentant $\alpha_0$ de $[*]_0$ tel que
        \[
          p_v i^*\alpha_0+d_v\gamma_1+d_vp_v\xi=0.
        \]
      \item l'op\'{e}rateur diff\'{e}rentiel $S_2$ soit de la forme g\'{e}n\'{e}rale
       \bea
         S_2 & = &
         \frac{1}{2}(\Delta+\gamma_1^v)^2+(\Delta+\gamma_1^v)\tilde{S}_1
           +\frac{1}{2}Ric_v^{E^* E^*}-\frac{1}{2}(\nabla\gamma_1)^{E^* E^*}
           \nonumber
           \\
             &   & -\alpha_0^{E^*E}-
             (\alpha_0+d\tau^*\gamma_1+d\tau^*\xi)^{E^*F}+\gamma_2^v
             +\tilde{S}_2
       \eea
       o\`{u} $\gamma_2\in\Ginf(C,E^*)$ et $\tilde{S}_2$ est un op\'{e}rateur
       diff\'{e}rentiel arbitraire de
     $\CinfK{\check{M}}$ dans $\CinfK{\check{M}}$ qui s'annule sur les
     constantes et soit adapt\'{e} \`{a} $C$.
    \een
  \een
\esat
\bbew
 Ce th\'{e}or\`{e}me est une v\'{e}rification directe, mais tr\`{e}s longue et
 fastidieuse. Une identit\'{e} tr\`{e}s utile est
 \[
    Ric_v(V,H)= (d\xi)(V,H)~~\mathrm{et}~~Ric_v(H,V)=0
    ~~\forall~V\in\Ginf(C,E),~H\in\Ginf(C,F),
 \]
 qui est une cons\'{e}quence de la premi\`{e}re identit\'{e} de Bianchi.
\ebew

\noindent Pour un champ vertical $V$ sur $C$ on d\'{e}finit la fonction
\[
    f_V:\check{M}\ra\korps:\varphi\mapsto \langle \varphi ,V_{\tau(\varphi)}
                            \rangle.
\]
Il est clair que tout $f_V$ est un \'{e}l\'{e}ment de l'id\'{e}al annulateur
$\mathcal{I}$. D'apr\`{e}s le th\'{e}or\`{e}me \ref{TObstructionsRecurrentes},
les obstructions pour la repr\'{e}sentabilit\'{e} d'un star-produit \`{a} l'ordre
$3$ s'obtiennent en calculant $i^*\mathsf{C}_3(f_V,f_W)-i^*\mathsf{C}_3(f_W,f_V)$
(pour tous $V,W\in\Ginf(C,E)$) qui
est un $2$-cocycle de la cohomologie longitudinale si le star-produit $*$ est
adapt\'{e} \`{a} $C$ jusqu'\`{a} l'ordre $2$. D'apr\`{e}s le lemme
\ref{LObstructionsRecursives},
\'{e}quation (\ref{EqLOb3}) il suffit de consid\'{e}rer ces fonctions $f_V$ qui
sont en bijection avec le module conormal
$\mathcal{I}/\mathcal{I}^2\cong \Ginf(C,E)$. Au bout d'un tr\`{e}s long
calcul direct on arrive au r\'{e}sultat (\ref{EqObstructionImportante}).
Le terme avec la classe d'Atiyah-Molino provient du terme
$S^3_{\tilde{\Gamma}}$. Ceci donne la d\'{e}monstration du th\'{e}or\`{e}me
\ref{TObstructionsOrdreTrois}.

\section{Existence des morphismes et des re\-pr\'{e}\-sen\-ta\-tions lorsque la
       classe d'Atiyah-Molino s'annule}
       \label{SecMorphRepAtiyahMolinoNulle}

\subsection{Quantification des morphismes de Poisson entre deux vari\'{e}t\'{e}s
            symplectiques}
  \label{SubSecQuantifMorphVarSymp}

Soient $(M,\omega)$ et $(M',\omega')$ deux vari\'{e}t\'{e}s symplectiques
et soit $\phi:M\ra M'$ un morphisme de Poisson. D'apr\`{e}s la proposition
\ref{PMorphPoissonSymp} l'application $\phi$ est une submersion dont
l'image $\check{M}'$ est un ouvert de $M'$. On rappelle les sous-fibr\'{e}s
int\'{e}grables $F:=\mathrm{Ker}~T\phi$ et $E:=F^\omega$ de $TM$. Soient
$\mathsf{P}_F$ et $\mathsf{P}_E$ les projections sur $F$ et $E$. On
d\'{e}finit
$\omega^E:=\omega(\mathsf{P}_E\underline{~},\mathsf{P}_E\underline{~})$
et $\omega^F:=\omega(\mathsf{P}_F\underline{~},\mathsf{P}_F\underline{~})$.
Alors, d'apr\`{e}s la proposition \ref{PMorphPoissonSymp} il vient
$\omega=\omega^E+\omega^F$ et $\phi^*\omega'=\omega^E$, et par cons\'{e}quent,
$d\omega^E=0=d\omega^F$. Le couple $(M,\omega^F)$ est donc pr\'{e}symplectique
\`{a} fibr\'{e} caract\'{e}ristique $E$. Soit $(F,\hat{\nabla})$ une connexion sans torsion
adapt\'{e}e \`{a} $(M,\omega^F)$: on a toutes les propri\'{e}t\'{e}s $i)$ -- $v)$ de la
proposition \ref{PExistConnAdaptees} pour $\hat{\nabla}$ avec
$R^F=0$ gr\^{a}ce \`{a} l'int\'{e}grabilit\'{e} du fibr\'{e} $F$. Soit $\nabla'$ une connexion
sans torsion dans le fibr\'{e} tangent de $M'$ qui soit symplectique.
Puisque $E$ est isomorphe au fibr\'{e} retir\'{e} $\phi^*TM'$ la connexion retir\'{e}e
$\phi^*\nabla'$ induit  une connexion $\nabla''$ dans le fibr\'{e} $E$ sur $M$
qui pr\'{e}serve $\omega^E=\phi^*\omega'$. On a la proposition suivante:
\bprop\label{PApplPoissonConnCompatibles}
 Avec les notations mentionn\'{e}es ci-dessus, pour toute con\-ne\-xion sans torsion
 symplectique $\nabla'$ dans le fibr\'{e} tangent de $M'$ il existe une connexion
 sans torsion
 symplectique $\tilde{\nabla}$ dans le fibr\'{e} tangent de $M$
 avec les propri\'{e}t\'{e}s
 suivantes pour tous $X\in\Ginf(C,TC)$, $V,W\in\Ginf(C,E)$ et
 $H,H_1,H_2,H_3\in\Ginf(C,F)$:
 \begin{itemize}
  \item[i)] $\tilde{\nabla}_XV\in\Ginf(C,E)$,
  \item[ii)] $\tilde{\nabla}_VH=\mathsf{P}_F[V,H]\in\Ginf(C,F)$,
  \item[iii)] $\tilde{\nabla}_HV=\mathsf{P}_E[H,V]\in\Ginf(C,E)$,
  \item[iv)] $\tilde{\nabla}_{X}H \in\Ginf(C,F)$,
  \item[v)] $\tilde{\nabla}\omega^E=0=\tilde{\nabla}\omega^F$,
  \item[vi)] $\phi$ est une application affine.
 \end{itemize}
 Pour toute connexion sans torsion $\tilde{\nabla}$ satisfaisant les
 conditions ci-dessus, les seules composantes en g\'{e}n\'{e}ral non nulles de son
 tenseur de courbure $\tilde{R}$ sont les suivantes:
 \bea
     R^E_{EE}(V_1,V_2;V,W) & := & \omega\big(V_1,\tilde{R}(V,W)V_2\big)
                                 \label{EqCourbureEEE} \\
     R^F_{FF}(H_1,H_2;H,H') & := & \omega\big(H_1,\tilde{R}(H,H')H_2\big)
                                 \label{EqCourbureFFF} \\
     R^F_{EF}(H_1,H_3;V,H_2) & := & \omega\big(H_1,\tilde{R}(V,H_2)H_3\big)
                                 =-\rho_{AME}(V;H_1,H_2,H_3) \nonumber \\
                              & &   \label{EqCourbureFEF}
 \eea
\eprop
\bbew
 On d\'{e}finit $\tilde{\nabla}$ par
 \[
   \tilde{\nabla}_XY := \hat{\nabla}_{\mathsf{P}_FX}(\mathsf{P}_FY)
                 +\mathsf{P}_F[\mathsf{P}_EX,\mathsf{P}_FY]
                 +\mathsf{P}_E[\mathsf{P}_FX,\mathsf{P}_EY]
                 +\nabla''_{\mathsf{P}_EX}\mathsf{P}_EY
 \]
 On v\'{e}rifie toutes les propri\'{e}t\'{e}s \'{e}nonc\'{e}es sur des champs de vecteurs
 de la forme $X=X^{\prime h}+H$ avec $X'\in \Ginf(M',TM')$ et
 $H\in\Ginf(M,F)$ o\`{u} $(~)^h$ d\'{e}signe le rel\`{e}vement horizontal de $X'$
 \`{a} $M$. On note que $(X,X')$ est $\phi$-li\'{e}. Puisque $E$ et $F$ sont
 int\'{e}grables, les tenseurs de courbure $\tilde{R}(V,W)H$
 et $\tilde{R}(H,H')V$ s'annulent (propri\'{e}t\'{e} de Bott). Il restent les
 termes purement dans $E$ ou dans $F$, et le seul terme mixte est
 donn\'{e} par la classe d'Atiyah-Molino du fibr\'{e} $E$ qui est proportionnel
 \`{a} $R^F_{EF}$. Puisque $F$ est le
 fibr\'{e} noyau d'une submersion,
alors son $1$-cocycle d'Atiyah-Molino s'annule.
\ebew

On rappelle maintenant les structures n\'{e}cessaires pour la {\em construction
de Fedosov partielle}: soit $N$ une vari\'{e}t\'{e} diff\'{e}rentiable avec un sous-fibr\'{e}
int\'{e}grable $K$ de son fibr\'{e} tangent muni d'une forme
symplectique $\omega^K$ sur les fibres. $\omega^K$ est une $2$ forme ferm\'{e}e
longitudinale. Soit $(A_u)_{u\in\nat}$ une famille de
fibr\'{e} vectoriels sur $M$. On suppose qu'il existe une famille
$(\mu_{r,uvw})_{r,u,v,w\in\nat}$ d'homomorphismes
dans $\Ginf\big(N,\mathrm{Hom}(A_u\otimes A_v,A_w)\big)$ de sorte que pour
chaque
$r\in\nat$ il existe $N_r\in\entier$ tel que $\mu_{r,uvw}=0$ lorsque
$w-u-v<N_r$ et que $\sum_{r,u,v,w=0}^\infty\nu^r\mu_{r,uvw}$ d\'{e}finisse
une multiplication associative `fibre-par-fibre' sur le $\korps[[\nu]]$-
module
\[
   \mathcal{A}:=\times_{u=0}^\infty \Ginf(A_u)[[\nu]].
\]
Soit $P\in\Ginf(N,\Lambda^2K)$ le champs de
bivecteurs correspondant \`{a} $\omega$. Il vient que $P$ est une structure
de Poisson r\'{e}guli\`{e}re sur $M$. On note
\[
    \mathcal{W}(K^*):=\Ginf(N,\overline{SK^*})[[\nu]]
         :=\big(\times_{k=0}^\infty \Ginf(N,S^k K^*)\big)[[\nu]]
\]
et pour $l\in\nat$
\[
    \mathcal{W}\otimes \Lambda^l(K^*)
    :=\big(\times_{k=0}^\infty \Ginf(N,S^k K^*\otimes \Lambda^lK^*)
             \big)[[\nu]]
\]
avec $\mathcal{W}\otimes \Lambda(K^*)
:=\oplus_{l\in\nat}\mathcal{W}\otimes \Lambda^l(K^*)$. Soit
$\mathcal{A}:=\Ginf(N,A)[[\nu]]$ et
\[
  \mathcal{A}\otimes \mathcal{W}\otimes \Lambda(K^*)
   :=\big(\times_{k,u=0}^\infty \Ginf(N,A_u\otimes S^k K^*\otimes \Lambda^lK^*)
             \big)[[\nu]].
\]
 On d\'{e}finit la multiplication fibre par fibre
non d\'{e}form\'{e}e entre sections factoris\'{e}es $\alpha_1\otimes f_1\otimes
\gamma_1$ et $\alpha_2\otimes f_2\otimes\gamma_2$ (avec
$\alpha_1,\alpha_2\in \Ginf(N,A)$, $f_1,f_2\in\mathcal{W}(K^*)$ et
$\gamma_1,\gamma_2\in\Lambda(K^*)$ par
\[
    (\alpha_1\otimes f_1\otimes\gamma_1)
           (\alpha_2\otimes f_2\otimes\gamma_2)
              :=\alpha_1\alpha_2\otimes f_1f_2\otimes
                      \gamma_1\wedge\gamma_2.
\]
Pour $X\in\Ginf(N,K)$ soit $i_s(X)(\alpha\otimes f\otimes \gamma):=
\alpha\otimes \big(i(X)f\big) \otimes \gamma$ et
$i_a(X)(\alpha\otimes f\otimes \gamma):=
\big(\alpha\big)\otimes f \otimes i(X)\gamma$.
Soit $e_1,\ldots,e_{2m}$ une base locale de $K$ et $e^1,\ldots,e^{2m}$
la base locale de $K^*$. Le champs de bivecteurs $P$ s'\'{e}crit
$P=\sum_{i,j=1}^{2m}P^{ij}e_i\otimes e_j$ pour certaines fonctions
locales $P^{ij}$
On d\'{e}finit la multiplication
d\'{e}form\'{e}e fibre par fibre $\circ_K$:
\bea
   \lefteqn{\xi\circ_K\eta:=\sum_{r=0}^\infty\frac{\nu^r}{r!}
             \sum_{i_1,\ldots,i_r,j_1,\ldots,j_r=1}^{2m}
             P^{i_1j_1}\cdots P^{i_rj_r}} \nonumber \\
           & &~~~~~~~~~~~~~~~~~~~~~~~~~~~~~~~
             \big(i_s(e_{i_1})\cdots i_s(e_{i_r})\xi\big)
             \big(i_s(e_{j_1})\cdots i_s(e_{j_r})\eta\big)
               \label{EqMultiplicationFibreparfibre}
\eea
qui ne d\'{e}pend \'{e}videmment pas de la base choisie. On va \'{e}crire
\beq \label{EqStructurePoissonFibreparfibre}
 \{\xi,\eta\}_K:=P^{ij}\big(i_s(e_i)\xi\big)\big(i_s(e_j)\eta\big)
\end{equation}
pour la structure de Poisson fibre-par-fibre d\'{e}finie sur
$\mathcal{A}\otimes \mathcal{W}\otimes\Lambda(K^*)$.
Muni de la multiplication $\circ_K$,
l'espace $\mathcal{A}\otimes \mathcal{W}\otimes\Lambda(K^*)$ est une alg\`{e}bre
associative. On note
$deg_{sK}\xi:=\sum_{i=1 }^{2m}\big(1\otimes e^i\otimes 1\big)
\big(i_s(e_i)\xi\big)$
l'endomorphisme
$\korps[[\nu]]$-lin\'{e}aire de $\mathcal{A}\otimes\mathcal{W}\otimes\Lambda(K^*)$
dont les vecteurs propres \`{a} valeur propre $k$
sont des sections
de $\Ginf(N,A\otimes S^kK^*\otimes \Lambda K^*)[[\nu]]$ quel que soit
l'entier positif $k$. On note \'{e}galement
$deg_{aK}\xi:=\sum_{i=1 }^{2m}\big(1\otimes 1\otimes e^i\big)
\big(i_a(e_i)\xi\big)$ l'endomorphisme
$\korps[[\nu]]$-lin\'{e}aire de $\mathcal{A}\otimes\mathcal{W}\otimes\Lambda(K^*)$
dont les vecteurs propres \`{a} valeur propre $l$
sont des sections
de $\Ginf(N,A\otimes \overline{SK^*}\otimes \Lambda^l K^*)[[\nu]]$ quel que
soit l'entier positif $l$. Soit $Deg_K$ l'endomorphisme $\korps$-lin\'{e}aire
de $\mathcal{A}\otimes\mathcal{W}\otimes\Lambda(K^*)$ donn\'{e} par
$Deg_K:=2\nu\frac{\partial}{\partial \nu}+deg_{sK}$. Il vient que $Deg_K$ et
$deg_{aK}$ sont des d\'{e}rivations de l'alg\`{e}bre
$\big(\mathcal{A}\otimes\mathcal{W}\otimes\Lambda(K^*),\circ_K\big)$ qui est
donc
une alg\`{e}bre associative $\entier$-gradu\'{e}e par rapport \`{a} la d\'{e}rivation
$deg_{aK}$. Un endomorphisme $\korps$-lin\'{e}aire $\psi$ de
$\mathcal{A}\otimes\mathcal{W}\otimes\Lambda(K^*)$ est dit {\em de degr\'{e}
antisym\'{e}trique $l$} lorsque $[deg_{aK},\psi]=l\psi$ et {\em de degr\'{e} total
$r$} lorsque $[Deg_{K},\psi]=r\psi$. On introduit le commutateur gradu\'{e}
$[f,h]:=f\circ_K h-(-1)^{l_1l_2}h\circ_K f:=ad_{\circ_K}(f)(h)$ quel que
soient les \'{e}l\'{e}ments $f,g\in \mathcal{A}\otimes\mathcal{W}\otimes\Lambda(K^*)$
avec $deg_{aK}f=l_1f$ et $deg_{aK}h=l_2h$. Il vient que
$ad_{\circ_K}(f)$ est une d\'{e}rivation gradu\'{e}e de degr\'{e} $l_1$ de
$\mathcal{A}\otimes\mathcal{W}\otimes\Lambda(K^*)$.
En outre, tout \'{e}l\'{e}ment
$\xi$ de
$\mathcal{A}\otimes\mathcal{W}\otimes\Lambda(K^*)$ s'\'{e}crit comme s\'{e}rie formelle
$\xi=\sum_{r=0}^\infty \xi^{(r)}$ avec $Deg_K(\xi^{(r)})=r\xi^{(r)}$.
Par cons\'{e}quent, $\mathcal{A}\otimes\mathcal{W}\otimes\Lambda(K^*)$ n'est pas
gradu\'{e} par rapport \`{a} $Deg_K$ (parce que
$\mathcal{A}\otimes\mathcal{W}\otimes\Lambda(K^*)$ n'est pas somme directe des
espaces propres de $Deg_K$), mais une alg\`{e}bre associative filtr\'{e}e
par
\[
  \mathcal{A}\otimes\mathcal{W}_r\otimes\Lambda(K^*)
  := \{\xi\in \mathcal{A}\otimes\mathcal{W}\otimes\Lambda(K^*)~|~
        \xi^{(0)}=0=\cdots=\xi^{(r-1)}\}
\]
o\`{u} $r\in\nat$ et $\mathcal{A}\otimes\mathcal{W}_0\otimes\Lambda(K^*)
:=\mathcal{A}\otimes\mathcal{W}\otimes\Lambda(K^*)$. Evidemment,
$\mathcal{A}\otimes\mathcal{W}_r\otimes\Lambda(K^*)\supset
\mathcal{A}\otimes\mathcal{W}_{r+1}\otimes\Lambda(K^*)$ et l'intersection
de tous les $\mathcal{A}\otimes\mathcal{W}_r\otimes\Lambda(K^*)$ est \'{e}gale
\`{a} $\{0\}$. Soient $\delta_K$ et $\delta_K^*$ les endomorphismes
$\korps[[\nu]]$-lin\'{e}aires de
$\mathcal{A}\otimes\mathcal{W}\otimes\Lambda(K^*)$ d\'{e}finis par
$\delta_K\xi:=\sum_{i=1}^{2m}\big(1\otimes 1\otimes e^i\big)
\big(i_s(e_i)\xi\big)$ et
$\delta_K^*\xi:=\sum_{i=1}^{2m}\big(1\otimes e^i\otimes 1\big)
\big(i_a(e_i)\xi\big)$. Il vient que $\delta_K$ est une d\'{e}rivation gradu\'{e}e
de degr\'{e} antisym\'{e}trique $1$ et de degr\'{e} total $-1$ de
$\mathcal{A}\otimes\mathcal{W}\otimes\Lambda(K^*)$ et que
$\delta_K^*$ est un endomorphisme $\korps[[\nu]]$ lin\'{e}aire de
$\mathcal{A}\otimes\mathcal{W}\otimes\Lambda(K^*)$ qui est de degr\'{e}
antisym\'{e}trique $-1$ et de degr\'{e} total $+1$. $\delta_K^*$ n'est pas en g\'{e}n\'{e}ral
une d\'{e}rivation gradu\'{e}e de
$\mathcal{A}\otimes\mathcal{W}\otimes\Lambda(K^*)$. De plus, on a
\beas
   \delta_K^2 & = & 0 \\
   {\delta_K^*}^2 & = & 0 \\
   \delta_K\delta_K^*+\delta_K^*\delta_K & := & deg_{sK}+ deg_{aK}
\eeas
donc $\big(\mathcal{A}\otimes\mathcal{W}\otimes\Lambda(K^*),
\circ_K,\delta\big)$
est une alg\`{e}bre diff\'{e}rentielle gradu\'{e}e dont la cohomologie est concentr\'{e}
en degr\'{e} antisym\'{e}trique nul et est isomorphe \`{a} l'espace
$\mathcal{A}$. Soit $\sigma_K$ la projection de
$\mathcal{A}\otimes\mathcal{W}\otimes\Lambda(K^*)$ sur $\mathcal{A}$
dont le noyau est \'{e}gal \`{a}
$\mathcal{A}\otimes\mathcal{W}\otimes\Lambda(K^*)^+$,
le $\korps[[\nu]]$-sous-module de
$\mathcal{A}\otimes\mathcal{W}\otimes\Lambda(K^*)$ qui consiste en
toutes les s\'{e}ries formelles  termes dont la somme de degr\'{e} sym\'{e}trique et
du degr\'{e} antisym\'{e}trique est strictement positive. Alors
$deg_{sK}+ deg_{aK}$ est inversible sur
$\mathcal{A}\otimes\mathcal{W}\otimes\Lambda(K^*)^+$ et on note
$(deg_{sK}+ deg_{aK})^{-1}$ son application r\'{e}ciproque. On d\'{e}finit
\[
    \delta_K^{-1}\xi:=\left\{\begin{array}{cl}
                              0 & \mathrm{si}~~\sigma_K\xi=0 \\
                  (deg_{sK}+ deg_{aK})^{-1}\delta_K^* &
                               \mathrm{si}~~\xi\in
                  \mathcal{A}\otimes\mathcal{W}\otimes\Lambda(K^*)^+
                             \end{array}\right.
\]
et il vient
\[
   \delta_K\delta_K^{-1}+\delta_K^{-1}\delta_K=id-\sigma_K.
\]
On rappelle qu'une {\em connexion partielle le long de $K$}
dans un fibr\'{e} vectoriel $L$ sur $M$ d'apr\`{e}s Ehresmann est une application
$\korps$-bilin\'{e}aire
$\nabla^K:\Ginf(N,K)\times \Ginf(N,L)\ra \Ginf(N,L):
(X,\psi)\mapsto \nabla^K_X\psi$ telle que
$\nabla^K_{fX}\psi=f\nabla^K_X\psi$ et $\nabla^K_X(f\psi)=(Xf)\psi
+f\nabla^K_X\psi$ quel que soit $f\in\CinfK{N}$. Pour le cas $L=K$,
la connexion partielle $\nabla^K$ est dite {\em sans torsion} lorsque
$\nabla^K_XY-\nabla^K_YX=[X,Y]$ quels que soient $X,Y\in\Ginf(N,K)$.
Soit $\nabla^{KK}$ une connexion partielle sans torsion le long de $K$ dans
$K$ qui soit symplectique dans le sens que la forme symplectique $\omega^K$
sur $K$ soit pr\'{e}serv\'{e}e. Une telle connexion partielle existe toujours
gr\^{a}ce \`{a} l'astuce de
Tondeur, Lichnerowicz et Hess, voir par exemple (\ref{EqAstuceTondeur}).
La connexion partielle $\nabla^{KK}$ induit une connexion partielle le long
de $K$ dans tous les fibr\'{e}s $S^kK^*\otimes\Lambda K^*$. Pour chaque entier
positif $u$, soit $\nabla^{A_uK}$
une connexion partielle le long de $K$ dans le fibr\'{e} $A_u$. La s\'{e}rie
formelle des $\nabla^{A_uK}$ d\'{e}finit une connexion partielle $\nabla^{AK}$
le long de $K$ dans $\mathcal{A}$ (o\`{u} les axiomes d'Ehresmann se g\'{e}n\'{e}ralisent
de fa\c{c}on \'{e}vidente au cas $\korps$ est remplac\'{e} par $\korps[[\nu]]$).
On suppose que $\nabla^{AK}_X$ soit une d\'{e}rivation de l'alg\`{e}bre
associative $\mathcal{A}$ quel que soit le champ de vecteurs
$X\in\Ginf(N,K)$.
 Les deux
connexions partielles $\nabla^{AK}$ et $\nabla^{KK}$ induisent une
connexion partielle, appel\'{e}e $\nabla^K$, dans tous les fibr\'{e}s
$A_u\otimes \oplus_{l=0}^kS^lK^*\otimes \Lambda K^*$. Ainsi la d\'{e}riv\'{e}e
covariante ext\'{e}rieure $\partial_K$ est bien d\'{e}finie dans
$\mathcal{A}\otimes\mathcal{W}\otimes\Lambda (K^*)$: \`{a} l'aide des
bases locales $e_1,\ldots,e_{2m}$ de $K$ et $e^1,\ldots,e^{2m}$ de $K^*$
mentionn\'{e}es ci-dessus on pose
\[
      \partial_KF:=\sum_{i=1}^{2m}e^i\wedge \nabla^{K}_{e_i}F
\]
quel que soit $F\in\mathcal{A}\otimes\mathcal{W}\otimes\Lambda (K^*)$.
La d\'{e}finition ne d\'{e}pend pas des bases choisies. Il vient que
$\partial_K$ est de degr\'{e} antisym\'{e}trique $1$ et de degr\'{e} total $0$.
Le th\'{e}or\`{e}me de Fedosov suivant permet de construire des d\'{e}formations
de l'alg\`{e}bre associative $\mathcal{A}$:

\bsat[Fedosov]\label{TFedosovConstructionPartielle}
  Avec les structures $N,K,A,\delta_K,\partial_K$ mentionn\'{e}s
  ci-dessus, on suppose qu'il existe un \'{e}l\'{e}ment
  $R^A_{KK}\in\mathcal{A}\otimes\Lambda^2 (K^*)$ de degr\'{e} total $0$
  tel que
  \[
     \partial_K^2=-\frac{1}{2\nu}ad_{\circ_K}(R^A_{KK}+R^K_{KK})
  \]
  (o\`{u} $R^K_{KK}\in\Ginf(S^2K^*\otimes \Lambda^2K^*)$ est
  li\'{e} au tenseur de courbure  $R^K$ de
  $\nabla^{KK}$par $R^K_{KK}(V,W;X,Y):=\omega^K(V,R^K(X,Y)W)$ quels que
  soient les champs de vecteurs
  $V,W,X,Y\in\Ginf(N,K)$) et tel que
  \[
     \partial_K(R^A_{KK}+R^K_{KK})=0.
  \]
  Soit
  $\Omega=\sum_{r=1}^\infty\nu^r \Omega_r\in \Ginf(N,\Lambda^2K^*)[[\nu]]$
  une s\'{e}rie formelle de
  $2$-formes ferm\'{e}es (le long de $K$) et
  $\mathsf{s}\in\mathcal{A}\otimes \mathcal{W}_3(K^*)$ avec
  $\sigma_K(\mathsf{s})=0$ donn\'{e}s.
  \ben
   \item Il existe un unique \'{e}l\'{e}ment $\mathsf{r}\in \mathcal{A}\otimes
   \mathcal{W}_2(K^*)$ tel que {\em l'\'{e}quation de courbure} suivante
   \beq \label{EqFedCourbureGenerale}
     0=-\delta_K\mathsf{r}+\partial_K\mathsf{r}
          -\frac{1}{2\nu}\mathsf{r}\circ_K\mathsf{r} + R^A_{KK}+R^K_{KK}
          +1\otimes \Omega
   \end{equation}
   et {\em l'\'{e}quation de normalisation}
   \beq \label{EqFedNormalisation}
      \delta_K^{-1}\mathsf{r}=\mathsf{s}.
   \end{equation}
   soient satisfaites.
   \item L'application
   \beq
     D_K:=-\delta_K+\partial_K-\frac{1}{2\nu}ad_{\circ_K}(\mathsf{r})
   \end{equation}
   est une d\'{e}rivation gradu\'{e}e de degr\'{e} antisymetrique $1$ de
   $\big(\mathcal{A}\otimes\mathcal{W}\otimes\Lambda (K^*),\circ_K\big)$,
   dite la {\em d\'{e}rivation de Fedosov}
   telle que
   \beq
    D_K^2=0.
   \end{equation}
   \item Les groupes de cohomologie de $D_K$ sont isomorphes \`{a}
   $\mathcal{A}$ (pour le degr\'{e} antisym\'{e}trique $0$) et \`{a} $\{0\}$
   (pour le degr\'{e} antisym\'{e}trique strictement positif).
   \item Il existe une unique injection $\korps[[\nu]]$-lin\'{e}aire
   $\tau_K:\mathcal{A}=\Ginf(N,A)[[\nu]]\ra
   \mathcal{A}\otimes\mathcal{W}(K^*)$, dite la {\em s\'{e}rie de Fedosov-Taylor}
    telle que
   \bea
     \sigma_K\big(\tau_K(a)\big) & = & a ~~~\forall a\in\mathcal{A}, \\
     \tau_K(\mathcal{A}) & = & Ker
     D_K|_{\mathcal{A}\otimes\mathcal{W}(K^*)}.
   \eea
   \item Il existe une unique d\'{e}formation associative formelle
   bidiff\'{e}rentielle $*$ de la multiplication point-par-point dans
   $\mathcal{A}$ telle que
   \beq
    \tau_K(a*b)=\tau_K(a)\circ_K\tau_K(b)~~~\mathrm{donc}~~~
                             a*b=\sigma_K\big(\tau_K(a)\circ_K\tau_K(b)\big)
   \end{equation}
   quels que soient $a,b\in\mathcal{A}$. Les op\'{e}rateurs bidiff\'{e}rentiels
   de $*$ sont `le long de $K$', c.-\`{a}-d. dans toute carte distingu\'{e}e
   de la vari\'{e}t\'{e} feuillet\'{e}e $N$ associ\'{e}e au fibr\'{e} int\'{e}grable $K$, ces
   op\'{e}rateurs ne contiennent que des d\'{e}riv\'{e}es partielles param\'{e}trant
   les feuilles.
  \een
\esat
Voir par exemple \cite{Fed96} pour une d\'{e}monstration. Pour calculer
la s\'{e}rie de Fedosov-Taylor, la formule de Neumaier suivante est tr\`{e}s utile
(voir \cite{Neu01}):
\beq
  \tau_K(a)=\frac{id}
   {id-[\delta_K^{-1},\partial_K-\frac{1}{2\nu}ad_{\circ_K}(\mathsf{r})]}(a)
\end{equation}
Pour le cas $K=TM$ et $\mathcal{A}=\CinfK{M}[[\nu]]$ (munie de la
multiplication point-par-point), cette construction
donne un star-produit $*$ sur la vari\'{e}t\'{e} symplectique $(M,\omega)$ dont la
classe de Deligne est donn\'{e}e par
\beq\label{EqClasseDeligneNeumaier}
    [*]=\frac{[\omega]}{\nu}+\frac{[\Omega]}{\nu},
\end{equation}
voir \cite{Neu02}. Dans ce
cas, on peut toujours modifier $\mathsf{r}$ en
$\hat{\mathsf{r}}:=\mathsf{r}-\sigma_{TM}(\mathsf{r})$ (si bien que
$\Omega$ sera remplac\'{e}e par $\Omega + d\sigma_{TM}(\mathsf{r})$ et
$s$ par $s-\delta_{TM}^{-1}\sigma_{TM}(\mathsf{r})$).

On va d'abord appliquer le th\'{e}or\`{e}me pr\'{e}cedent au cas o\`{u}
$(N,K,\mathcal{A})$ est donn\'{e} par $(M,TM,\CinfK{M}[[\nu]])$ et
$(N,K,\mathcal{A})=(M',TM',\CinfK{M'}[[\nu]])$: Gr\^{a}ce \`{a} la proposition
\ref{PApplPoissonConnCompatibles} les d\'{e}riv\'{e}es covariantes
$\tilde{\nabla}$ sur $M$ et $\nabla'$ sur $M'$ sont telles que
l'application de Poisson $\phi$ est une application affine. Alors
en \'{e}crivant $\delta_{TM}=:\delta$, $\partial_{TM}=:\partial$,
$\delta_{TM'}=:\delta'$ et $\partial_{TM'}=:\partial'$ les
deux \'{e}quations suivantes sont faciles \`{a} d\'{e}duire:
\bea
   \delta \phi^* & = & \phi^*\delta' \label{EqCompaDeltaDeltaPrime}\\
   \partial \phi^* & = & \phi^*\partial' \label{EqCompaPartialPartialPrime}
\eea
o\`{u} $\phi^*:\mathcal{W}\otimes\Lambda(T^*M')\ra\mathcal{W}\otimes\Lambda(T^*M)$
est le pull-back usuel. On voit \'{e}galement que $\phi^*$ est un morphisme
d'alg\`{e}bres associatives gradu\'{e}es par rapport aux multiplications
fibre-par-fibre
$\circ':=\circ_{TM'}$ et $\circ:=\circ_{TM}$. On a d'abord le th\'{e}or\`{e}me
suivant:
\bsat \label{TFedQuantApplPois}
 Soit $\phi:M\ra M'$ une application de Poisson entre deux vari\'{e}t\'{e}s
 symplectiques $(M,\omega)$ et $(M',\omega')$. Alors les deux \'{e}nonc\'{e}s
 suivants sont \'{e}quivalents:
 \ben
  \item Il existent des star-produits $*$ sur $M$ et $*'$ sur $M'$ tels
    que $\phi$ admette une quantification diff\'{e}rentielle
    en tant que morphisme d'alg\`{e}bres
    associatives.
  \item Il existent des d\'{e}rivations de Fedosov $D$ dans
    $\mathcal{W}\otimes\Lambda(T^*M)$ et $D'$ dans
    $\mathcal{W}\otimes\Lambda(T^*M')$ telles que
    \[
        D\phi^*=\phi^*D'.
    \]
 \een
 De plus, si un de ces \'{e}nonc\'{e}s est satisfait, on peut trouver des
 star-produits
 $*_F$ (\'{e}quivalent \`{a} $*$) et $*'_F$ (\'{e}quivalent \`{a} $*'$) de
 telle mani\`{e}re que le pull-back
 $\phi^*:(\CinfK{M'}[[\nu]],*'_F)\ra (\CinfK{M}[[\nu]],*_F)$
 est d\'{e}j\`{a} un morphisme d'alg\`{e}bres associatives.
\esat
\bbew
 ``$1) ~\Leftarrow~ 2)$'': on suppose que $D\phi^*=\phi^*D'$. Soient
 $\tau$ et $\tau'$ les s\'{e}ries de Fedosov-Taylor associ\'{e}es \`{a} $D$ et \`{a} $D'$,
 soient $\sigma:=\sigma_{TM}$ et $\sigma'=\sigma_{TM'}$ et soient
 $*$ et $*'$ les star-produits r\'{e}sultant de la construction de Fedosov.
 Soient $f',h'\in\CinfK{M'}[[\nu]]$.
 Il vient $D\phi^*\tau'(f')=\phi^*D'\tau'(f')=0$, alors
 \[
   \phi^*\tau'(f')=\tau\big(\sigma(\phi^*\tau'(f'))\big)
 =\tau\big(\phi^*\sigma'\tau'(f')\big)=\tau(\phi^*f'),
 \]
 donc
 \beas
   \phi^*(f'*'h')& = & \sigma\tau\big(\phi^*(f'*'h')\big)
                      = \sigma\phi^*\tau'(f'*'h')
                      =\sigma \phi^*\big(\tau'(f')\circ'\tau'(h')\big) \\
                 & = & \sigma\big((\phi^*\tau'(f'))
                              \circ (\phi^*\tau'(h'))\big)
                 = \sigma\big((\tau(\phi^*f'))\circ(\tau(\phi^*h'))\big)
                 \\
                 & = & \sigma \tau\big((\phi^*f')*(\phi^*h')\big)
                  =  (\phi^*f')*(\phi^*h')
 \eeas
 et $\phi^*$ est un morphisme d'alg\`{e}bres associatives, ce qui d\'{e}montre
 $1)$ et en m\^{e}me temps l'\'{e}nonc\'{e} apr\`{e}s $2)$.\\
 ``$1) ~\Rightarrow~ 2)$'': Soient $*$ et $*'$ des star-produits sur
 $M$ et $M'$ et soit
 $\Phi=\phi^*+\sum_{r=1}^\infty\nu^r\Phi_r$ une s\'{e}rie d'op\'{e}rateurs
 diff\'{e}rentiels dans $\Dop\big(\CinfK{M'},\CinfK{M}\big)$ le long de $\phi$
 tels que $\Phi:(\CinfK{M'}[[\nu]],*)\ra(\CinfK{M}[[\nu]],*')$ soit un
 homomorphisme d'alg\`{e}bres associatives. Puisque tout star-produit est
 \'{e}quivalent \`{a} un star-produit construit \`{a} l'aide de la m\'{e}thode de Fedosov
 \`{a} partir d'une connexion symplectique arbitraire (voir \cite{NT95a},
 \cite{NT95b}, \cite{Neu02}), il existe des transformations d'\'{e}quivalence
 $\hat{S}=id+\sum_{r=1}^\infty\nu^r\hat{S}_r$ et
 $S'=id+\sum_{r=1}^\infty \nu^r S'_r$
 (avec $\hat{S}_r\in \Dop^1\big(\CinfK{M};\CinfK{M}\big)$ et
 $S'_r\in\Dop^1\big(\CinfK{M'};\CinfK{M'}\big)$, s'annulant sur les constantes
  quel que soit
 l'entier strictement positif $r$) telles que $*_F:=\hat{S}(*)$
 et $*'_F:=S'(*')$ sont des star-produits de Fedosov sur $M$ et $M'$.
 Visiblement, $\hat{S}\Phi {S'}^{-1}$ est un homomorphisme d'alg\`{e}bres
 associatives $(\CinfK{M'}[[\nu]],*_F)\ra(\CinfK{M}[[\nu]],*'_F)$.
 On consid\`{e}re l'application $\Psi$ de l'espace $\CinfK{M'}[[\nu]]$
 dans $\mathcal{W}(T^*M)$ d\'{e}finie par
 \[
    f'\mapsto \tau\big(\hat{S}\Phi {S'}^{-1}(f')\big)~~
      \forall f'\in\CinfK{M'}[[\nu]]
 \]
 $\Psi$ est une s\'{e}rie d'op\'{e}rateurs diff\'{e}rentiels le long de $\phi$. Donc
 il existe un morphisme de fibr\'{e}s vectoriels $\hat{\Psi}$ du fibr\'{e}
 en jets
 $\phi^*J^\infty(M',M'\times\korps)$ dans $\times_{k=0}S^kT^*M$, alors
 $\Psi(f')=\hat{\Psi}\big(j^\infty(f')\big)$. Puisque les fibr\'{e}s
 $J^r(M',M'\times\korps)$ et $\oplus_{k=0}^rS^kT^*M'$ sont isomorphes
 (en utilisant la d\'{e}riv\'{e}e covariante $\nabla'$) quel que soit l'entier
 positif $r$, il existe une s\'{e}rie d'op\'{e}rateurs diff\'{e}rentiels $Q$ d'ordre
 $0$, telle que $j^\infty(f')=Q\big(\tau'(f')\big)$. Donc il existe
 une s\'{e}rie d'op\'{e}rateurs diff\'{e}rentiels d'ordre $0$, $\hat{\Phi}$, le long
 de $\phi$, de $\mathcal{W}(T^*M')$ dans $\mathcal{W}(T^*M)$ telle que
 \[
   \tau\big(\hat{S}\Phi {S'}^{-1}(f')\big)=\hat{\Phi}\big(\tau'(f')\big)
        ~~~~~~~~\forall f'\in\CinfK{M'}[[\nu]].
 \]
 Il vient
 \beas
  \hat{\Phi}\big(\tau'(f')\circ'\tau'(h')\big) & = &
          \hat{\Phi}\big(\tau'(f'*'_F h')\big) =
          \tau\big(\hat{S}\Phi {S'}^{-1}(f'*'_F h') \big)\\
          & = & \tau\big((\hat{S}\Phi {S'}^{-1}(f'))*_F
                     (\hat{S}\Phi {S'}^{-1}(h'))\big) \\
          & = &
              \tau\big(\hat{S}\Phi {S'}^{-1}(f')\big)\circ
               \tau\big(\hat{S}\Phi {S'}^{-1}(h')\big) \\
          & = &
              \hat{\Phi}\big(\tau'(f')\big)\circ
              \hat{\Phi}\big(\tau'(h')\big),
 \eeas
 et puisque $\hat{\Phi}$ est `tensoriel', c.-\`{a}-d. d'ordre $0$, et puisque
 pour tout $p'\in M'$ la fibre $\mathcal{W}(T^*M')_{p'}:=\times_{k=0}^\infty
 (S^kT^*M')_{p'}[[\nu]]$ est donn\'{e}e par
 $\{\tau'(f')_{p'}~|~f'\in\CinfK{M'}[[\nu]]\}$ gr\^{a}ce au lemme classique
 d'Emile Borel, il s'ensuit que $\hat{\Phi}$ est un homomorphisme
 d'alg\`{e}bres associatives tensoriel $(\mathcal{W}(T^*M'),\circ')$ dans
 $(\mathcal{W}(T^*M),\circ)$ qui d\'{e}forme $\phi^*$. On peut prolonger
 $\hat{\Phi}$ de mani\`{e}re canonique sur $\mathcal{W}\otimes \Lambda
 (T^*M')$ par $\hat{\Phi}(A'\otimes \gamma'):=\hat{\Phi}(A')\otimes
 \phi^*\gamma'$ avec $A'\in\mathcal{W}(T^*M')$ et $\gamma'\in\Ginf(\Lambda
 T^*M')$, et $\hat{\Phi}$ restera un homomorphisme d'alg\`{e}bres associatives
 gradu\'{e}es.\\
 On va montrer que $\hat{\Phi}$ est de la forme $A~\phi^*$ o\`{u} $A$ est
 un automorphisme int\'{e}rieur de
 $\big(\mathcal{W}\otimes\Lambda(T^*M),\circ\big)$: soit
 $\hat{\Phi}=\sum_{r=0}^\infty\nu^r\hat{\Phi}_r$. Chaque $\hat{\Phi}_r$
 est un op\'{e}rateur diff\'{e}rentiel d'ordre $0$, donc un \'{e}l\'{e}ment de
 $\times_{k',k\in\nat}\Ginf\big(\mathrm{Hom}(\phi^*S^{k'}T^*M',S^kT^*M)\big)$
 d'apr\`{e}s l'\'{e}q. (\ref{EqOpMultiGlobalUndSo}). Puisque le fibr\'{e} retir\'{e}
 $\phi^*T^*M'$ est isomorphe au fibr\'{e} $E^*$, tout $\hat{\Phi}_r$ est
 de la forme $\check{\Phi}_r~\phi^*$ avec
 $\check{\Phi}_r\in
 \times_{k',k\in\nat}\Ginf\big(\mathrm{Hom}(S^{k'}E^*,S^kT^*M)\big)$.
 Il s'ensuit
 que $\check{\Phi}_0$ est donn\'{e}e par l'injection canonique de
 $\mathcal{W}\otimes\Lambda(E^*)$ dans $\mathcal{W}\otimes\Lambda(T^*M)$.
 Soit $s$ le plus petit entier
 strictement positif tel que $\check{\Phi}_s\neq 0$. Puisque
 $\check{\Phi}_0$ et
 $\check{\Phi}=\sum_{r=0}^\infty\nu^r\check{\Phi}_r$ sont des
 homomorphismes
 d'alg\`{e}bres associatives de
 $\big(\mathcal{W}\otimes \Lambda(E^*),\circ_E\big)$ dans
 $\big(\mathcal{W}\otimes \Lambda(T^*M),\circ\big)$, la propri\'{e}t\'{e}
 d'homomorphisme \`{a} l'ordre $s$ montre que
 $\check{\Phi}_s$ est une d\'{e}rivation de l'alg\`{e}bre commutative gradu\'{e}e
 $\mathcal{W}\otimes \Lambda(E^*)$ dans $\mathcal{W}\otimes \Lambda(T^*M)$
 (munies de la multiplication fibre-par-fibre non d\'{e}form\'{e}e) le long de
 $\check{\Phi}_0$. Il existe donc un \'{e}l\'{e}ment
 $L=\sum_{i=1}^{2m}L^i\otimes e_i$ de
 $\times_{k=0}^\infty\Ginf\big(S^kT^*M\otimes E\big)$ avec
 $L^i\in \times_{k=0}^\infty\Ginf\big(S^kT^*M\big)$ tel que
 $\check{\Phi}_s(A)=i_s(L)\big(A\big)
 :=\sum_{i=1}^{2m}L^i~i_s(e_i)\big(A\big)$ quel que soit
 $A\in\mathcal{W}\otimes\Lambda(E^*)$. La propri\'{e}t\'{e}
 d'homomorphisme de $\check{\Phi}$ \`{a} l'ordre $s+1$ montre que $i_s(L)$
 est une d\'{e}rivation d'alg\`{e}bres de Poisson le long de $\check{\Phi}_0$.
 Un petit calcul montre que ceci est \'{e}quivalent \`{a} dire que
 $L^\flat:=(id\otimes \omega_E^\flat)(L)\in
 \times_{k=0}^\infty\Ginf\big(S^kT^*M\otimes E^*\big)$ est
 $\delta_E$-ferm\'{e}, i.e.~$\delta_E L^\flat=0$. Ici on a regard\'{e} $L^\flat$
 comme
 \'{e}l\'{e}ment de $\mathcal{A}\otimes\mathcal{W}\otimes \Lambda^1(E^*)$ avec
 $\mathcal{A}:=\big(\times_{k=0}^\infty\Ginf\big(S^k F^*\big),\circ_F\big)$.
 Alors il existe un \'{e}l\'{e}ment
 $\mathsf{C}_s\in\mathcal{W}(T^*M)$ tel que $L^\flat=\delta_E(\mathsf{C}_s)$
 et par
 cons\'{e}quent: $\check{\Phi}_s(A)=-\{\mathsf{C}_s,A\}$, avec la structure de
 Poisson fibre-par-fibre associ\'{e}e \`{a} $\omega$, voir \'{e}qn
 (\ref{EqStructurePoissonFibreparfibre}). Si l'on d\'{e}finit
 \[
    T^{[s]}:\mathcal{W}\otimes\Lambda(T^*M)\ra
             \mathcal{W}\otimes\Lambda(T^*M):
             A\mapsto e^{\nu^{s-1}ad_\circ(\mathsf{C}_s/2)}A
 \]
 alors $T^{[s]}$ est un homomorphisme d'alg\`{e}bres associatives et
 \[
     T^{[s]}~\hat{\Phi}=\phi^*~~\mathrm{mod}~~\nu^{s+1}.
 \]
 Par r\'{e}currence on montre l'existence d'un \'{e}l\'{e}ment
 $C=\nu^{s-1}\mathsf{C}_s/2+\cdots \in\mathcal{W}_1(T^*M)$ tel que
 \[
     \hat{\Phi} = e^{-ad_\circ(C)}~\phi^*.
 \]
 Il est \'{e}vident que $e^{-ad_\circ(C)}$ est un automorphisme de l'alg\`{e}bre
 $\big(\mathcal{W}\otimes\Lambda(T^*M),\circ\big)$. Soient
 $D'=-\delta'+\partial'-\frac{1}{2\nu}ad_{\circ'}(\mathsf{r}')$ et
 $D=-\delta+\partial-\frac{1}{2\nu}ad_{\circ}(\mathsf{r})$ les
 d\'{e}rivations de Fedosov \`{a} l'aide desquelles on a construit $*_F$ et
 $*'_F$. On peut supposer que les d\'{e}riv\'{e}es covariantes $\nabla'$ (dans $TM'$)
 et $\tilde{\nabla}$ (dans $TM$) sont choisies de telle fa\c{c}on que
 l'\'{e}quation (\ref{EqCompaPartialPartialPrime}) soit satisfaite.
 Soit
 \[
   \tilde{D}:=e^{ad_\circ(C)}~D~e^{-ad_\circ(C)}
             =:-\delta+\partial-\frac{1}{2\nu}ad_{\circ}(\tilde{\mathsf{r}}).
 \]
 Alors le membre de droite est bien d\'{e}fini (c.-\`{a}-d.
 $\tilde{\mathsf{r}}\in\mathcal{W}_3\otimes \Lambda^1(T^*M)$ car
 $C\in\mathcal{W}_1(T^*M)$) et d\'{e}fini une d\'{e}rivation gradu\'{e}e de degr\'{e}
 antisym\'{e}trique $+1$ de carr\'{e} nul. On peut montrer que $\tilde{D}$ est
 toujours une d\'{e}rivation de Fedosov d\'{e}finie par une s\'{e}rie de $2$-formes
 ferm\'{e}es cohomologue \`{a} celle qui d\'{e}finit $D$ et par une normalisation
 $\delta^{-1}\tilde{\mathsf{r}}$ diff\'{e}rente de celle de
 $\mathsf{r}$, et le star-produit $\tilde{*}$ qui en r\'{e}sulte est \'{e}quivalent
 \`{a} $*_F$ (alors il existe une transformation d'\'{e}quivalence $S$ telle que
 $\tilde{*}=S(*)$), voir par exemple \cite{Fed96}, \cite{Neu01}. Soit
 $f'\in\CinfK{M'}[[\nu]]$. Alors
 \beas
   \tilde{D}\phi^*\tau'(f')
   & = & e^{ad_\circ(C)}~D~e^{-ad_\circ(C)}~\phi^*\tau'(f')
   =e^{ad_\circ(C)}~D~\hat{\Phi}\tau'(f')\\
   & = & e^{ad_\circ(C)}~D~\tau\big(\hat{S}\Phi {S'}^{-1}(f')\big)
   =0=\phi^*~D'\tau'(f'),
 \eeas
 alors, gr\^{a}ce \`{a} l'\'{e}qn (\ref{EqCompaPartialPartialPrime})
 \[
  0=\tilde{D}\phi^*\tau'(f')-\phi^*~D'\tau'(f')
    =-\frac{1}{2\nu}\big(ad_\circ(\tilde{\mathsf{r}})\phi^*
                         -\phi^*ad_{\circ'}(\mathsf{r}')\big)\tau'(f').
 \]
 Puisque les s\'{e}ries de Fedosov-Taylor engendrent chaque fibre de
 $\mathcal{W}(T^*M')$ et puisque la partie grassmannienne $\Lambda(T^*M')$
 est dans le noyau de tout $ad_{\circ'}(A')$ il s'ensuit finalement que
 \[
    ad_\circ(\tilde{\mathsf{r}})\phi^*=
    \phi^*ad_{\circ'}(\mathsf{r}'),~~~\mathrm{donc}~~~
    \tilde{D}\phi^*=\phi^*D',
 \]
 ce qui montre le th\'{e}or\`{e}me.
\ebew

\noindent On va maintenant \'{e}tudier l'existence des d\'{e}rivations de Fedosov
qui sont
`$\phi$-li\'{e}es' dans le sens de l'\'{e}nonc\'{e} 2) du Th\'{e}or\`{e}me pr\'{e}c\'{e}dent.
On rappelle que $TM=E\oplus F$, et on a la d\'{e}composition en composantes
de $\delta$ et de $\partial$: avec $(e_1,\ldots,e_{2m}$ une base locale
des sections $\Cinf$ de $E$ et $e^1,\ldots,e^{2m}$ sa base duale et
$f_1,\ldots,f_{2n}$ une base locale des sections
$\Cinf$ et $f^1,\ldots,f^{2n}$ sa base duale, on d\'{e}finit
\[
   \delta_EA:=\sum_{i=1}^{2m}(1\otimes e^i) i_s(e_i)A~~~
    \delta_FA:=\sum_{i=1}^{2n}(1\otimes f^i) i_s(f_i)A~~~
    \mathrm{donc}~~~\delta=\delta_E+\delta_F
\]
(et $\delta_E^{-1}$, $\delta_F^{-1}$  de mani\`{e}re analogue) ainsi que
\[
   \partial_EA:=\sum_{i=1}^{2m}(1\otimes e^i) \tilde{\nabla}_{e_i}A~~~
    \partial_FA:=\sum_{i=1}^{2n}(1\otimes f^i) \tilde{\nabla}_{f_i}A~~~
    \mathrm{donc}~~~\partial=\partial_E+\partial_F.
\]
Il vient
\bea
     \partial_E^2 & = & -\frac{1}{2\nu}ad_{\circ}(R^E_{EE})
          \label{EqIdUtilesPartialCarreEE}\\
     ~[\partial_E,\partial_F] & = &
          -\frac{1}{2\nu}ad_{\circ}(R^F_{EF})
          \label{EqIdUtilesPartialCarreEF}\\
     \partial_F^2 & = & -\frac{1}{2\nu}ad_{\circ}(R^F_{FF})
          \label{EqIdUtilesPartialCarreFF}
\eea
avec les composantes de la courbure (\ref{EqCourbureEEE}),
(\ref{EqCourbureFEF}) et (\ref{EqCourbureFFF}).

On voit rapidement que $\delta_E$ et $\partial_E$ pr\'{e}servent
$\mathcal{W}\otimes\Lambda(E^*)$ et $\mathcal{W}(F^*)\otimes
\Lambda(T^*M)$ ainsi que $\delta_F$ et $\partial_F$ pr\'{e}servent
$\mathcal{W}\otimes\Lambda(F^*)$ et $\mathcal{W}(E^*)\otimes
\Lambda(T^*M)$. Nous allons donc regarder $\delta_E$ et $\partial_E$
dans le cas $(N,K,\mathcal{A})=\big(M,E,\CinfK{M}[[\nu]]\big)$
et $\delta_F$ et $\partial_F$ dans les cas
$(N,K,\mathcal{A})=\big(M,F,\CinfK{M}[[\nu]]\big)$ et
$(N,K,\mathcal{A})
=\big(M,F,(\Ginf\big(M,\Lambda(E^*)\big)[[\nu]],\wedge)\big)$.

\bprop \label{PConditionsPourFedLiees}
 Avec les notations pr\'{e}c\'{e}dentes: soient $D'=-\delta'+\partial'-
 \frac{1}{2\nu}ad_{\circ'}(\mathsf{r}')$ et $D=-\delta+\partial
 -\frac{1}{2\nu}ad_{\circ}(\mathsf{r})$ deux d\'{e}rivations de Fedosov
 telles que
 \[
      D\phi^* = \phi^*D'.
 \]
 Alors on a les \'{e}nonc\'{e}s suivants:
 \ben
  \item $\mathsf{r}=\mathsf{r}^E_E+\mathsf{r}^F_E+\mathsf{r}^F_F$ avec
    $\mathsf{r}^E_E\in\mathcal{W}_3\otimes\Lambda^1(E^*)$,
    $\mathsf{r}^F_F\in\mathcal{W}_3\otimes\Lambda^1(F^*)$ et
    $\mathsf{r}^F_E\in\mathcal{W}_3(F^*)\otimes\Lambda^1(E^*):=
    \times_{k=0}^\infty\Ginf(M,S^kF^*\otimes E^*)[[\nu]]$. De plus
    $\mathsf{r}^E_E=\phi^*\mathsf{r}'$.
  \item Il vient
    \[
         R^E_{EE}=\phi^*R'
    \]

  \item Soit $\Omega =\Omega_{EE} + \Omega_{EF} +\Omega_{FF}$ une s\'{e}rie
    formelle de $2$-formes ferm\'{e}es sur $M$ avec
    $\Omega_{EE}\in\nu\Ginf\big(M,\Lambda^2(E^*)\big)[[\nu]]$,
    $\Omega_{EF}\in\nu\Ginf\big(M,E^*\wedge F^*\big)[[\nu]]$ et
    $\Omega_{FF}\in\nu\Ginf\big(M,\Lambda^2(F^*)\big)[[\nu]]$.
    Soit $\Omega'\in\nu\Ginf\big(M',\Lambda^2(T^*M')\big)[[\nu]]$
    On a les quatre \'{e}quations suivantes pour les composantes de
    $\mathsf{r}$ et pour $\mathsf{r}'$:
    \bea
      0 & =& -\delta'\mathsf{r}'+\partial'\mathsf{r}'-\frac{1}{2\nu}
                \mathsf{r}'\circ'\mathsf{r}'+R'+\Omega'
                       \label{EqFedCourbureUnten}\\
      0 & = & \partial_E\mathsf{r}^F_E-\frac{1}{2\nu}\mathsf{r}^F_E
               \circ_F\mathsf{r}^F_E+\Omega_{EE}-\phi^*\Omega'
                       \label{EqFedKruemmungFENullCrucial}\\
      0 & = & \big(-\delta_F+\partial_F-\frac{1}{2\nu}ad_{\circ_F}\big)
              (\mathsf{r}^F_E) + R^F_{EF} +\partial_E\mathsf{r}^F_F
              +\Omega_{EF}
                       \label{EqFedrFEalsFunktionvonAtiyahMolino} \\
      0 & = & -\delta_F\mathsf{r}^F_F+\partial_F\mathsf{r}^F_F-\frac{1}{2\nu}
                \mathsf{r}^F_F\circ_F\mathsf{r}^F_F+R^F_{FF}+\Omega_{FF}
                       \label{EqFedCourbureObenlelongdeF}
    \eea
    avec les normalisations ${\delta'}^{-1}(\mathsf{r}')=\mathsf{s}'$,
    $\delta_E^{-1}(\mathsf{r}^F_E)=\mathsf{s}^{EF}$ et
    $\delta_F^{-1}(\mathsf{r}^F_F)=\mathsf{s}^F$ o\`{u}
    $\mathsf{s}'\in\mathcal{W}_3(T^*M')$,
    $\mathsf{s}^{EF}\in\times_{k=1}^\infty\Ginf(M,E^*\otimes S^kF^*)[[\nu]]$
    et $\mathsf{s}^F\in\mathcal{W}_3(F^*)$.
 \een
\eprop
\bbew
 1. L'\'{e}quation $D\phi^* = \phi^*D'$ implique que
 \[
  [\mathsf{r},\phi^*A']_\circ=\phi^*[\mathsf{r}',A']_{\circ'}
                             = [\phi^*\mathsf{r}',\phi^*A']_\circ
 \]
 ce qui appartient \`{a}
 $\mathcal{W}\otimes\Lambda(E^*)$ quel que soit
 $A'\in\mathcal{W}\otimes\Lambda(T^*M')$.
 Puisque  le fibr\'{e} $E$ est isomorphe au sous-fibr\'{e} retir\'{e} $\phi^*T^*M'$ et
 puisque l'alg\`{e}bre $\big(\mathcal{W}\otimes\Lambda(T*^M),\circ\big)$ est
 le produit tensoriel compl\'{e}t\'{e} (par rapport \`{a} la filtration mentionn\'{e}e
 ci-dessus) de ses sous-alg\`{e}bres
 $\big(\mathcal{W}\otimes\Lambda(E^*),\circ_E\big)$
 et $\big(\mathcal{W}\otimes\Lambda(F^*),\circ_F\big)$, il n'est pas
 difficile de montrer que $\mathsf{r}$ ne peut avoir que les trois
 composantes $\mathsf{r}^E_E$, $\mathsf{r}^F_E$ et $\mathsf{r}^F_F$
 donn\'{e}es. De plus, puisque $[\mathsf{r}^F_E,\phi^*A']_{\circ}=0=
 [\mathsf{r}^F_F,\phi^*A']$ il vient que
 $[\mathsf{r}^E_E-\phi^*\mathsf{r}',\phi^*A']_{\circ_E}=0$, et le fait
 que $\mathsf{r}^E_E\in\mathcal{W}_2\otimes\Lambda^1(E^*)$ et
 $\mathsf{r}'\in\mathcal{W}_2\otimes\Lambda^1(T^*M')$ entra\^{\i}ne que
 $\mathsf{r}^E_E-\phi^*\mathsf{r}'=0$.\\
 2. Calcul direct en utilisant le fait que $\partial_E$ est la connexion
 retir\'{e}e de $\partial'$.\\
 3. On fait la d\'{e}composition de l'\'{e}quation
 (\ref{EqFedCourbureGenerale}) selon les composantes en $\Lambda^2E^*$,
 $E^*\wedge F^*$ et $\Lambda^2 F^*$. Pour la composante $\Lambda^2E^*$,
 on utilise l'\'{e}quation (\ref{EqFedCourbureGenerale}) pour $\mathsf{r}'$
 et le fait que $\mathsf{r}^E_E=\phi^*\mathsf{r}'$ (qui implique
 $\partial_F\mathsf{r}^E_E=0$). On traite de mani\`{e}re analogue les
 \'{e}quations pour les normalisations.
\ebew

\noindent On va \'{e}tudier la question de savoir sous quelles conditions les
quatre \'{e}quations (\ref{EqFedCourbureUnten}),
(\ref{EqFedKruemmungFENullCrucial}),
(\ref{EqFedrFEalsFunktionvonAtiyahMolino})
et (\ref{EqFedCourbureObenlelongdeF}) soient r\'{e}solubles:

\bprop \label{PFedosovEinfuehrungGamma}
 Soient $(M,\omega)$ et $(M',\omega')$ deux vari\'{e}t\'{e}s symplectiques
 et $\phi:M\ra M'$ une application de Poisson. Soient $E$ et $F$ les
 sous-fibr\'{e}s int\'{e}grables de $TM$ et $\delta,\delta_E,\delta_F,\partial,
 \partial_E,\partial_F,\delta',\partial'$ les structures mentionn\'{e}es
 ci-dessus avec $\delta\phi^*=\phi^*\delta$ et
 $\partial\phi^*=\phi^*\partial'$.
 %Soient $\mathsf{r}'$,
 %$\mathsf{r}^F_{F}$ et $\mathsf{r}^F_E$ comme dans la proposition
 %\ref{PConditionsPourFedLiees}.
 Pour tous les choix $\Omega'$, $\Omega=\Omega_{EE}+\Omega_{EF}+\Omega_{FF}$
 et $\mathsf{s}',\mathsf{s}^F$ (voir la proposition
 \ref{PConditionsPourFedLiees}, partie 3.), il vient
 \ben
  \item Il existe un unique \'{e}l\'{e}ment
  $\mathsf{r}'\in\mathcal{W_3}\otimes\Lambda^1(T^*M')$ satisfaisant
   (\ref{EqFedCourbureUnten}) avec
   ${\delta'}^{-1}\mathsf{r}'=\mathsf{s}'$.
  \item Il existe un unique \'{e}l\'{e}ment
  $\mathsf{r}^F_F\in\mathcal{W_3}\otimes\Lambda^1(F^*)$ satisfaisant
   (\ref{EqFedCourbureObenlelongdeF}) avec
   ${\delta_F}^{-1}\mathsf{r}^F=\mathsf{s}^F$.
  \item Il existe un unique \'{e}l\'{e}ment
   $\mathsf{r}^F_E\in\mathcal{W_3}(F^*)\otimes\Lambda^1 E^*$ satisfaisant
   (\ref{EqFedrFEalsFunktionvonAtiyahMolino}) qui d\'{e}pend de
   $\mathsf{r}^F_F$, de $\Omega_{FF}$, et de fa\c{c}on au plus lin\'{e}aire
   de $\Omega_{EF}$, du
   $1$-cocycle d'Atiyah-Molino $\rho_{AME}$ du fibr\'{e} $E$ et d'une
   s\'{e}rie de $1$-formes,
   $\sigma_F(\mathsf{r}^F_E)\in\nu\Ginf(M,\Lambda^1 E^*)[[\nu]]$.
  \item Il existe une d\'{e}formation formelle associative diff\'{e}rentielle
    $\bullet$ de la
    multiplication ext\'{e}rieure dans $\Ginf(M,\Lambda E^*)$ qui consiste en
   des op\'{e}rateurs bidiff\'{e}rentiels dont les d\'{e}riv\'{e}es partielles sont
   le long des feuilles de $F$ telle que
   \ben
    \item $\partial_E-\frac{1}{2\nu}ad_{\circ_F}(\mathsf{r^F_E})$ induit une
   d\'{e}formation diff\'{e}rentielle $\mathsf{d}_E$
   de la diff\'{e}rentielle de Cartan le long de
   $E$, $d_E$, qui soit une d\'{e}rivation gradu\'{e}e de degr\'{e} $+1$ de
   $\big(\Ginf(M,\Lambda E^*)[[\nu]],\bullet\big)$,
    \item le carr\'{e} $\mathsf{d}_E^2$ de $\mathsf{d}_E$ est donn\'{e} par
      \[
        \mathsf{d}_E^2\beta=
            -\frac{1}{2\nu}(\gamma\bullet\beta-\beta\bullet\gamma)
      \]
       avec
      \beq\label{EqLaDeuxFormeGamma}
        \gamma:=\sigma_F(\partial_E\mathsf{r}^F_E)
             -\frac{1}{2\nu}\sigma_F(\mathsf{r}^F_E
               \circ_F\mathsf{r}^F_E)+\Omega_{EE}-\phi^*\Omega'
      \end{equation}
    \item et le membre de droite de l'\'{e}qn (\ref{EqFedKruemmungFENullCrucial})
      s'annule si et seulement si $\gamma=0$.
   \een
 \een
\eprop
\bbew
 $1.$ et $2.$ se d\'{e}duisent directement du th\'{e}or\`{e}me de Fedosov
 \ref{TFedosovConstructionPartielle}. L'alg\`{e}bre $\mathcal{A}$ pour $2.$ est
 $\mathcal{W}\otimes\Lambda(E^*)$.\\
 $3.$ Les identit\'{e}s suivantes sont utiles et faciles \`{a} calculer:
  \beq \label{EqIdentitesutilesDeltaRVV}
   \delta_ER^E_{EE}=0,~~\delta_FR^E_{EE}=0,~~\delta_ER^F_{EF}=0,
   ~~\delta_FR^F_{EF}=0,~~\delta_ER^F_{FF}=0,~~\delta_FR^F_{FF}=0,
 \end{equation}
 \beq \label{EqIdentitesutilesPartialRVV}
  \partial_ER^E_{EE}=0,~~\partial_FR^E_{EE}=0,~~\partial_ER^F_{EF}=0,
   ~~\partial_FR^F_{EF}+\partial_ER^F_{FF}=0,~~\partial_FR^F_{FF}=0.
 \end{equation}
 Avec $D^F_F:=-\delta_F+\partial_F-\frac{1}{2\nu}ad_\circ{\mathsf{r}^F_F}$
 on voit que $(D^F_F)^2=0$ et on calcule que
 \beq \label{EqDFFRFEFpluspErFFplusOmegaEF}
    D^F_F(R^F_{EF} +\partial_E\mathsf{r}^F_F
              +\Omega_{EF})=0
 \end{equation}
 en utilisant les identit\'{e}s pr\'{e}c\'{e}dentes
 (\ref{EqIdentitesutilesDeltaRVV}),
 (\ref{EqIdentitesutilesPartialRVV}),
 (\ref{EqIdUtilesPartialCarreEF}) et (\ref{EqIdUtilesPartialCarreFF}).
 Parce que le degr\'{e} antisym\'{e}trique par rapport au fibr\'{e} $F$ de
 $R^F_{EF} +\partial_E\mathsf{r}^F_F+\Omega_{EF}$ est \'{e}gal \`{a} $1$,
 le fait que la cohomologie d\'{e}finie par $D^F_F$ sur
 $\mathcal{W}\otimes\Lambda(T^*M)$ est concentr\'{e} dans le noyau de
 $deg_{aF}$ implique que $R^F_{EF} +\partial_E\mathsf{r}^F_F
 +\Omega_{EF}$ est dans l'image de $D^F_F$: par cons\'{e}quent, l'\'{e}quation
 (\ref{EqFedrFEalsFunktionvonAtiyahMolino}) est r\'{e}soluble. En applicant
 $\delta_F^{-1}$ \`{a} l'\'{e}quation (\ref{EqFedrFEalsFunktionvonAtiyahMolino})
 et en observant que $\delta_F^{-1}\mathsf{r}^F_E=0$ on obtient, \`{a} l'aide
 de
 \beq \label{EqClasseAtiyahMolinoEtCourbureFedosov}
    \delta_F^{-1}R^F_{EF}=-\rho_{AME}.
 \end{equation}
 l'\'{e}quation suivante:
 \beq \label{EqFormuleExplicitefFE}
   \mathsf{r}^F_E=\sum_{k=0}^\infty
         [\delta_F^{-1},\partial_F-\frac{1}{2\nu}ad_\circ(\mathsf{r}^F_F)]^k
         \big(\sigma_F\mathsf{r}^F_E
           -\rho_{AME}-\partial_E\mathsf{s}^F
           +\delta_F^{-1}\Omega_{EF}\big).
 \end{equation}
 On rappelle $\mathsf{s}^F=\delta_F^{-1}\mathsf{r}^F_F$.

 4. On regarde le cas $(N,K,\mathcal{A})=(M,F,\Ginf(M,\Lambda E^*)[[\nu]])$.
 Il est est facile de voir que la partie de la courbure $R^A_{FF}$ de
 $\partial_F$ dans $\Lambda E^*$ s'annule. Donc le th\'{e}or\`{e}me de
 Fedosov \ref{TFedosovConstructionPartielle} s'applique et nous donne
 une d\'{e}formation associative $\bullet$ de la multiplication ext\'{e}rieure
 $\wedge$ de $\Ginf(M,\Lambda E^*)[[\nu]]$. On voit d'ailleurs que le
 terme $\sum_{k=0}^\infty
         [\delta_F^{-1},\partial_F-\frac{1}{2\nu}ad_\circ(\mathsf{r}^F_F)]^k
         \big(\sigma_F\mathsf{r}^F_E\big)$ dans la formule de
         $\mathsf{r}^F_E$,
 (\ref{EqFormuleExplicitefFE}) est la s\'{e}rie de Fedosov-Taylor
 de la s\'{e}rie de $1$-formes
 $\sigma_F(\mathsf{r}^F_E)\in\Ginf(M,\Lambda^1E^*)$ et se trouve dans le
 noyau de $D^F_F$. De plus, l'op\'{e}rateur
 $D^F_E:=\partial_E-\frac{1}{2\nu}ad(\mathsf{r}^F_E)$ pr\'{e}serve l'alg\`{e}bre
 $\mathcal{A}\otimes\mathcal{W}\otimes\Lambda(F^*)$ et est une d\'{e}rivation
 gradu\'{e}e de degr\'{e} $1$ par rapport \`{a} $deg_{aE}$ et de degr\'{e} $0$ par rapport \`{a}
 $deg_{aF}$ $0$. Un calcul direct
 montre que $D_F^F=-\delta_F+\partial_F-\frac{1}{2\nu}ad(r^F_F)$ et $D^F_E$
 anticommutent. Il vient que $D^F_E$ pr\'{e}serve le noyau de $D^F_F$, donc
 les s\'{e}ries de Fedosov-Taylor,
 alors  il existe
 une unique application $\korps[[\nu]]$-lin\'{e}aire $\mathsf{d}_E$ de
 $\Ginf(M,\Lambda E^*)[[\nu]]$ telle que
 \[
    D^F_E\tau_F(\beta)=\tau_F(\mathsf{d}_E\beta).
 \]
 quel que soit $\beta\in \Ginf(M,\Lambda E^*)[[\nu]]$.
 Il est clair que
 $\mathsf{d}_E$ est une d\'{e}rivation gradu\'{e}e de degr\'{e} $1$ par rapport \`{a}
 $deg_{aE}$ de
 l'alg\`{e}bre d\'{e}form\'{e}e $\big(\Ginf(M,\Lambda E^*)[[\nu]],\bullet\big)$, car
 $D^F_E$ est une d\'{e}rivation gradu\'{e}e.
 %Soit $\beta'\in
 %\Ginf(M',\Lambda T^*M')[[\nu]]$. Puisque $\partial_F\phi^*\beta'=0$
 %il vient $D^F_F\phi^*\beta'=0$, donc $\phi^*\beta'=\tau_F(\phi^*\beta')$.
 %Par cons\'{e}quent,
 %\[
 %  \mathsf{d}_E\phi^*\beta' =\partial_E\phi*\beta'=d_E\phi^*\beta'
 %\]
 Puisque $r^F_E\in\mathcal{W}_2(F^*)\otimes \Lambda^1(E^*)$ il vient
 que $(\mathsf{d}_E)_0=d_E$.
 Finalement, le carr\'{e} de $D^E_F$
 (bien s\^{u}r restreint \`{a} $\mathcal{A}\otimes \mathcal{W}\otimes
 \Lambda(F^*)$) est \'{e}gal \`{a}
 \[
    (D^F_E)^2=-\frac{1}{2\nu}ad_{\circ_F}\big(
                 \partial_E\mathsf{r}^F_E
             -\frac{1}{2\nu}\mathsf{r}^F_E\circ_F\mathsf{r}^F_E
              +\Omega_{EE}-\phi^*\Omega'\big)
 \]
 et un calcul direct montre que
 \[
    D^F_F\big(\partial_E\mathsf{r}^F_E
             -\frac{1}{2\nu}\mathsf{r}^F_E\circ_F\mathsf{r}^F_E
              +\Omega_{EE}-\phi^*\Omega'\big)=0.
 \]
 Alors il existe une unique s\'{e}rie formelle de $2$-formes
 $\gamma\in \nu\Ginf(M,\Lambda^2 E^*)[[\nu]]$ donn\'{e}e
 par $\sigma_F\big(\partial_E\mathsf{r}^F_E
             -\frac{1}{2\nu}\mathsf{r}^F_E\circ_F\mathsf{r}^F_E
              +\Omega_{EE}-\phi^*\Omega'\big)$ telle que
 \[
   \partial_E\mathsf{r}^F_E
             -\frac{1}{2\nu}\mathsf{r}^F_E\circ_F\mathsf{r}^F_E
              +\Omega_{EE}-\phi^*\Omega'=\tau_F(\gamma).
 \]
 Ceci montre la proposition.
 \ebew

\noindent Les premiers deux ordres de la s\'{e}rie formelle de $2$-formes
$\gamma$ donnent des obstructions enti\`{e}rement analogues \`{a} celles du
th\'{e}or\`{e}me \ref{TObstructionsOrdreTrois}:

\bsat \label{TQuantApplPoissCondNec}
 Soit $\phi:M\ra M'$ une application de Poisson entre deux vari\'{e}t\'{e}s
 symplectiques $(M,\omega)$ et $(M',\omega')$. Soit $F:=\mathrm{Ker}~T\phi$
 et $E$ le sous-fibr\'{e} $F^\omega$ (voir la proposition \ref{PMorphPoissonSymp}).
 Soit $*$ un star-produit
 sur $M$ et $*'$ un star-poduit sur $M'$. Pour que $\phi$ soit
 quantifiable il est n\'{e}cessaire que les deux conditions soient
 satisfaites:
 \ben
  \item $p_v \big([*]_0-\phi^*[*']_0\big)=0$, c.-\`{a}-d. il existe un
    repr\'{e}sentant $\alpha_0$ de
    $[*]_0-\phi^*[*']_0$ telle que $\alpha_0$ est une $2$-forme relative.
  \item De plus, $\alpha_0$ a la propri\'{e}t\'{e} suivante:
    \bea
      \lefteqn{\frac{1}{12}\overline{P_F}^{(3)}
            \big(\kappa_{AM}(M,E),\kappa_{AM}(M,E)\big)}\nonumber \\
        & & ~~~~~~~
      +\frac{1}{2}\overline{P_F}^{(1)}([\alpha_0]_{(1,1)},[\alpha_0]_{(1,1)})
        -p_v[*]_1=0 \label{EqObstructionImportantePoisson}
    \eea
 \een
 o\`{u} $\kappa_{AM}(M,E)$ d\'{e}signe la classe d'Atiyah-Molino de la vari\'{e}t\'{e}
 pr\'{e}symplectique $(M,\varpi:=\omega_F)$ et $[~~]_{(1,1)}$ d\'{e}signe la
 classe de cohomologie d'une $2$-forme relative vue comme $1$-forme
 logitudinale \`{a} valeurs dans $Q^*:=(TM/E)^*\cong F^*$.
\esat
\bbew
 D'apr\`{e}s les propositions \ref{PConditionsPourFedLiees} et
 \ref{TQuantApplPoissCondNec} il faut calculer la s\'{e}rie formelle de
 $2$-formes $\gamma=\sum_{r=1}^\infty\nu^r\gamma_r$,
  voir (\ref{EqLaDeuxFormeGamma}): soit
 $\chi=\sum_{r=0}^\infty\nu^r\chi_r$ la s\'{e}rie d'$1$-formes
 $\sigma_F(\mathsf{r}^F_E)$. En g\'{e}n\'{e}ral, soit
 $\mathsf{r}_{k,r}\in \Ginf(M,S^kF^*\otimes \Lambda^1 E^*)$ la composante
 de $\mathsf{r}^F_E$ de degr\'{e} $k$ par rapport \`{a} $deg_{sF}$ et de l'ordre
 $r$ (par rapport au param\`{e}tre formel $\nu$). Soit
 $\mathsf{s}_{k,r}\in \Ginf(M,S^kF^*)$ form\'{e}e de mani\`{e}re analogue.
 Puisque $\mathsf{r}^F_E$
 est de degr\'{e} total au moins $2$ il vient que $\mathsf{s}_{1,0}=0$.
 \[
   \gamma_1 = \partial_E \chi_1
                   +\Omega_{EE1}-\phi^*\Omega'_1.
 \]
 On calcule que $\mathsf{r}_{1,1}=\delta_F^{-1}(\Omega_{EF1}+d\chi_1)
 -\partial_E \mathsf{s}^F_{1,1}$ et $\mathsf{r}_{3,0}=-\rho_{AME}-\partial
 s^F_{3,0}$. Il vient
 \beas
   \gamma_2 & = & \partial_E\chi_2 -\frac{1}{2}
                    P_F^{(1)}\big(\delta_F^{-1}(\Omega_{EF1}+d\chi_1)
                                    -\partial_E \mathsf{s}^F_{1,1},
                                  \delta_F^{-1}(\Omega_{EF1}+d\chi_1)
                                    -\partial_E \mathsf{s}^F_{1,1}\big) \\
            &   & ~~~~~~~~-\frac{1}{12}
                      P_F^{(3)}\big(\rho_{AME}+\partial_E s^F_{3,0},
                                     \rho_{AME}+\partial_E s^F_{3,0}\big)
                     +\Omega_{EE2}-\phi^*\Omega'_2.
 \eeas
 En passant aux classes de cohomologie, on obtient le r\'{e}sultat.
\ebew

De la proposition \ref{TQuantApplPoissCondNec} on peut tirer le cas
particulier suivant:

\bsat \label{TQuantApplPoissonAtiyahMolinoNulle}
 Soit $\phi:M\ra M'$ une application de Poisson entre deux vari\'{e}t\'{e}s
 symplectiques $(M,\omega)$ et $(M',\omega')$. Soit $F:=Ker T\phi$
 et $E$ le sous-fibr\'{e} $F^\omega$
 (voir la proposition \ref{PMorphPoissonSymp}).
 Soit $*'$ un star-produit
 sur $M'$ et soit $*$ un star-poduit sur $M$. On suppose que les
 deux conditions suivantes soient satisfaites:
 \ben
  \item La classe d'Atiyah-Molino $\kappa_{AM}(M,E)$ du sous-fibr\'{e} $E$
  s'annule.
  \item Il existe une s\'{e}rie $\Omega_{FF}=\sum_{r=1}^\infty\nu^r\Omega_{FFr}$
  de $2$-formes ferm\'{e}es sur $M$ telle que $i_V\Omega_{FF}=0$ quel que soit
  le champ de vecteurs $V\in \Ginf(M,E)$ et telle que la classe de Deligne
  $[*]$ de $*$ est \'{e}gale \`{a}
  \[
      [*](\nu)= \frac{[\omega]}{\nu} +
             \phi^*\left([*'](\nu)-\frac{[\omega']}{\nu}\right)
                  + \frac{[\Omega_{FF}](\nu)}{\nu}.
  \]
 \een
 Alors $\phi$ est quantifiable.
\esat
\bbew
 Puisque la classe d'Atiyah-Molino s'annule, on peut choisir les
 connexions symplectiques $\tilde{\nabla}$ sur $M$ et $\nabla'$ sur
 $M'$ de telle fa\c{c}on  que $\phi$ est toujours une application
 affine et le $1$-cocycle d'Atiyah-Molino
 $\rho_{AME}$ du fibr\'{e} s'annule. Dans la construction de Fedosov
 d\'{e}crite ci-dessus, il vient que la courbure $R^F_{EF}=-\delta_F(\rho_{AME})=0$
 voir (\ref{EqClasseAtiyahMolinoEtCourbureFedosov}). Grace \`{a} l'hypoth\`{e}se
 2), on peut choisir $\Omega_{EF}=0$ et
 $\Omega_{EE}=\phi^*\Omega'$. On choisit la normalisation $\mathsf{s}^F=0$
 et $\chi=\sigma(\mathsf{r}^F_E)=0$. La formule (\ref{EqFormuleExplicitefFE})
 montre que $\mathsf{r}^F_E=0$. Il s'ensuit que $\gamma=0$, alors il
 existe deux d\'{e}rivations de Fedosov $D$ et $D'$ qui sont $\phi$-li\'{e}es, et
 $\phi$ est quantifiable d'apr\`{e}s le th\'{e}or\`{e}me \ref{TFedQuantApplPois}.
\ebew

\noindent Dans le cas o\`{u} le fibr\'{e} $F$ s'annule on obtient le r\'{e}sultat suivant
dont le cas particulier d'un symplectomorphisme est bien connu, voir
par exemple \cite{GR99}, Cor. 9.7.:
\bcor
 Soit $\phi:M\ra M'$ une symplectomorphisme local entre deux vari\'{e}t\'{e}s
 symplectiques $(M,\omega)$ et $(M',\omega')$.
 Soit $*'$ un star-produit
 sur $M'$ et soit $*$ un star-poduit sur $M$.
 Alors $\phi$ est quantifiable si et seulement si
 \[
      \phi^*[*']=[*].
 \]
 En particulier, si $(M,\omega)=(M',\omega')$, $*=*'$ et si $\phi$ est un
 diff\'{e}omorphisme, alors $\phi$ est quantifiable en tant qu'automorphisme
 si et seulement si
 \[
       \phi^*[*]=[*].
 \]
\ecor
\bbew
 Puisque les deux vari\'{e}t\'{e}s ont la m\^{e}me dimension, alors le sous-fibr\'{e} $F$
 s'annule, et
 les hypoth\`{e}ses du th\'{e}or\`{e}me pr\'{e}c\'{e}dent
 \ref{TQuantApplPoissonAtiyahMolinoNulle} sont satisfaites, alors
 la condition $\phi^*[*']=[*]$ est suffisante pour la quantification
 de $\phi$. R\'{e}ciproquement, les propositions \ref{PConditionsPourFedLiees}
 et \ref{PFedosovEinfuehrungGamma}, surtout l'\'{e}qn (\ref{EqLaDeuxFormeGamma})
 montrent la n\'{e}cessit\'{e} de cette condition.\\
 Le cas particulier est \'{e}vident.
\ebew

\subsection{La fibration en espaces r\'{e}duits locaux}
   \label{SubSecFibEspRedLoc}

Soit $(M,\omega)$ une vari\'{e}t\'{e} symplectique et $i:C\ra M$ une sous-vari\'{e}t\'{e}
co\"{\i}sotrope ferm\'{e}e. Soit $\tau:\check{M}\ra C$ un voisinage tubulaire.
Dans ce paragraphe on va construire une vari\'{e}t\'{e} fibr\'{e}e $\phi :E^C\ra\check{M}$
sur $\check{M}$ qui g\'{e}n\'{e}ralisera le produit cart\'{e}sien
$M\times M_{\mathrm{red}}$ utilis\'{e} pour la r\'{e}duction des star-produits en
paragraphe \ref{SubSecEspRedReduction}.

\bsat \label{TVarieteEspacesReduitsLocaux}
 Soit $(M,\omega)$ une vari\'{e}t\'{e} symplectique, $i:C\ra M$ une sous-vari\'{e}t\'{e}
co\"{\i}sotrope ferm\'{e}e et $\tau:\check{M}\ra C$ un voisinage tubulaire de
$C$ dans $M$. Alors il existe une vari\'{e}t\'{e} fibr\'{e}e $\phi:E^C \ra \check{M}$ sur
$\check{M}$ telle que:
\ben
 \item $E^C$ admet une forme symplectique $\omega^C$ telle que la
 submersion $\phi$ soit une application de Poisson.
 \item Il existe un plongement $j:C\ra E^C$ tel que $\phi\circ j=i$ et
 que $j(C)$  soit une sous-vari\'{e}t\'{e} lagrangienne de $E^C$.
 \item Soit $E$ le sous-fibr\'{e} caract\'{e}ristique de la vari\'{e}t\'{e}
 pr\'{e}symplectique $(C,\varpi:=i^*\omega)$. Si la classe d'Atiyah-Molino de
 $E$ s'annule, alors la classe d'Atiyah-Molino du fibr\'{e}
 $(Ker~T\phi)^{\omega^E}$ s'annule.
\een
\esat
\bbew
 Soit $F$ un sous-fibr\'{e} de $TC$ qui soit compl\'{e}mentaire au fibr\'{e} $E$
 (i.e. $TC=F\oplus E$). Soit $\nabla$ une connexion sans torsion dans le
 fibr\'{e} tangent de $C$ telle que $(F,\nabla)$ soit adapt\'{e}e dans le sens
 du paragraphe \ref{SubSecConnSympAdaptCois}. Soit $Exp$ l'application
 exponentielle de $\nabla$. D'apr\`{e}s la th\'{e}orie g\'{e}n\'{e}rale
 des connexions il existe un voisinage ouvert $W^{(1)}$ de la section nulle de
 $TC$ tel que l'application $x\mapsto
 \big(\tau_C(x),Exp_{\tau_C(x)}(c)\big)$ soit un diff\'{e}omorphisme sur un
 voisinage ouvert de la diagonale $\Delta C:=
 \{(c,c)~|~c\in C\}$ dans $C\times C$.
 Soit $\hat{W}^{(1)}$ l'ouvert $(\tau\times id)^{-1}W^{(1)}$ du fibr\'{e}
 retir\'{e} $\tau^*TC$ sur $\check{M}$. Pour une courbe $\gamma:]a,b[\ra C$
 (avec $a<0<b$ et $\gamma(0)=c$) soit
 $\mathsf{Tp}_\gamma(t): T_cC\ra T_{\gamma(t)}C$ le transport parall\`{e}le
 le long de la courbe $\gamma$. Puisque la connexion $\nabla$ pr\'{e}serve
 les champs de vecteurs verticaux (i.e. $\nabla_XV\in\Ginf(C,E)$ quels que
 soient $X\in\Ginf(C,TC)$ et $V\in\Ginf(C,E)$), le transport parall\`{e}le
 pr\'{e}serve le sous-fibr\'{e} $E$. Soit $W^{(2)}$ l'ouvert dans $W^{(1)}$
 d\'{e}fini par
 \[
  W^{(2)}:=
   \{w+v\in F\oplus E~|~c=\tau_C(w+v)~\mathrm{et}
        ~\mathsf{Tp}_{t\mapsto Exp_{c}(tw)}(1)\big(v\big)\in W^{(1)}\cap E\}
             \cap W^{(1)},
 \]
 et soit $\hat{W}^{(2)}$ l'ouvert $(\tau\times id)^{-1}W^{(2)}$ du fibr\'{e}
 retir\'{e} $\tau^*TC$ sur $\check{M}$. Evidemment, $\check{M}\subset
 \hat{W}^{(2)}$.
 Alors l'application
 \beq\label{EqHatApplFibreVecVarFeuilletee}
    \hat{\Phi}:\hat{W}^{(2)}\ra \check{M}\times C:
        (m,w+v)\mapsto \Big(m,Exp_{Exp_{\tau(m)}(w)}
              \big(\mathsf{Tp}_{t\mapsto Exp_{\tau(m)}(tw)}(1)\big(v\big)\big)
              \Big)
 \end{equation}
 est bien d\'{e}finie et de classe $\Cinf$. Pour $(m,0_c)\in \tau^*TC$,
 $y\in T_mM$
 et $w+v\in F_c\oplus E_c$ on calcule
 \[
   T_{(m,0_c)}\hat{\Phi}(y,w+v)=(y,T_c\tau y + w + v),
 \]
 et $\hat{\Phi}$ est donc un diff\'{e}omorphisme local dans un
 ouvert de $\hat{W}^{(2)}$ qui contient $\check{M}$. De la
 m\^{e}me fa\c{c}on utilis\'{e}e pour la construction des voisinages tubulaires
 (voir \cite{Lan95}), on peut r\'{e}tr\'{e}cir $\hat{W}^{(2)}$ jusqu'\`{a} ce qu'on
 obtienne un ouvert $\hat{W}^{(3)}$
 avec $\check{M}\subset \hat{W}^{(3)}
 \subset \hat{W}^{(2)}$ tel que la restriction de $\hat{\Phi}$ \`{a}
 $\hat{W}^{(3)}$ est un diff\'{e}omorphisme sur un ouvert $\hat{U}^{(3)}$ de
 $\check{M}\times C$. Soit $g^E$ (resp. $g^F$) une m\'{e}trique d\'{e}finie positive
 dans le fibr\'{e} $\tau^*E$ (resp. $\tau^*F$). D'apr\`{e}s le lemme
 \ref{LOuvertBeauSectionNulle} il existe $h_1,h_2\in\CinfR{\check{M}}$
 \`{a} valeurs strictement positives telles qu'on a
 \beas
  \check{M} & \subset & \hat{W} \\
      &:= &\{v=(m,w+v)\in\tau^*(F\oplus E)~|~
                       g^E(v,v)<h_1(m)~\mathrm{et}~
                                     g^F(w,w)<h_2(m)\} \\
            & \subset &
              \hat{W}^{(3)}.
 \eeas
 Alors la restriction de l'application $\hat{\Phi}$ (voir l'\'{e}q.
 (\ref{EqHatApplFibreVecVarFeuilletee})) \`{a} $\hat{W}$ est un
 diff\'{e}omorphisme sur un ouvert $\hat{U}$ de $\check{M}\times C$
 qui contient $(\tau\times id)^{-1}(\Delta C)$.
 On regarde la submersion surjective
 $\tau^*TC\ra \tau^*F$ qui consiste en projetant fibre-par-fibre le fibr\'{e}
 vectoriel $\tau^*(F\oplus E)$ sur le sous-fibr\'{e} $\tau^*F$ le long du
 sous-fibr\'{e} $E$. Soit $\psi$
 la restriction de cette projection \`{a} $\hat{W}$. Par construction il vient
 que l'ouvert $\psi(\hat{W})$ de $\tau^*(F)$ est \'{e}gal \`{a}
 \[
   \psi(\hat{W})=\{(m,w)\in\tau^*F~|~g^F(w,w)<h_2(m)\}=:W=\hat{W}\cap
   \tau^*F.
 \]
 Puisque $\tau^*F$ est une sous-vari\'{e}t\'{e} ferm\'{e}e de $\tau^*TC$, il vient que
 $W$ est une sous-vari\'{e}t\'{e} ferm\'{e}e de $\hat{W}$.
 Soit $\Phi:W\ra \check{M}\times C$ l'application
 \beq
   \Phi(m,w)=\big(m,Exp_{\tau(m)}(w)\big).
 \end{equation}
 Il s'ensuit que la restriction du diff\'{e}omorphisme $\hat{\Phi}:\hat{W}\ra
 \hat{U}$ \`{a} $W$ est \'{e}gale \`{a} $\Phi$ (car c'est le cas o\`{u} la composante
 $v\in E$ est toujours nulle en (\ref{EqHatApplFibreVecVarFeuilletee})).
 Alors
 \[
    \hat{E}^C:=\Phi(W)=\hat{\Phi}(W)
 \]
  est une sous-vari\'{e}t\'{e} ferm\'{e}e de l'ouvert $\hat{U}$ de $\check{M}\times
  C$. Soit $p:\hat{U}\ra \hat{E}^C$ la submersion surjective
  \[
     p:=\Phi\circ\psi\circ\hat{\Phi}^{-1}.
  \]
  Pour tout $m\in\check{M}$, soit
  $\tau^*(F\oplus E)_m$ la fibre sur $m$ du fibr\'{e} $\tau^*(F\oplus E)$.
  Alors l'intersection $\hat{W}\cap \tau^*(F\oplus E)_m$ est donn\'{e}e
  par le produit cart\'{e}sien $W_m\times W'_m$ avec
  $W_m:=\{w\in \tau^*(F)_m~|~g^F(w,w)<h_2(m)\}$ et
  $W'_m:=\{v\in \tau^*(E)_m~|~g^E(v,v)<h_1(m)\}$. Evidemment, $W_m$ et
  $W'_m$ sont des ouverts dans $\tau^*(F)_m$ et $\tau^*(E)_m$,
  respectivement. Soit $\hat{\Phi}_m:W_m\times W'_m\ra C$ la restriction de
  $pr_2\circ\hat{\Phi}$
  \`{a} $W_m\times W'_m$ et $U_m:=\hat{\Phi}_m(W_m\times W'_m)\subset C$.
  Puisque $pr_1\circ\hat{\Phi}$ est \'{e}gale \`{a} la projection canonique
  $\tau^*TC\ra\check{M}$ et puisque tous les espaces tangents \`{a} la
  sous-vari\'{e}t\'{e} $W_m\times W'_m$ sont dans le noyau de l'application
  tangente de $pr_1\circ\hat{\Phi}$ il vient que $\hat{\Phi}_m$ est
  un diff\'{e}omorphisme, donc tous les $U_m$ sont des ouverts de $C$.
  Puisque $U_m$ contient $\tau(m)\in C$ il vient que les $U_m$ recouvrent
  $C$. Donc $(U_m,\hat{\Phi}_m^{-1})_{m\in\check{M}}$ est un atlas de
  $C$. Puisque $\nabla_VV'\in\Ginf(C,E)$ quels que soient
  $V,V'\in\Ginf(C,E)$ il vient que pour tout $c\in C$ la partie
  $P_c:=\{Exp_c(v)~|~v\in Z_c\}$ (o\`{u} $Z_c$ est une boule ouverte de centre
  $0\in
  E_c$ et de rayon strictement positif suffisamment petit) donne une
  sous-vari\'{e}t\'{e} int\'{e}grale
  de $C$ dont le fibr\'{e} tangent est donn\'{e} par la restriction de $E$ \`{a}
  $P_c$. Par cons\'{e}quent, $(U_m,\hat{\Phi}_m^{-1})_{m\in\check{M}}$
  consiste en cartes distingu\'{e}es de la vari\'{e}t\'{e} feuillet\'{e}e $C$: la plaque
  de $(U_m,\hat{\Phi}_m^{-1})$ passant par $c=\hat{\Phi}_m(w,v)\in U_m$ est
  donn\'{e} par $\{\hat{\Phi}_m(w,v')~|~v'\in W'_m\subset E_m\}$. On voit
  rapidement que les fibres de la submersion $p$ sont donn\'{e}es par
  \beas
     \lefteqn{p^{-1}\big(p(m,c)\big)=}\\
            & &\{(m',c')\in\hat{U}~|~m'=m~\mathrm{et}~
                  c,c'\mathrm{~appartiennent~\grave{a}~la~m\hat{e}me~
                                               plaque~de~}U_m\}.
  \eeas
  quel que soit $(m,c)\in\hat{U}$. Ceci montre que la relation
  d'\'{e}quivalence $\sim$ dans $\hat{U}$ d\'{e}finie par ``$(m,c)\sim (m',c')$
  lorsque $m=m'$ et $c,c'$ appartiennent \`{a} la m\^{e}me plaque de $U_m$''
  a un espace quotient $E^C$ qui s'identifie \`{a} $\hat{E}^C$ en choisissant
  l'unique \'{e}l\'{e}ment de l'intersection de la classe d'\'{e}quivalence $[m,c]$
  de $(m,c)$ avec la sous-vari\'{e}t\'{e} $\hat{E}^C$. On va d\'{e}signer la
  projection canonique $\hat{U}\ra E^C$ \'{e}galement par $p$.
  De plus, puisque la restriction de la projection $\tau^{*}F\ra \check{M}$
  \`{a} $W\subset \tau^{*}F$ est une submersion $\phi'$ sur $\check{M}$, il
  en est de m\^{e}me avec $\hat{\phi}:=\phi'\circ\Phi^{-1}:\hat{E}^C\ra \check{M}$.
  La
  vari\'{e}t\'{e} fibr\'{e}e $\hat{\phi}:\hat{E}^C\ra\check{M}$ peut donc \^{e}tre identifi\'{e}e
  avec $\phi:E^C\ra\check{M}$ o\`{u} $\phi([m,c]):=m$ quel que soit
  $m\in\check{M}$, autrement dit
  \[
        \phi\circ p = pr_1.
  \]
  Ensuite on regarde la restriction de la $2$-forme ferm\'{e}e
  $pr_1^*\omega-pr_2^*\varpi$ \`{a} l'ouvert $\hat{U}$ de $M\times C$.
  Visiblement, cette $2$-forme est pr\'{e}symplectique et son fibr\'{e}
  caract\'{e}ristique consiste en des espaces tangents aux plaques
  dans tous les $U_m$.
  Par cons\'{e}quent, la vari\'{e}t\'{e} quotient $E^C$ est munie d'une forme
  symplectique $\omega^C$ telle que
  \[
     p^*\omega^C= \big(pr_1^*\omega-pr_2^*\varpi\big)|_{\hat{U}}.
  \]
  Soit $P$ la structure de Poisson associ\'{e}e \`{a} $\omega$ et soit
  $\tilde{P}\in\Ginf(C,\Lambda^2F)$ l'unique champ de bivecteurs tel que
  \[
   \varpi^\flat~\tilde{P}^\sharp = \mathsf{P}_{F^*}~~~
   \mathrm{et}~~~
   \tilde{P}^\sharp~\varpi^\flat = \mathsf{P}_{F}.
  \]
  o\`{u} $\mathsf{P}_{F}$ d\'{e}signe la projection $TC\ra F$ le long de $E$ et
  $\mathsf{P}_{F^*}$ d\'{e}signe la projection $T^*C\ra F^*$ le long de $E^*$.
  Alors la structure de Poisson $P^C$ associ\'{e}e \`{a} $\omega^C$ est donn\'{e}e
  par $P^C\circ p = Tp\otimes Tp\big(P_{(1)}-\tilde{P}_{(2)}\big)$, alors
  \beas
     T\phi\otimes T\phi\big(P^C\circ p\big)
        & = & T(\phi\circ p)\otimes T(\phi\circ p)
             \big(P_{(1)}-\tilde{P}_{(2)}\big) \\
        & = & T pr_1\otimes T pr_1\big(P_{(1)}-\tilde{P}_{(2)}\big) \\
        & = & P\circ pr_1 = P\circ \phi \circ p,
  \eeas
  et $\phi$ est une application de Poisson car $p$ est surjective.

  2. Soit $\hat{\j}:C\ra \check{M}\times C$ le plongement $c\mapsto
  (i(c),c)$. Alors l'image de $\hat{\j}$ est une sous-vari\'{e}t\'{e} ferm\'{e}e de
  $\hat{U}$ et il vient $pr_1\circ \hat{\j}=i$.
  On d\'{e}finit $j:C\ra E^C$ par $j(c):=p\big(\hat{\j}(c)\big)$. Alors
  $j$ est de classe $\Cinf$ et
  \[
       \phi\circ j= \phi \circ p\circ \hat{\j} =pr_1\circ \hat{\j} = i,
  \]
  donc $j$ est un plongement de $C$ dans $E^C$. Il est \'{e}vident que
  $\hat{\j}(C)$ est une sous-vari\'{e}t\'{e} isotrope de la sous-vari\'{e}t\'{e}
  pr\'{e}symplectique
  $\big(\hat{U},(pr_1^*\omega-pr_2^*\varpi)|_{\hat{U}}\big)$, alors
  $0=\hat{\j}^*(pr_1^*\omega-pr_2^*\varpi)|_{\hat{U}}$. Par cons\'{e}quent,
  \[
    j^*\omega^C=\hat{\j}^*p^*\omega^C
     =\hat{\j}^*(pr_1^*\omega-pr_2^*\varpi)|_{\hat{U}}=0
  \]
  donc $j(C)$ est une sous-vari\'{e}t\'{e} isotrope de $(E^C,\omega^C)$ et puisque
  $\dim E^C=\dim M + \dim C -(\dim M - \dim C)=2\dim C$ il s'ensuit que
  $j(C)$ est une sous-vari\'{e}t\'{e} lagrangienne de $(E^C,\omega^C)$.

  3. Supposons que $\kappa_{AM}(C,E)=0$. Alors il existe une connexion
  adapt\'{e}e $(F,\nabla)$ sur $C$ telle que $\mathsf{P}_F\big(R(V,X)Y\big)=0$
  quels que soient $V\in\Ginf(C,E)$ et $X,Y\in\Ginf(C,TC)$
  pour le tenseur de courbure $R$ de $\nabla$, voir l'\'{e}qn
  (\ref{EqCocycleAtiyahMolinoPresymp}) et le th\'{e}or\`{e}me
  \ref{TMolinoPresymp}. En outre, $\nabla\varpi=0$. Soit $\tilde{\nabla}$
  la connexion symplectique sur $(\check{M},\omega|_{\check{M}})$
  induite par $\nabla$, voir le th\'{e}or\`{e}me \ref{TConnSympAutourCoisotrope}.
  Alors le produit cart\'{e}sien $\check{M}\times C$ est muni d'une connexion
  $\hat{\nabla}$ d\'{e}finie par $\hat{\nabla}_{\tilde{X}+X}(\tilde{Y}+Y):=
  \tilde{\nabla}_{\tilde{X}}\tilde{Y}+\nabla_XY$ quels que soient
  $\tilde{X},\tilde{Y}\in\Ginf(\check{M},T\check{M})$ et
  $X,Y\in\Ginf(C,TC)$. On va noter la restriction de $\hat{\nabla}$
  \`{a} l'ouvert $\hat{U}$ de $\check{M}\times C$ par le m\^{e}me symbole
  $\hat{\nabla}$. Evidemment,
  $\hat{\nabla}(pr_1^*\omega-pr_2^*\varpi)|_{\hat{U}}=0$. De plus, le
  sous-fibr\'{e} $Ker T\phi$ de $T\hat{U}$ est visiblement donn\'{e} par
  $pr_2^*E$, i.e. par $\{(0_m,v_c)~|~(m,c)\in\hat{U}, ~v_c\in E_c\}$.
  On peut toujours repr\'{e}senter un vecteur tangent $v_c\in pr_2^*E$ par
  la valeur en $c$ d'un champs de vecteurs $V\in\Ginf(C,E)$. Le fibr\'{e}
  tangent de $\hat{U}$ se d\'{e}compose comme
  $T\hat{U}=pr_1^*T\check{M}\oplus pr_2^*F \oplus pr_2^*E$. Soit
  $\mathsf{P}_1$ la projection sur le sous-fibr\'{e}
  $pr_1^*T\check{M}\oplus pr_2^*F$ le long du sous-fibr\'{e} $pr_2^*E$.
  Avec $\tilde{X},\tilde{Y}\in\Ginf(\check{M},T\check{M})$ et
  $X,Y\in\Ginf(C,TC)$ il vient que
  \[
    \mathsf{P}_1\big(\hat{R}(V,\tilde{X}+X)(\tilde{Y}+Y)\big)
               =\mathsf{P}_1\big(R(V,X)Y\big)
               =\mathsf{P}_F\big(R(V,X)Y\big)=0.
  \]
  D'apr\`{e}s le lemme \ref{LSubmersionAffine} on peut conclure qu'il existe
  une unique connexion $\nabla'$ sur $E^C$ telle que la
  submersion $p$ est une application affine (car les fibres de $p$ sont
  connexes). Soit $B$ le sous-fibr\'{e} $Ker T\phi$ de $TE^C$ et $A\subset TE^C$
  le fibr\'{e} orthogonal \`{a} $B$ par rapport \`{a} $\omega^C$. L'\'{e}quation
  $\phi\circ p=pr_1$ implique que $Tp(pr_2^*TC)=Tp(pr_2^*F)=B$. Puisque
  $p$ est une submersion on a $Tp(T\hat{U})=TE^C$. Soit
  $(m,c)\in\hat{U}$ et $(z,x)\in T_{(m,c)}\hat{U}=T_mM\times T_cC$.
  Soit $T_{(m,c)}p(z,x)\in A_{p(m,c)}$. Alors quel que soit $(0,v)\in
  (pr_2^*E)_{(m,c)}$ il vient
  \beas
    0& = & \omega^\mathsf{C}_{p(m,c)}\big(T_{(m,c)}p(z,x),T_{(m,c)}p(0,v)\big)
      =(p^*\omega^C)_{(m,c)}\big((z,x),(0,v)\big)\\
     & = & \omega_m(z,0)-\varpi_c(x,v)=-\varpi_c(x,v),
  \eeas
  alors $x\in E_c$, donc $T_{(m,c)}p(x)=0$. Alors
  \[
      Tp(pr_1^*TM)=A.
  \]
  La d\'{e}composition $T\hat{U}=(pr_1^*TM\oplus pr_2^*F)\oplus pr_2^*E$
  d\'{e}finit un rel\`{e}vement horizontal $X'\mapsto X^{\prime h}$ des champs de
  vecteurs sur $E^C$ aux champs de vecteurs dans
  $\Ginf\big(\hat{U},(pr_1^*TM\oplus pr_2^*F)|_{\hat{U}}\big)$. Puisque
  $p$ est affine on a $Tp~\hat{\nabla}_{X^{\prime h}}Y^{\prime h}=
  \nabla'_{X'}Y'\circ p$. On calcule que $\nabla'$ pr\'{e}serve $A$ et $B$ car
  $\hat{\nabla}$ pr\'{e}serve $pr_1^*TM$ et $pr_2^*TC$. En outre, il vient que
  les tenseurs de courbure
  $\hat{R}$ de $\hat{\nabla}$ et $R'$ de $\nabla'$ sont $p$-li\'{e}s, i.e.
  $Tp~\hat{R}(X^{\prime h},Y^{\prime h})Z^{\prime h}=R'(X',Y')Z'\circ p$
  quels que soient $X',Y',Z'\in\Ginf(E^C,TE^C)$. Soit $(m,c)\in \hat{U}$.
  Soit $x'\in A_{p(m,c)}$
  et $y',z'\in B_{p(m,c)}$ et soient $\tilde{x}\in T_mM$, $y,z\in T_cC$ tels
  que $T_{(m,c)}p(\tilde{x})=x'$, $T_{(m,c)}p(y)=y'$ et
  $T_{(m,c)}p(z)=z'$. Soit $\mathsf{P}_B$ la projection de $TE^C$
  sur $B$ le long du fibr\'{e} $A$. Alors
  \[
     \mathsf{P}_B\big(R'_{p(m,c)}(x',y')z'\big)
     =\mathsf{P}_B\big(T_{(m,c)}p\big(\hat{R}_{(m,c)}(\tilde{x},y)z\big)\big)
     =0
  \]
  car $\hat{R}_{(m,c)}(\tilde{x},y)z=\tilde{R}_m(\tilde{x},0)0+
                                       R_c(0,y)z=0$. Alors le $1$-cocycle
  d'Atiyah-Molino de la vari\'{e}t\'{e} pr\'{e}symplectique $(E^C,\omega^C|_{B\times B})$
  relatif \`{a} $\nabla'$ s'annule, donc la classe d'Atiyah-Molino
  $\kappa_{AM}(E^C,A)$ est nulle.
\ebew

\noindent La vari\'{e}t\'{e} fibr\'{e}e $E^C$ construite dans la d\'{e}monstration du
th\'{e}or\`{e}me
pr\'{e}c\'{e}dent a comme fibres des `morceaux locaux des espaces r\'{e}duits' et
est diff\'{e}omorphe \`{a} un voisinage ouvert de la section nulle
du fibr\'{e} vectoriel $\tau^*F$. Cette construction sera appel\'{e}e
{\em la fibration
en espaces r\'{e}duits locaux sur $\check{M}$}.

\subsection{Quantification des sous-vari\'{e}t\'{e}s co\"{\i}sotropes \`{a} classe
            d'Atiyah-Molino nulle}
   \label{SubSecQuantSousVarCoisClasAMNulle}

Les consid\'{e}rations pr\'{e}c\'{e}dentes permettent de montrer l'existence des
repr\'{e}\-sentations des star-produits dans le cas particulier important
suivant:

\bsat\label{TQuantSousVarCoisClasAMNulle}
 Soit $(M,\omega)$ une vari\'{e}t\'{e} symplectique et $i:C\ra M$ une sous-vari\'{e}t\'{e}
 co\"{\i}sotrope ferm\'{e}e de $M$. Soit $E\subset TC$ le sous-fibr\'{e}
 caract\'{e}ristique de la vari\'{e}t\'{e} pr\'{e}symplectique $(C,\varpi:=i^*\omega)$.
 De plus, soit $*$ un star-produit sur $M$ et $[*]$ sa classe de Deligne.
 On suppose que les \'{e}nonc\'{e}s suivants
 soient vrais:
 \ben
  \item La classe d'Atiyah-Molino $\kappa_{AM}(C,E)$ s'annule.
  \item $i^*[*]=0$.
 \een
 Alors il existe une repr\'{e}sentation diff\'{e}rentielle de $*$ sur $C$.
\esat
\bbew
 On consid\`{e}re un voisinage tubulaire $\tau:\check{M}\ra C$ de $C$ (o\`{u}
 $\check{M}$ est un ouvert de $M$ qui contient $C$). De plus, on consid\`{e}re
 la fibration en espaces r\'{e}duits locaux $\phi:E^C\ra\check{M}$, voir
 le th\'{e}or\`{e}me \ref{TVarieteEspacesReduitsLocaux}. Puisque par construction
 la sous-vari\'{e}t\'{e} $C$ est un r\'{e}tracte par d\'{e}formation de $\check{M}$ et de
 $E^C$, la cohomologie de de Rham de $E^C$ et de $\check{M}$
 est isomorphe \`{a} celle de $C$. Alors
 la restriction de la classe de Deligne de $*$ \`{a} $\check{M}$ s'annule.
 D'apr\`{e}s le th\'{e}or\`{e}me \ref{TVarieteEspacesReduitsLocaux}, la classe
 d'Atiyah-Molino du sous-fibr\'{e} $(Ker T\phi)^{\omega^C}$ de $TE^C$ s'annule
 car $\kappa_{AM}(C,E)$ s'annule. On peut choisir un star-produit $*^C$
 sur $E^C$ tel que
 \[
      [*^C]=\frac{[\omega^C]}{\nu}
               + \phi^*\left([*]-\frac{[\omega]}{\nu}\right).
 \]
 D'apr\`{e}s le th\'{e}or\`{e}me \ref{TQuantApplPoissonAtiyahMolinoNulle}, il existe
 une quantification $\Phi$ de l'application de Poisson $\phi$, c.-\`{a}-d.
 un homomorphisme diff\'{e}rentiel d'alg\`{e}bres associatives sur $\korps[[\nu]]$ de
 $\big(\CinfK{\check{M}}[[\nu]],*\big)$ dans
 $\big(\CinfK{E^C}[[\nu]],*^C\big)$. De plus, d'apr\`{e}s le th\'{e}or\`{e}me
 \ref{TVarieteEspacesReduitsLocaux} il existe un plongement $j$ de
 la sous-vari\'{e}t\'{e} co\"{\i}sotrope $C$ en tant que sous-vari\'{e}t\'{e} lagrangienne
 de $E^C$ avec $\phi\circ j=i$. Il vient que
 \[
   j^*[*^C]=j^*\frac{[\omega^C]}{\nu}
             +j^*\phi^*\left([*]-\frac{[\omega]}{\nu}\right)
           =  0 + i^*[*]-i^*\frac{[\omega]}{\nu}
           =  0  + 0 -0 = 0
 \]
 car $j(C)$ est une sous-vari\'{e}t\'{e} lagrangienne ($j^*\omega^C=0$) et,
 par hypoth\`{e}se, $i^*[*]=0$, donc $i^*\frac{[\omega]}{\nu}=0$. Alors,
 d'apr\`{e}s le corollaire \ref{CorRepLagArbi} il existe une repr\'{e}sentation
 diff\'{e}rentielle $\rho^C$ de $\big(\CinfK{E^C}[[\nu]],*^C\big)$ sur
 la sous-vari\'{e}t\'{e} lagrangienne $j(C)$ de $E^C$. Puisque la restriction
 $f\mapsto f|_{\check{M}}$ de $\CinfK{M}$ dans $\CinfK{\check{M}}$ induit
 une homomorphisme diff\'{e}rentiel de star-produits, l'application
 $\rho:\CinfK{M}[[\nu]]\ra \Dop\big(\CinfK{C},\CinfK{C}\big)[[\nu]]$,
 d\'{e}finie par
 \[
    \rho(f):=\rho^C\big(\Phi(f|_{\check{M}})\big)
 \]
 est une repr\'{e}sentation diff\'{e}rentielle de $*$ sur $j(C)$, donc sur $C$.
\ebew
On y voit le cas particulier suivant, d\'{e}j\`{a} trait\'{e} implicitement
dans \cite{BHW00}, Lemma 27 et Theorem 32 (o\`{u} la transformation d'\'{e}quivalence
$S$ appara\^{\i}t comme d\'{e}formation de l'application de restriction $i$):
\bcor
 Soit $(M,\omega)$ une vari\'{e}t\'{e} symplectique et $\Phi:G\times M\ra M$
 l'action propre et symplectique ($\Phi_g^*\omega =\omega~~\forall g\in G$)
 du groupe de Lie connexe sur $M$. Soit $\mathfrak{g}$ l'alg\`{e}bre de Lie de
 $G$ et soit $J:M\ra \mathfrak{g}^*$ une application moment pour l'action
 de $G$. Soit $0$ une valeur r\'{e}guli\`{e}re de $J$ et soit la sous-vari\'{e}t\'{e}
 co\"{\i}sotrope $C:=J^{-1}(0)$ de $M$ non vide.
 Alors il existe un star-produit repr\'{e}sentable sur $M$.
\ecor
\bbew
 Il est clair que l'action de $G$ pr\'{e}serve $C$ et y est propre. D'apr\`{e}s le
 corollaire \ref{CActionsPropresClasseAMNulle} la classe d'Atiyah-Molino de
 la vari\'{e}t\'{e} feuillet\'{e}e
 $C$ s'annule, et le th\'{e}or\`{e}me pr\'{e}c\'{e}dent \ref{TQuantSousVarCoisClasAMNulle}
 s'applique.
\ebew

\section{Exemples}
 \label{SecExemples}

 \subsection{Le fibr\'{e} cotangent du $2$-tore}
    \label{SubSecFibCotDeuxTore}

L'exemple qu'on va discuter dans ce paragraphe se trouve plus ou moins
dans \cite[p.135-137]{BHW00}: Soit $M:=T^*T^2$ donn\'{e}e par
$\{(e^{i\varphi},e^{i\psi},p,J)\in S^1\times S^1 \times \real^2
~|~\varphi,\psi,p,J\in\real\}$ munie de
la forme symplectique $d\varphi\wedge dp + d\psi\wedge dJ$.
On consid\`{e}re la sous-vari\'{e}t\'{e} de codimension $1$ (alors co\"{\i}sotrope)
\beq
 C:=\{(e^{i\varphi},e^{i\psi},p,J)\in M~|~J=0\}.
\end{equation}
L'id\'{e}al annulateur $\mathcal{I}$ est visiblement donn\'{e} par
\beq
  \mathcal{I}=\CinfK{M}J=\{fJ\in\CinfK{M}~|~f\in\CinfK{M} \},
\end{equation}
son normalisateur $\mathcal{N}(\mathcal{I})$ par
\beq
  \mathcal{N}(\mathcal{I}) = \left\{f\in\CinfK{M}~\Big|~
    \frac{\partial f}{\partial \psi}=0=\frac{\partial f}{\partial J}\right\}
     ~+~\mathcal{I},
\end{equation}
et l'espace r\'{e}duit $M_r$ est symplectomorphe \`{a} $T^*S^1$ donn\'{e} par
$\{(e^{i\varphi},p)\in S^1\times \real~|~\varphi,p\in\real \}$ et muni
de la forme symplectique canonique $d\varphi\wedge dp$.
Soit $\mu_M:\CinfK{M}\otimes \CinfK{M}\ra\CinfK{M}$ la multiplication
point-par point, et on consid\`{e}re d'abord le star-produit
\beq
  * = \mu_M\circ e^{-2\nu
       \big(\frac{\partial}{\partial p}\otimes
             \frac{\partial}{\partial \varphi}
             +\frac{\partial}{\partial J} \otimes
              \frac{\partial}{\partial \psi}\big)}.
\end{equation}
Ceci est le star-produit sur $T^*T^2$ de type standard, voir
l'\'{e}qn (\ref{EqProdStarStandardRk}), donc la classe de Deligne de
$*$ s'annule.
On voit toute suite que $*$ est un star-produit adapt\'{e} \`{a} $C$ qui est
projetable, donc le commutant de la repr\'{e}sentation $\rho(f)i^*g=i^*(f*g)$
est isomorphe \`{a} $\big(\CinfK{M_r}[[\nu]],*_r\big)$ avec
\beq
  *_r = \mu_{M_r}\circ e^{-2\nu
       \big(\frac{\partial}{\partial p}\otimes
             \frac{\partial}{\partial \varphi}\big)}.
\end{equation}
o\`{u} $\mu_{M_r}$ d\'{e}signe la multiplication point-par-point dans
$\CinfK{M_r}$.

On consid\`{e}re maintenant la transformation d'\'{e}quivalence $S$
\beq
  S:=e^{2i\nu p\frac{\partial}{\partial J}}
%         +4i\nu^2\frac{\partial^2}{\partial \varphi \partial J}
\end{equation}
En utilisant le fait que $p\frac{\partial}{\partial J}$ est une d\'{e}rivation
de la multiplication point-par-point et les identit\'{e}s
\beas
  \left[ p\frac{\partial}{\partial J}\otimes id
            + id \otimes p\frac{\partial}{\partial J}~,~
         \frac{\partial}{\partial p}\otimes
             \frac{\partial}{\partial \varphi}
             +\frac{\partial}{\partial J} \otimes
              \frac{\partial}{\partial \psi}\right]
              & = & -\frac{\partial}{\partial J}\otimes
                      \frac{\partial}{\partial \varphi} \\
  \left[ p\frac{\partial}{\partial J}\otimes id
            + id \otimes p\frac{\partial}{\partial J}~,~
         \frac{\partial}{\partial J}\otimes
             \frac{\partial}{\partial \varphi}\right] & = & 0
\eeas
on montre que le star-produit $*':=S(*)$ est donn\'{e} par
\beq
   *' = \mu_M\circ e^{-2\nu
       \big(\frac{\partial}{\partial p}\otimes
             \frac{\partial}{\partial \varphi}
             +\frac{\partial}{\partial J} \otimes
              (\frac{\partial}{\partial \psi}-2i\nu
               \frac{\partial}{\partial \varphi})\big)}.
\end{equation}
Visiblement, $*'$ est \'{e}galement adapt\'{e} \`{a} $C$.
En choisissant $g(e^{i\varphi},e^{i\psi},p,J)=Je^{i\psi}$ et
$h(e^{i\varphi},e^{i\psi},p,J)=e^{i\varphi}$ on a $g\in\mathcal{I}\subset
\mathcal{N}(\mathcal{I})$ et $h\in \mathcal{N}(\mathcal{I})$, mais
\[
   (g*'h)(\varphi,\psi,p,J) = (J-4\nu^2)e^{i(\varphi+\psi)} \not\in
                             \mathcal{N}(\mathcal{I}),
\]
donc $*'$ n'est pas projetable. En outre, on voit que la repr\'{e}sentation
$\rho'$ de $*'$ d\'{e}finie par $\rho'(f)i^*g=i^*(f*'g)$ a un commutant
tr\`{e}s petit: soit $h$ un \'{e}l\'{e}ment de
l'id\'{e}alisateur de $\mathcal{I}[[\nu]]$ par rapport \`{a} $*'$,
$\mathcal{N}_{*'}(\mathcal{I})$ (voir l'\'{e}qn (\ref{EqDefIdIStar})). En
particulier on a
\[
  0 = i^*(J*'h)=i^*\left(Jh-2\nu\left(\frac{\partial h}{\partial \psi}-
                                   2i\nu\frac{\partial h}{\partial \varphi}
                                   \right)\right),~~\mathrm{alors~~}
  \frac{\partial i^*h}{\partial \psi}=
       2i\nu\frac{\partial i^*h}{\partial \varphi}.
\]
Ceci donne
\bea
    \frac{\partial i^*h_0}{\partial \psi} & = & 0 \label{EqGegenH0}\\
     \frac{\partial i^*h_1}{\partial \psi}  & = &
       2i\frac{\partial i^*h_0}{\partial \varphi} \label{EqGegenH1}\\
         \vdots  & \vdots & \vdots  \nonumber \\
    \frac{\partial i^*h_{r+1}}{\partial \psi}  & = &
       2i\frac{\partial i^*h_r}{\partial \varphi} \label{EqGegenHrplus1}
\eea
Evidemment, l'\'{e}qn (\ref{EqGegenH0}) implique que $i^*h_0$ ne d\'{e}pend pas de
$\psi$. Dans l'\'{e}qn (\ref{EqGegenH1}),
l'int\'{e}gration sur $S^1$ (param\'{e}tr\'{e} par $\psi$) fait dispara\^{\i}tre le membre
de gauche, tandisque le membre de droite est multipli\'{e} par $2\pi$. Il
s'ensuit que $i^*h_0$ ne d\'{e}pend pas non plus de $\varphi$. Alors le membre
de droite de l'\'{e}qn (\ref{EqGegenH1}) s'annule, alors $i^*h_1$ ne d\'{e}pend
pas de $\psi$. Par r\'{e}currence, il s'ensuit que toutes les fonctions
$i^*h_r$ ne d\'{e}pendent ni de $\psi$ ni de $\varphi$. Il est clair que
si $i^*h=P^*\tilde{h}$ (o\`{u} $h\in\CinfK{\real}[[\nu]]$ et $P:C\ra\real$
est la projection sur la variable $p$ le long de $\phi$ et $\psi$),
alors $i^*J*'h=0$.\\
D'un autre c\^{o}t\'{e}, puisque $*'$ est adapt\'{e} les \'{e}l\'{e}ments de
$\mathcal{I}[[\nu]]$ sont de la forme $gJ=g*'J$. Donc si
$i^*h=P^*\tilde{h}$, alors $J*'h\in\mathcal{I}[[\nu]]$, et puisque
$\mathcal{I}[[\nu]]$ est un id\'{e}al \`{a} gauche il vient que
$(gJ)*'h=(g*'J)*'h=g*'(J*'h)\in\mathcal{I}[[\nu]]$, donc $h\in
\mathcal{N}_{*'}(\mathcal{I})$. Par cons\'{e}quent
\[
   \mathrm{Hom}_{(*',\rho')}(C,C)\cong
    \big(\mathcal{N}_{*'}(\mathcal{I})/\mathcal{I}[[\nu]]\big)^{\mathrm{opp}}
    \cong \CinfK{\real}[[\nu]],
\]
et ceci peut \^{e}tre regard\'{e}e comme une sous-alg\`{e}bre commutative propre de
$\big(\CinfK{T^*S^1}[[\nu]],*_r\big)$ qui n'est donc pas isomorphe \`{a}
une d\'{e}formation de l'alg\`{e}bre r\'{e}duite.

Finalement, le m\^{e}me ph\'{e}nom\`{e}ne du trop petit commutant arrive quand
on choisit un star-produit avec une classe de Deligne non nulle, par
exemple
\beq
   *'' = \mu_M\circ e^{-2\nu
       \big(\frac{\partial}{\partial p}\otimes
             \frac{\partial}{\partial \varphi}
             +\frac{\partial}{\partial J} \otimes
              (\frac{\partial}{\partial \psi}-2i\nu
               \frac{\partial}{\partial p})\big)}.
\end{equation}
La classe de Deligne de $*''$ est un multiple non nul dans $\korps[[\nu]]$
de la classe $[d\varphi\wedge d\psi]$. On obtient les m\^{e}mes \'{e}quations
pour le calcul du commutant que pour $*'$, en rempla\c{c}ant $\varphi$ par
$p$ dans les \'{e}quations (\ref{EqGegenH0}), (\ref{EqGegenH1}) et
(\ref{EqGegenHrplus1}), alors on reste avec les fonctions de la variable
$\varphi$, donc
\[
   \mathrm{Hom}_{(*'',\rho'')}(C,C)\cong
    \big(\mathcal{N}_{*''}(\mathcal{I})/\mathcal{I}[[\nu]]\big)^{\mathrm{opp}}
    \cong \CinfK{S^1}[[\nu]].
\]

 \subsection{Encore l'espace projectif complexe comme r\'{e}d\-uc\-tion quantique
           $\ldots$}
    \label{SubSecEspProj}

Dans ce paragraphe, je vais rediscuter l'exemple de l'espace projectif
complexe
en tant que vari\'{e}t\'{e} r\'{e}duite pour illustrer des star-produits r\'{e}ductibles.
La construction a \'{e}t\'{e} d\'{e}j\`{a} faite dans \cite{BBEW96a}, \cite{BBEW96b} et
\cite{Wal98}; ici j'ajoute la construction des star-produits adapt\'{e}s.

On consid\`{e}re la vari\'{e}t\'{e} $M:=\complex^{n+1}\setminus \{0\}$ (munie des
coordonn\'{e}es complexes $z:=(z^1,\ldots,z^{n+1})$ et la forme symplectique
$\omega:=\frac{i}{2}\sum_{k=1}^{n+1}dz^k\wedge d\bar{z}^k$). Soit $\pi$
la projection canonique de $M$ sur l'espace projectif complexe
$\complex P(n)$ comme espace quotient de l'action du groupe
$\complex\setminus \{0\}$. Soit
$x:M\ra\real$ la fonction $z\mapsto |z|^2:=\sum_{k=1}^{n+1}\bar{z}^kz^k$. Soit
$C$ la sous-vari\'{e}t\'{e} $\{z\in M~|~x(z)=1\}$ qui est co\"{\i}sotrope parce
qu'elle est de codimension $1$. Il est bien connu que la vari\'{e}t\'{e}
symplectique r\'{e}duite
$M_{\mathrm{red}}=C/\mathcal{F}$ est donn\'{e} par l'espace projectif complexe,
et les fibres de la projection canonique
$p:C\ra M_{\mathrm{red}}$ co\"{\i}ncident avec les orbites de
l'action du sous-groupe $U(1)$ de $\complex\setminus \{0\}$ sur $C$. On
consid\`{e}re le star-produit suivant sur $M$ pour
$F,G\in\CinfK{M}[[\lambda]]$:
\beq
  F*G := \sum_{r=0}^\infty\frac{\lambda^r}{r!}
                    \sum_{i_1,\ldots,i_r=1}^{n+1}
                    \frac{\partial^r F}{\partial z^{i_1}\cdots\partial z^{i_r}}
           \frac{\partial^r G}{\partial \bar{z}^{i_1}\cdots\partial
                                  \bar{z}^{i_r}}
\end{equation}
o\`{u} on a utilis\'{e} le param\`{e}tre formel $\lambda :=-4i\nu$ de \cite{BBEW96a}.
Une fonction $R$ dans $\CinfK{M}[[\lambda]]$ est dite {\em radiale} lorsqu'il
existe une fonction $\varphi\in\CinfK{\real^+\setminus\{0\}}[[\lambda]]$
telle que $R=\varphi\circ x$. Soient
$E:=\sum_{k=1}^{n+1}z^k\frac{\partial}{\partial z^k}$ et
$\bar{E}:=\sum_{k=1}^{n+1}\bar{z}^k\frac{\partial}{\partial \bar{z}^k}$
deux {\em op\'{e}rateurs d'Euler}. On voit facilement que $i(E-\bar{E})$ est
le champ de vecteurs fondamental de l'action de $U(1)$ sur $M$ et que
$\frac{1}{2}(E+\bar{E})$ s'identifie avec $x\frac{\partial}{\partial x}$
quand on \'{e}crit $M\cong (\real^+\setminus\{0\})\times C$ et utilise $x$
comme variable le long du premier facteur.
On calcule
\bea
    (F*R)(z) & = & \sum_{r=0}^\infty\frac{\lambda^r}{r!}
                    \sum_{i_1,\ldots,i_r=1}^{n+1}
                   z^{i_1}\cdots z^{i_r}~
    \frac{\partial^r F}{\partial z^{i_1}\cdots\partial z^{i_r}}(z)
                    ~\frac{\partial^r \varphi}{\partial x^r}\big((x(z)\big)
                    \nonumber \\
             & = & \sum_{r=0}^\infty\frac{\lambda^r}{r!}
                             x(z)^r ~\Big(\left(\frac{1}{x}E\right)^rF\Big)(z)
             ~ \frac{\partial^r \varphi}{\partial x^r}\big((x(z)\big).
\eea
En particulier, pour deux fonctions radiales $R_1=\varphi_1\circ x$ et
$R_2=\varphi_2\circ x$ on obtient
\beq
    R_1*R_2=\sum_{r=0}^\infty\frac{\lambda^r x^r}{r!}~
                             \frac{\partial^r \varphi_1}{\partial x^r}(x)~
                             \frac{\partial^r \varphi_2}{\partial x^r}(x)
           =:\big(\varphi_1\star\varphi_2\big)(x)
\end{equation}
avec la notation de \cite[p.361]{BBEW96a}. Soit $D=1+\sum_{r=1}\lambda^r d_r$
une s\'{e}rie formelle dans $\korps[[\lambda]]$.
Dans le m\^{e}me travail \cite{BBEW96a}
on avait montr\'{e} l'existence d'une s\'{e}rie $S_D$ d'op\'{e}rateurs diff\'{e}rentiels sur
$\real^+\setminus\{0\}$ qui \'{e}tablie un isomomorphisme entre $\star$ et la
multiplication point-par-point, i.e. qui est telle que
\beq \label{EqSDcommeequivalence}
     S_D(\varphi_1\star \varphi_2)=(S_D\varphi_1)(S_D\varphi_2).
\end{equation}
$S_D$ est d\'{e}finie par son symbole, c.-\`{a}-d. par ses valeurs sur des
fonctions exponentielles: soit $\alpha\in\korps$ et
$e_\alpha(x):=e^{\alpha x}$, alors
\beq
  S_D(e_\alpha)=e_{D\ln(1+\lambda\alpha)/\lambda}
    \mathrm{~~et~~}
  S_D^{-1}(e_\alpha)= e_{(\exp(\lambda\alpha/D)-1)/\lambda},
\end{equation}
voir \cite{BBEW96a} et \cite{BBEW96b} pour des preuves (dans \cite{BBEW96a}
la s\'{e}rie $D$ pouvait d\'{e}pendre de $x$, ce qu'on ne fait pas ici pour avoir une
formule explicite pour $S_D^{-1}$). Si l'on remplace dans l'\'{e}qn
(\ref{EqSDcommeequivalence}) $\varphi_1$ par
$S_D^{-1}(\varphi)$ et $\varphi_2$ par $e_\alpha$ on obtient l'\'{e}quation
de s\'{e}ries d'op\'{e}rateurs diff\'{e}rentiels sur $\real^+\setminus \{0\}$
\beq \label{EqSDundpositvreell}
  (e_{D\ln(1+\lambda\alpha)/\lambda}(x))^{-1}
    S_D\circ \sum_{r=0}^\infty \frac{\lambda^r \alpha^r x^r}{r!}
              \frac{\partial^r}{\partial x^r} \circ S_D^{-1} =
              \mathrm{id}.
\end{equation}
Si l'on regarde $\alpha$ comme un param\`{e}tre formel dans cette \'{e}quation,
alors les coefficients de $\alpha^m\lambda^n$ sont des op\'{e}rateurs
diff\'{e}rentiels sur $\real^+\setminus\{0\}$ dont les coefficients sont
des fonctions rationnelles en $x$. Cette sous-alg\`{e}bre d'op\'{e}rateurs
diff\'{e}rentiels est isomorphe \`{a} l'alg\`{e}bre
$\mathcal{A}$
engendr\'{e}e par des mots $\xi,\xi^{-1}$ et $\eta$ avec les relations
$\xi\xi^{-1}=1=\xi^{-1}\xi$, $\eta\xi-\xi\eta=1$ et
$\eta\xi^{-1}-\xi^{-1}\eta=-\xi^{-2}$ o\`{u} on a identifi\'{e}e $\xi$ avec $x$ et
$\eta$ avec $\frac{\partial}{\partial x}$. Mais on peut \'{e}galement envoyer les
g\'{e}n\'{e}rateurs $\xi,\xi^{-1}$ et $\eta$ sur $x$, $x^{-1}$ et l'op\'{e}rateur
$\frac{1}{x}E$ dans l'alg\`{e}bre des op\'{e}rateurs diff\'{e}rentiels sur $M$ parce
qu'on a les m\^{e}mes relations, surtout
$\frac{1}{x}E(xF)-x(\frac{1}{x}E(F))=F$: par cons\'{e}quent, l'\'{e}quation
(\ref{EqSDundpositvreell}) est toujours satisfaite pour l'op\'{e}rateur
$\hat{S}_D$ o\`{u} on a remplac\'{e} $\frac{\partial}{\partial x}$ par
$\frac{1}{x}E$ (et non pas par $\frac{1}{2x}(E+\bar{E})$ comme dans
\cite{BBEW96a}). Il vient
\beq \label{EqSDtildeinM}
  (e_{D\ln(1+\lambda\alpha)/\lambda}(x))^{-1}
    \hat{S}_D\circ \sum_{r=0}^\infty \frac{\lambda^r \alpha^r x^r}{r!}
              \left(\frac{1}{x}E\right)^r \circ \hat{S}_D^{-1} =
              \mathrm{id},
\end{equation}
et donc
\beq
  \hat{S}_D\big(F*(e_\alpha\circ x)\big)=
        \big(\hat{S}_DF\big)\big(\hat{S}_D(e_\alpha\circ x)\big)
\end{equation}
quels que soient $F\in \CinfK{M}[[\lambda]]$ et $\alpha\in\korps$. Il
s'ensuit que pour le star-produit $\hat{*}_D$ d\'{e}fini par $\hat{S}_D(*)$
on a
\[
     F\hat{*}_D R~=~FR
\]
quels que soient $F\in \CinfK{M}[[\lambda]]$ et $R$ radiale. En
particulier
\[
    F\hat{*}_D(x-1)~=~F~(x-1),
\]
alors l'id\'{e}al annulateur $\mathcal{I}[[\lambda]]$ de $C$ (dont les
\'{e}l\'{e}ments sont des multiples de $(x-1)$ (voir par exemple \cite{BBEW96a}))
est un id\'{e}al \`{a} gauche pour $\hat{*}_D$, donc {\em $\hat{*}_D$ est un
star-produit adapt\'{e} \`{a} $C$.}\\
De plus, puisque $x$ est $U(1)$-invariant et $\frac{1}{x}E$ commute avec
$Y=i(E-\bar{E})$ il vient que la transformation d'\'{e}quivalence $\hat{S}_D$
est $U(1)$-invariante. Il s'ensuit que $\hat{*}_D$ est $U(1)$-invariant
parce que $*$ l'est. Par cons\'{e}quent, le sous-espace
$\CinfK{M}^{U(1)}[[\lambda]]$ de toutes les s\'{e}ries formelles \`{a} coefficients
dans l'espace des fonctions $U(1)$-invariantes sur $M$
est une sous-alg\`{e}bre par rapport \`{a} $\hat{*}_D$. En outre, puisque
$\frac{1}{x}Ef=\frac{1}{2x}(E+\bar{E})f$ pour toutes fonction
$U(1)$-invariante $f$, il s'ensuit que le star-produit $\hat{*}_D$
restreint \`{a} $\CinfK{M}^{U(1)}[[\lambda]]$ co\"{\i}ncide avec le
star-produit $\tilde{*}$ d\'{e}fini dans \cite{BBEW96a}. Il est clair que
l'id\'{e}alisateur de Lie de $\mathcal{I}$, $\mathcal{N}(\mathcal{I})$,
co\"{\i}cide avec
\[
  \mathcal{N}(\mathcal{I})=\CinfK{M}^{U(1)}+\mathcal{I}
\]
et puisque pour tout $f\in\CinfK{M}^{U(1)}[[\lambda]]$ on a d'apr\`{e}s
\cite{BBEW96a}
\[
   (x-1)\hat{*}_D f = (x-1)\tilde{*} f = (x-1)f \in \mathcal{I}[[\lambda]]
\]
il vient que $\mathcal{N}(\mathcal{I})[[\lambda]]$ est une sous-alg\`{e}bre
de $\big(\CinfK{M}[[\lambda]],\hat{*}_D\big)$, alors {\em $\hat{*}_D$ est
un star-produit projetable}. Puisque
\[
  \mathcal{N}(\mathcal{I})/\mathcal{I}=
       \CinfK{M}^{U(1)}/\big(\CinfK{M}^{U(1)}\cap\mathcal{I}\big)
       \cong \CinfK{\complex P^n}
\]
et la restriction du star-produit $\hat{*}_D$ \`{a}
$\CinfK{M}^{U(1)}[[\lambda]]$ co\"{\i}ncide avec $\tilde{*}$ d\'{e}fini dans
\cite{BBEW96a} il s'ensuit que
l'alg\`{e}bre quotient $\CinfK{\complex P^n}[[\lambda]]$ est munie d'un
star-produit r\'{e}duit $*^D$ qui est \'{e}gal au star-produit $*^D_\mu$ de
\cite{Wal98} (pour $D$ arbitraire et $\mu=-1$). Pour le cas $D=1$
ce star-produit r\'{e}duit appara\^{\i}t \'{e}galement dans
\cite{BBEW96a}. Si $D\neq D'$, alors les star-produits r\'{e}duits $*^D$ et
$*^{D'}$ sont non \'{e}quivalents (voir par exemple \cite{Wal98} pour une
d\'{e}monstration). Puisque $H^2_{dR}(\complex P^n,\korps)\cong\korps$ il
vient que l'ensemble de toutes les classes de Deligne des star-produits
sur $\complex P^n$ est en bijection avec toutes les s\'{e}ries formelles
$D$.

On a donc le ph\'{e}nom\`{e}ne suivant: bien que --pour $D\neq D'$-- les deux
star-produits projetables $*_{D}$ et $*_{D'}$ soient \'{e}quivalents, ils ne
le sont pas par rapport \`{a} une transformation d'\'{e}quivalence adapt\'{e}e
selon la proposition \ref{PCouplesEquiSiTrafoAdap}. Pour le changement
de transformations d'\'{e}quivalence $S_{D'}S_D^{-1}$ on voit explicitement
qu'elle ne pr\'{e}serve pas l'id\'{e}al annulateur: on
calcule
\beq
  S_{D'}S_D^{-1}e_\alpha = e_{\frac{D'}{D}\alpha},
\end{equation}
alors, avec
$\hat{*}_{D'}=\hat{S}_{D'}\big(\hat{S}_D^{-1}(\hat{*}_D)\big)$, on voit
ais\'{e}ment que
\beq \label{EqTrafodEqtransitoire}
  \hat{S}_{D'}\hat{S}_D^{-1}=e^{\big(\ln(D'/D)\big)E},
\end{equation}
alors (puisque $Ex=x$)
\beq
  \hat{S}_{D'}\hat{S}_D^{-1}(x)=\frac{D'}{D}x,
\end{equation}
et finalement pour tout $F\in\CinfK{M}[[\lambda]]$ on a
\beq
  \hat{S}_{D'}\hat{S}_D^{-1}\big(F~(x-1)\big)
   =\big(\hat{S}_{D'}\hat{S}_D^{-1}(x)(F)\big)\big(\frac{D'}{D}x-1\big)
\end{equation}
gr\^{a}ce \`{a} l'\'{e}qn (\ref{EqTrafodEqtransitoire}) et le fait que $E$ est une
d\'{e}rivation pour la multiplication point-par-point. Ceci montre pour ce
changement de transformations d'\'{e}quivalence que l'id\'{e}al annulateur n'est
pas pr\'{e}serv\'{e}, car $x-1$ est envoy\'{e} sur $(D'x/D)-1$ qui n'appartient pas
\`{a} $\mathcal{I}[[\lambda]]$ pour $D\neq D'$.

On voit aussi que le th\'{e}or\`{e}me \ref{Tdaserstemalfalsch} et son corollaire
\ref{Cdaserstemalfalsch} ne s'appliquent pas \`{a} cet exemple car
\[
   H^1_v(S^{2n+1},\korps)\cong H^1_{dR}(S^1,\korps)\cong \korps \neq
   \{0\}.
\]

\section{Probl\`{e}mes ouverts}
  \label{SecRemProblemOuverts}

\noindent Les probl\`{e}mes suivants me semblent int\'{e}ressants:

\ben
 \item Chaque morphisme de star-produits est-il diff\'{e}rentiel (compare
   la d\'{e}\-fi\-ni\-tion \ref{DDefHomStarGenDiff}) ? J'ai
   l'impression qu'il devrait y avoir des r\'{e}ponses positives au moins
   pour des applications de Poisson $M'\ra M$ `bien comport\'{e}es' par
   une construction ordre-par-ordre en utilisant l'\'{e}quation de morphisme
   et
   le th\'{e}or\`{e}me de Hochschild-Kostant-Rosenberg (voir par exemple
   \cite{HKR62}, \cite{CDG80}, \cite{Pfl98a}, \cite{Nad99}, \cite{GR99})
   pour la cohomologie de
   Hochschild de $\CinfK{M}$ \`{a} valeurs dans $\CinfK{M'}$.
 \item Chaque repr\'{e}sentation de star-produits sur $\CinfK{C}[[\nu]]$
    est-elle dif\-f\'{e}\-ren\-tielle (compare la d\'{e}finition \ref{DRepDiff}) ?
    J'ai la m\^{e}me impression que pour le probl\`{e}me
    pr\'{e}c\'{e}dent: l'identit\'{e} de r\'{e}pr\'{e}sentation pose un probl\`{e}me
    cohomologique
    ordre-par-ordre; et il s'agit de calculer la cohomologie de Hochschild
    de $\CinfK{M}$ \`{a} valeur dans $\Dop^1\big(\CinfK{C};\CinfK{C}\big)$.
 \item Pour le probl\`{e}me g\'{e}n\'{e}ral pers\'{e}v\'{e}rant de quantifier les
   applications moment,
  je ne connais pas d'exemple pour la situation assez \'{e}l\'{e}mentaire suivante:
  soient $f_0$ et
  $g_0$ deux polyn\^{o}mes sur $\real^{2n}$ (muni de la structure symplectique
  canonique) dont les diff\'{e}rentielles $df_0$ et $dg_0$ soient ind\'{e}pendants
  sur un ouvert et qui commutent par rapport au crochet de Poisson, i.e.
  $\{f_0,g_0\}=0$. Trouve-t-on toujours un star-produit $*$ sur
  $\real^{2n}$ et des termes d'ordre sup\'{e}rieur $f_1,\ldots$ et
  $g_1,\ldots$ tels que $f*g-g*f=0$ pour $f:=\sum_{r=0}^\infty\nu^r f_r$ et
  $g:=\sum_{r=0}^\infty\nu^r g_r$, ou y a-t-il un contre-exemple?
 \item Je ne connais pas d'exemple concret non plus pour une sous-vari\'{e}t\'{e}
  co\"{\i}sotrope ferm\'{e}e $C$ d'une vari\'{e}t\'{e} symplectique $(M,\omega)$
  telle que l'obstruction (\ref{EqObstructionImportante}) du th\'{e}or\`{e}me
  \ref{TObstructionsOrdreTrois} --qui implique la classe d'Atiyah-Molino
  du feuilletage de $C$-- ne soit pas nulle pour tout choix de
  $2$-formes ferm\'{e}es $\alpha_0$ et $\alpha_1$.
 \item M\^{e}me probl\`{e}me pour la quantification des morphismes de Poisson
   entre deux vari\'{e}t\'{e}s symplectiques, comparer l'\'{e}quation
   (\ref{EqObstructionImportantePoisson}) du th\'{e}or\`{e}me
   \ref{TQuantApplPoissCondNec}: il faudrait construire un exemple
   concret.
\een


\begin{thebibliography}{99}
\bibitem{AM85} Abraham, R., Marsden, J.~E.:
         {\em Foundations of Mechanics}, second edition.
         Addison Wesley Publishing Company, Inc., Reading Mass. 1985.

\bibitem {ACMP83}
 Arnal, D., Cortet, J.~C., Molin, P., Pinczon, G.:{\em
  Covariance and Geometrical Invariance in {$*$-Quantization}}.
  J. Math. Phys.  {\bf 24}.2 (1983), 276--283.

\bibitem {AMM99} Arnal, D., Manchon, D., Masmoudi, M.:
     {\em Choix des signes pour la formalit\'{e} de Kontsevitch}.
     Pacific J. Math. {\bf 203} (2002), 23-66,
     {\tt math.QA/0003003}.


\bibitem{BDM03}
  Baklouti, A., Dhieb, S., Manchon, D.:
  {\em Orbites coadjointes et vari\'{e}t\'{e}s caract\'{e}ristiques}.
 {\tt math.RT/0302171}, 2003.


\bibitem{BW95}
        Bates, S., Weinstein, A.: {\em Lectures on the Geometry on
        Quantization.} Berkeley Mathematics Lecture Notes, Volume 8,
        1995.

\bibitem{BFFLS78}
  Bayen, F., Flato, M., Fr{\o}nsdal, C., Lichnerowicz, A., Sternheimer, D.:
  {\em Deformation Theory and Quantization.}
  Annals of Physics {\bf 111} (1978), part I: 61-110, part II: 111-151.

\bibitem{BCG97} Bertelson, M., Cahen, M., Gutt, S.: {\em Equivalence of
       Star-Products.} Class.Quant.Grav. {\bf 14} (1997), A93-A107.

\bibitem{Bon98} Bonneau, P.: {\em Fedosov Star-Products and
       $1$-Differentiable Deformations.} {\tt math.QA/9809032}, septembre
       1998.

\bibitem{BGHHW03p} M.~Bordemann, G.~Ginot, G.~Halbout, H.-C.~Herbig,
       S.~Waldmann: {\em Star-repr\'{e}sentations sur des sous-vari\'{e}t\'{e}s
       co-isotropes}, {\tt math.QA/0309321}, septembre 2003.

\bibitem{BNPW03} Bordemann, M., Neumaier, N., Pflaum, M., Waldmann, S.:
      {\em On representations
      of star product algebras over cotangent spaces on Hermitian line
      bundles}, J.~Funct.~Analysis {\bf 199} (2003), 1-47,
      {\tt math.QA/9811055}.

\bibitem{Bor02p} Bordemann, M.: {\em Sur l'existence d'une prescription
    d'ordre
    naturelle projectivement invariante},
     {\tt math.DG/0208171}.

\bibitem{Bor00} Bordemann, M.: {\em The deformation quantization of
      certain super-Poisson brackets and BRST cohomology}. Dans:
      Dito, G., Sternheimer, D.: {\em Conf\'{e}rence Mosh\'{e} Flato 1999. Volume
      II}. Kluwer, Dordrecht, 2000, 45-68.

\bibitem{BHW00}
     Bordemann, M., Herbig, H.-C., Waldmann, S.: {\em BRST cohomology and
     Phase Space
     Reduction in Deformation Quantisation}, Commun.Math.Phys. {\bf 210}
     (2000), 107-144.

\bibitem{Bor99v} Bordemann, M.: {\em (GNS) Representations and (KMS) States
        in Deformation Quantization}, Contribution au `XVII-th Workshop on
        Geometric Methods in Physics' \`{a} Bia{\l}owie\.{z}a, Pologne,
        3.7.-9.7.1998:
        Schlichenmaier, Martin, Ali, S.T., Strasburger, A., Odzijewicz, A.
        (Ed.):
        {\em Coherent States, Quantization and Gravity, Proceedings of the
        XVII-th Workshop
        On Geometrc Methods in Physics}, Warsaw University Press,
        Varsovie, 2001.

\bibitem{BW99}
      Bordemann, M., Walter, M.: {\em Quantum Integrable Toda-Like Systems
      in Deformation Quantization.} Lett.~Math.~Phys. {\bf 48} (1999),
      123-133.

\bibitem{BNW99}
      Bordemann, M., Neumaier, N., Waldmann, S.: {\em Homogeneous Fedosov
      Star Products on Cotangent Bundles II: GNS Representations, the WKB
      Expansion, traces, and applications},
      J.~Geom. Phys.~{\bf 29} (1999), 199-234.

\bibitem{BNW98}
      Bordemann, M., Neumaier, N., Waldmann, S.: {\em Homogeneous Fedosov
      Star Products on Cotangent Bundles I: Weyl and Standard Ordering
      with Differential Operator Representation},
      Comm.~Math. Phys.~{\bf 198} (1998), 363-396.

\bibitem {BRW98}
       M.~Bordemann, H.~R\"omer, S.~Waldmann: {\em A Remark on Formal
       KMS States in Deformation Quantization},
       Lett.~Math.~Phys.~{\bf 45} (1998), 49-61.

\bibitem {BW98}
     M.~Bordemann, S.~Waldmann:
     {\em Formal GNS Construction and States in Deformation Quantization},
     Comm.~Math.~Phys.~{\bf 195} (1998), 549-583.


\bibitem {BBEW96b}
  M.~Bordemann, M.~Brischle, C.~Emmrich, S.~Waldmann:
 {\em Subalgebras with Converging Star Products in Deformation Quantization:
 An Algebraic Construction for $\complex P^n$}, J.~Math.~Phys. {\bf 37}
  (1996), 6311-6323.

\bibitem{BBEW96a}
    M.~Bordemann, M.~Brischle, C.~Emmrich, S.~Waldmann:
  {\em Phase Space Reduction for Star Products: An Explicit Construction
  for $\complex P^n$}, Lett.~Math. ~Phys. {\bf 36} (1996), 357-371.

\bibitem{Bot72} Bott, R.: {\em Lectures on characteristic classes and
      foliations}. Dans Bott, R., Gitler, S., James, I.M.: {\em Lectures on
      Algebraic and Differential Topology}, LNM {\bf 279}, Springer,
      Berlin, 1972, 1-94.

\bibitem{BJ73} Br\"{o}cker, T., J\"{a}nich, K.: {\em Einf\"{u}hrung in die
    Differentialtopologie.} Springer Verlag, Berlin 1973.

\bibitem {bursztyn.waldmann:2000b}
Bursztyn, H., Waldmann, S.: {\em Deformation Quantization of
  Hermitian Vector Bundles}. Lett. Math. Phys.  {\bf 53} (2000), 349--365.

\bibitem {bursztyn.waldmann:2000a}
Bursztyn, H., Waldmann, S.: {\em On Positive Deformations of
  {$^*$}-Algebras}. In: Dito, G., Sternheimer, D. (Ed.): {\em
  Conf\'{e}rence Mosh\'{e} Flato 1999. Quantization, Deformations, and
  Symmetries}. \cite{DS00},   69--80.

\bibitem {bursztyn.waldmann:2001b}
 Bursztyn, H., Waldmann, S.: {\em {$^*$}-Ideals and Formal
  Morita Equivalence of {$^*$}-Algebras}.
 Int. J. Math.  {\bf 12}.5 (2001), 555--577.

\bibitem {bursztyn.waldmann:2001a}
  Bursztyn, H., Waldmann, S.: {\em Algebraic Rieffel Induction,
  Formal Morita Equivalence and Applications to Deformation Quantization}.
 J. Geom. Phys.  {\bf 37} (2001), 307--364.

\bibitem {bursztyn.waldmann:2002a:pre}
 Bursztyn, H., Waldmann, S.: {\em Bimodule deformations,
  {P}icard groups and contravariant connections}.
  Pr\'{e}publication (Freiburg FR-THEP 2002/10)  {\tt math.QA/0207255}
  (juillet 2002), 32 pages, \`{a} para\^{\i}tre dans K-Theory.

\bibitem {bursztyn.waldmann:2002a}
 Bursztyn, H., Waldmann, S.: {\em The characteristic classes of
  {M}orita equivalent star products on symplectic manifolds}.
  Commun. Math. Phys.  {\bf 228} (2002), 103--121.

\bibitem{CDG80} Cahen, M., DeWilde, M., Gutt, S.: {\em Local cohomology
        of the algebra of $\Cinf$-functions on a connected manifold.}
        Lett.~Math.~Phys. {\bf 4} (1980), 157-167.

\bibitem{CG83} Cahen, M., Gutt, S.: {\em Regular $*$-representations of
      Lie Algebras}. Lett.Math.Phys. {\bf 6} (1983), 395-404.

\bibitem{CRM75} Calogero, P., Ragnisco, O, Marchioro, C: {\em Exact
      solutions of the classical and quantal one-dimensional many-body
      problems with the two-body potential $V_a(x)=g^2a^2/\sinh^2(ax)$}.
      Lett.~Nuovo Cimento {\bf 13} (1975), 383-387.

\bibitem{CHW98} Cannas da Silva, A., Hartshorn, K., Weinstein, A.:
      {\em Lectures on Geometric Models for Noncommutative Algebras.}
      University of Berkeley, 1998.

\bibitem{CF03p} Cattaneo, A., Felder, G.: {\em Coisotropic submanifolds in
      Poisson Geometry and Branes in the Poisson Sigma Model}
      {\tt math.QA/0309180}, septembre 2003.

\bibitem{CF99} Cattaneo, A., Felder, G.: {\em A path integral approach to the
      Kontsevitch quantization formula}. Commun.Math.Phys. {\bf 212}
      (2000), 591-611.

\bibitem{CFT00} Cattaneo, A., Felder, G., Tomassini, L.:
      {\em From local to global deformation quantization of Poisson
      manifolds}. Duke Math. Journal {\bf 115} (2002), 329-352,
      {\tt math.QA/0012228}.

\bibitem{CE56} Cartan, H., Eilenberg, S.: {\em Homological Algebra}.
        Princeton University Press, Princeton, 1956.

\bibitem{Con94} Connes, A.: {\em Noncommutative Differential Geometry.}
       Academic Press, San Diego, 1994.

\bibitem {CFS92} Connes, A., Flato, M., Sternheimer, D.:
      {\em Closed Star Products and Cyclic Cohomology.} Lett.Math.Phys.
      {\bf 24} (1992), 1-12.

\bibitem {Del95} Deligne, P.: {\em D\'{e}formations de l'alg\`{e}bre des fonctions
        d'une vari\'{e}t\'{e} symplectique: comparaison entre Fedosov et DeWilde,
        Lecomte.} Sel.Math., New series {\bf 1 (4)} (1995), 667-697.

\bibitem{DL83a} DeWilde, M., Lecomte, P.B.A.:
   {\em Star-products on cotangent bundles.} Lett.~Math.~Phys. {\bf 7}
   (1983), 235-241.

\bibitem{DL83} DeWilde, M., Lecomte, P.B.A.:
  {\em Existence of star-products and of formal deformations
   of the Poisson Lie Algebra of arbitrary symplectic manifolds.}
  Lett.~Math.~Phys. {\bf 7} (1983), 487-49.

\bibitem{DL88} DeWilde, M., Lecomte, P.B.A.: {\em Formal Deformations of the
   Poisson Lie Algebra of a Symplectic Manifold and Star Products. Existence,
   Equivalence, Derivations.} Dans: Hazewinkel, M., Gerstenhaber, M. (eds.)
   {\em Deformation Theory of Algebras and Structures and Applications.}
   Dordrecht, Kluwer, 1988.

\bibitem{DM93} Dimakis, A., M\"{u}ller-Hoissen, F.:
   {\em Stochastic differential calculus, the Moyal $*$-product, and
    noncommutative differential geometry}. Lett.Math.Phys.
    {\bf 28} (1993), p. 123-137.

\bibitem{DS00}  Dito, G., Sternheimer, D.: {\em Conf\'{e}rence Mosh\'{e} Flato 1999.
    Volume I, II}. Kluwer, Dordrecht, 2000, 45-68.

\bibitem{Dri83} Drinfel'd, V.: {\em On constant quasiclassical solutions of
        the
        Yang-Baxter Quantum Equation.} Sov.Math.Doklady {\bf 28} (1983),
        667-671.

\bibitem{Dub87} Dubois-Violette, M.: {\em Syst\`{e}mes dynamiques contraints:
           l'approche
             homologique.} Ann.~Inst.~Fourier {\bf 37}, No.4, 45-57 (1987).

\bibitem{DLO99} Duval, C., Lecomte, P.B.A., Ovsienko, V.:
  {\em Conformally equivariant quantization: existence and uniqueness.}
  Ann. Inst. Fourier, Grenoble {\bf 49}, 6 (1999), 1999-2029.

\bibitem{Fed85} Fedosov, B.: {\em Formal Quantization.}
  Some Topics of Modern Mathematics and Their Applications
   to Problems of Mathematical Physics, Moscow (1985), 129-136.

\bibitem{Fed94} Fedosov, B.:
   {\em A Simple Geometrical Construction of Deformation Quantization.}
    J. of Diff.~Geom. {\bf 40} (1994), 213-238.

\bibitem{Fed94b} Fedosov, B.:
   {\em Reduction and eigenstates in deformation quantization}. Dans:
   Demuth, M., Schrohe, E., Schulze, B.-W. (\'{e}d.): {\em Pseudodifferential
   calculus and
    mathematical physics}, vol 5, in {\em Advances in Partial Differential
    Equations},
    277-297. Akademie Verlag, Berlin, 1994.

\bibitem{Fed96} Fedosov, B.: {\em Deformation Quantization and Index Theory.}
           Akademie Verlag, Berlin, 1996.

\bibitem{Fed98} Fedosov, B.: {\em Nonabelian reduction in Deformation
             Quantization}.
             Lett.~Math.~Phys. {\bf 43} (1998), 137-154.


\bibitem {FHST89} Fish, J., Henneaux, M., Stasheff, J., Teitelboim, C.:
     {\em Existence, uniqueness and cohomology of the classical BRST charge
     with ghosts of ghosts}. Commun. Math. Phys.  {\bf 120} (1989), 379--407.

\bibitem {Ger63} Gerstenhaber, M.: {\em The Cohomology Structure of an
    Associative Ring}. Ann. Math. {\bf 78} (1963), 267-288.

\bibitem{Ger64} Gerstenhaber, M.: {\em On the deformation of rings and
algebras}. Ann.Maths {\bf 79} (1964), 59-103.

\bibitem{Ger66} Gerstenhaber, M.: {\em On the deformation of rings and
algebras II}. Ann.Maths {\bf 84} (1966), 1-19.

\bibitem{Ger68} Gerstenhaber, M.: {\em On the deformation of rings and
algebras III}. Ann.Maths {\bf 88} (1968), 1-34.

\bibitem{Ger74} Gerstenhaber, M.: {\em On the deformation of rings and
algebras IV}. Ann.Maths {\bf 99} (1974), 257-276.

\bibitem{Glo98} Gl\"{o}{\ss}ner, P.: {\em Star-Product Reduction for Coisotropic
     Submanifolds of Codimension 1}. Pr\'{e}publication Facult\'{e} de Physique de
     l'Universit\'{e} de Freiburg FR-THEP-98/10, {\tt math.QA/9805049}, mai
     1998.

\bibitem{Got82} Gotay, M.: {\em On coisotropic imbeddings of presysmplectic
       manifolds.} Proc.~Am.~Math.~Soc. {\bf 84} (1982), 111-114.

\bibitem{Gut83} Gutt, S.: {\em An explicit $*$-product on the cotangent
         bundle of a Lie group.} Lett.Math.Phys. {\bf 7} (1983), 249-258.

\bibitem{GR99} Gutt, S., Rawnsley, J.: {\em Equivalence of
    star products on a symplectic manifold; an introduction to Deligne's
    \v{C}ech cohomology classes}. J.Geom.Phys. {\bf 29} (1999),
    347-399.

\bibitem{GR03} Gutt, S., Rawnsley, J.: {\em Natural star products on
       symplectic
       manifolds and quantum moment maps.}
       {\tt math.SG/0304498v2}.

\bibitem{Haa92} Haag, R.: {\em Local Quantum Physics}. Springer Verlag,
       Berlin, 1992.

\bibitem {halbout:2002a}
  Halbout, G. (Ed.): {\em Deformation Quantization}, tome~1 in
  {\em IRMA Lectures in Mathematics and Theoretical Physics}.
  Walter de Gruyter, Berlin, New York, 2002.


\bibitem {HT92}
  Henneaux, M., Teitelboim, C.: {\em Quantization of Gauge
  Systems}. Princeton University Press, New Jersey, 1992.

\bibitem{Hes81} He{\ss}, H.: {\em Symplectic connections in geometric
                 quantization and factor orderings}. Dissertation
                 (Fachbereich Physik, Freie Universit\"{a}t Berlin, F.R.G.,
                 1981).

\bibitem{HKR62} Hochschild, G., Kostant, B., Rosenberg, A.:
      {\em Differential forms on regular affine algebras.}
      Trans.~Am.~Math.~Soc. {\bf 102} (1962), 383-408.

\bibitem{KT75} Kamber, F., Tondeur, P.:
          {\em Foliated Bundles and Characteristic
          Classes.} Lecture Notes in Mathematics {\bf 493}, Springer,
          Berlin 1975.

\bibitem{Kim92} Kimura, T.: {\em Prequantum BRST Cohomology.}
      Contemp.~Math. {\bf 132} (1992), 439-457.

%\bibitem {KN63} Kobayashi, S. \& Nomizu, K. Fondations of Differential
%       Geometry,
%        Vol I, Insterscience Publishers, Wiley and Sons, 1963.

\bibitem {KMS93} Kol\'{a}\v{r}, I., Michor, P., Slov\'{a}k, J.: {\em Natural
     Operations in Differential Geometry.} Springer, Berlin, 1993.

\bibitem{Kon97b} Kontsevitch, M.: {\em Deformation Quantization of Poisson
           Manifolds. I.}
             Pr\'{e}publication IHES, {\tt q-alg}/9709040, septembre 1997.

\bibitem{KNP04p} Kowalzig, N., Neumaier, N., Pflaum, M.:
        {\em Phase space reduction of star-products on cotangent bundles.}
        Pr\'{e}publication {\tt arXiv:math.SG/0403239v1}, mars 2004.

\bibitem{Lam99} Lam, T.Y.: {\em Lectures on Modules and Rings}.
   Springer Verlag, New York, 1999.

\bibitem{Lan95} Lang, S.: {\em Differential and Riemannian Manifolds.}
      New York - Berlin - Heidelberg, Springer Verlag, 1995.

\bibitem{LPS01} Landsman, N.P., Pflaum, M., Schlichenmaier, M. (\'{e}diteurs):
    {\em Quantization of Singular Symplectic Quotients.} Birkh\"{a}user,
    Basel, 2001.

\bibitem{LO99} Lecomte, P.B.A., Ovsienko, V.Y.: {\em Projectively Equivariant
   Symbol Calculus}. Lett.Math,Phys. {\bf 49} (1999), 173-196.


\bibitem{Lic82} Lichnerowicz, A.: {\em D\'{e}formations d'alg\`{e}bres associ\'{e}es
    \`{a} une vari\'{e}t\'{e} symplectique (les $*_\nu$-produits}. Ann.~Inst.~Fourier
    {\bf 32} (1982), 157-209.

\bibitem{Lod92} Loday, J.-L.: {\em Cyclic Homology.} Springer, Berlin,
       1992.

\bibitem{Lu93} Lu, J.-H.: {\em Moment Maps at the Quantum Level.}
     Comm.~Math.~Phys.~{\bf 157} (1993), 389-404.

\bibitem{MR86} Marsden, J.E., Ra\c{t}iu, T.: {\em Reduction of Poisson
     manifolds.} Lett.Math.Phys. {\bf 11} (1986), 161-169.

\bibitem{MW74} Marsden, J.E., Weinstein, Alan: {\em Reduction of symplectic
           manifolds
            with symmetry}. Rep. on Math.~Phys. {\bf 5} (1974), 121-130.


\bibitem{Mol71a} Molino, P.: {\em Classe d'Atiyah d'un feuilletage et
      connexions transversales projetables}. C.R.Acad.Sc. Paris
      {\bf 272} (1971), 779-781.

\bibitem{Mol71b} Molino, P.: {\em Classes caract\'{e}ristiques et obstruction
      d'Atiyah pour les fibr\'{e}s principaux feuillet\'{e}s}.
      C.R.Acad.Sc. Paris {\bf 272} (1971), 1376-1378.

\bibitem{Mol73} Molino, P.: {\em Propri\'{e}t\'{e}s cohomologiques et propri\'{e}t\'{e}s
      topologiques des feuilletages \`{a} connexion transverse projetable}.
      Topology {\bf 12} (1973), 317-325.

\bibitem{Mol74} Molino, P.: {\em La classe d'Atiyah d'un feuilletage
      comme cocycle de d\'{e}formation infinit\'{e}simale}. C.R.Acad.Sc. Paris
      {\bf 278} (1974), 719-721.

\bibitem{Mol88} Molino, P.: {\em Riemannian foliations.} Birkh\"{a}user,
      Boston, Basel 1988.

\bibitem{Nad99} Nadaud, F.: {\em On continuous and differential Hochschild
         cohomology.} Lett.~Math.~Phys. {\bf 47} (1999), 85-95.

\bibitem{NV81} Neroslavski, O., Vlassov, A.T.: {\em Sur les d\'{e}formations
     de l'alg\`{e}bre des fonctions d'une vari\'{e}t\'{e} symplectique.} C.R.Acad.Sc.
     Paris I {\bf 292} (1981), 71-73.

\bibitem{NT95a} Nest, R., Tsygan, B.: {\em Algebraic Index Theorem.}
     Commun.Math.Phys. {\bf 172} (1995), 223-262.

\bibitem{NT95b} Nest, R., Tsygan, B.: {\em Algebraic Index Theorem for
    Families.} Adv.Math. {\bf 113} (1995), 151-205.

\bibitem{Neu01} Neumaier, N.: {\em Klassifikationsergebnisse in der
    Deformationsquantisierung.} Th\`{e}se de doctorat \`{a} la facult\'{e} de physique
    de l'Universit\'{e} de Fribourg-en-Brisgau, RFA, octobre 2001.

\bibitem{Neu02} Neumaier, N.: {\em Local $\nu$-Euler Derivations and Deligne's
    Characteristic Class for Fedosov Star-Products and Star-Products of
    Special Type.} Commun.~Math.~Phys.~{\bf 230} (2002), 271-288, voir
    aussi {\tt math.QA/9905176v2}.

\bibitem{OP83} Olshanetzky, M.~A., Perelomov, A.~M.:{\em
        Quantum integrable systems related to Lie algebras}.
        Phys.~Rep. {\bf 94} (1983), 313-404.

\bibitem{OMY91} Omori, H., Maeda, Y., Yoshioka, A.: {\em Weyl manifolds and
    Deformation Quantization}. Adv.Math. {\bf 85} (1991), 224-255.

\bibitem{OMY92} Omori, H., Maeda, Y., Yoshioka, A.: {\em Existence of a
Closed Star-Product.} Lett.Math.Phys. {\bf 26} (1992), 285-294.

\bibitem{Pal61} Palais, R.: {\em On the Existence of Slices for Actions
      of non-compact Lie Groups.} Ann.Math. {\bf 73} (1961), 295-323.

\bibitem{Pfl01} Pflaum, M.: {\em Analytic and Geometric Study of Stratified
      Spaces.} LNM 1768, Springer, Berlin 2001.

\bibitem{Pfl98b} Pflaum, M.: {\em The normal symbol of Riemannian manifolds.}
        New York J.~Math. {\bf 4} (1998), 97-125.

\bibitem{Pfl98a} Pflaum, M.: {\em On Continuous Hochschild Homology and
        Cohomology Groups.} Lett.~Math.~Phys. {\bf 44} (1998), 43-51.

\bibitem{Sch97} Schirmer, J.: {\em A star-product for Grassmann manifolds}.
        Pr\'{e}publication Facult\'{e} de physique de l'universit\'{e} de Freiburg,
        {\tt q-alg/9709021}, septembre 1997.

\bibitem{SL91} Sjamaar, R., Lerman, E.: {\em Stratified symplectic
       spaces and reduction}. Ann. Math. {\bf 134} (1991), 375-422.

\bibitem{Ste98} Sternheimer, D.: {\em Deformation Quantization: Twenty Years
       After.} Dans: Rembieli\'{n}ski, J. (\'{e}d.) {\em Particles, Fields,
       and Gravitation. (\L \'{o}d\'{z} 1998)}. AIP Press, New York, 1998, p.
       107-145; voir aussi {\tt math.QA/9809056}.

\bibitem{Ton88} Tondeur, P.: {\em Foliations on Riemannian Manifolds.}
      Springer, New York, 1988.

\bibitem{Ton97} Tondeur, P.: {\em Geometry of Foliations.} Birkh\"{a}user,
Basel, 1997.

\bibitem{Vai94} Vaisman, I.: {\em Lectures on the Geometry of Poisson
              Manifolds}. Birkh\"{a}user, Basel, 1994.

\bibitem{Wal98} Waldmann, S.: {\em A Remark on Nonequivalent Star Products
   via Reduction for $\complex \mathbb{P}^n$}. Lett.Math.Phys.
   {\bf 44} (1998), 331-338.

\bibitem{Wal00} Waldmann, S.:
   {\em Locality in GNS Representations of Deformation Quantization.}
   Commun.Math.Phys. {\bf 210} (2000), 467-495.

\bibitem {waldmann:2002b}
 Waldmann, S.: {\em {M}orita equivalence of {F}edosov star
  products and deformed {H}ermitian vector bundles}.
  Lett. Math. Phys.  {\bf 60} (2002), 157--170.

\bibitem {waldmann:2002a}
 Waldmann, S.: {\em On the representation theory of deformation
  quantization}. In: Halbout, G. (Ed): {\em Deformation
  quantization}. \cite{halbout:2002a},   107--133.

\bibitem{Wei71} Weinstein, A.: {\em Symplectic manifolds and their
   Lagrangian submanifolds.} Adv.~Math.~{\bf 6} (1971), 329-346.

\bibitem{Wei77b} Weinstein, A.: {\em Lectures on symplectic manifolds.}
    Dans C.B.M.S.~Conf.~Series, Am.~Math.~Soc. no. {\bf 29}, Providence,
    R.I. 1977.

\bibitem{Wei88} Weinstein, A.: {\em Coisotropic calculus and Poisson
   groupoids.} J.Math.Soc.Japan {\bf 40} (1988), 705-727.

\bibitem{Wei95} Weinstein, A.: {\em Deformation Quantization, S\'{e}minaire
      Bourbaki.} Vol. 1993/94. Ast\'{e}risque {\bf 227} (1995), Exp. No. 789,
      5, p. 389-409.

\bibitem {weinstein.xu:1998a}
Weinstein, A., Xu, P.: {\em Hochschild cohomology and
  characteristic classes for star-products}.
 Dans: Khovanskij, A., Varchenko, A., Vassiliev, V. (ed.):
  {\em Geometry of differential equations. Dedicated to V. I. Arnold
  on the occasion of his 60th birthday},   177--194. Am. Math.
 Soc., Providence, R.I. 1998.

\bibitem {widom:1980a}
  Widom, H.: {\em A Complete Symbolic Calculus for
  Pseudodifferential Operators}.
  Bull. Sc. Math.  {\bf 104} (1980), 19--63.


\bibitem {Xu98}
Xu, P.: {\em Fedosov $*$-Products and Quantum Momentum Maps}.
     Commun. Math. Phys.  {\bf 197} (1998), 167--197.

\end{thebibliography}
\end{document}